\documentclass[10pt,a4paper,reqno,oneside]{amsart}
\usepackage[T1]{fontenc}
\usepackage{lmodern,amssymb,mathtools,booktabs,longtable,array,microtype}
\usepackage{enumitem,needspace,etoolbox}
\usepackage[hidelinks]{hyperref}
\usepackage{url}
\usepackage[textwidth=184mm,textheight=268mm,centering,footskip=18pt]{geometry}
\setlist[enumerate]{itemsep=0pt,parsep=0pt,topsep=2pt}
\newtheoremstyle{compactplain}{3pt}{3pt}{\itshape}{}{\bfseries}{.}{.5em}{}
\theoremstyle{compactplain}
\newtheorem{theorem}{Theorem}[section]
\newtheorem{lemma}[theorem]{Lemma}
\newtheorem{proposition}[theorem]{Proposition}
\newtheorem{corollary}[theorem]{Corollary}
\theoremstyle{definition}

\theoremstyle{remark}
\newtheorem{remark}[theorem]{Remark}
\numberwithin{equation}{section}
\newcommand{\eq}{\approx}
\newcommand{\SR}{\mathsf{SR}}
\newcommand{\Var}{\operatorname{Var}}
\newcommand{\Id}{\operatorname{Id}}
\newcommand{\HSP}{\mathsf{HSP}}
\newcommand{\supp}{\operatorname{supp}}
\newcommand{\head}{\operatorname{hd}}
\newcommand{\tail}{\operatorname{tl}}
\newcommand{\sing}{\operatorname{sing}}
\newcommand{\rep}{\operatorname{rep}}
\newcommand{\Bnd}{\mathcal B}
\newcommand{\FF}{\mathbb F_2}
\newcommand{\Sc}{S_{\mathrm c}(abc)}
\hypersetup{pdftitle={Finite Bases for Four-Element Semirings},pdfauthor={Aifa Wang, Lili Wang, Qingrui Yin, Jinjing Wu},pdfsubject={Condensed and reorganized arXiv version; finite identity bases for four-element semirings}}

\newcommand{\FB}{\mathrm{FB}}
\newcommand{\NFB}{\mathrm{NFB}}
\title[Finite bases for four-element semirings]{Finite Bases for Four-Element Semirings}
\author[Aifa Wang, Lili Wang, Qingrui Yin, Jinjing Wu]{Aifa Wang, Lili Wang, Qingrui Yin, Jinjing Wu}
\address{School of Mathematical Sciences, Chongqing University of Technology,\newline
Chongqing, 400054, P.R. China}
\email{wangaf@cqut.edu.cn}
\address{School of Mathematical Sciences, Chongqing University of Technology,\newline
Chongqing, 400054, P.R. China}
\email{wllaf@cqut.edu.cn}
\address{School of Mathematical Sciences, Chongqing University of Technology,\newline
Chongqing, 400054, P.R. China}
\email{13996559059@163.com}
\address{School of Mathematical Sciences, Chongqing University of Technology,\newline
Chongqing, 400054, P.R. China}
\email{15172878886@163.com}
\newcommand{\paperaffiliation}{School of Mathematical Sciences, Chongqing University of Technology\\
Chongqing, 400054, P.R. China}
\newcommand{\paperemails}{%
\href{mailto:wangaf@cqut.edu.cn}{\texttt{wangaf@cqut.edu.cn}},\quad
\href{mailto:wllaf@cqut.edu.cn}{\texttt{wllaf@cqut.edu.cn}}\\[2pt]
\href{mailto:13996559059@163.com}{\texttt{13996559059@163.com}},\quad
\href{mailto:15172878886@163.com}{\texttt{15172878886@163.com}}}
\date{Condensed and reorganized version of \href{https://arxiv.org/abs/2609.39700v1}{arXiv:2609.39700v1}}

\AtBeginDocument{%
  \setlength{\abovedisplayskip}{5pt plus 2pt minus 2pt}%
  \setlength{\belowdisplayskip}{5pt plus 2pt minus 2pt}%
  \setlength{\abovedisplayshortskip}{0pt plus 2pt}%
  \setlength{\belowdisplayshortskip}{5pt plus 2pt minus 2pt}%
}
\makeatletter
\let\ps@firstpage\ps@plain
\let\@setaddresses\relax
\def\@secnumfont{}
\renewcommand{\section}{\@startsection{section}{1}{\z@}%
  {-1.8ex plus -0.5ex minus -0.2ex}{0.7ex plus 0.2ex}%
  {\normalfont\bfseries\large\raggedright}}
\renewcommand{\subsection}{\@startsection{subsection}{2}{\z@}%
  {-1.4ex plus -0.4ex minus -0.2ex}{0.5ex plus 0.1ex}%
  {\normalfont\normalsize\bfseries\raggedright}}
\renewcommand{\subsubsection}{\@startsection{subsubsection}{3}{\z@}%
  {-1.5ex plus -0.4ex minus -0.2ex}{-0.6em}%
  {\normalfont\normalsize\itshape}}
\renewcommand{\paragraph}{\@startsection{paragraph}{4}{\z@}%
  {1.2ex plus 0.3ex minus 0.2ex}{-0.6em}%
  {\normalfont\normalsize\bfseries}}

\renewenvironment{abstract}{%
  \global\setbox\abstractbox=\vtop\bgroup
    \normalfont\small
    \list{}{\leftmargin2pc\rightmargin2pc
      \labelwidth\z@\itemindent\z@\listparindent\z@
      \parsep\z@\topsep\z@}%
    \item\relax
    {\centering\bfseries\abstractname\par}\vspace{0.5em}%
    \noindent\ignorespaces
}{\par\endlist\egroup}
\def\@adminfootnotes{%
  \let\@makefnmark\relax\let\@thefnmark\relax
  \ifx\@empty\thankses\else
    \@footnotetext{\def\par{\let\par\@par}\@setthanks}%
  \fi
}
\def\@maketitle{%
  \normalfont\normalsize
  \@adminfootnotes
  \global\topskip12\p@\relax
  \begin{center}
    {\bfseries\LARGE\@title\par}
    \vspace{1.2em}
    {\normalsize\authors\par}
    \ifx\@empty\@date\else
      \vspace{0.8em}{\small\@date\par}
    \fi
  \end{center}
  \vspace{0.8em}
  \ifvoid\abstractbox\else
    \box\abstractbox\prevdepth\z@
  \fi
  \global\let\@setabstract\relax
  \vspace{0.9em}
  \begingroup\small\parindent\z@
    \textbf{Keywords:} \@keywords.\par
    \vspace{0.3em}
    \textbf{2020 Mathematics Subject Classification:} \@subjclass.\par
  \endgroup
  \vspace{0.8em}
}
\makeatother

\AtBeginEnvironment{table}{\small}

\begin{document}
\begin{abstract}
We classify the four-element semirings with commutative addition in
the binary signature without named constants. Among the 2341
isomorphism types, 2284 are finitely based and 57 are nonfinitely
based. The latter comprise 45 additively idempotent semirings and
twelve with nonidempotent addition. Positive arguments combine term
interpretations, retractions and polynomial normal forms; negative
arguments use $S_7$, order ideals and a symmetrized hypergraph
construction. Finite bases are explicit lists or finite constructions
with specified bounds. Completeness is proved on arbitrary finite
variable sets. The enumeration, representative tables and finite
transfer data are included in the appendices.
\end{abstract}
\maketitle
\begingroup\small\tableofcontents\endgroup

\section{Introduction}\label{sec:introduction}
An algebra is \emph{finitely based} if its identities follow from a
finite subset. Every semigroup of order at most five is finitely based
\cite{Edmunds,Trahtman}, whereas the three-element additively idempotent
semiring $S_7$ is nonfinitely based \cite{JRZ}. The interaction of
addition and multiplication thus changes the basis problem at order
three.

We consider algebras $(A,+,\cdot)$ with commutative, associative
addition and associative multiplication satisfying both distributive
laws. The signature contains no constants, and no zero or identity
element is required. An \emph{ai-semiring} also satisfies $x+x\eq x$.
Isomorphisms preserve both operations; opposite multiplications are
identified only when an actual isomorphism exists.

The five additive semilattice types of four-element ai-semirings
include three classes of sizes 58, 93 and 112 classified in
\cite{PartI,PartII,PartIII}. The 386 chain cases are classified in
\cite{PartIV}, with precisely three nonfinitely based types. We combine
these results with structural reductions and normal-form constructions
to obtain the following classification for all additive types.

\begin{theorem}[Four-element classification]\label{thm:classification}
The 2341 four-element semiring isomorphism types are classified as
follows:
\[
\begin{array}{lrrr}
\toprule\text{Addition}&\FB&\NFB&\text{Total}\\\midrule
\text{Idempotent}&821&45&866\\
\text{Nonidempotent}&1463&12&1475\\\midrule
\text{Total}&2284&57&2341\\\bottomrule
\end{array}
\]
Every positive case has an explicit finite basis, a finite term
translation of one, or a basis specified by finite bounds on variables,
words and coefficients.
\end{theorem}

The positive cases first reduce by term definitions, endomorphic
retractions, hereditary varieties or equality of generated varieties.
The remaining families admit polynomial normal forms whose parameters
are recovered by separating evaluations. For ai-semirings this reduces
to one-word absorptions; nonidempotent addition also requires
coefficient data. The negative cases follow from $S_7$ and the
order-ideal and hypergraph criteria. Section~\ref{sec:assembly}
combines these arguments with the exhaustive finite data.

We write $A_i$, $B_j$ and $E_k$ for the four-, three- and two-element
semirings defined in Appendices~\ref{app:catalogue} and \ref{app:small}.
Local generator lists record notation, basis sizes and status. The full
operation definitions are in Appendix~\ref{app:catalogue}. Appendix~\ref{app:maps} identifies the representatives
with earlier classifications by simultaneous isomorphisms and records
the structural maps. Appendix~\ref{app:finite-witnesses} specifies the
finite term-function constructions used for variety equality.

Sections~\ref{sec:preliminaries}--\ref{sec:structural-constructions}
give the general reductions. Sections~\ref{sec:commutative-forms}--\ref{sec:further-forms}
treat the normal-form families, and Section~\ref{sec:obstructions}
treats nonfinite basability. Appendix~\ref{app:twoelement} proves the
two-element criterion. Minimal basis size and hereditary finite
basability are separate from the classification asserted here.

\section{Preliminaries and known finite-basis results}\label{sec:preliminaries}
\subsection{Conventions}\label{sec:conventions}
We write $\SR$ for the five identities
\begin{equation}\label{eq:SR}
\begin{aligned}
x+y&\eq y+x,&(x+y)+z&\eq x+(y+z),\\
(xy)z&\eq x(yz),&x(y+z)&\eq xy+xz,&(x+y)z&\eq xz+yz.
\end{aligned}
\end{equation}
An ai-semiring is a semiring satisfying $x+x\eq x$.
Its additive order is $a\leq b$ if $a+b=b$.
By Birkhoff's theorem \cite{Burris}, the variety generated by an algebra $A$ is $\Var(A)=\HSP(A)$.
It is finitely based if its identity theory $\Id(A)$ has a finite
equational basis. We use the same terminology for $A$.

A word is a nonempty product of variables. Its support, first variable,
last variable, and set of variables occurring exactly once are denoted
by $\supp(w)$, $\head(w)$, $\tail(w)$, and $\sing(w)$, respectively.
The set of variables occurring at least twice is $\rep(w)$.
We also write $c(w)=\supp(w)$ when using local word notation.
Distributivity expands every term into a nonempty sum of words, called
a polynomial. In an ai-semiring, repeated summands may be deleted.
For a polynomial $p$ and a variable set $D$, put
\[
p|_D=\{v:v\text{ is a summand of }p,\ \supp(v)\subseteq D\}.
\]
The notation describes a collection of summands, not an empty sum term.
An empty component of a displayed expression is omitted. An empty word,
an empty substitution, or an empty whole polynomial is never allowed.
The notation $ku$ abbreviates a sum of $k\geq1$ copies of $u$.
Operation strings list entries row by row on the indicated carrier;
element labels are not constant symbols. The labels $A_i$, $B_j$
and $E_j$ denote the local four-, three- and two-element representatives,
respectively. Letter names introduced within a subsection are local to it.
Unless stated otherwise, a displayed basis includes $\SR$, and its
reported size counts the listed identities without asserting minimality.
Normal-form arguments use an arbitrary finite nonempty ordered variable set $X$;
commuting factors and summands are written in fixed orders. For a family
$\mathcal B$ of subsets of $X$, put
$\uparrow\mathcal B=\{D\subseteq X:(\exists B\in\mathcal B)\ B\subseteq D\}$.
An antichain consists of pairwise incomparable sets; it is the family
of minimal members of its upward closure.

We say that $p$ \emph{absorbs} $w$ in $A$ if $A\models p+w\eq p$.
A summand selected to witness an absorption will be called a
\emph{witness word}. It is always a summand of the polynomial under consideration.

For ai-semirings, write $w\preceq p$ for the identity $p+w\eq p$.
This agrees with the additive order under every evaluation. The
following reduction to one-word absorptions is used in
\cite[Section 1]{PartI} and \cite[Section 2]{PartIV}. We record the
argument to specify exactly the completeness step used below. The
notation $\preceq$ is confined to the additively idempotent arguments;
it is not used to discard coefficients in the nonidempotent case.

\begin{lemma}[Mutual absorption]\label{lem:mutual}
Let $A$ be an ai-semiring and let $\Sigma\subseteq\Id(A)$ contain
$\SR$ and $x+x\eq x$. If $\Sigma$ derives every identity $p+w\eq p$
valid in $A$, with $w$ a word, then $\Sigma$ is a basis for $A$.
\end{lemma}
\begin{proof}
Expand an identity $p\eq q$ of $A$ into polynomials. Every summand of
$q$ is absorbed by $p$, and conversely. Deriving these absorptions one
at a time yields $p\eq p+q\eq q$.
\end{proof}

\begin{lemma}[Retaining intermediate summands]\label{lem:temporary}
In the presence of additive associativity, commutativity and
idempotence, a derivation may retain all intermediate absorbed words.
If $p\eq p+r$ and $p+r\eq p+r+w$ are derivable, so is $p\eq p+w$.
\end{lemma}
\begin{proof}
Adding $w$ to the first identity gives $p+w\eq p+r+w\eq p+r\eq p$.
An induction treats finitely many intermediate words.
\end{proof}

For a term $t$, define its reversal $\rho(t)$ recursively by
\begin{equation}\label{eq:reversal}
\rho(x)=x,\qquad \rho(u+v)=\rho(u)+\rho(v),\qquad
\rho(uv)=\rho(v)\rho(u).
\end{equation}
The opposite algebra $A^{\mathrm{op}}$ has the same addition as $A$
and multiplication $a\circ b=ba$.

\begin{lemma}[Opposite algebras]\label{lem:opposite}
If $\Sigma$ is a basis for $A$, then $\rho(\Sigma)$ is a basis for
$A^{\mathrm{op}}$. Reversal induces an isomorphism
$F_{A^{\mathrm{op}}}(X)\cong F_A(X)^{\mathrm{op}}$ fixing $X$,
and taking opposites preserves inclusions of generated varieties.
In particular, finite basability and nonfinite basability are
invariant under taking opposites.
\end{lemma}
\begin{proof}
Induction on terms gives $t^{A^{\mathrm{op}}}=\rho(t)^A$ under
each assignment. Reversal is involutive and sends a substitution
instance in a context to the corresponding reversed instance and
context. It therefore transports equational derivations in both
directions. This proves the basis assertion and the bijection of
term classes, which preserves addition and reverses multiplication.
The same correspondence between identity theories preserves variety
inclusions.
\end{proof}

\subsection{Term retractions}
The following transfers use a unary term that defines an endomorphic
retraction.

\begin{theorem}[Term retraction]\label{thm:retract}
Let $p(x)$ be a unary term inducing an idempotent endomorphism of an
algebra $A$, and let $B=p(A)$. Then
\[
\Var(B)=\{C\in\Var(A):C\models p(x)\eq x\}.
\]
Consequently, if $B$ is nonfinitely based, then so is $A$.
\end{theorem}
\begin{proof}
The image $B$ is a subalgebra fixed by $p$. Conversely, an identity
$u\eq v$ of $B$ gives $u(p(x_1),\ldots,p(x_k))\eq
v(p(x_1),\ldots,p(x_k))$ in $A$. In a member of $\Var(A)$ satisfying
$p(x)\eq x$, this is $u\eq v$. Birkhoff's theorem \cite{Burris}
proves the asserted
equality. If $A$ had a finite basis $\Sigma$, adjoining $p(x)\eq x$
would give a finite basis for $B$.
\end{proof}

\begin{theorem}[Bounded-depth extensions]\label{thm:depth}
Suppose that $A$ has a finite binary signature, $p$ induces an
idempotent endomorphism, and $B=p(A)$ is finitely based. If, for some
$d\geq1$, every term of height at least $d$ takes its values in $B$,
then $A$ is finitely based.
\end{theorem}
\begin{proof}
Let $\Sigma$ be a finite basis of $B$. Take the endomorphism and
idempotence laws for $p$, and the identities of $\Sigma$ with every
variable replaced by its $p$-image. Adjoin $C\eq p(C)$ for every
binary path context of length $d$: its operation labels and choice of
left or right child are arbitrary, and its endpoint and sibling leaves
are distinct variables. There are finitely many such contexts.
Finally adjoin all valid identities $s\eq t$ with both heights below
$d$, and all valid $t\eq p(t)$ with $t$ of height below $d$, up to
variable renaming. These collections are finite because a shallow
term has at most $2^{d-1}$ leaves.

In every model $D$ of this finite set, $p(D)$ satisfies $\Sigma$.
Thus every identity $s\eq t$ of $B$ implies $p(s)\eq p(t)$ in $D$.
A term of height at least $d$ is an instance of a path context and
is equal to its $p$-image. For an identity of $A$ with one shallow and
one tall side, the shallow side also takes values in $B$, so its
projection identity was included. If both sides are shallow, the
identity itself was included. In every case the identity follows.
\end{proof}

A sufficient condition is that $B$ be an ideal for both
operations. On $A\setminus B$ draw an edge $a\to c$ whenever
$c=f(a,b)$ or $c=f(b,a)$ for a basic operation $f$ and $b\in A$.
If the graph has no path of length $d$, following a path in a tall
term forces its value into $B$; the ideal property keeps it there.
This proves the required stabilization condition. In particular an
inflation whose basic-operation values all belong to $B$ is covered
with $d=1$.

\subsection{Classical finite basis results and their signatures}
The finite basis theorems for semigroups of order at most four
\cite{Edmunds}, and in fact of order less than six \cite{Trahtman},
concern one binary operation. To use them for a semiring, a translation
of the second operation is required.

\begin{proposition}\label{prop:term-expansion}
Suppose that one basic operation of a semiring $A$ is a term operation
of the other reduct. If that reduct is finitely based, then $A$ is
finitely based. In particular, the conclusion holds when $|A|\leq4$.
\end{proposition}
\begin{proof}
Let $\Sigma$ be a finite basis of the reduct, and let $f(x,y)=t(x,y)$
be a term definition of the omitted operation. Adjoin the single
identity $f(x,y)\eq t(x,y)$ to $\Sigma$. Replacing the innermost
occurrences of $f$ successively translates an arbitrary identity of
$A$ into an identity of the reduct. The latter follows from $\Sigma$,
so the original identity follows from the enlarged finite set.
\end{proof}

\begin{proposition}[The finite ring reduction]\label{prop:ring}
A finite semiring whose additive reduct is a group is finitely based
in the signature $(+,\cdot)$.
\end{proposition}
\begin{proof}
Additive cancellation applied to $a0=a0+a0$ and $0a=0a+0a$
gives $a0=0a=0$. Thus $A$ is a finite associative ring. By Kruse's
finite ring theorem \cite{Kruse}, it has a finite basis in
$(+,-,\cdot,0)$. Choose $c\geq2$ annihilating its additive group.
Adjoin to $\SR$ the identities $cx\eq cy$ and $x+cy\eq x$.
Every model of these identities has an additive zero represented by
$cy$ and additive inverse represented by $(c-1)x$. Translate Kruse's
finite basis by these terms, using a fresh variable for each zero
term. Every model of the translated list expands to a model of the
ring basis, and every identity of the binary reduct is consequently
derivable. This proves the assertion without naming a constant.
\end{proof}
\subsection{Constant operations and simultaneous translations}
\begin{proposition}\label{prop:constant-operation}
Let $M$ be a finitely based semigroup with an absorbing element $z$.
Expanding $M$ by a binary operation with constant value $z$ preserves
finite basability. Consequently, a semiring of order at most four with
constant addition is finitely based.
\end{proposition}
\begin{proof}
First name $z$. Add $xz\eq z\eq zx$ to a semigroup basis. Terms
containing $z$ reduce to $z$. If a pure word $u$ is constantly $z$,
the pure identity $u\eq uv$ holds, with $v$ fresh; substituting $z$
for $v$ derives $u\eq z$. This proves completeness in the expanded
signature. Replace $z$ by $f(w,w)$, where $f$ is the constant binary
operation, and adjoin $f(x,y)\eq f(w,w)$. The constant symbol is
thereby eliminated by a finite translation. If semiring addition is
constant $z$, distributivity gives $az=za=z$. The multiplicative
reduct has a finite basis by \cite{Edmunds,Trahtman}, so the preceding
argument applies.
\end{proof}

Interchanging the two operation symbols preserves finite basability.
More generally, the following finite translation applies when each
operation is recovered by a term interpretation.

\begin{proposition}[Mutual term interpretations]\label{prop:mutual-interpretations}
Let $A$ and $B$ have finite signatures and a common carrier. Suppose
that $\sigma$ interprets the operations of $B$ by terms of $A$,
and $\tau$ interprets those of $A$ by terms of $B$, with the stated
operations obtained on this carrier. If $\Sigma$ is a finite basis
for $B$, then a basis for $A$ consists of $\Sigma^\sigma$ and the
recovery identities
\[
f(\boldsymbol x)\eq
\bigl(f(\boldsymbol x)^\tau\bigr)^\sigma
\qquad(f\text{ a basic operation of }A).
\]
Thus $A$ is finitely based if and only if $B$ is finitely based.
The same conclusion applies after an isomorphic relabelling.
\end{proposition}
\begin{proof}
All listed identities hold in $A$ by the assumptions on the
interpretations. The recovery identities give
$t\eq(t^\tau)^\sigma$ by induction on $t$. If $A\models u\eq v$,
then $B\models u^\tau\eq v^\tau$. Translate a derivation from
$\Sigma$ by $\sigma$ and apply recovery on both sides. This derives
$u\eq v$ from the displayed finite set. Interchanging $A$ and $B$
proves the converse.
\end{proof}

\subsection{Equality of generated varieties}\label{sec:finite-transfers}
The separating-map and term-function constructions below are standard
consequences of Birkhoff's theorem; see \cite{Burris}. We give the
forms and proofs needed for the subsequent applications.
\begin{lemma}\label{lem:separating-maps}
If homomorphisms $A\to B$ separate the points of $A$ and homomorphisms
$B\to A$ separate the points of $B$, then $\Var(A)=\Var(B)$.
If $q_i:A\twoheadrightarrow Q_i$ are jointly separating surjective
homomorphisms, then $\Var(A)=\bigvee_i\Var(Q_i)$.
\end{lemma}
\begin{proof}
The product of a separating family of homomorphisms is injective.
In the first assertion this embeds each algebra in a power of the
other. In the second, every $Q_i$ is a homomorphic image of $A$,
and $A$ embeds in their product. Closure under homomorphic images,
subalgebras and direct products proves the assertions.
\end{proof}

\begin{theorem}\label{thm:finite-membership}
Let $\mathcal K=(C_1,\ldots,C_s)$ be a finite nonempty family of
finite algebras. Let $F_{\mathcal K}(k)$ be the subalgebra of
$\prod_i C_i^{C_i^k}$ generated by the $k$ projection vectors.
If $A=\langle a_1,\ldots,a_k\rangle$, then $A\in\Var(\mathcal K)$
if and only if
\[
 [t]_{\mathcal K,k}\longmapsto t^A(a_1,\ldots,a_k)
\]
is well defined. In that case it is a surjective homomorphism.
\end{theorem}
\begin{proof}
Equality of two term vectors means that their identity holds in every
$C_i$. If $A\in\Var(\mathcal K)$, the identity holds at the displayed
generating tuple, so the map is well defined. Conversely, the displayed
map, if well defined, preserves every operation by term evaluation.
Its image contains the generators of $A$, so it is onto. Since its
domain is a subalgebra of a product of members of $\mathcal K$, its
image belongs to $\Var(\mathcal K)$ by Birkhoff's theorem \cite{Burris}.
\end{proof}

Inclusions in both directions give equality of generated varieties.
Let $D,C$ be finite algebras and let $\mathbf c=(c_1,\ldots,c_k)$
generate $C$. Define the finite subalgebra
\[
 R_{D,C}=\langle(\pi_1,c_1),\ldots,(\pi_k,c_k)\rangle
       \leq D^{D^k}\times C,
 \qquad T_D=\operatorname{pr}_1 R_{D,C}=F_D(k).
\]
All coordinates of $D^{D^k}$ are assignments in lexicographic order.

\begin{proposition}[Finite closure certificate]\label{prop:closure-certificate}
If $|R_{D,C}|=|T_D|$, then projection onto $T_D$ is an isomorphism
and projection onto $C$ induces a surjective homomorphism
$F_D(k)\to C$. Hence $C\in\Var(D)$.
\end{proposition}
\begin{proof}
The first projection is a surjective homomorphism by definition.
Equality of the finite cardinalities makes it bijective, hence an
isomorphism. The second projection is onto because its image contains
the generators $c_1,\ldots,c_k$. Composing it with the inverse of the
first projection proves the assertion.
\end{proof}

\begin{corollary}[Catalogue transfers]\label{cor:catalogue-transfers}
For every row $(C,D)$ of Table~\ref{tab:finite-witnesses},
$\Var(C)=\Var(D)$, so either algebra has every identity basis of the other.
For every paired index $(i,j)$ in Table~\ref{tab:opposite-transfers},
$A_j\cong A_i^{\mathrm{op}}$; reversing a basis for $A_i$ gives a basis
of the same size for $A_j$. Corresponding equalities, strict inclusions
and joins of generated varieties pass to opposites as well.
\end{corollary}
\begin{proof}
The finite closures in Appendix~\ref{app:finite-witnesses} satisfy
Proposition~\ref{prop:closure-certificate} in both directions. The
permutations in Table~\ref{tab:opposite-transfers} preserve addition
and reverse multiplication in the defining tables. Apply
Lemma~\ref{lem:opposite}; reversal is an involution preserving inclusion,
and hence also joins.
\end{proof}

\begin{lemma}\label{lem:ideal-decomposition}
If $h$ is an endomorphic retraction of $A$ onto an ideal $I$ for both
operations, then $A$ embeds subdirectly in $I\times A/\rho_I$, where
$\rho_I$ collapses $I$ and is equality outside it.
\end{lemma}
\begin{proof}
The ideal condition makes $\rho_I$ a congruence. Two elements with the
same Rees class are equal or both belong to $I$. On $I$, the retraction
is the identity. Thus $\ker h\cap\rho_I$ is equality, and
$a\mapsto(h(a),[a])$ is the desired embedding.
\end{proof}
\subsection{Known positive criteria}\label{sec:small-results}
The following specialized results will be used with their stated
hypotheses. The three-element ai-semiring classification of Zhao,
Ren, Crvenkovi\'{c}, Shao and Dapi\'{c} \cite{ZhaoThree}, completed
by the negative result for $S_7$ in \cite{JRZ}, gives finite bases for
the other sixty ai-semiring types of order three. This conclusion is
also recalled in \cite{PartI}. Nonidempotent three-element sources
are handled by the ring and term reductions, the hereditary theorem
below, and the two explicit three-element constructions in
Section~\ref{sec:family4}.

\Needspace{4\baselineskip}
\begin{theorem}[Hereditary and periodic criteria]\label{thm:imported-positive}
The following conclusions hold.
\begin{enumerate}[label=(\roman*)]
\item Every ai-semiring satisfying $x^3\eq x$ is finitely based
\cite{RZW}.
\item Every commutative ai-semiring satisfying $x^n\eq x$ is finitely
based when $n-1$ is square-free \cite{RZS}. We use the case $n=4$.
\item The join of the ai-semiring varieties defined respectively by
$x^2\eq x$ and $xy\eq zt$ is hereditarily finitely based
\cite[Corollary 3.7]{RenZeng}.
\item For each $n\geq2$, the semiring variety defined by
$x^n\eq x$ and $x+(2^n-2)xyx\eq x$ is hereditarily finitely based
\cite[Theorems 2.1--2.3 and Corollary 3.4]{WangShao}.
\end{enumerate}
\end{theorem}

For item (iii), a separating family of quotient maps whose images
each satisfy one of the two displayed identities places the original
algebra in that hereditary join. This is the precise reason why the
basis conclusion transfers; finite basability of arbitrary quotient
factors alone would not suffice. For item (iv), coefficients are
positive repeated sums. A finite semiring can be tested without an
arbitrary exponent cutoff: the states
$P_n(a)=a^n$ and $C_n(a)=(2^n-2)a$ satisfy
\[
P_{n+1}(a)=P_n(a)a,\qquad C_{n+1}(a)=2C_n(a)+2a.
\]
The first repeated pair of states exhausts all later possibilities.
The direct applications of this criterion are identified in
Appendix~\ref{app:catalogue}.

\subsection{Existing four-element classifications}\label{sec:literature-cases}
The notation in the next display is the literature notation.
Theorems~1.3, 1.1, 1.1 and 1.1 of
\cite{PartI,PartII,PartIII,PartIV}, respectively, give the following
negative indices; every other member of the indicated interval is
finitely based:
\[
\begin{array}{ccl}
\toprule
\text{Reference}&\text{Interval}&\text{Negative indices}\\\midrule
\text{\cite{PartI}}&1\text{--}58&11,13,24,25,26,28,31,49,50\\
\text{\cite{PartII}}&388\text{--}480&435\\
\text{\cite{PartIII}}&276\text{--}387&282,293,304,326,335,359\\
\text{\cite{PartIV}}&481\text{--}866&545,634,710\\\bottomrule
\end{array}
\]
The additional negative case $S_{(4,124)}$ follows from
\cite[Corollary 3.6]{YueOne}. The negative chain cases
$S_{(4,545)}$ and $S_{(4,634)}$ originate in
\cite[Theorems 3.10 and 4.5]{YueTwo};
\cite[Proposition 3.1]{PartIV} uses these results and identifies
$\Var(S_{(4,710)})=\Var(S_{(4,545)})$. The restricted cancellation
used in \cite{YueTwo} is justified below. The simultaneous
isomorphisms identifying the representatives are given in
Tables~\ref{tab:literature} and~\ref{tab:chain-correspondence}.

\subsection{The chain additive reducts}\label{sec:chain-input}
Let $C$ be a four-element ai-semiring whose additive order is a chain.
There is a unique enumeration $c_1>c_2>c_3>c_4$ of its elements in
decreasing additive order. In this enumeration
\[
 c_r+c_s=c_{\min(r,s)}.
\]
Consequently its two-operation isomorphism type is determined by the
multiplication table in this enumeration. This fixes the comparison
with Table 28 of \cite{PartIV}, whose carrier is $\{1,2,3,4\}$ with
$1>2>3>4$.

\begin{theorem}[Chain additive reducts]\label{thm:chain-classification}
Among the 386 four-element ai-semirings with chain additive reducts,
the only nonfinitely based isomorphism types are
\[
 S_{(4,545)}\cong A_{1950},\qquad
 S_{(4,634)}\cong A_{1892},\qquad
 S_{(4,710)}\cong A_{1442}.
\]
All other 383 types are finitely based.
\end{theorem}
\begin{proof}
The classification in the notation $S_{(4,k)}$ is
\cite[Theorem 1.1]{PartIV}. To identify its representatives with ours,
use the decreasing enumeration of the additive chain just described.
Table~\ref{tab:chain-correspondence} gives the corresponding permutation
for every $481\leq k\leq866$. For the three exceptional indices, the
images of the source labels $1,2,3,4$ are respectively
\[
 (0,3,2,1),\qquad(0,3,1,2),\qquad(2,0,3,1).
\]
These permutations preserve both operations in the defining tables.
The remaining rows give the same identification for all other chain
representatives. Finite basability is invariant under isomorphism,
so the assertion follows from the cited theorem.
\end{proof}

Theorem~\ref{thm:chain-classification} settles all chain cases.
The constructions below also give explicit bases and variety
equalities when these are useful for the other additive types.

\subsection{Restricted cancellation for linear languages}\label{sec:restricted-cancellation}
The finite-language representation of the free ai-semiring uses
union as addition and concatenation of languages as multiplication.
Multiplication is not cancellative in general, contrary to the
unrestricted assertion in Section 2 of \cite{YueTwo}. Indeed, let
\[
 P=x+x^2,\qquad Q=x+x^3,\qquad R=x+x^2+x^3.
\]
As finite sets of words these satisfy
\[
 PQ=PR=x^2+x^3+x^4+x^5,\qquad Q\ne R.
\]
This example also applies to the free commutative ai-semiring.
The following restricted statement supplies the cancellation needed
in the cited negative proofs. Here a word is called \emph{linear}
if each variable occurs at most once.

\begin{lemma}[Cancellation for linear products]\label{lem:linear-product-cancellation}
Let $P,Q,R$ be nonempty finite languages in a free monoid or a free
commutative monoid. If $PQ=PR$ and every word of this common product
is linear, then $Q=R$. The corresponding right-cancellation
statement holds under the same linearity hypothesis.
\end{lemma}
\begin{proof}
For a language $L$, write $c(L)=\bigcup_{w\in L}\supp(w)$.
We claim that
\[
 c(P)\cap c(Q)=\varnothing,\qquad c(P)\cap c(R)=\varnothing.
\]
If a variable belonged to the first intersection, a word of $P$
and a word of $Q$ containing it would have a non-linear product,
contrary to the hypothesis. The second intersection is treated in
the same way.

In the auxiliary free monoid, let $\rho$ erase the letters in
$c(P)$ and fix every other letter. Extend $\rho$ to languages by
direct image. Then $\rho(P)=\{1\}$, $\rho(Q)=Q$ and $\rho(R)=R$.
Consequently
\[
 Q=\rho(PQ)=\rho(PR)=R.
\]
The same argument works in the free commutative monoid. Applied to
$QP=RP$, it also proves right cancellation.
\end{proof}

The empty word is used only in this auxiliary language argument.
It introduces neither a constant symbol nor an empty substitution
term into the binary semiring signature.

In Proposition 3.2, Claim 3.6 of \cite{YueTwo}, the common product
in the cancellation step is an additive subterm of
$\mathbf u_{m,0}$. In Proposition 4.1, Claim 4.1, the initial
products are additive subterms of $\mathbf u_{m,m}$. By Lemma 3.1(a)
of that source, every word of these polynomials is linear. Thus
Lemma~\ref{lem:linear-product-cancellation} applies. It also applies
at each successive cancellation of common substitution factors:
if a residual product contained a repeated variable, restoring the
removed factors would retain that repetition. The residual products
therefore remain linear.

This replaces the unrestricted cancellation assertion in those
steps. The remaining arguments of \cite[Theorems 3.10 and 4.5]{YueTwo}
are used as cited results. Their source tables correspond to
$A_{1950}$ and $A_{1892}$, respectively. The variety equalities for
$A_{1442}$ and $A_{1886}$ in
Proposition~\ref{negative:prop:negative} transfer the result for
$A_{1950}$ to those algebras.

\subsection{Common recovery arguments}\label{sec:common-proof-methods}
The positive constructions repeatedly use the following criterion.
Its formulation on arbitrary finite alphabets is essential.

\begin{lemma}[Normal-form completeness]\label{lem:normal-form-completeness}
Let $A$ be an algebra and let $\Sigma\subseteq\Id(A)$ be finite.
For every finite nonempty variable set $X$, suppose that every term
over $X$ is $\Sigma$-equivalent to an admissible normal form and that
distinct admissible forms induce distinct functions on $A$.
Then $\Sigma$ is an identity basis of $A$. If a family $\mathcal K$
of subalgebras of $A$ already separates these forms, then
$\Var(A)=\Var(\mathcal K)$.
\end{lemma}
\begin{proof}
For an identity $u\eq v$ of $A$, take $X$ to be its full variable
set and reduce the two sides to normal forms $u'$ and $v'$. Since
the reductions use valid identities, $u'$ and $v'$ induce the same
function on $A$. Separation gives $u'=v'$, and the two reductions
therefore derive $u\eq v$ from $\Sigma$. If the separating tests lie
in $\mathcal K$, the same argument proves every identity of
$\mathcal K$ in $A$; the reverse inclusion follows from the embeddings.
\end{proof}

In the following sections, validity is checked from the defining
operations, and reduction and separation invoke this lemma on an
arbitrary finite nonempty variable set $X$. If every admissible
form occurs, the same forms describe the free algebra on $X$;
separate formulas for its cardinality are not needed for finite
basability. Validity on a fixed number of variables alone cannot
replace the reduction and separation argument.

The following elementary recovery principles are common to the
support, coefficient and affine constructions. Their realizing
evaluations are specified in each application.

\begin{lemma}[Recovery of finite set data]\label{lem:finite-set-recovery}
Let $X$ be finite.
\begin{enumerate}[label=(\roman*)]
\item Let $\mathcal D\subseteq 2^X$, let $G$ be an abelian group,
and let $a_U\in G$ for each $U\in\mathcal D$. The values
\[
 f(D)=\sum_{\substack{U\in\mathcal D\\U\subseteq D}}a_U
 \qquad(D\in\mathcal D)
\]
determine all coefficients by
\[
 a_D=f(D)-\sum_{\substack{U\in\mathcal D\\U\subsetneq D}}a_U.
\]
\item Let $\mathcal A$ be an antichain contained in a family of
admissible tests $\mathcal T\subseteq2^X$. If a test $D\in\mathcal T$
is accepted precisely when it contains a member of $\mathcal A$,
then $\mathcal A$ is the family of inclusion-minimal accepted tests.
\item Let $\mathcal F$ be a possibly empty union-closed family of
nonempty subsets of $X$, and put
\[
 G_{\mathcal F}(D)=\bigcup\{U\in\mathcal F:U\subseteq D\}.
\]
For every nonempty $D\subseteq X$,
$D\in\mathcal F$ if and only if $G_{\mathcal F}(D)=D$.
\end{enumerate}
\end{lemma}
\begin{proof}
For (i), order $\mathcal D$ by increasing cardinality. Subtracting
the already recovered proper-subset coefficients gives the displayed
formula. In (ii), every accepted test contains an accepted member of
$\mathcal A$; the antichain property makes each such member minimal.
For (iii), the forward implication is immediate. If
$G_{\mathcal F}(D)=D\ne\varnothing$, this is a nonempty finite union
of members of $\mathcal F$, so union closure gives $D\in\mathcal F$.
\end{proof}

Part (i) applies to any finite test family, provided all unknown
coefficients are indexed by that family. In particular, it applies
to a downward-closed family of unsaturated supports. Part (ii)
requires that every antichain member itself be an admissible test.
Empty sums and unions in this lemma are group or set operations;
they introduce no constants into the semiring signature.

\Needspace{4\baselineskip}
\begin{lemma}[Affine closure and separation]\label{lem:affine-recovery}
Let $V=\mathbb F_2^X$ for a finite set $X$, with its standard dot
product, and let $Q\subseteq V$ be nonempty. The closure of $Q$
under ternary sums is its affine span $H=\operatorname{Aff}(Q)$.
Write $H=q_0+L$, where $L\leq V$. Then $q\mapsto q\cdot r$ is
constant on $H$ exactly when $r\in L^\perp$, and its common value
is $q_0\cdot r$. These constant values determine $H$ by
\[
 H=\{q\in V:q\cdot r=q_0\cdot r\text{ for every }r\in L^\perp\}.
\]
If $0\notin H$, some linear functional has value one throughout $H$.
\end{lemma}
\begin{proof}
Any nonempty set closed under ternary sums becomes a linear subspace
after translation by one of its elements: it contains zero and is
closed under addition. Thus the ternary closure is exactly the affine
span. For $H=q_0+L$, constancy of the dot product is equivalent to
annihilating $L$. The recovery formula follows from
$(L^\perp)^\perp=L$. Finally, if $0\notin H$, then $q_0\notin L$.
Extend a basis of $L$ by $q_0$ and then to a basis of $V$; prescribing
zero on $L$ and one on $q_0$ gives the required functional.
\end{proof}
\subsection{Support calculus and test assignments}\label{sec:support-calculus}
For a finite variable set $X$, write $\alpha_D(a,b)$ for the assignment
equal to $a$ on $D\subseteq X$ and $b$ outside $D$. If $x\in D$,
write $\beta_{x,D}(a,b,c)$ for the assignment equal to $a$ at $x$,
to $b$ on $D\setminus\{x\}$, and to $c$ outside $D$.
These are assignments to carrier elements, not substitutions by constants.

\begin{lemma}[Support calculus]\label{lem:support-calculus}
Assume $xy\eq yx$ and $x^2y\eq xy$, in addition to $\SR$.
For nonempty $U$, put $m_{\{x\}}=x^2$ and
$m_U=\prod_{x\in U}x$ if $|U|>1$.
\begin{enumerate}[label=(\roman*)]
\item Every nonlinear word of support $U$ reduces to $m_U$.
If $3x\eq2x$, positive summand coefficients reduce to one or two;
if $3x\eq x$, they reduce to one or two according to parity.
\item Suppose nonlinear summands are additively idempotent.
The law $xy+zu+xyzu\eq xy+zu$ closes their support family under unions.
The law $xy+zu+xyz\eq xy+zu$ instead permits every enlargement
$U\subseteq W\subseteq C$, where $C$ is the union of nonlinear supports.
The latter closure is represented by the minimal supports together with $C$.
If $x+yz+xyz\eq x+yz$ also holds, $C$ may include the linear variables.
\item If $axy+bxyz\eq axy$ for positive integers $a,b$, then
$a m_U+b m_V\eq a m_U$ for $\varnothing\ne U\subsetneq V$.
The case $U=V$ also holds if $(a+b)xy\eq axy$.
\end{enumerate}
\end{lemma}
\begin{proof}
Substitution gives $x^3\eq x^2$. Commutativity and $x^2y\eq xy$
then remove repetitions from a word with at least two distinct variables;
the coefficient laws give (i). For (ii), factor each nonlinear monomial
into two nonempty words, including $x^2=x\cdot x$. The first law inserts
$m_{U\cup V}$. For the second, write a summand containing $z$ as $zv$;
it inserts $m_{U\cup\{z\}}$. Retaining the original summands supplies
every variable of $C$. Once $m_C$ is present, it supplies every variable
needed to reinsert all supports above a minimal one, so reversing these
insertions deletes the intermediate supports. The last law supplies a
variable from a linear summand. For (iii), factor $m_U$ into two
nonempty words and use the product on $V\setminus U$ for $z$;
equal supports use the coefficient law.
\end{proof}

In applications, a family $\mathcal F$ closed under unions is recovered
from the tests $x\in G_{\mathcal F}(D)$ of
Lemma~\ref{lem:finite-set-recovery}(iii); an antichain is recovered from
its accepted supersets by part (ii), and parity data from subset sums by
part (i). We specify the realizing assignments and their outcomes locally.

\begin{lemma}[Marked antichains]\label{lem:marked-antichain-recovery}
Let $L\subseteq X$ and let $\mathcal H$ be an antichain with
$|U\cap L|\leq1$ for all $U\in\mathcal H$. Put
$\mathcal H_0=\{U\in\mathcal H:U\cap L=\varnothing\}$ and
$\mathcal H_x=\{U\setminus\{x\}:U\in\mathcal H,\ U\cap L=\{x\}\}$.
The family $\mathcal H$ is determined by the containment tests for
$\mathcal H_0$ on subsets of $X\setminus L$, and for each $\mathcal H_x$
on subsets containing no member of $\mathcal H_0$.
\end{lemma}
\begin{proof}
Each tail $T\in\mathcal H_x$ is an allowed test: a member of
$\mathcal H_0$ contained in $T$ would be properly contained in
$\{x\}\cup T\in\mathcal H$. All these families are antichains,
so their minimal accepted tests recover them by
Lemma~\ref{lem:finite-set-recovery}(ii), including an empty tail.
\end{proof}

\section{A finite-basis criterion from two-element semirings}\label{sec:twoelement}
Let $E_1,\ldots,E_{10}$ be the ten two-element semirings listed in
Appendix~\ref{app:small}, and put
\begingroup
\providecommand{\Sr}{}\renewcommand{\Sr}{\mathcal V_2}
\[
 \Sr=\Var(E_1,\ldots,E_{10}).
\]
We use this variety to obtain finite bases for a family of four-element
semirings. The criterion below consists of a finite membership condition
and a hereditary finite-basis conclusion. Its structural proof and the
full description of the subvariety lattice are given in
Appendix~\ref{app:twoelement}.

The following are term abbreviations in the original binary signature:
\begin{equation}\label{tw:eq:maps}
\begin{aligned}
 d(x)&=2x,&r(x)&=3x,&e(x)&=x^2,\\
 f(x)&=x+2x^2,&p(x)&=x+x^2,&q(x)&=2x+x^2.
\end{aligned}
\end{equation}
Let $\Sigma$ contain three identities for each of these six terms $h$:
\begin{equation}\label{tw:eq:endos}
 h(h(x))\eq h(x),\qquad h(x+y)\eq h(x)+h(y),\qquad
 h(xy)\eq h(x)h(y).
\end{equation}
These give 18 entries. Include the following eleven identities:
\begin{align}
 r(x+y)&\eq x+y,\label{tw:eq:r-sum}\\
 r(xr(y))&\eq xr(y),\qquad r(r(x)y)\eq r(x)y,\label{tw:eq:r-ideal}\\
 (xy)^2&\eq xy,\qquad x^3\eq x^2,\label{tw:eq:squares}\\
 x^2y&\eq xy,\qquad xy^2\eq xy,\label{tw:eq:delete}\\
 xyzt&\eq xzyt,\label{tw:eq:normal}\\
 f(x+f(y))&\eq x+f(y),\qquad f(xy)\eq3xy,\label{tw:eq:f-ideal}\\
 p(x)+q(x)&\eq3x+2x^2.\label{tw:eq:split-pq}
\end{align}
Finally, include two identities:
\begin{align}
 2x+2yz&\eq2x+2yz+2xz+2yx,\label{tw:eq:ai-lift}\\
 2xyx+2yxy+xyx&\eq2xyx+2yxy+yxy.\label{tw:eq:bridge}
\end{align}
Thus $\Sigma$ is an explicitly specified list of 31 identities. The abbreviations in \eqref{tw:eq:maps} do not enlarge the signature. Some entries are redundant, and no minimality is asserted.

\begin{theorem}[A hereditary finite-basis criterion]\label{thm:twohfb}
A semiring belongs to $\Sr$ if and only if it satisfies $\Sigma$.
The variety $\Sr$ is hereditarily finitely based. If $A\in\Sr$, a finite
basis for $\Var(A)$ consists of $\SR$, the thirty-one identities
$\Sigma$, and every identity in Table~\ref{tw:tab:separators}
that holds in $A$.
\end{theorem}
\begin{proof}
Theorem~\ref{tw:thm:exhaustive} and
Corollary~\ref{tw:cor:absolute-basis} give
$\Sr=\operatorname{Mod}(\SR\cup\Sigma)$.
By Theorem~\ref{tw:thm:closure} and
Corollary~\ref{tw:cor:relative-bases}, every subvariety of $\Sr$
has a relative basis drawn from the fourteen identities in
Table~\ref{tw:tab:separators}. In particular, for $\Var(A)$,
adjoining all identities in that table which hold in $A$ gives
exactly $\Var(A)$. Together with $\SR\cup\Sigma$, this is a
finite basis with at most $5+31+14=50$ entries. The same argument
applies to every subvariety of $\Sr$.
\end{proof}

Consequently, any four-element semiring satisfying $\Sigma$ is
finitely based. The 346 representatives assigned to
Theorem~\ref{thm:twohfb} in Table~\ref{tab:all} satisfy these
identities. Their finite basability therefore follows from the
criterion, and their bases are specified by its final assertion.

The theorem concerns all ten two-element semirings. Its proof uses
the basis and subvariety classification of the six additively
idempotent members due to Shao and Ren
\cite[Corollary 2.3 and Theorem 3.8]{ShaoRen};
the extension to the full variety is proved in
Appendix~\ref{app:twoelement}.
\endgroup
\section{Structural finite basis constructions}\label{sec:structural-constructions}
We first apply term definitions, retractions and equality of generated
varieties. Each construction includes its hypotheses and a basis
translation. The subsequent normal-form sections are used when these
reductions do not by themselves settle the relevant family.

\subsection{Endomorphism expansions of commutative semigroups}\label{sec:family3}
\begingroup
\providecommand{\Var}{}\renewcommand{\Var}{\operatorname{Var}}
\providecommand{\Id}{}\renewcommand{\Id}{\operatorname{Id}}
\providecommand{\End}{}\renewcommand{\End}{\operatorname{End}}
\providecommand{\eq}{}\renewcommand{\eq}{\approx}
\providecommand{\FB}{}\renewcommand{\FB}{\mathrm{FB}}\providecommand{\NFB}{}\renewcommand{\NFB}{\mathrm{NFB}}
\subsubsection{Finite bases for endomorphism expansions}
Write the operation of a commutative semigroup additively. Endomorphisms
need not commute under composition. Their pointwise sum is again an
endomorphism: if $e,f\in\End(S,+)$, then
\[
 (e+f)(a+b)=e(a)+e(b)+f(a)+f(b)
           =(e(a)+f(a))+(e(b)+f(b)).
\]
This is where commutativity is essential.

\begin{lemma}[Stabilization of powers of a finite subsemigroup]\label{f3:lem:powers}
Let $T$ be a nonempty finite subsemigroup of a semigroup, with $|T|\leq n$.
For setwise products, $T^k=T^n$ for every $k\geq n$.
\end{lemma}
\begin{proof}
Closure gives $T^2\subseteq T$. Multiplying an inclusion by $T$ gives
$T\supseteq T^2\supseteq T^3\supseteq\cdots$. Once equality occurs, all
subsequent terms are equal. If no equality occurred among the first $n$
inclusions, the nonempty set $T^{n+1}$ would have size at most $|T|-n\leq0$.
Thus equality occurs at or before $T^n=T^{n+1}$, as required.
\end{proof}

\begin{theorem}\label{f3:thm:main}
Let $S$ be a finite commutative semigroup of order $n$. Expand $S$ by
finitely many unary endomorphisms and finitely many distinguished constants,
either collection possibly empty. The resulting algebra is finitely based.
\end{theorem}
\begin{proof}
\emph{Step 1: a finite definitional expansion.}
Let $E$ be the set of endomorphisms generated by the identity map and the
given unary maps under composition and pointwise addition. It is finite,
with $|E|\leq n^n$, and every member has a representing unary term in the
original expanded signature. Add a unary symbol for each member of $E$.
If constants are present, let $K$ be their closure under addition and
the maps in $E$, and add a constant symbol for each element of $K$.
Each added constant is represented by a ground term. If there are no
constants, put $K=\varnothing$. These are finite definitional expansions;
in particular, no arbitrary, non-term-definable endomorphism is added.

\emph{Step 2: normalization.}
Take associativity and commutativity of addition, together with the finite
tables of valid identities
\begin{align*}
 \mathrm{id}(x)&\eq x,&
 e(x+y)&\eq e(x)+e(y),\\
 e(f(x))&\eq(e\circ f)(x),&
 e(x)+f(x)&\eq(e+f)(x),
\end{align*}
for $e,f\in E$, and the constant addition and action tables for $K$.
Identify each original unary or constant symbol with its chosen representative
by an equation if necessary.
The notation on the right denotes a unary symbol from $E$.
Every term reduces to a sum with exactly one summand $e_x(x)$ for each
occurring variable $x$, and at most one constant from $K$. Terms are
nonempty: a sum with neither variables nor a constant is not allowed.
For example, repeated occurrences of a variable combine using the
pointwise-addition table. The representative terms defining the added
symbols reduce to those symbols by these same tables.

\emph{Step 3: compression of variable types.}
Consider a normalized identity $u\eq v$. Use a formal symbol $\bot$ for
an absent coefficient. Each variable has a type
\[
 \tau=(e,f)\in(E\cup\{\bot\})^2\setminus\{(\bot,\bot)\},
\]
according to its coefficient in $u$ and $v$.
Adjoin a \emph{fresh} identity element $1$ to $S$, even if $S$ already
has an identity, and denote the resulting monoid by $S^1$.
Interpret $\bot$ as the constant map from $S$ to this new identity.
For each type put
\[
 T_\tau=\{(e(a),f(a)):a\in S\}\subseteq (S^1)^2.
\]
Each $T_\tau$ is a subsemigroup: the map $a\mapsto(e(a),f(a))$ is a
homomorphism, including when one component is absent. Also
$|T_\tau|\leq n$.

If $k_\tau$ variables have type $\tau$, their independent contributions
to the pair of term values form precisely $T_\tau^{k_\tau}$. More
explicitly, the set of all pairs of values of $(u,v)$ is
\begin{equation}\label{f3:eq:range}
 (c_u,c_v)+\sum_{\tau:k_\tau>0}T_\tau^{k_\tau},
\end{equation}
where a missing constant is interpreted as the fresh identity $1$;
the sums in this display are setwise sums in $(S^1)^2$.
This formula follows because variables of distinct types can be assigned
independently. Actual term values remain in $S$.
By Lemma~\ref{f3:lem:powers}, replacing every $k_\tau>n$ by $n$ leaves
\eqref{f3:eq:range} unchanged. Consequently the identity with these reduced
multiplicities is valid in $S$ if and only if the original identity is valid:
both assertions say that the same set of value pairs lies on the diagonal.

\emph{Step 4: a finite complete set of identities.}
Put
\[
 N=n\bigl((|E|+1)^2-1\bigr).
\]
On a fixed alphabet of $N$ variables there are only finitely many normalized
terms: for each variable choose an element of $E$ or absence, and choose
at most one member of $K$. Add every valid identity between these finitely
many terms to the normalization identities. Ground identities are included
when $K\ne\varnothing$. This gives a finite set $\Sigma$.

Every valid identity normalizes and compresses to a valid identity on at
most $N$ variables, hence to an instance of one in $\Sigma$ after renaming.
To recover the uncompressed identity, for a type with original multiplicity
$k>n$, keep $n-1$ of its variables and substitute the sum of $k-n+1$ fresh
variables for its last retained variable. Distribute its coefficient on
each side. If that variable is absent from one side, it stays absent there.
Perform these substitutions with disjoint variable sets for all types.
Associativity and commutativity then recover the original normalized
identity. Thus $\Sigma$ is complete.

Finally eliminate the finitely many added symbols using their representing
terms. Their defining equations follow from the normalization tables, so
this is a finite definitional reduction of a complete basis. It gives a
finite basis in the originally specified expanded signature.
\end{proof}

\begin{remark}
The bound is an existence bound, not a practical claim that every resulting
basis is short. The argument allows constants with no absorption or identity
property. It proves ordinary finite basability. It does not establish
hereditary finite basability of the generated variety, since its arbitrary
subvarieties need not be generated by a finite algebra of the displayed form.
\end{remark}

\subsubsection{A term-definition criterion for semirings}
\begin{corollary}[Linear reduction]\label{f3:cor:linear}
Let $A=(S,+,\cdot)$ be a finite semiring. Retain either binary operation
$\circ$, provided it is commutative. Suppose the other operation $*$ has
the form
\begin{equation}\label{f3:eq:linear}
 x*y=f(x)\circ g(y)\circ c,
\end{equation}
where some of the three displayed contributions may be absent, but not all.
Assume every present $f,g$ is a unary term operation of $A$ and an
endomorphism of $(S,\circ)$. Assume a present $c$ is the constant value of
a unary term $p(x)$ of $A$. Then $A$ is finitely based.
\end{corollary}
\begin{proof}
The expansion $C=(S,\circ,f,g,c)$, omitting absent symbols, is finitely
based by Theorem~\ref{f3:thm:main}. Let $\Sigma$ be a finite basis for it.
Replace $f,g$ by their original unary terms and, when needed, replace
$c$ by $p(z)$ with a fresh variable $z$. Add \eqref{f3:eq:linear}, with these
same replacements, and $p(x)\eq p(y)$ when $c$ is present.
All these finitely many identities hold in $A$.

Any nonempty model $D$ of them expands to a $C$-type algebra by interpreting
$f,g$ by the given terms and $c$ by the common value of $p$. The translated
identities ensure that this expansion satisfies $\Sigma$, so it belongs
to $\Var(C)$. Conversely, the $*$ operation of $D$ is defined by
\eqref{f3:eq:linear}. Every identity of $A$, after eliminating $*$, is an
identity of $C$ and therefore holds in $D$. This proves completeness in
the original binary signature. The fresh-variable translation explicitly
handles the absence of nullary symbols in that signature.
\end{proof}

\begin{corollary}
Every finite semiring with constant multiplication is finitely based.
\end{corollary}
\begin{proof}
Retain addition and take only the constant contribution in
\eqref{f3:eq:linear}. Its value is defined by the unary term $p(x)=x^2$.
\end{proof}

For a less degenerate example, the specified $A_{140}$ has tables
\[
 +=\begin{pmatrix}0&0&0&0\\0&1&0&3\\0&0&2&0\\0&3&0&0\end{pmatrix},
 \qquad
 \cdot=\begin{pmatrix}0&0&0&0\\0&1&0&0\\0&0&0&0\\0&0&0&0\end{pmatrix}.
\]
The map $f(x)=x^2$ has values $(0,1,0,0)$, satisfies
$f(x+y)=f(x)+f(y)$, and gives
\[
 xy=x^2+y^2.
\]
Corollary~\ref{f3:cor:linear} proves finite basability. This conclusion uses
both the endomorphism theorem and a two-way term translation; it does
not follow merely from finite basability of the additive reduct.

\endgroup
\subsection{Capped and additively nilpotent semirings}\label{sec:family4}
\begingroup
\providecommand{\Var}{}\renewcommand{\Var}{\operatorname{Var}}
\providecommand{\eq}{}\renewcommand{\eq}{\approx}
\providecommand{\SR}{}\renewcommand{\SR}{\mathrm{SR}}
\providecommand{\cS}{}\renewcommand{\cS}{\mathcal S}
\providecommand{\FB}{}\renewcommand{\FB}{\mathrm{FB}}\providecommand{\NFB}{}\renewcommand{\NFB}{\mathrm{NFB}}
\subsubsection{The capped natural-number semiring}
Let $T=\{0,1,2\}$, with
\[
 a+b=\min(2,a+b\text{ in }\mathbb N),\qquad
 ab=\min(2,ab\text{ in }\mathbb N).
\]
It is the semiring obtained by identifying all natural numbers at least
two, then forgetting the constant symbols. Its tables are
\[
 +=\begin{pmatrix}0&1&2\\1&2&2\\2&2&2\end{pmatrix},\qquad
 \cdot=\begin{pmatrix}0&0&0\\0&1&2\\0&2&2\end{pmatrix}.
\]
For a nonempty finite set $D$ of variables, write $m_D=\prod_{x\in D}x$.
An empty product is never used as a term.

\begin{theorem}\label{f4:thm:T}
The following four identities, together with $\SR$, form a finite basis
for $T$:
\begin{align}
 xy&\eq yx, & x^2&\eq x,\label{f4:eq:Tband}\\
 x+y+xy&\eq x+y,\label{f4:eq:Tabs2}\\
 x+y+xyz&\eq x+y.\label{f4:eq:Tabs3}
\end{align}
Thus $T$ has a basis of nine identities.
\end{theorem}
\begin{proof}
All four identities hold in $T$. For the absorption identities, if both
of the first two summands are nonzero their sum is already $2$; if either
is zero, the additional product is zero.

Distributivity and \eqref{f4:eq:Tband} reduce every term to a nonempty sum
of square-free monomials. Substituting $y=x$ in \eqref{f4:eq:Tabs2} gives
$3x\eq2x$, so each monomial need occur at most twice.
Fix the finite alphabet $X$ consisting of all variables on the two sides
of an identity. For a polynomial $t$, let $c_t(D)\in\{0,1,2\}$ be the
multiplicity of $m_D$. For $Y\subseteq X$ define
\[
 b_t(Y)=\min\left(2,\sum_{\varnothing\ne D\subseteq Y}c_t(D)\right).
\]
This is exactly the value of $t$ on the assignment taking variables in
$Y$ to $1$ and all other variables to $0$.
Call $Y$ saturated if $b_t(Y)=2$.

If $Y$ is saturated, there are two summand occurrences $m_D,m_E$ of
$t$, possibly equal, with $D\cup E\subseteq Y$. In their presence one
can add $m_Y$: use \eqref{f4:eq:Tabs2} with $x=m_D,y=m_E$ when
$Y=D\cup E$, and \eqref{f4:eq:Tabs3} with the extra factor
$z=m_{Y\setminus(D\cup E)}$ otherwise. Multiplicative idempotence
reduces the resulting product to $m_Y$. These equations can be used
twice, so two copies can be added. Adding a monomial on a saturated
support creates no new saturated support: every superset of that
support was already saturated.

Consequently $t$ is provably equal to the following maximal polynomial:
put coefficient two on every saturated support; leave coefficient one
on each originally occurring unsaturated support; and put coefficient
zero on all remaining supports. There are only finitely many supports
in $X$, so this is a finite equational derivation for each $t$.

The maximal polynomial is determined by $b_t$. Indeed, an originally
occurring unsaturated support $D$ has $c_t(D)=1$ and contains no
proper occurring support; otherwise $b_t(D)$ would be two. Thus such
supports are precisely the inclusion-minimal positive supports $D$
with $b_t(D)=1$. Conversely, at an inclusion-minimal positive support
the unique active summand with value one must have that exact support.
The saturated supports are read directly from $b_t$.

If $t\eq u$ holds in $T$, their values on all the displayed $0,1$
assignments agree, hence $b_t=b_u$. Their maximal polynomials are
identical, so the proposed identities derive $t\eq u$. This proves
completeness for identities with arbitrarily many variables.
\end{proof}

\subsubsection{A semiring with constant nonzero products}
Let $U=\{0,1,2\}$ have the same addition as $T$, and multiplication
\[
 a\mathbin{\cdot_U}b=
 \begin{cases}0,&a=0\text{ or }b=0,\\2,&a\ne0\text{ and }b\ne0.
 \end{cases}
 \qquad
 \cdot_U=\begin{pmatrix}0&0&0\\0&2&2\\0&2&2\end{pmatrix}.
\]
This operation is associative and distributes over the capped addition.
For example, $b+c$ is nonzero exactly when at least one of $b,c$ is
nonzero, which verifies distributivity directly. All products of at
least two nonzero arguments are two.

\begin{theorem}\label{f4:thm:U}
The following six identities, together with $\SR$, form a finite basis
for $U$:
\begin{align}
 xy&\eq yx,&2x&\eq x^2,&3x&\eq2x,\label{f4:eq:Ubasic}\\
 x^2y&\eq xy,\label{f4:eq:Urepeat}\\
 xy+xyz&\eq xy,\label{f4:eq:Uabsprod}\\
 x+y+xy&\eq x+y.\label{f4:eq:Uabslinear}
\end{align}
Thus $U$ has a basis of eleven identities.
\end{theorem}
\begin{proof}
Each displayed identity holds in $U$. A product of length at least two
is either zero or two, which verifies \eqref{f4:eq:Urepeat} and
\eqref{f4:eq:Uabsprod}; \eqref{f4:eq:Uabslinear} is checked as in $T$.

For a nonempty support $D$, let $u_D=x^2$ when $D=\{x\}$, and let
$u_D=m_D$ when $|D|\geq2$. After distributive expansion, commutativity
and \eqref{f4:eq:Urepeat} reduce every nonlinear word to its $u_D$.
Moreover,
\[
 2(xy)=(2x)y=x^2y=xy,
\]
so nonlinear summands are additively idempotent. Pairs of equal linear
summands become $x^2$ by \eqref{f4:eq:Ubasic}. We obtain a polynomial
described by a set $L\subseteq X$ of linear variables and a family
$\mathcal F$ of nonempty supports for nonlinear summands. Neither need
be nonempty separately, but the polynomial itself is nonempty.

On a $0,1$ assignment supported on $Y\subseteq X$, its value is
\[
 b_t(Y)=\min\bigl(2,\ |L\cap Y|
           +2|\{D\in\mathcal F:D\subseteq Y\}|\bigr).
\]
Again call the supports of value two saturated.
If $Y$ contains a member $D$ of $\mathcal F$, use
\eqref{f4:eq:Uabsprod} to add $u_Y$ in the presence of $u_D$:
factor $u_D$ into two nonempty words and use the additional variables
as the third factor when $Y\ne D$.
If $Y$ is saturated only because it contains two linear variables,
first add their product by \eqref{f4:eq:Uabslinear}, then enlarge its
support using \eqref{f4:eq:Uabsprod}. Thus $u_Y$ can be inserted for
every saturated $Y$. Such insertions do not create new saturated
supports, because supersets of a saturated support are saturated.

After all insertions, a linear variable $x$ with saturated singleton
$\{x\}$ can be removed. The term $x^2$ is present and
\[
 x+x^2=3x=2x=x^2.
\]
The resulting polynomial has one nonlinear summand $u_Y$ for every
saturated support $Y$, and has precisely the linear variables $x$
for which $b_t(\{x\})=1$. These data depend only on $b_t$.
Every transformation just used follows from the proposed finite basis.

If $t\eq v$ holds in $U$, their $0,1$ evaluations agree, so they reduce
to the same polynomial. Hence the displayed identities derive every
identity of $U$.
\end{proof}

The local representatives and the isomorphisms to these displayed tables
are as follows. A map lists images in the local order $0,1,2$.
\begin{center}
\begin{tabular}{llll}
\toprule Representative&Addition&Multiplication&Map to standard form\\\midrule
$B_{113}\cong T$&\texttt{012111211}&\texttt{000011012}&$(0,2,1)$\\
$B_{69}\cong U$&\texttt{000012020}&\texttt{010111010}&$(2,0,1)$\\
\bottomrule
\end{tabular}
\end{center}
The numerical labels in the standard forms are used to explain their
operations; zero and one are still not named symbols of their signatures.

\subsubsection{Commutative multiplication and nilpotent addition}
We call addition nilpotent of degree $d\geq2$ if there is an element $z$
such that every sum of $d$ elements equals $z$. This refers to sums
of arbitrary, not necessarily equal, elements.

\begin{lemma}\label{f4:lem:zero}
If addition is nilpotent of degree $d$ in a semiring, its constant
$d$-fold sum $z$ is absorbing for both operations.
\end{lemma}
\begin{proof}
Associativity gives
\[
 z+a=x_1+\cdots+x_{d-1}+(x_d+a)=z.
\]
Distributivity gives $az=ax_1+\cdots+ax_d=z$, and similarly $za=z$.
Thus all longer sums also equal $z$.
\end{proof}

\begin{theorem}\label{f4:thm:nil}
Every finite semiring with commutative multiplication and nilpotent
addition is finitely based.
\end{theorem}
\begin{proof}
Let $A$ have $n$ elements and additive nilpotence degree $d$. Temporarily
name the absorbing element $z$ of Lemma~\ref{f4:lem:zero}. Choose positive
integers $r,q$ such that the finite multiplicative semigroup satisfies
$x^{r+q}\eq x^r$, and put $L=r+q-1$.
Include $\SR$, commutativity of multiplication, this power identity,
the absorption identities for $z$, and
\[
 x_1+\cdots+x_d\eq z.
\]
These are finitely many valid identities. Distributivity normalizes any
term to $z$ or a nonempty sum of fewer than $d$ commutative monomials.
Every positive variable exponent in a monomial can be reduced to one
of $1,\ldots,L$.

Consider an identity between two such normalized terms, and list all
monomial slots on both sides. There are $m\leq2(d-1)$ slots; a side
equal to $z$ has no slots. For each occurring variable record the vector
of its exponents in the slots, using exponent zero for absence. Thus
its profile lies in
\[
 \{0,1,\ldots,L\}^{m}\setminus\{(0,\ldots,0)\}.
\]
If $m=0$, the identity is just $z\eq z$. Otherwise adjoin a fresh formal
identity to the multiplicative semigroup, obtaining $M^1$, even if the
original semigroup already had an identity. For a profile
$\alpha=(\alpha_1,\ldots,\alpha_m)$, put
\[
 T_\alpha=\{(a^{\alpha_1},\ldots,a^{\alpha_m}):a\in A\}
       \subseteq(M^1)^m,
\]
interpreting exponent zero as the fresh identity. Because multiplication
is commutative, this map from $A$ is a semigroup homomorphism.
Consequently $T_\alpha$ is a subsemigroup of size at most $n$.

Setwise powers of a finite subsemigroup form a descending chain:
$T_\alpha\supseteq T_\alpha^2\supseteq\cdots$. A nonempty set with
at most $n$ elements cannot strictly decrease $n$ times; once equality
occurs it persists. Hence $T_\alpha^k=T_\alpha^n$ for $k\geq n$.
Variables of a fixed profile contribute independently to the vector
of all monomial values. Replacing more than $n$ variables of that
profile by exactly $n$ variables therefore preserves the set of
possible vectors. Doing this for each profile preserves that set
globally, since different profiles use independent variables.
Taking the sums of the indicated slots on the two sides then preserves
the set of possible pairs of term values. Thus compression preserves
validity of the identity.

Every normalized identity compresses to at most
\[
 N=n\bigl((L+1)^{2(d-1)}-1\bigr)
\]
variables. There are finitely many normalized terms on a fixed alphabet
of $N$ variables: finitely many monomials with bounded exponents, fewer
than $d$ summands, and the one constant term $z$. Add all valid identities
between them to the structural identities above.

To recover an original group of $k>n$ variables with a given profile,
replace one of its retained variables by a product of $k-n+1$ fresh
variables and leave the other $n-1$ alone. Commutativity distributes
each exponent across that product, restoring exactly the original
exponent profile in every slot. Different profiles use disjoint fresh
variables. Hence the finite set just constructed is complete.

Finally replace $z$ by the unary term $dx$ in a fresh variable and add
$dx\eq dy$. Any nonempty model then has a uniquely determined value
for $z$ and expands to a model of the named-constant basis. Conversely
all translated identities hold in $A$. This eliminates the constant
and supplies a finite basis in the original two-binary-operation signature.
\end{proof}

The commutativity assumption concerns multiplication, and is used in
the homomorphism property of each exponent profile. Nilpotence concerns
addition, and bounds the number of monomial slots. Neither assumption
can simply be omitted from this proof. The bound is constructive but
potentially large; no minimal basis claim is intended.

Theorem~\ref{f4:thm:nil} directly adds the following 25 four-element classes:
\[
\begin{gathered}
122,128,435,518,524,610,645,646,1345,1402,1403,1404,1405,\\
1428,1457,1458,1459,1472,1473,1474,1766,1831,1833,1834,1934.
\end{gathered}
\]
The tables give commutativity and the stated additive nilpotence degrees.
Most have degree three; $A_{1404}$ has degree four. The criterion also
applies to the three-element source $B_{86}$.

\endgroup
\subsection{Commuting interiors and mixed annihilation}\label{sec:family5}
\begingroup
\providecommand{\eq}{}\renewcommand{\eq}{\approx}
\providecommand{\Var}{}\renewcommand{\Var}{\operatorname{Var}}
\providecommand{\SR}{}\renewcommand{\SR}{\mathrm{SR}}
\subsubsection{A commutative quotient for interior factors}
Let $M$ be a semigroup satisfying
\begin{equation}\label{f5:eq:perm}
 xyzt\eq xzyt.
\end{equation}
Adjoin a fresh identity $1$ to $M$, even when $M$ already has an
identity, and denote the resulting monoid by $M^1$. Define a relation
on $M^1$ by
\begin{equation}\label{f5:eq:congruence}
 u\sim v \quad\Longleftrightarrow\quad
 aub=avb\quad\text{for every }a,b\in M.
\end{equation}
The endpoints belong to $M$, whereas $u,v$ may be the new identity.

\begin{lemma}\label{f5:lem:quotient}
The relation $\sim$ is a monoid congruence and $C=M^1/{\sim}$ is a
commutative monoid. If $M$ has $n$ elements, then $|C|\leq n+1$.
\end{lemma}
\begin{proof}
Reflexivity, symmetry and transitivity follow from equality. If
$u\sim v$ and $s\in M^1$, then $as,sb\in M$ for all $a,b\in M$.
Thus
\[
 a(su)b=(as)ub=(as)vb=a(sv)b,
 \qquad a(us)b=au(sb)=av(sb)=a(vs)b.
\]
This proves compatibility on both sides. For $u,v\in M$, identity
\eqref{f5:eq:perm} gives $a(uv)b=a(vu)b$ for every $a,b\in M$, so
$uv\sim vu$. If either $u$ or $v$ is the new identity the same
conclusion is immediate. The quotient is therefore commutative.
The size bound follows from $|M^1|=n+1$.
\end{proof}

The quotient need not distinguish the original elements of $M$.
Its role is precisely to record everything an interior product can
contribute to a word whose two endpoints are fixed. In particular,
for each $a,b\in M$ the map
\[
 F_{a,b}:C\longrightarrow M,\qquad F_{a,b}([u])=aub
\]
is well defined. No homomorphism property of $F_{a,b}$ is needed.

\begin{lemma}\label{f5:lem:powers}
If $T$ is a nonempty finite subsemigroup of any semigroup and
$|T|\leq n$, then its setwise powers satisfy $T^k=T^n$ for $k\geq n$.
\end{lemma}
\begin{proof}
This is Lemma~\ref{f3:lem:powers}.
\end{proof}

\subsubsection{Nilpotent addition with commuting interiors}
Addition is nilpotent of degree $d\geq2$ if every sum of $d$ elements,
not necessarily equal, has a fixed value $o$.

\begin{lemma}\label{f5:lem:zero}
This element $o$ is absorbing for both operations.
\end{lemma}
\begin{proof}
Apply Lemma~\ref{f4:lem:zero}.
\end{proof}

\begin{theorem}\label{f5:thm:interior}
Every finite semiring with nilpotent addition and multiplication
satisfying \eqref{f5:eq:perm} is finitely based. More precisely, let its
order be $n$, let $d\geq2$ be an additive nilpotence degree, and choose
positive $r,q$ such that $x^{r+q}\eq x^r$ holds multiplicatively.
Set $L=r+q-1$. In the signature with the absorbing element temporarily
named, a complete finite basis can be constructed using normalized
identities on at most
\begin{equation}\label{f5:eq:bound}
 N=4(d-1)+n\bigl((L+1)^{2(d-1)}-1\bigr)
\end{equation}
variables, in addition to finitely many structural identities.
\end{theorem}
\begin{proof}
\emph{Normalization.} Name the absorbing element $o$. Include
$\SR$, \eqref{f5:eq:perm}, $x^{r+q}\eq x^r$, the absorption identities
for $o$, and $x_1+\cdots+x_d\eq o$.
Finite multiplicative semigroups always admit such positive $r,q$.
Distributive expansion reduces a term to $o$ or a nonempty sum of
fewer than $d$ nonempty words. Words of length one remain single
variables. In every longer word retain its first and last letters.
All interior letters may be permuted: an adjacent interior swap is an
instance of \eqref{f5:eq:perm}, since its prefix and suffix are both
nonempty. Group equal interior letters, and use the power identity
to reduce each positive interior exponent to $1,\ldots,L$.
Exponents zero record absence. Reductions are performed only inside
the retained endpoints; words of length two have empty interiors.

\emph{Protected variables.} In a normalized identity, list the word
slots on both sides. Their total number is $m\leq2(d-1)$, with a
side $o$ contributing no slots. If $m=0$ the identity is trivial.
Otherwise protect every variable appearing as either endpoint of any
slot, including the variable in each length-one slot. There are at
most $2m$ protected variables. Keep every occurrence of these variables,
including their interior occurrences. Every unprotected variable
appears only in interiors and has a nonzero profile
\[
 \alpha=(\alpha_1,\ldots,\alpha_m)\in
 \{0,\ldots,L\}^{m}\setminus\{0\},
\]
where $\alpha_j$ is its interior exponent in slot $j$.

\emph{Preservation of evaluation ranges.} Apply
Lemma~\ref{f5:lem:quotient} to the multiplicative reduct $M$. For a profile
$\alpha$ put
\[
 T_\alpha=\{([a]^{\alpha_1},\ldots,[a]^{\alpha_m}):a\in M\}
 \subseteq C^m,
\]
using the identity of $C$ for exponent zero. The displayed map from
$M$ to $C^m$ is a semigroup homomorphism: $C$ is commutative and
$[ab]=[a][b]$. Thus $T_\alpha$ is a subsemigroup of size at most $n$.

Fix an assignment to all protected variables. It determines each
endpoint and the contribution of all protected interior letters in
each slot. A length-one slot is now fixed. In every other slot, the
quotient class of the remaining interior product, together with these
fixed data, determines its value by the map $F_{a,b}$. A group of
$k$ unprotected variables of profile $\alpha$ contributes precisely
$T_\alpha^k$ to the vector of interior classes. Distinct groups use
independent variables, so their contributions multiply setwise in
$C^m$. Lemma~\ref{f5:lem:powers} therefore permits replacing $k>n$
variables of any profile by $n$ variables without changing the set of
possible word-value vectors, conditional on the protected assignment.
Taking the prescribed sums on both sides preserves the possible
pairs of term values as well. This holds for every protected assignment.
Consequently the compressed identity is valid if and only if the
original identity is valid.

\emph{A finite set of identities.} The compressed identity has at most
$2m+n((L+1)^m-1)\leq N$ variables. On a fixed alphabet of $N$
variables there are only finitely many normalized words: there are
finitely many choices of endpoints and of bounded interior exponents,
as well as finitely many single letters. Hence there are finitely many
normalized terms with fewer than $d$ slots, together with $o$.
Add all valid identities between these terms to the structural
identities already listed. This gives a finite set $\Sigma$.

\emph{Recovery by substitution.} For a profile originally represented
by $k>n$ variables, replace one retained variable by a product of
$k-n+1$ fresh variables and replace the other $n-1$ retained variables
by the remaining single variables. Every occurrence being replaced
lies in an interior. An occurrence of exponent $j$ becomes $j$ copies
of that product. Interior swaps rearrange these copies into exponent
$j$ for each fresh variable. Thus the substitution recovers exactly
the original profile in every slot, without changing any endpoint.
Different profiles use disjoint variables. Renaming and these
substitutions recover the entire original normalized identity from
one in $\Sigma$. Normalization then proves completeness for arbitrary
identities, not merely those of bounded length or variable number.

\emph{Removing the constant.} The unary term $h(x)=dx$ has constant
value $o$. Replace $o$ in each member of $\Sigma$ by $h(v)$, with $v$
fresh, and add $h(x)\eq h(y)$. Any nonempty model has a uniquely
determined interpretation of $o$ and expands to a model of $\Sigma$.
It therefore satisfies every identity of the original semiring.
Conversely the translated finite set holds in that semiring. This
gives a complete finite basis in the original signature.
\end{proof}

\begin{remark}
Commutativity of multiplication implies \eqref{f5:eq:perm}, so
Theorem~\ref{f5:thm:interior} contains the earlier commutative result
Theorem~\ref{f4:thm:nil}. The present proof does not assert that the multiplicative
reduct itself is commutative. Protecting endpoints is essential both
for the sandwich quotient and for the recovery substitution.
\end{remark}

\subsubsection{Mixed annihilation and reduct bases}\label{f5:sec:mixed}
We next treat a different class of semirings. The condition that
$xy+z$ is constant can be written without naming that constant as
\begin{equation}\label{f5:eq:mixed}
 xy+z\eq uv+w.
\end{equation}
Denote its common value by $o$. Associativity gives $o+a=o$, while
distributivity gives
$ao=a(xy+z)=(ax)y+az=o$ and $oa=(xy+z)a=x(ya)+za=o$.
Thus $o$ is again absorbing for both operations.

\begin{lemma}\label{f5:lem:namedzero}
If a semigroup with an absorbing zero is finitely based, its expansion
by a symbol naming that zero is finitely based.
\end{lemma}
\begin{proof}
Take a finite semigroup basis and add $ox\eq xo\eq o$.
Every term containing $o$ reduces to $o$. An identity without $o$
is already covered by the original basis. If a constant-free word
$w$ has value always $o$, the identity $w\eq wy$, with $y$ fresh,
holds in the reduct. Derive it from the reduct basis and substitute
$o$ for $y$ to obtain $w\eq o$. This handles all remaining identities.
The same argument applies to an additively written semigroup.
\end{proof}

\begin{lemma}\label{f5:lem:nilsemigroup}
A finite semigroup of order $n$ with zero as its only idempotent is
nilpotent. Every product of $n$ elements equals zero.
\end{lemma}
\begin{proof}
Every element has an idempotent power, hence a power equal to zero.
Suppose a product of $n$ elements were nonzero. All its nonempty
prefix products would be nonzero. Among $n$ prefixes, two must agree,
since there are only $n-1$ nonzero elements. Writing the shorter
prefix as $a$ and the intervening nonempty block as $b$ gives $ab=a$.
Then $ab^k=a$ for all $k$, contrary to a power of $b$ being zero.
\end{proof}

\begin{theorem}\label{f5:thm:mixed}
Let $A$ be a finite semiring satisfying \eqref{f5:eq:mixed}. If its
multiplicative semigroup reduct is finitely based, then $A$ is
finitely based. In particular, every such semiring of order at most
four is finitely based.
\end{theorem}
\begin{proof}
Temporarily name $o$. Every finite commutative semigroup is finitely
based, as also follows from the general expansion theorem in
Section~\ref{sec:family3}. Lemma~\ref{f5:lem:namedzero} supplies finite bases
$\Sigma_+$ and $\Sigma_\cdot$ for both reducts with $o$ named.
Include these bases, $\SR$, $xy+z\eq o$, and absorption by $o$.

Any term containing both operations reduces to $o$. Indeed a sum
with a product summand is covered by $xy+z\eq o$. A product with
a sum factor becomes a sum of products by distributivity, and is
therefore also $o$. Applying this at a first change of operation in
a term, then absorption in its surrounding context, proves the
assertion. Normal terms are consequently $o$, pure additive sums
of variables, or pure multiplicative words.

Suppose first that the multiplicative reduct has an idempotent
$e\ne o$. Then $e+e=e^2+e=o$. Under the assignment taking every
variable to $e$, every multiplicative word, including a single
variable, has value $e$, whereas every additive sum of length at
least two has value $o$. No identity between these two families
can hold. Every valid identity between normal terms therefore
belongs to one of the two reduct languages (with $o$ allowed),
and follows from $\Sigma_+$ or $\Sigma_\cdot$.

Suppose instead that the additive reduct has an idempotent $f\ne o$.
Distributivity gives $f^2+f^2=(f+f)f=f^2$, while mixed annihilation
gives $f^2+f^2=o$, so $f^2=o$. Assigning every variable to $f$
makes every additive sum, including a single variable, equal $f$,
and every multiplicative word of length at least two equal $o$.
This again separates the two families, and the same completeness
argument applies.

Finally, if neither reduct has a nonzero idempotent, both are
nilpotent by Lemma~\ref{f5:lem:nilsemigroup}. Put $d=k=\max(2,|A|)$;
all $d$-fold sums and all $k$-fold products equal $o$. Add these
two identities. Every normal term is then $o$, a pure sum of fewer
than $d$ variables, or a pure word of fewer than $k$ variables.
There are only finitely many such terms on a fixed alphabet of
$2\max(d-1,k-1)$ variables. Adding all valid identities between
them is complete by normalization and variable renaming.

In every case eliminate the named constant by the unary term
$h(x)=x^2+x$, adding $h(x)\eq h(y)$ as in
Theorem~\ref{f5:thm:interior}. This proves the general assertion.
For order at most four, the multiplicative reduct is finitely based
by Edmunds' theorem \cite{Edmunds}, which proves the final assertion.
\end{proof}

This criterion uses a normal-form separation argument. It does not
infer finite basability from the finite basability of arbitrary
semiring factors, nor from an arbitrary join of finitely based varieties.

\endgroup
\subsection{Regularization by a semilattice}\label{sec:family7}
\begingroup
\providecommand{\Var}{}\renewcommand{\Var}{\operatorname{Var}}
\providecommand{\Id}{}\renewcommand{\Id}{\operatorname{Id}}
\providecommand{\Reg}{}\renewcommand{\Reg}{\operatorname{Reg}}
\providecommand{\SL}{}\renewcommand{\SL}{\mathcal{SL}}
\providecommand{\V}{}\renewcommand{\V}{\mathcal V}
\providecommand{\eq}{}\renewcommand{\eq}{\approx}
\subsubsection{A finite basis for regularization}
We recall the partition construction of P\l onka and give the
finite-basis translation explicitly.
Let $E$ be the two-element algebra in which both operations are the
same semilattice operation, and put $\SL=\Var(E)$. A term has
\emph{content} equal to its set of variables. An identity is
\emph{regular} if its two terms have the same content. Exactly the
regular identities hold in $E$. Consequently, for any variety $\V$
in the present signature,
\[
 \Reg(\V):=\V\vee\SL
 \quad\text{satisfies exactly the regular identities of }\V.
\]
We call $\V$ \emph{strongly irregular} if a binary term $q(x,y)$
containing both variables satisfies $q(x,y)\eq x$ in $\V$.
The term $q$ is an abbreviation, not an additional basic operation.

For the moment write $x\circ y=q(x,y)$. Consider the following
partition identities:
\begin{align}
 x\circ x&\eq x, \tag{P1}\\
 (x\circ y)\circ z&\eq x\circ(y\circ z), \tag{P2}\\
 x\circ(y\circ z)&\eq x\circ(z\circ y), \tag{P3}\\
 f(x,y)\circ z&\eq f(x\circ z,y\circ z),\qquad f\in\{+,\cdot\}, \tag{P4}\\
 z\circ f(x,y)&\eq(z\circ x)\circ y,\qquad f\in\{+,\cdot\}. \tag{P5}
\end{align}
There are seven identities in this list. They are regular and hold
in both $\V$ and $E$: in $\V$, $\circ$ is the first projection;
in $E$, it is the semilattice operation.

The following is the form of P\l onka's partition decomposition
\cite{Plonka} for the present binary signature.

\begin{lemma}[Partition decomposition]\label{f7:lem:partition}
An algebra satisfying (P1)--(P5) is a P\l onka sum of subalgebras
on which $x\circ y=x$. Its index algebra is a semilattice, with both
basic operations inducing the same semilattice operation.
\end{lemma}
\begin{proof}
Define $a\sim b$ by $a\circ b=a$ and $b\circ a=b$.
Reflexivity and symmetry are immediate. If $a\sim b\sim c$, then
\[
 a\circ c=(a\circ b)\circ c
 =a\circ(b\circ c)=a\circ b=a;
\]
the reverse equation follows in the same way. Hence $\sim$ is an
equivalence relation. If $x,y$ belong to the class of $a$, then
(P4) and (P5) give
\[
 f(x,y)\circ a=f(x,y),\qquad
 a\circ f(x,y)=(a\circ x)\circ y=a.
\]
Thus every class is a subalgebra.

By (P1)--(P3), an iterated $\circ$-product depends only on its first
element and the set of subsequent elements: subsequent factors may
be permuted and repetitions deleted, including a subsequent copy of
the first factor. These rules show that
\[
 [a]\vee[b]=[a\circ b]
\]
is well-defined, associative, commutative, and idempotent. For
example, replacing $a$ by an equivalent element does not change the
class of $a\circ b$, since multiplication of the two prospective
representatives in either order absorbs the equivalent factors.
They also show that, for $u=f(x,y)$ and $v=x\circ y$,
\[
 u\circ v=u\circ f(x,y)=u\circ u=u,
 \qquad v\circ u=v\circ x\circ y=v.
\]
Here the first equality uses (P5), and the last uses the preceding
permutation and deletion rules. Thus $[f(x,y)]=[x]\vee[y]$.

Let $A_i$ denote class $i$. For $i\leq j$, define
$\phi_{ij}:A_i\to A_j$ by $\phi_{ij}(x)=x\circ b$, where $b\in A_j$.
If $b\sim c$, then
\[
 x\circ b=x\circ b\circ c=x\circ c\circ b=x\circ c,
\]
so this definition is independent of $b$. Its image lies in $A_j$.
It is a homomorphism by (P4), is the identity when $i=j$, and satisfies
$\phi_{jk}\phi_{ij}=\phi_{ik}$ by (P2). For $x\in A_i$, $y\in A_j$,
and $k=i\vee j$, choose $b\in A_k$. Then
\[
 f(x,y)=f(x,y)\circ b
        =f(\phi_{ik}(x),\phi_{jk}(y)).
\]
This is the P\l onka sum formula.
\end{proof}

For a finite nonempty variable set $X=\{x_1,\ldots,x_m\}$, define
\[
 R_X(t)=((t\circ x_1)\circ\cdots)\circ x_m.
\]
The ordering of $X$ can be fixed arbitrarily. If $s\eq t$ is an
identity of $\V$ and $X$ is the union of its two contents, then
$R_X(s)\eq R_X(t)$ is a regular identity of $\V$.

\begin{theorem}[Finite-basis transfer]\label{f7:thm:regularization}
Let $\V$ be strongly irregular, witnessed by $q(x,y)\eq x$.
If $\Sigma$ is a finite basis for $\V$, then a finite basis for
$\Reg(\V)$ consists of (P1)--(P5) and the regularized identities
$R_X(s)\eq R_X(t)$ for $s\eq t$ in $\Sigma$.
Identities in $\Sigma$ that are already regular may be retained
unchanged. Moreover,
\[
 \V=\Reg(\V)\cap\operatorname{Mod}(q(x,y)\eq x),
\]
so $\V$ is FB if and only if $\Reg(\V)$ is FB.
\end{theorem}
\begin{proof}
Every listed identity is regular and valid in $\V$, so it holds in
$\V\vee\SL$. Conversely, let $A$ satisfy the listed basis. By
Lemma~\ref{f7:lem:partition}, it is a P\l onka sum of subalgebras $A_i$
on which $\circ$ is the first projection. Each regularized basis
identity therefore reduces to its original identity within $A_i$.
Consequently $A_i\in\V$.

Term evaluation in a P\l onka sum first sends all variables to the
component indexed by the join of their indices, and then evaluates
there. This follows by induction from the sum formula and the
homomorphism property of the connecting maps. If $s\eq t$ is regular,
both sides use the same join. Every regular identity of $\V$ thus
holds in $A$, proving $A\in\Reg(\V)$.

Finally, if $q(x,y)=x$ throughout such an algebra, there is only one
partition class, hence the algebra belongs to $\V$. Alternatively,
regularizing any identity of $\V$ and then imposing $q(x,y)\eq x$
recovers that identity. This proves the displayed equality. A finite
basis of $\Reg(\V)$ together with this single identity is then a
finite basis of $\V$, proving the converse.
\end{proof}

This theorem does not assert that an arbitrary join of FB varieties
is FB. The shared semilattice operation and the specified strongly
irregular term are essential to the proof.

\subsubsection{Applications and explicit maps}
The five source semirings are listed in Table~\ref{f7:tab:sources}.
Their carriers are $\{0,1,2\}$. The last column records a term $q$
that contains both variables and is identically the first projection.
The four sources other than $B_{78}$ are ai-semirings distinct from
$S_7$, and are FB by \cite{ZhaoThree}; see also
Section~\ref{sec:small-results}. The algebra $B_{78}$ is the field $\mathbb F_3$ in our
signature, under the label map $0\mapsto1$, $1\mapsto0$, $2\mapsto2$.
A direct finite basis for it is proved in Section~\ref{f7:sec:field}.

\begin{table}[ht]
\centering
\begin{tabular}{@{}llll@{}}
\toprule Source & Addition & Multiplication & $q(x,y)$\\\midrule
$B_{40}$ & \texttt{012111212} & \texttt{000020000} & $x+xy$ \\
$B_{54}$ & \texttt{012111212} & \texttt{000012000} & $x+xy$ \\
$B_{83}$ & \texttt{012111212} & \texttt{000010020} & $x+xy$ \\
$B_{88}$ & \texttt{012111212} & \texttt{000012020} & $x+xy$ \\
$B_{78}$ & \texttt{201012120} & \texttt{012111210} & $x+3y$ \\
\bottomrule
\end{tabular}
\caption{Strongly irregular finite generators; operation strings are row-major.}
\label{f7:tab:sources}
\end{table}

For the following maps, $E=E_7$ has
both tables \texttt{0001}. A two-digit string $bc$ denotes the
element $(b,c)$ of $C_j\times E$. Each row in Table~\ref{f7:tab:forward}
lists a subalgebra $D\leq C_j\times E$ and the images, in the same
order, of a surjective homomorphism $D\to A_i$.

\begin{longtable}{@{}llll@{}}
\caption{Forward witnesses $A_i\in\Var(C_j,E)$.}\label{f7:tab:forward}\\
\toprule Target & Source & Elements of $D$ & Images in $A_i$\\\midrule
\endfirsthead
\toprule Target & Source & Elements of $D$ & Images in $A_i$\\\midrule
\endhead
$A_{665}$ & $B_{40}$ & \texttt{00,01,10,20} & \texttt{0213} \\
$A_{714}$ & $B_{54}$ & \texttt{00,01,10,20} & \texttt{0213} \\
$A_{998}$ & $B_{54}$ & \texttt{00,10,11,20} & \texttt{0123} \\
$A_{1267}$ & $B_{40}$ & \texttt{00,01,10,11,20,21} & \texttt{101213} \\
$A_{1316}$ & $B_{54}$ & \texttt{00,01,10,11,20,21} & \texttt{101213} \\
$A_{1426}$ & $B_{83}$ & \texttt{00,01,10,20} & \texttt{0213} \\
$A_{1444}$ & $B_{88}$ & \texttt{00,01,10,20} & \texttt{0213} \\
$A_{1660}$ & $B_{78}$ & \texttt{00,10,11,20} & \texttt{0123} \\
$A_{1835}$ & $B_{78}$ & \texttt{00,01,10,11,20,21} & \texttt{101213} \\
$A_{1872}$ & $B_{83}$ & \texttt{00,01,10,11,20,21} & \texttt{202123} \\
$A_{1888}$ & $B_{88}$ & \texttt{00,01,10,11,20,21} & \texttt{202123} \\
$A_{1908}$ & $B_{83}$ & \texttt{00,10,11,20} & \texttt{0123} \\
$A_{1941}$ & $B_{88}$ & \texttt{00,10,11,20} & \texttt{0123} \\
\bottomrule
\end{longtable}

Both $C_j$ and $E$ also occur as subalgebras of $A_i$.
Table~\ref{f7:tab:backward} lists their element sets and the images of
the corresponding isomorphisms. Each displayed subset is closed
under both operations, and each displayed map preserves them.

\begin{longtable}{@{}lllll@{}}
\caption{Reverse witnesses inside each $A_i$.}\label{f7:tab:backward}\\
\toprule Target & Copy of $C_j$ & Images in $C_j$ & Copy of $E$ & Images in $E$\\\midrule
\endfirsthead
\toprule Target & Copy of $C_j$ & Images in $C_j$ & Copy of $E$ & Images in $E$\\\midrule
\endhead
$A_{665}$ & \texttt{013} & \texttt{012} & \texttt{02} & \texttt{01} \\
$A_{714}$ & \texttt{013} & \texttt{012} & \texttt{02} & \texttt{01} \\
$A_{998}$ & \texttt{013} & \texttt{012} & \texttt{12} & \texttt{01} \\
$A_{1267}$ & \texttt{023} & \texttt{012} & \texttt{01} & \texttt{10} \\
$A_{1316}$ & \texttt{023} & \texttt{012} & \texttt{01} & \texttt{10} \\
$A_{1426}$ & \texttt{013} & \texttt{012} & \texttt{02} & \texttt{01} \\
$A_{1444}$ & \texttt{013} & \texttt{012} & \texttt{02} & \texttt{01} \\
$A_{1660}$ & \texttt{013} & \texttt{012} & \texttt{12} & \texttt{01} \\
$A_{1835}$ & \texttt{023} & \texttt{012} & \texttt{12} & \texttt{01} \\
$A_{1872}$ & \texttt{013} & \texttt{012} & \texttt{02} & \texttt{10} \\
$A_{1888}$ & \texttt{013} & \texttt{012} & \texttt{02} & \texttt{10} \\
$A_{1908}$ & \texttt{013} & \texttt{012} & \texttt{12} & \texttt{01} \\
$A_{1941}$ & \texttt{013} & \texttt{012} & \texttt{12} & \texttt{01} \\
\bottomrule
\end{longtable}

\begin{theorem}\label{f7:thm:thirteen}
Every algebra in Table~\ref{f7:tab:forward} is finitely based. More
precisely, its generated variety is
\[
 \Var(A_i)=\Var(C_j,E)=\Reg(\Var(C_j)),
\]
with $C_j$ as specified in that table.
\end{theorem}
\begin{proof}
The forward map gives $A_i\in\operatorname{HS}(C_j\times E)$,
hence the first inclusion of generated varieties. The two reverse
subalgebras give the opposite inclusion. The source $C_j$ is FB
and strongly irregular, so Theorem~\ref{f7:thm:regularization} applies.
\end{proof}

The subalgebras in the forward maps have at most six elements.
Together with the reverse embeddings, they establish both inclusions
directly from the displayed operation tables.

\subsubsection{The field of three elements}
\label{f7:sec:field}
Put $\mathcal G_3=\Var(\mathbb F_3)$ and
$\mathcal R_3=\mathcal G_3\vee\SL$.

\begin{lemma}\label{f7:lem:field}
In the present signature, a basis for $\mathcal G_3$ consists of
the five semiring axioms together with
\begin{equation}\label{f7:eq:fieldbasis}
 xy\eq yx,\qquad x^3\eq x,\qquad
 3x\eq3y,\qquad x+3y\eq x.
\end{equation}
This variety is minimal nontrivial.
\end{lemma}
\begin{proof}
The last two identities make $0=3x$ a well-defined additive identity,
and $2x$ is an additive inverse of $x$. Distributivity gives
$0y=y0=0$. The models are therefore commutative rings of
characteristic dividing three satisfying $x^3=x$, viewed in a
signature without named constants.

Using these identities, each term becomes a polynomial over
$\mathbb F_3$ with zero constant coefficient, with each variable
exponent in $\{0,1,2\}$. A polynomial of degree at most two in each
variable that vanishes on all of $\mathbb F_3^k$ has all coefficients
zero. Indeed, in one variable a nonzero polynomial of degree at most
two cannot have three roots; induction applies this fact to its
coefficient polynomials in the remaining variables. Thus every
identity of $\mathbb F_3$ reduces, by the displayed axioms, to
equality of identical reduced polynomials. This proves completeness
for arbitrarily many variables. The zero polynomial is represented
by $3x$, whose value does not depend on $x$.

In any nontrivial model choose $a\ne0$ and put $e=a^2$. Then
$e\ne0$, since $a^3=a$, and $e^2=e$. The subalgebra
$\{0,e,2e\}$ is isomorphic to $\mathbb F_3$; its three elements
are distinct because the additive group has exponent three.
Every nontrivial subvariety consequently contains $\mathbb F_3$
and equals $\mathcal G_3$.
\end{proof}

\begin{theorem}\label{f7:thm:sixteen}
The algebras $A_{1660}$ and $A_{1835}$ generate $\mathcal R_3$.
A complete basis for their common variety has sixteen identities.
\end{theorem}
\begin{proof}
The equality of varieties follows from Theorem~\ref{f7:thm:thirteen}
with source $B_{78}\cong\mathbb F_3$. Here choose
\[
 q(x,y)=x+3y,\qquad R(t)=(t\circ x)\circ y.
\]
The sixteen identities are the following:
\begin{enumerate}
\item the five semiring axioms;
\item $xy\eq yx$ and $x^3\eq x$;
\item $R(3x)\eq R(3y)$ and $R(x+3y)\eq R(x)$;
\item the seven identities (P1)--(P5), with $u\circ v=u+3v$.
\end{enumerate}
The first seven identities are regular; the next two regularize the
nonregular identities in Lemma~\ref{f7:lem:field}. Completeness is
therefore exactly Theorem~\ref{f7:thm:regularization}. The notation
$\circ$ and $R$ merely abbreviates terms in $+$ and $\cdot$.
The displayed substitutions specify all sixteen identities in the original binary signature.
No assertion of minimal basis size is intended.
\end{proof}

The operation tables are included in Appendix~\ref{app:catalogue}.

\endgroup
\subsection{Cyclic zero-product joins}\label{sec:family9}
\begingroup
\providecommand{\eq}{}\renewcommand{\eq}{\approx}
\providecommand{\Var}{}\renewcommand{\Var}{\operatorname{Var}}
\providecommand{\Mod}{}\renewcommand{\Mod}{\operatorname{Mod}}
\providecommand{\Id}{}\renewcommand{\Id}{\operatorname{Id}}
\providecommand{\D}{}\renewcommand{\D}{\mathcal D}
\providecommand{\K}{}\renewcommand{\K}{\mathcal K}
\providecommand{\Z}{}\renewcommand{\Z}[1]{Z_{#1}^{\,0}}
\providecommand{\W}{}\renewcommand{\W}[1]{\mathcal W_{#1}}
\providecommand{\SR}{}\renewcommand{\SR}{\mathsf{SR}}
For $m\geq 2$, let $\Z{m}$ have carrier $\mathbb Z/m\mathbb Z$, addition
modulo $m$, and multiplication identically zero. It is a zero-product
ring, not the field of order $m$. The zero in this description is an
element of the algebra, not a symbol in our signature. Write $mx$ for
the sum of $m$ copies of $x$.

\subsubsection{A basis transformation for cyclic zero-product joins}
\begin{lemma}\label{f9:lem:normal}
Let $A$ satisfy $\SR$ and
\begin{equation}\label{f9:eq:ambient}
 (m+1)x\eq x,\qquad xy+xy\eq xy.
\end{equation}
Then $h(a)=ma$ is a semiring retraction onto an ai-subsemiring $h(A)$.
For every term $P$, expand $P$ distributively into a nonempty sum of
nonempty words. Let $L_P$ be the sum of its length-one monomials with
coefficients reduced modulo $m$, omitting zero coefficients. Then
\begin{equation}\label{f9:eq:normal}
 P\eq L_P+h(P),
\end{equation}
where an empty $L_P$ is omitted, and $h(P)$ denotes $mP$.
\end{lemma}
\begin{proof}
The first identity in~\eqref{f9:eq:ambient} gives $kx\eq(k+m)x$ for every
positive integer $k$. Consequently
\[
 h(h(x))\eq h(x),\qquad h(x)+h(x)\eq h(x),\qquad x+h(x)\eq x.
\]
Commutativity and associativity of addition give $h(x+y)\eq h(x)+h(y)$.
Since every product is additively idempotent, distributivity gives
\[
 h(xy)\eq xy\eq h(x)h(y).
\]
These equations establish the retraction claim and show that every
word of length at least two is fixed by $h$ and is additively idempotent.

Consider a length-one monomial $x$ with positive coefficient $c$ in
the expansion of $P$. If $r\in\{0,\ldots,m-1\}$ is its residue modulo
$m$, then $cx\eq rx+mx$ when $r>0$, and $cx\eq mx$ when $r=0$.
The contribution of that monomial to $h(P)$ is always $mx$.
Each nonlinear summand already agrees with its contribution to $h(P)$.
Adding these equalities proves~\eqref{f9:eq:normal}. Empty sums have not
been introduced as terms; omission is only a convention for this
description of a nonempty polynomial.
\end{proof}

\begin{theorem}\label{f9:thm:join}
Let $\K$ be any ai-semiring variety and let $\Sigma$ be any identity
basis for $\K$. For each term $P$, let $P^{[m]}$ be obtained by replacing
every variable $x$ by $mx$. Then
\begin{equation}\label{f9:eq:basis}
\begin{aligned}
\K\vee\Var(\Z{m})
 &=\Mod\bigl(\SR\cup\{(m+1)x\eq x,\ xy+xy\eq xy\}\\
 &\qquad\cup\{P^{[m]}\eq Q^{[m]}:(P\eq Q)\in\Sigma\}\bigr).
\end{aligned}
\end{equation}
In particular, $\K\vee\Var(\Z{m})$ is FB if and only if $\K$ is FB.
\end{theorem}
\begin{proof}
Every displayed identity holds in $\K$: there $mx\eq x$. It holds
in $\Z{m}$ as well, since all substituted variables evaluate to zero.
This proves one inclusion in~\eqref{f9:eq:basis}.

Conversely, suppose that $A$ satisfies the displayed identities.
By Lemma~\ref{f9:lem:normal}, $h$ is a retraction and the lifted identities
ensure $h(A)\in\K$. For every term, its evaluation on the tuple of
images agrees with the image of its evaluation; thus
$P^{[m]}\eq h(P)$ holds in $A$.

Take any identity $P\eq Q$ true in both $\K$ and $\Z{m}$.
In $\Z{m}$ every word of length at least two is zero. Therefore an
identity holds there precisely when the coefficients of every
length-one monomial agree modulo $m$. Necessity follows by assigning
one variable the value $1$ and all others $0$; sufficiency is immediate
from the evaluation formula. Choose the same variable order in $L_P$
and $L_Q$; then they are the same formal sum, or both are absent.
Since $h(A)\in\K$, we also have $h(P)\eq h(Q)$ in $A$.
Equation~\eqref{f9:eq:normal} now gives $P\eq Q$ in $A$. Thus $A$ belongs
to the join, proving equality and the forward finite-basis implication.

For the converse, intersect the join with the ai-semiring variety.
In an ai-semiring $h$ is the identity map, so~\eqref{f9:eq:basis} shows
that this intersection is exactly $\K$. A finite basis of the join,
together with $x+x\eq x$, is consequently a finite basis of $\K$.
\end{proof}

\begin{corollary}\label{f9:cor:equality}
Suppose $A$ satisfies~\eqref{f9:eq:ambient}, set $B=h(A)$, and suppose
$\Z{m}\in\Var(A)$. Then
\[
 \Var(A)=\Var(B)\vee\Var(\Z{m}).
\]
Hence $A$ is FB if and only if $B$ is FB.
\end{corollary}
\begin{proof}
Apply Theorem~\ref{f9:thm:join} with $\Sigma=\Id(B)$. The retraction
ensures that $A$ satisfies all lifted identities, giving containment
in the displayed join. The reverse containment follows from
$B\leq A$ and the assumed membership of $\Z{m}$.
\end{proof}
The equality in this corollary is essential. Containment in an
arbitrary FB variety alone would not justify an FB classification.

\subsubsection{Applications and explicit maps}
Let $D_2$ be the two-element distributive lattice with addition join
and multiplication meet. Put $\D=\Var(D_2)$. The following table records
each equality $\Var(A_j)=\Var(B)\vee\Var(\Z{m})$ proved by
Corollary~\ref{f9:cor:equality}.
\begin{center}
\begin{tabular}{crl}\toprule
$B$ & $m$ & Local four-element indices $j$\\\midrule
$D_2$ & 3 & 165, 1114 \\
$B_{40}$ & 2 & 226, 245, 871 \\
$B_{53}$ & 2 & 284, 976 \\
$B_{54}$ & 2 & 280, 391, 984 \\
$B_{76}$ & 2 & 542, 1110 \\
$B_{80}$ & 2 & 600, 1752 \\
$B_{82}$ & 2 & 420, 1616 \\
$B_{83}$ & 2 & 416, 1334, 1620 \\
$B_{88}$ & 2 & 441, 1352, 1764 \\
$B_{89}$ & 2 & 447, 1762 \\
$B_{118}$ & 2 & 732, 1217, 1784 \\
$B_{121}$ & 2 & 1182, 1184, 1814 \\
$B_{122}$ & 2 & 792, 794, 1820 \\
$B_{125}$ & 2 & 777, 1710, 1714 \\
$B_{130}$ & 2 & 805, 1279, 1283 \\
\bottomrule\end{tabular}
\end{center}
All fourteen three-element inputs are ai-semirings and are not
isomorphic to $S_7$; their FB status follows from
\cite{ZhaoThree}, as recalled in Section~\ref{sec:small-results}. The other input is $D_2$, whose variety has the usual
finite distributive-lattice basis. Thus all 39 classes in the table
are FB. The displayed equalities group the applications by their
chosen small input; no claim that all fifteen generated varieties
are pairwise distinct is needed.

For every row, the displayed maps give $h(x)=mx$, the image
$h(A_j)$, and an isomorphism from the indicated $B$ onto this image.
The tables also satisfy~\eqref{f9:eq:ambient}. For 34 rows the maps give an embedding
of $\Z{m}$ into $A_j$. The remaining five rows use the following
explicit homomorphisms. Write $F_j$ for the indicated
subsemiring of the power $A_j^{A_j}$ consisting of unary term functions:
\begin{center}
\begin{tabular}{rlc}\toprule
$j$ & \begin{tabular}[c]{@{}l@{}}Distinct functions generating\\the whole displayed $F_j$\end{tabular} & Images in $\Z{2}$\\\midrule
245, 871 & $x,\ x^2,\ x^3,\ 2x$ & $1,0,0,0$\\
391, 1334, 1352 & $x,\ x^2,\ 2x$ & $1,0,0$\\\bottomrule
\end{tabular}
\end{center}
Pointwise use of the operation tables shows that each displayed
set is closed under both operations and that its image assignment
is a surjective homomorphism. Hence $\Z{2}\in\operatorname{HSP}(A_j)$ in these cases.
The table entries can be reconstructed directly from the four-element
tables in Subsection~\ref{f9:app:tables}.

For each application, substituting a finite basis of the indicated
source algebra into \eqref{f9:eq:basis} gives a finite basis of $A_j$.

\subsubsection{Source tables and embeddings}\label{f9:app:tables}
All strings are row-major tables on $\{0,\ldots,n-1\}$. Digit $0$ is
an element label. The following small representatives identify every
input used in this cyclic-join construction.
\begin{center}\small
\begin{tabular}{cll}\toprule
Input & Addition & Multiplication\\\midrule
$D_2$ & \texttt{0111} & \texttt{0001} \\
$B_{40}$ & \texttt{012111212} & \texttt{000020000} \\
$B_{53}$ & \texttt{000010002} & \texttt{000012000} \\
$B_{54}$ & \texttt{012111212} & \texttt{000012000} \\
$B_{76}$ & \texttt{000012022} & \texttt{010111010} \\
$B_{80}$ & \texttt{011111112} & \texttt{012111210} \\
$B_{82}$ & \texttt{000010002} & \texttt{000010020} \\
$B_{83}$ & \texttt{012111212} & \texttt{000010020} \\
$B_{88}$ & \texttt{012111212} & \texttt{000012020} \\
$B_{89}$ & \texttt{000012022} & \texttt{000012020} \\
$B_{118}$ & \texttt{012111212} & \texttt{000011012} \\
$B_{121}$ & \texttt{012111212} & \texttt{010011012} \\
$B_{122}$ & \texttt{012111212} & \texttt{000111012} \\
$B_{125}$ & \texttt{000012022} & \texttt{012111012} \\
$B_{130}$ & \texttt{000012022} & \texttt{010111212} \\
\bottomrule\end{tabular}
\end{center}
In the next table, $E(a_0,\ldots,a_{m-1})$ denotes the embedding
$i\mapsto a_i$ of $\Z{m}$; $H_s$ denotes the $s$-element subsemiring
of unary term functions and its homomorphism displayed above.
The term $h(x)=mx$ specifies its image on each carrier; comparison
of the two induced operations gives the stated small isomorphism type.

\begingroup\small\setlength{\tabcolsep}{2pt}
\begin{longtable}{@{}llll@{\hspace{2em}}llll@{}}
\caption{The 39 additional FB representatives and reverse-inclusion witnesses.}\\
\toprule $j$ & $h(A_j)\cong$ & $m$ & Witness & $j$ & $h(A_j)\cong$ & $m$ & Witness\\
\midrule\endfirsthead
\toprule $j$ & $h(A_j)\cong$ & $m$ & Witness & $j$ & $h(A_j)\cong$ & $m$ & Witness\\
\midrule\endhead
\bottomrule\endlastfoot
165 & $D_2$ & 3 & $E(0,2,3)$ & 226 & $B_{40}$ & 2 & $E(0,3)$\\
245 & $B_{40}$ & 2 & $H_{4}$ & 280 & $B_{54}$ & 2 & $E(0,3)$\\
284 & $B_{53}$ & 2 & $E(0,3)$ & 391 & $B_{54}$ & 2 & $H_{3}$\\
416 & $B_{83}$ & 2 & $E(0,3)$ & 420 & $B_{82}$ & 2 & $E(0,3)$\\
441 & $B_{88}$ & 2 & $E(0,3)$ & 447 & $B_{89}$ & 2 & $E(0,3)$\\
542 & $B_{76}$ & 2 & $E(0,3)$ & 600 & $B_{80}$ & 2 & $E(0,3)$\\
732 & $B_{118}$ & 2 & $E(0,3)$ & 777 & $B_{125}$ & 2 & $E(0,3)$\\
792 & $B_{122}$ & 2 & $E(0,3)$ & 794 & $B_{122}$ & 2 & $E(0,3)$\\
805 & $B_{130}$ & 2 & $E(0,3)$ & 871 & $B_{40}$ & 2 & $H_{4}$\\
976 & $B_{53}$ & 2 & $E(1,2)$ & 984 & $B_{54}$ & 2 & $E(1,2)$\\
1110 & $B_{76}$ & 2 & $E(0,3)$ & 1114 & $D_2$ & 3 & $E(0,2,3)$\\
1182 & $B_{121}$ & 2 & $E(0,3)$ & 1184 & $B_{121}$ & 2 & $E(0,3)$\\
1217 & $B_{118}$ & 2 & $E(0,3)$ & 1279 & $B_{130}$ & 2 & $E(0,3)$\\
1283 & $B_{130}$ & 2 & $E(0,3)$ & 1334 & $B_{83}$ & 2 & $H_{3}$\\
1352 & $B_{88}$ & 2 & $H_{3}$ & 1616 & $B_{82}$ & 2 & $E(0,3)$\\
1620 & $B_{83}$ & 2 & $E(0,3)$ & 1710 & $B_{125}$ & 2 & $E(0,3)$\\
1714 & $B_{125}$ & 2 & $E(0,3)$ & 1752 & $B_{80}$ & 2 & $E(0,3)$\\
1762 & $B_{89}$ & 2 & $E(0,3)$ & 1764 & $B_{88}$ & 2 & $E(0,3)$\\
1784 & $B_{118}$ & 2 & $E(0,3)$ & 1814 & $B_{121}$ & 2 & $E(0,3)$\\
1820 & $B_{122}$ & 2 & $E(0,3)$ &  &  &  & \\
\end{longtable}\endgroup

\endgroup
\subsection{Joins with an absorbing operation}\label{sec:family10}
\begingroup
\providecommand{\eq}{}\renewcommand{\eq}{\approx}
\providecommand{\Var}{}\renewcommand{\Var}{\operatorname{Var}}
\providecommand{\Id}{}\renewcommand{\Id}{\operatorname{Id}}
\providecommand{\K}{}\renewcommand{\K}{\mathcal K}
\providecommand{\SR}{}\renewcommand{\SR}{\mathsf{SR}}
\providecommand{\SL}{}\renewcommand{\SL}{\mathsf{SL}}
Consider the two-element semirings $M$ and $P$, on $\{0,1\}$:
\[
\begin{array}{c|cc}
 & + & \cdot\\\hline
M & \min & \text{constant }0\\
P & \text{constant }0 & \min
\end{array}
\]
Their local representatives are $E_3$ and $E_5$. In $M$, any term
containing a multiplication symbol evaluates constantly to $0$;
terms containing only addition evaluate as the meet of their variables.
In $P$ the analogous statement holds with the operations interchanged.
Thus the two syntactic classes of terms are separated by assigning
all variables the value $1$. These observations concern evaluations
of terms, and do not introduce a constant into the signature.

For a choice of $N\in\{M,P\}$, call its constant operation $\sigma$
and its semilattice operation $\tau$. A \emph{pure} term uses only
$\tau$; a \emph{marked} term contains $\sigma$ at least once. Equality
of two pure terms holds in $N$ exactly when they have the same set
of variables. Every two marked terms agree in $N$.

\subsubsection{A complete finite-basis construction}
\begin{theorem}\label{f10:thm:construction}
Let $\K$ be a semiring variety with finite basis $\Sigma$.
Suppose there is a unary marked term $h(x)$ such that $\K$ satisfies
$h(x)\eq x$. Suppose also that the pure identities common to $\K$
and $\SL$ have a finite semigroup basis $\Theta$, in the operation
$\tau$. Then $\K\vee\Var(N)$ is finitely based.

More precisely, a basis consists of $\SR$, $\Theta$, the identities
obtained from $\Sigma$ by replacing each variable $x$ with $h(x)$,
and the following structural identities:
\begin{align}
h(h(x))&\eq h(x),\label{f10:eq:retract}\\
h(x\mathbin\circ y)&\eq h(x)\mathbin\circ h(y)
                    &&(\circ\in\{+,\cdot\}),\label{f10:eq:hom}\\
h(\sigma(x,y))&\eq\sigma(x,y),\label{f10:eq:marked}\\
\tau(x,h(y))&\eq\tau(h(x),h(y)),\qquad
\tau(h(x),y)\eq\tau(h(x),h(y)).\label{f10:eq:absorb}
\end{align}
\end{theorem}
\begin{proof}
All displayed identities hold in $\K$, since $h$ is the identity
there. In $N$, $h$ is the constant zero map. Its value is absorbing
for both operations, so the identities hold in $N$ as well. The
lifted members of $\Sigma$ hold in $N$ because every substituted
variable is zero. The members of $\Theta$ hold in both generators
by hypothesis. This proves soundness of the proposed basis.

Conversely, let $A$ satisfy that basis. Equations
\eqref{f10:eq:retract}--\eqref{f10:eq:hom} make $h$ an endomorphic retraction,
and the lifted basis $\Sigma$ gives $h(A)\in\K$. For any term $T$,
write $T^h$ for the simultaneous substitution of $h(x)$ for each
variable $x$. Then $h(T)\eq T^h$ follows from~\eqref{f10:eq:hom}.
We claim that every marked term satisfies $T\eq h(T)$ in $A$.
If its root is $\sigma$, equation~\eqref{f10:eq:marked} proves this.
Otherwise its root is $\tau$, and at least one child is marked.
Apply induction to that child and then~\eqref{f10:eq:absorb}; finally
use~\eqref{f10:eq:hom}. This is induction on term structure, with no
restriction on depth, number of variables, or length.

Now take an arbitrary identity $U\eq V$ of the join. It holds in
$N$, so either both terms are marked, or both are pure with the
same variable set. In the first case,
\[
 U\eq h(U)=U^h\eq V^h=h(V)\eq V,
\]
because $h(A)\in\K$. In the second case, $U\eq V$ is a pure
identity common to $\K$ and $\SL$, hence follows from $\Theta$.
This proves completeness of the finite basis.
\end{proof}

For a finite source $B$, the pure identity theory in the theorem is
the theory of the semigroup direct product $B_\tau\times\SL_2$.
If $\SL_2\in\Var(B_\tau)$, this is just $\Id(B_\tau)$.
In our applications $|B_\tau|=3$, so the required finite basis follows
from the finite-basis theorem for semigroups of order at most four
\cite{Edmunds}; see also~\cite{Trahtman}. We verify the containment
by finite separating semigroup homomorphisms. For the few cases
without that containment, explicit bases are given next.

\begin{lemma}\label{f10:lem:pure}
Let $S$ be one of the following semigroups, using multiplicative
notation solely within this statement. Each listed pair is a basis
for $\Var(S)\vee\SL$, relative to associativity:
\begin{center}
\begin{tabular}{ll}\toprule
$S$ & Basis\\\midrule
A nontrivial constant-product semigroup & $xy\eq yx,\ x^2y\eq xy$\\
$\{0,a,a^2\}$, with $a^3=0$ & $xy\eq yx,\ xyz\eq x^2yz$\\
A cyclic group of order three & $xy\eq yx,\ x^4\eq x$\\\bottomrule
\end{tabular}
\end{center}
\end{lemma}
\begin{proof}
In every case, $\SL_2$ requires the two words to have the same
support (set of occurring variables). In a constant-product
semigroup, every word of length at least two has the same value;
a one-letter word is separated from these by a nonconstant element.
The first pair of identities identifies all nonlinear words with
the same support, by inserting or deleting repetitions, and identifies
no additional words.

For the second semigroup, all words of length at least three are
zero. Distinct commutative words of length at most two are separated
by assigning the value $a$ to a suitable variable and zero to an
extra variable, or by using $a\ne a^2\ne0$ for singleton supports.
They are also separated from all longer words. The identity
$xyz\eq x^2yz$ permits duplication of any chosen letter in a word
of length at least three. Two such words with the same support
can each be enlarged to their componentwise maximum exponent vector.
Thus the listed pair is complete.

For the group of order three, identities are equality of exponent
vectors modulo three. Together with equality of support this is
exactly the equivalence induced by commutativity and $x^4\eq x$:
every positive exponent reduces to one of $1,2,3$. This proves the
last assertion.
\end{proof}

\subsubsection{Applications and explicit maps}
The following table gives $B$, the two-element seed, the marked unary
identity term in $B$, and the local target indices. For each target
the asserted equality is
\begin{equation}\label{f10:eq:equality}
 \Var(A_j)=\Var(B)\vee\Var(N).
\end{equation}
\begin{center}\small
\begin{tabular}{cccl}\toprule
$B$ & $N$ & $h(x)$ & Four-element indices $j$\\\midrule
$B_{10}$ & $P$ & $4x$ & 166 \\
$B_{13}$ & $P$ & $x^2+x$ & 179 \\
$B_{34}$ & $P$ & $2x$ & 633, 764 \\
$B_{39}$ & $P$ & $2x$ & 655 \\
$B_{40}$ & $P$ & $2x$ & 664 \\
$B_{41}$ & $P$ & $2x$ & 657 \\
$B_{53}$ & $P$ & $2x$ & 705, 993 \\
$B_{54}$ & $P$ & $2x$ & 713, 997 \\
$B_{55}$ & $P$ & $2x$ & 707, 1003 \\
$B_{76}$ & $P$ & $2x$ & 1226, 1609 \\
$B_{78}$ & $M$ & $x^3$ & 604 \\
$B_{78}$ & $P$ & $4x$ & 1659, 1896 \\
$B_{80}$ & $P$ & $2x$ & 1665, 1898 \\
$B_{82}$ & $P$ & $2x$ & 1417, 1903 \\
$B_{83}$ & $P$ & $2x$ & 1425, 1907 \\
$B_{84}$ & $P$ & $2x$ & 1419, 1913 \\
$B_{88}$ & $P$ & $2x$ & 1443, 1939 \\
$B_{89}$ & $P$ & $2x$ & 1437, 1947 \\
$B_{113}$ & $M$ & $x^2$ & 755, 1221 \\
$B_{118}$ & $P$ & $2x$ & 1991, 2109, 2254 \\
$B_{121}$ & $P$ & $2x$ & 2100, 2101, 2280 \\
$B_{122}$ & $P$ & $2x$ & 2025, 2026, 2286 \\
$B_{125}$ & $P$ & $2x$ & 2017, 2186, 2196 \\
$B_{130}$ & $P$ & $2x$ & 2033, 2141, 2144 \\
\bottomrule\end{tabular}
\end{center}
Every source has order three and is FB. The ai sources use the
three-element classification of \cite{ZhaoThree}; none is $S_7$.
The sources $B_{10}$ and $B_{78}$ are, respectively, the cyclic
zero-product ring of order three and the field of order three,
covered by the finite-ring identity theorem~\cite{Kruse} in the
constant-free signature by Proposition~\ref{prop:ring}.
The finite basability of $B_{13}$ and $B_{113}$ follows, respectively, from
the commutative-semigroup expansion proof and the explicit basis
in Sections~\ref{sec:family3} and~\ref{sec:family4}.

In forty-one applications one has
$\SL_2\in\Var(B_\tau)$. The remaining six use
Lemma~\ref{f10:lem:pure}: $A_{166},A_{179}$ use constant-product pure
reducts; $A_{655},A_{664},A_{657}$ use the three-element nilpotent
semigroup; and $A_{604}$ uses the additive cyclic group of order three.
The indicated unary term $h$ satisfies $h(b)=b$ for every element of its source algebra.
Thus Theorem~\ref{f10:thm:construction} proves the right-hand side of
\eqref{f10:eq:equality} is FB in all 47 applications.

For each equality~\eqref{f10:eq:equality}, exactly two maps are needed.
The maps $q_B:A_j\twoheadrightarrow B$ and
$q_N:A_j\twoheadrightarrow N$ are homomorphisms, and the map
$(q_B,q_N)$ is injective. Consequently $A_j$ embeds into $B\times N$,
giving one containment in~\eqref{f10:eq:equality}. The surjectivity of
both maps gives the reverse containment. Subsection~\ref{f10:app:tables}
specifies every map and target table. The asserted finite bases
then follow from Theorem~\ref{f10:thm:construction}.

For example, the maps for $A_{604}$ are
\[
 q_{B_{78}}=(1,0,2,1),\qquad q_M=(0,0,0,1).
\]
Their four ordered pairs are distinct, and both maps are onto.
The source $B_{78}\cong\mathbb F_3$ satisfies the marked identity
$x^3\eq x$ for type $M$. This proves
$\Var(A_{604})=\Var(\mathbb F_3)\vee\Var(M)$ and provides a finite
basis through Theorem~\ref{f10:thm:construction}.

\subsubsection{Source tables and quotient maps}\label{f10:app:tables}
Tables are row-major strings on $\{0,\ldots,n-1\}$. The notation
$q_B/q_N$ records the four values of each map in argument order
$0,1,2,3$. The digit $0$ is an element label, not a named constant.
The two-element target tables are $M:(0001,0000)$ and
$P:(0000,0001)$, listed as (addition, multiplication).
\begin{center}\small
\begin{tabular}{cll}\toprule
Source & Addition & Multiplication\\\midrule
$B_{10}$ & \texttt{012120201} & \texttt{000000000} \\
$B_{13}$ & \texttt{012111211} & \texttt{000000000} \\
$B_{34}$ & \texttt{000012022} & \texttt{000010000} \\
$B_{39}$ & \texttt{000010002} & \texttt{000020000} \\
$B_{40}$ & \texttt{012111212} & \texttt{000020000} \\
$B_{41}$ & \texttt{000012022} & \texttt{000020000} \\
$B_{53}$ & \texttt{000010002} & \texttt{000012000} \\
$B_{54}$ & \texttt{012111212} & \texttt{000012000} \\
$B_{55}$ & \texttt{000012022} & \texttt{000012000} \\
$B_{76}$ & \texttt{000012022} & \texttt{010111010} \\
$B_{78}$ & \texttt{201012120} & \texttt{012111210} \\
$B_{80}$ & \texttt{011111112} & \texttt{012111210} \\
$B_{82}$ & \texttt{000010002} & \texttt{000010020} \\
$B_{83}$ & \texttt{012111212} & \texttt{000010020} \\
$B_{84}$ & \texttt{000012022} & \texttt{000010020} \\
$B_{88}$ & \texttt{012111212} & \texttt{000012020} \\
$B_{89}$ & \texttt{000012022} & \texttt{000012020} \\
$B_{113}$ & \texttt{012111211} & \texttt{000011012} \\
$B_{118}$ & \texttt{012111212} & \texttt{000011012} \\
$B_{121}$ & \texttt{012111212} & \texttt{010011012} \\
$B_{122}$ & \texttt{012111212} & \texttt{000111012} \\
$B_{125}$ & \texttt{000012022} & \texttt{012111012} \\
$B_{130}$ & \texttt{000012022} & \texttt{010111212} \\
\bottomrule\end{tabular}
\end{center}
\begingroup\small\setlength{\tabcolsep}{3.5pt}
\begin{longtable}{@{}llll@{\hspace{2em}}llll@{}}
\caption{The 47 classes and their two surjective quotient maps.}\\
\toprule $j$ & $B$ & $N$ & $q_B/q_N$ & $j$ & $B$ & $N$ & $q_B/q_N$\\
\midrule\endfirsthead
\toprule $j$ & $B$ & $N$ & $q_B/q_N$ & $j$ & $B$ & $N$ & $q_B/q_N$\\
\midrule\endhead
\bottomrule\endlastfoot
166 & $B_{10}$ & $P$ & \texttt{0012/0100} & 179 & $B_{13}$ & $P$ & \texttt{0012/0100}\\
604 & $B_{78}$ & $M$ & \texttt{1021/0001} & 633 & $B_{34}$ & $P$ & \texttt{0102/0010}\\
655 & $B_{39}$ & $P$ & \texttt{0102/0010} & 657 & $B_{41}$ & $P$ & \texttt{0102/0010}\\
664 & $B_{40}$ & $P$ & \texttt{0102/0010} & 705 & $B_{53}$ & $P$ & \texttt{0102/0010}\\
707 & $B_{55}$ & $P$ & \texttt{0102/0010} & 713 & $B_{54}$ & $P$ & \texttt{0102/0010}\\
755 & $B_{113}$ & $M$ & \texttt{0120/0001} & 764 & $B_{34}$ & $P$ & \texttt{0112/0010}\\
993 & $B_{53}$ & $P$ & \texttt{0112/0010} & 997 & $B_{54}$ & $P$ & \texttt{0112/0010}\\
1003 & $B_{55}$ & $P$ & \texttt{0112/0010} & 1221 & $B_{113}$ & $M$ & \texttt{1021/0001}\\
1226 & $B_{76}$ & $P$ & \texttt{0102/0010} & 1417 & $B_{82}$ & $P$ & \texttt{0102/0010}\\
1419 & $B_{84}$ & $P$ & \texttt{0102/0010} & 1425 & $B_{83}$ & $P$ & \texttt{0102/0010}\\
1437 & $B_{89}$ & $P$ & \texttt{0102/0010} & 1443 & $B_{88}$ & $P$ & \texttt{0102/0010}\\
1609 & $B_{76}$ & $P$ & \texttt{0112/0010} & 1659 & $B_{78}$ & $P$ & \texttt{0112/0010}\\
1665 & $B_{80}$ & $P$ & \texttt{0112/0010} & 1896 & $B_{78}$ & $P$ & \texttt{0012/0100}\\
1898 & $B_{80}$ & $P$ & \texttt{0012/0100} & 1903 & $B_{82}$ & $P$ & \texttt{0112/0010}\\
1907 & $B_{83}$ & $P$ & \texttt{0112/0010} & 1913 & $B_{84}$ & $P$ & \texttt{0112/0010}\\
1939 & $B_{88}$ & $P$ & \texttt{0112/0010} & 1947 & $B_{89}$ & $P$ & \texttt{0112/0010}\\
1991 & $B_{118}$ & $P$ & \texttt{0120/0001} & 2017 & $B_{125}$ & $P$ & \texttt{1021/0001}\\
2025 & $B_{122}$ & $P$ & \texttt{1021/0001} & 2026 & $B_{122}$ & $P$ & \texttt{0120/0001}\\
2033 & $B_{130}$ & $P$ & \texttt{1021/0001} & 2100 & $B_{121}$ & $P$ & \texttt{1021/0001}\\
2101 & $B_{121}$ & $P$ & \texttt{0120/0001} & 2109 & $B_{118}$ & $P$ & \texttt{1021/0001}\\
2141 & $B_{130}$ & $P$ & \texttt{0120/0001} & 2144 & $B_{130}$ & $P$ & \texttt{2102/0001}\\
2186 & $B_{125}$ & $P$ & \texttt{0102/0010} & 2196 & $B_{125}$ & $P$ & \texttt{2120/0010}\\
2254 & $B_{118}$ & $P$ & \texttt{0122/0001} & 2280 & $B_{121}$ & $P$ & \texttt{0122/0001}\\
2286 & $B_{122}$ & $P$ & \texttt{0122/0001} &  &  &  & \\
\end{longtable}\endgroup

\endgroup
\subsection{Ideal retractions and pure words}\label{sec:family11}
\begingroup
\providecommand{\eq}{}\renewcommand{\eq}{\approx}
\providecommand{\Var}{}\renewcommand{\Var}{\operatorname{Var}}
\providecommand{\SR}{}\renewcommand{\SR}{\mathsf{SR}}
\subsubsection{The ideal-retraction criterion}
Select one basic operation $\sigma\in\{+,\cdot\}$ and let $\tau$
be the other. A term is \emph{pure} if it uses only $\tau$, and
\emph{marked} if it contains $\sigma$ at least once. Pure terms are
semigroup words, with $\tau$ as their operation.

\Needspace{4\baselineskip}
\begin{theorem}\label{f11:thm:main}
Let $A$ be a semiring and $h(x)$ a unary semiring term. Assume that:
\begin{enumerate}
\item $h$ is an endomorphic retraction onto $I=h(A)$;
\item $I$ is an ideal for both operations, and $\sigma(A,A)\subseteq I$;
\item some $e\in A\setminus I$ satisfies $\tau(e,e)=e$;
\item the semiring $I$ and the semigroup reduct $A_\tau=(A,\tau)$
are both FB.
\end{enumerate}
Then $A$ is FB. If $\Sigma$ is a finite basis for $I$ and $\Theta$
is a finite semigroup basis for $A_\tau$, a basis for $A$ consists
of $\SR$, $\Theta$, the simultaneous substitutions
\[
 \{U(h(x_1),\ldots,h(x_k))\eq V(h(x_1),\ldots,h(x_k)):
                  U\eq V\in\Sigma\},
\]
and the structural identities
\begin{align}
h(h(x))&\eq h(x),\label{f11:eq:idempotent}\\
h(x\mathbin\circ y)&\eq h(x)\mathbin\circ h(y)
                         &&(\circ\in\{+,\cdot\}),\label{f11:eq:endomorphism}\\
h(\sigma(x,y))&\eq\sigma(x,y),\label{f11:eq:sigma}\\
\tau(x,h(y))&\eq\tau(h(x),h(y)),\qquad
\tau(h(x),y)\eq\tau(h(x),h(y)).\label{f11:eq:ideal}
\end{align}
\end{theorem}
\begin{proof}
All these identities hold in $A$. For~\eqref{f11:eq:ideal}, note that
$\tau(x,h(y))$ belongs to $I$, so it equals its image under $h$;
endomorphism and idempotence of $h$ give the right-hand side.
The other equality is obtained in the same way. The lifted basis
$\Sigma$ is valid because all its variables take values in $I$.
This proves soundness.

Let $B$ satisfy the proposed finite basis. Equations
\eqref{f11:eq:idempotent}--\eqref{f11:eq:endomorphism} make $h$ an
endomorphic retraction in $B$, and the lifted basis puts $h(B)$ in
$\Var(I)$. Write $T^h$ for the term obtained from $T$ by replacing
each variable $x$ with $h(x)$. Endomorphism gives $h(T)\eq T^h$.

Every marked term $T$ satisfies $T\eq h(T)$ in $B$. This follows
by induction on its syntax. A root labelled $\sigma$ is covered
by~\eqref{f11:eq:sigma}. For a root labelled $\tau$, at least one child
is marked. Replace that child by its $h$-image using induction,
apply~\eqref{f11:eq:ideal}, and finish with~\eqref{f11:eq:endomorphism}.

Take an arbitrary identity $U\eq V$ of $A$. It cannot have exactly
one pure side: assign every variable the element $e$ from assumption
(3). The pure side evaluates to $e$, while the marked side evaluates
in $I$. The latter assertion follows directly from assumption (2),
starting at any occurrence of $\sigma$ and using the ideal property
along the path to the root. Since $e\notin I$, the values differ.

If both $U,V$ are pure, their equality is an identity of $A_\tau$,
so it follows from $\Theta$ and holds in $B$. If both are marked,
the identity is valid in the subsemiring $I$, and therefore in $h(B)$.
Consequently
\[
 U\eq h(U)=U^h\eq V^h=h(V)\eq V
\]
holds in $B$. These two cases exhaust the identities of $A$.
Hence every model of the proposed basis belongs to $\Var(A)$,
proving completeness.
\end{proof}

\begin{corollary}\label{f11:cor:small}
If $|A|\leq4$, assumptions (1)--(3) hold, and $I$ is FB, then $A$
is FB.
\end{corollary}
\begin{proof}
Every semigroup of order at most four is FB by the small-semigroup
finite-basis theorem~\cite{Edmunds}; see also~\cite{Trahtman}. Apply
this to $A_\tau$ and then use Theorem~\ref{f11:thm:main}.
\end{proof}

The element $e$ is a witness used in the completeness proof; it is
not named in the basis. The theorem does not assert that a finitely
based retract alone implies finite basability of its extension.
The operation-image condition, ideal condition, and separation of
the two kinds of terms are substantive hypotheses.

For example, $A_{158}$ has $h(x)=x^2$ with image map $(0,1,0,0)$.
Its image is the field of two elements, every product lies in that
image, and element $2$ is additively idempotent outside the image.
The theorem combines a lifted basis of $\mathbb F_2$ with a finite
basis of the four-element additive semigroup.

\subsubsection{Applications and explicit maps}\label{f11:app:tables}
Table strings are row-major on $\{0,\ldots,n-1\}$; digit $0$ is only
an element label. The source tables below identify the retracts.
Restricting both operations to the displayed image of $h$ identifies
the retract with the indicated source algebra.
\begin{center}\small
\begin{tabular}{cll}\toprule
Source & Addition & Multiplication\\\midrule
$E_{2}$ & \texttt{0110} & \texttt{0000} \\
$E_{4}$ & \texttt{0111} & \texttt{0000} \\
$E_{6}$ & \texttt{0110} & \texttt{0001} \\
$E_{7}$ & \texttt{0001} & \texttt{0001} \\
$E_{8}$ & \texttt{0111} & \texttt{0001} \\
$E_{9}$ & \texttt{0001} & \texttt{0101} \\
$E_{10}$ & \texttt{0001} & \texttt{0011} \\
$B_{50}$ & \texttt{000012022} & \texttt{000111000} \\
$B_{64}$ & \texttt{000012022} & \texttt{010010010} \\
\bottomrule\end{tabular}
\end{center}
In the main table, $h$ is encoded by $h(0)h(1)h(2)h(3)$ and $e$
is idempotent for $\tau$, the operation other than $\sigma$.
\begingroup\footnotesize\setlength{\tabcolsep}{3pt}
\begin{longtable}{@{}llllll@{\hspace{2em}}llllll@{}}
\caption{The 63 applications, ideal retractions, and outside idempotents.}\\
\toprule $j$ & $h(A_j)\cong$ & $\sigma$ & $h(x)$ & $h$ & $e$ & $j$ & $h(A_j)\cong$ & $\sigma$ & $h(x)$ & $h$ & $e$\\
\midrule\endfirsthead
\toprule $j$ & $h(A_j)\cong$ & $\sigma$ & $h(x)$ & $h$ & $e$ & $j$ & $h(A_j)\cong$ & $\sigma$ & $h(x)$ & $h$ & $e$\\
\midrule\endhead
\bottomrule\endlastfoot
158 & $E_{6}$ & $\cdot$ & $x^2$ & \texttt{0100} & 2 & 159 & $E_{8}$ & $\cdot$ & $x^2$ & \texttt{0100} & 2\\
283 & $E_{2}$ & $+$ & $x+x^2$ & \texttt{0003} & 1 & 295 & $E_{4}$ & $+$ & $2x$ & \texttt{0003} & 1\\
419 & $E_{2}$ & $+$ & $x+x^2$ & \texttt{0003} & 1 & 431 & $E_{4}$ & $+$ & $2x$ & \texttt{0003} & 1\\
444 & $E_{2}$ & $+$ & $x+x^2$ & \texttt{0003} & 1 & 458 & $E_{4}$ & $+$ & $2x$ & \texttt{0003} & 1\\
567 & $E_{6}$ & $\cdot$ & $x^2$ & \texttt{0110} & 3 & 568 & $E_{8}$ & $\cdot$ & $x^2$ & \texttt{0110} & 3\\
599 & $E_{2}$ & $+$ & $x+x^2$ & \texttt{0003} & 1 & 606 & $E_{4}$ & $+$ & $2x$ & \texttt{0003} & 1\\
688 & $B_{50}$ & $+$ & $2x$ & \texttt{0103} & 2 & 702 & $E_{7}$ & $+$ & $2x$ & \texttt{0020} & 1\\
703 & $E_{6}$ & $+$ & $3x$ & \texttt{0020} & 1 & 704 & $E_{8}$ & $+$ & $2x$ & \texttt{0020} & 1\\
818 & $E_{7}$ & $+$ & $2x$ & \texttt{0220} & 1 & 820 & $E_{6}$ & $+$ & $3x$ & \texttt{0220} & 1\\
821 & $E_{8}$ & $+$ & $2x$ & \texttt{0220} & 1 & 932 & $B_{50}$ & $+$ & $2x$ & \texttt{0113} & 2\\
964 & $E_{10}$ & $+$ & $2x$ & \texttt{0020} & 1 & 965 & $E_{10}$ & $+$ & $2x$ & \texttt{0020} & 1\\
1105 & $E_{8}$ & $\cdot$ & $x^2$ & \texttt{0100} & 2 & 1166 & $B_{64}$ & $+$ & $2x$ & \texttt{0103} & 2\\
1297 & $E_{9}$ & $+$ & $2x$ & \texttt{0100} & 2 & 1298 & $E_{9}$ & $+$ & $2x$ & \texttt{0100} & 2\\
1308 & $E_{8}$ & $+$ & $2x$ & \texttt{0100} & 2 & 1309 & $E_{7}$ & $+$ & $2x$ & \texttt{0100} & 2\\
1311 & $E_{6}$ & $+$ & $3x$ & \texttt{0100} & 2 & 1414 & $E_{7}$ & $+$ & $2x$ & \texttt{0020} & 1\\
1415 & $E_{6}$ & $+$ & $3x$ & \texttt{0020} & 1 & 1416 & $E_{8}$ & $+$ & $2x$ & \texttt{0020} & 1\\
1431 & $E_{6}$ & $+$ & $3x$ & \texttt{0020} & 1 & 1433 & $E_{8}$ & $+$ & $2x$ & \texttt{0020} & 1\\
1446 & $E_{7}$ & $+$ & $2x$ & \texttt{0220} & 1 & 1448 & $E_{6}$ & $+$ & $3x$ & \texttt{0220} & 1\\
1449 & $E_{8}$ & $+$ & $2x$ & \texttt{0220} & 1 & 1462 & $E_{6}$ & $+$ & $3x$ & \texttt{0220} & 1\\
1464 & $E_{8}$ & $+$ & $2x$ & \texttt{0220} & 1 & 1479 & $E_{10}$ & $+$ & $2x$ & \texttt{0020} & 1\\
1480 & $E_{10}$ & $+$ & $2x$ & \texttt{0020} & 1 & 1490 & $E_{10}$ & $+$ & $2x$ & \texttt{0220} & 1\\
1492 & $E_{10}$ & $+$ & $2x$ & \texttt{0220} & 1 & 1499 & $E_{10}$ & $+$ & $2x$ & \texttt{0020} & 1\\
1501 & $E_{10}$ & $+$ & $2x$ & \texttt{0020} & 1 & 1572 & $B_{64}$ & $+$ & $2x$ & \texttt{0113} & 2\\
1650 & $E_{9}$ & $+$ & $2x$ & \texttt{0110} & 2 & 1651 & $E_{9}$ & $+$ & $2x$ & \texttt{0110} & 2\\
1663 & $E_{6}$ & $+$ & $3x$ & \texttt{1121} & 0 & 1671 & $E_{10}$ & $+$ & $2x$ & \texttt{1121} & 0\\
1672 & $E_{10}$ & $+$ & $2x$ & \texttt{1121} & 0 & 1825 & $E_{9}$ & $+$ & $2x$ & \texttt{1121} & 0\\
1826 & $E_{9}$ & $+$ & $2x$ & \texttt{1121} & 0 & 1838 & $E_{6}$ & $+$ & $3x$ & \texttt{2122} & 0\\
1845 & $E_{9}$ & $+$ & $2x$ & \texttt{0020} & 1 & 1846 & $E_{9}$ & $+$ & $2x$ & \texttt{0020} & 1\\
1852 & $E_{9}$ & $+$ & $2x$ & \texttt{0020} & 1 & 1854 & $E_{9}$ & $+$ & $2x$ & \texttt{0020} & 1\\
1864 & $E_{8}$ & $+$ & $2x$ & \texttt{0020} & 1 & 1865 & $E_{7}$ & $+$ & $2x$ & \texttt{0020} & 1\\
1867 & $E_{6}$ & $+$ & $3x$ & \texttt{0020} & 1 & 1877 & $E_{8}$ & $+$ & $2x$ & \texttt{0020} & 1\\
1883 & $E_{6}$ & $+$ & $3x$ & \texttt{0020} & 1 &  &  &  &  &  & \\
\end{longtable}\endgroup

\endgroup
\subsection{Finite-sum thresholds}\label{sec:family12}
\begingroup
\providecommand{\eq}{}\renewcommand{\eq}{\approx}
\providecommand{\Var}{}\renewcommand{\Var}{\operatorname{Var}}
\providecommand{\SR}{}\renewcommand{\SR}{\mathsf{SR}}
\subsubsection{Compression and the threshold theorem}
\begin{lemma}\label{f12:lem:interior}
Let $M$ be a semigroup satisfying $xyzt\eq xzyt$. Adjoin a fresh
identity to obtain $M^1$, even if $M$ already has an identity. Then
\[
 u\sim v\quad\Longleftrightarrow\quad
 aub=avb\ \text{for all }a,b\in M
\]
is a monoid congruence on $M^1$, and $C=M^1/{\sim}$ is commutative.
For each $a,b\in M$ there is a well-defined map
$F_{a,b}:C\to M$, $F_{a,b}([u])=aub$.
\end{lemma}
\begin{proof}
Apply Lemma~\ref{f5:lem:quotient}; the maps $F_{a,b}$ are
well defined by the definition of the quotient congruence.
\end{proof}

Fix a finite semiring $A$ of order $n$ whose multiplicative reduct
satisfies $xyzt\eq xzyt$. Choose positive $r,q$ with
$x^{r+q}\eq x^r$ in $A$, and set $L=r+q-1$. A \emph{normalized
word} is either a single variable or a word in which the first and
last letters are retained and each positive exponent among its
interior letters lies in $\{1,\ldots,L\}$. Interior letters are put
in a fixed order. Both endpoints are allowed to be the same variable,
and endpoint variables may occur in the interior. Every word has this
form modulo the two displayed multiplicative identities: adjacent
interior letters can be exchanged using nonempty prefix and suffix
substitutions in $xyzt\eq xzyt$.

\begin{lemma}\label{f12:lem:compression}
Let $W_1,\ldots,W_m$ be normalized words. There is a list
$\widehat W_1,\ldots,\widehat W_m$ involving at most
\[
 2m+n\bigl((L+1)^m-1\bigr)
\]
variables with exactly the same set of possible evaluation vectors
in $A^m$. Moreover, the original list is obtained from the compressed
list by simultaneous substitution and the commuting-interiors identity.
The assertion about evaluation vectors also holds conditionally on
every fixed assignment of all endpoint variables.
\end{lemma}
\begin{proof}
Protect every variable that appears as the first or last letter of any
word; a length-one word protects its variable. There are at most $2m$
protected variables. Retain all their occurrences, including interior
occurrences. An unprotected variable has a nonzero exponent profile
\[
 \alpha=(\alpha_1,\ldots,\alpha_m)
 \in\{0,\ldots,L\}^m\setminus\{0\},
\]
where $\alpha_i$ counts its occurrences in the interior of $W_i$.
For the commutative monoid $C$ of Lemma~\ref{f12:lem:interior}, define
\[
 T_\alpha=\{([a]^{\alpha_1},\ldots,[a]^{\alpha_m}):a\in A\}
 \subseteq C^m.
\]
Exponent zero means the identity of $C$. The displayed map from the
multiplicative semigroup of $A$ is a homomorphism because $C$ is
commutative. Therefore $T_\alpha$ is a subsemigroup of size at most $n$.

By Lemma~\ref{f3:lem:powers}, $T_\alpha^k=T_\alpha^n$ for
$k\geq n$. A group of $k>n$ variables of profile
$\alpha$ contributes precisely $T_\alpha^k$ to the vector of
interior products; replace it by $n$ variables of that profile.
Groups of different profiles are independent. Conditional on any
assignment of the protected variables, their contributions multiply
in $C^m$, together with the fixed contributions of the protected
interior letters. Applying the maps $F_{a,b}$ in each nonsingleton
word recovers its value. The singleton words already have fixed
values. This proves equality of the full evaluation-vector sets.

There are at most $(L+1)^m-1$ nonzero profiles. Keeping at most $n$
variables of each profile proves the stated bound. To recover a group
originally containing $k>n$ variables, substitute a product of $k-n+1$
fresh variables for one of the retained variables and single fresh
variables for the other $n-1$. All substituted occurrences are interior.
Interior commutations rearrange the resulting powers into the original
exponent profile in every word. Different profiles use disjoint fresh
variables. All protected variables remain fixed.
\end{proof}

\begin{theorem}\label{f12:thm:threshold}
Let $A$ be a finite semiring of order $n$ satisfying $xyzt\eq xzyt$.
Suppose a unary term $h$ is an endomorphic retraction onto an FB
subsemiring $I$, $A+I\subseteq I$, and, for some $d\geq2$,
\[
 a_1+\cdots+a_d\in I\qquad(a_1,\ldots,a_d\in A).
\]
Then $A$ is FB. For any positive $r,q$ with $x^{r+q}\eq x^r$ in $A$,
a finite basis can be constructed using a basis of $I$ and finite
sets of short-polynomial identities on an alphabet of size
\begin{equation}\label{f12:eq:bound}
 N=4(d-1)+n\bigl((r+q)^{2(d-1)}-1\bigr).
\end{equation}
\end{theorem}
\begin{proof}
Let $\Sigma$ be a finite basis of $I$. A \emph{short polynomial} is a
nonempty sum of fewer than $d$ normalized words. On a fixed alphabet
$X_N$ of $N$ variables there are only finitely many such polynomials:
endpoints and exponents have finitely many choices, and the number
of summands is bounded. Repeated summands are allowed. Let
$\mathcal P_N$ be this finite set, using a fixed order of summands.

Take as a proposed basis $\mathcal B$ the following finite collection:
$\SR$, $xyzt\eq xzyt$, $x^{r+q}\eq x^r$, the lifted basis
\[
 \Sigma^h=\{U(h(\mathbf x))\eq V(h(\mathbf x)):
                         U\eq V\in\Sigma\},
\]
the structural identities
\begin{align}
 h(h(x))&\eq h(x),\label{f12:eq:retract}\\
 h(x\mathbin\circ y)&\eq h(x)\mathbin\circ h(y)
                         &&(\circ\in\{+,\cdot\}),\label{f12:eq:endo}\\
 x+h(y)&\eq h(x)+h(y),\label{f12:eq:absorb}\\
 x_1+\cdots+x_d&\eq h(x_1+\cdots+x_d),\label{f12:eq:threshold}
\end{align}
and the finite sets
\begin{align*}
 \Delta_=&\{P\eq Q:P,Q\in\mathcal P_N,\ A\models P\eq Q\},\\
 \Delta_h=&\{P\eq h(P):P\in\mathcal P_N,\ A\models P\eq h(P)\}.
\end{align*}
Here $h(P)$ means substitution of the polynomial $P$ into the fixed
unary term $h$, not a new operation symbol. All members are valid in
$A$. In particular, $x+h(y)$ lies in $I$ and is fixed by $h$, giving
\eqref{f12:eq:absorb} by endomorphism and retraction. The threshold
assumption gives \eqref{f12:eq:threshold}.

We prove completeness. First, every valid identity $P\eq Q$ between
short polynomials, on an arbitrary finite alphabet, follows from
$\mathcal B$. List the word slots on the two sides. Their total number
is at most $2(d-1)$. Lemma~\ref{f12:lem:compression} gives compressed
polynomials on at most $N$ variables. Since compression preserves
the possible vectors of all slot values, it preserves the pairs of
values of the two sums, and hence preserves validity. The compressed
identity belongs, up to renaming, to $\Delta_=$. The recovery
substitution in that lemma derives the original identity.

Second, every valid identity $P\eq h(P)$ with $P$ short follows from
$\mathcal B$. Compress the fewer than $d$ slots of $P$. The same
evaluation-vector argument preserves the range of $P$, and therefore
preserves the assertion that all its values are fixed by $h$.
The compressed absorption identity belongs to $\Delta_h$.
Recovery substitution commutes with application of the unary term
$h$, so it recovers $P\eq h(P)$ as well. This step is needed for
identities having only one side in the retract.

Let $B$ be any model of $\mathcal B$. Equations
\eqref{f12:eq:retract}--\eqref{f12:eq:endo} make $h$ an endomorphic retraction
of $B$, and $\Sigma^h$ implies $h(B)\in\Var(I)$. Distributivity
expands any term to a nonempty sum of words. Normalize its words as
above. If the sum has fewer than $d$ summands, it is short. If it
has at least $d$, apply \eqref{f12:eq:threshold} to the first $d$ summands
and then \eqref{f12:eq:absorb} to each remaining summand. Endomorphism
shows that the result is its own $h$-image. Thus every term is
provably either short or fixed by $h$.

Consider an arbitrary identity $U\eq V$ of $A$ after this normalization.
If both sides are short, the first compression argument proves it in
$B$. If both are fixed by $h$, its restriction to $I$ is valid and
therefore holds in $h(B)$; hence
\[
 U\eq h(U)\eq h(V)\eq V
\]
holds in $B$. If, say, $U=P$ is short and $V$ is fixed by $h$, then
$A\models P\eq h(P)$, because $A\models P\eq V$ and $h(V)=V$.
The second compression argument yields $B\models P\eq h(P)$.
The restricted identity is valid in $I$, so again
$P\eq h(P)\eq h(V)\eq V$ holds in $B$. The other mixed case is
symmetric. All identities of $A$ hold in $B$, which proves completeness.
\end{proof}

The bounds on variables, summands and exponents make
$\Delta_=$ and $\Delta_h$ finite. Their completeness follows from
the compression and recovery arguments above.

\subsubsection{Applications and explicit maps}
For each representative in Table~\ref{f12:tab:targets}, the term is
$h(x)=3x=x+x+x$, the threshold is $d=3$, and the multiplicative
identities are $x^3\eq x^2$ and $xyzt\eq xzyt$. The image is a
two-element semiring. Consequently Theorem~\ref{f12:thm:threshold}
applies with $n=4$, $r=2$, $q=1$, giving $N=328$.

\begin{corollary}
The thirteen four-element semirings listed in Table~\ref{f12:tab:targets}
are finitely based.
\end{corollary}
\begin{proof}
The displayed tables and maps establish all the algebraic hypotheses
of Theorem~\ref{f12:thm:threshold}. The images are among the two-element
semirings treated in Section~\ref{sec:twoelement}.
Apply the theorem to each row.
\end{proof}

Table strings are row-major on $\{0,\ldots,n-1\}$; digit $0$ is
only an element label. The image map is written $h(0)h(1)h(2)h(3)$.
The image isomorphism gives the images of elements $0,1$ of $E_j$.
All retraction maps have proper images, and none of these thirteen
objects has all binary sums in this image: the threshold three is
essential for their inclusion in this extension of the method.

\begin{table}[htbp]\centering\small
\caption{The thirteen applications and their principal witnesses.}
\label{f12:tab:targets}
\setlength{\tabcolsep}{2pt}
\begin{tabular}{rllcll}\toprule
$j$ & Addition & Multiplication & Image & $h$ & Isomorphism\\\midrule
445 & \texttt{0003020300033330} & \texttt{0000012002000000} & $E_{2}$ & \texttt{0003} & $(0,3)$ \\
449 & \texttt{0000020000000003} & \texttt{0000012002000000} & $E_{3}$ & \texttt{0003} & $(0,3)$ \\
459 & \texttt{0003020300033333} & \texttt{0000012002000000} & $E_{4}$ & \texttt{0003} & $(0,3)$ \\
1430 & \texttt{0000030000200000} & \texttt{0000010300200300} & $E_{7}$ & \texttt{0020} & $(0,2)$ \\
1432 & \texttt{0020032022020020} & \texttt{0000010300200300} & $E_{6}$ & \texttt{0020} & $(0,2)$ \\
1434 & \texttt{0020032022220020} & \texttt{0000010300200300} & $E_{8}$ & \texttt{0020} & $(0,2)$ \\
1463 & \texttt{0220230220020220} & \texttt{0000012302200300} & $E_{6}$ & \texttt{0220} & $(0,2)$ \\
1500 & \texttt{0000030000200000} & \texttt{0000010322220300} & $E_{10}$ & \texttt{0020} & $(0,2)$ \\
1502 & \texttt{0020032022220020} & \texttt{0000010322220300} & $E_{10}$ & \texttt{0020} & $(2,0)$ \\
1853 & \texttt{0000030000200000} & \texttt{0020012300200320} & $E_{9}$ & \texttt{0020} & $(0,2)$ \\
1855 & \texttt{0020032022220020} & \texttt{0020012300200320} & $E_{9}$ & \texttt{0020} & $(2,0)$ \\
1878 & \texttt{0000030000200000} & \texttt{0020012322220320} & $E_{8}$ & \texttt{0020} & $(2,0)$ \\
1880 & \texttt{0020032022220020} & \texttt{0020012322220320} & $E_{7}$ & \texttt{0020} & $(2,0)$ \\
\bottomrule\end{tabular}
\end{table}

\begin{table}[htbp]\centering\small
\caption{The two-element image semirings.}\label{f12:tab:images}
\begin{tabular}{cll}\toprule
Image & Addition & Multiplication\\\midrule
$E_{2}$ & \texttt{0110} & \texttt{0000} \\
$E_{3}$ & \texttt{0001} & \texttt{0000} \\
$E_{4}$ & \texttt{0111} & \texttt{0000} \\
$E_{6}$ & \texttt{0110} & \texttt{0001} \\
$E_{7}$ & \texttt{0001} & \texttt{0001} \\
$E_{8}$ & \texttt{0111} & \texttt{0001} \\
$E_{9}$ & \texttt{0001} & \texttt{0101} \\
$E_{10}$ & \texttt{0001} & \texttt{0011} \\
\bottomrule\end{tabular}
\end{table}

For example, $A_{1430}$ has $h(A_{1430})=\{0,2\}$, isomorphic to
$E_7$, where the two operations are the same semilattice operation.
Its binary sum $1+1=3$ is outside the image, while every sum of three
elements is in $\{0,2\}$. The multiplication satisfies the two
identities above, so Theorem~\ref{f12:thm:threshold} applies although the
binary operation-image hypothesis of Theorem~\ref{f11:thm:main} fails.

For completeness, the corresponding Rees quotients can also be
identified explicitly. Put
\[
\begin{split}
 J_0&=\{445,449,459\},\\
 J_1&=\{1430,1432,1434,1463,1500,1502,1853,1855,1878,1880\}.
\end{split}
\]
For $j\in J_0\cup J_1$, write $I_j=h(A_j)$ and define
$q_j:A_j\to B_{86}$, in argument order $0,1,2,3$, by
\[
 q_j=\begin{cases}
       (0,1,2,0),&j\in J_0,\\
       (0,1,0,2),&j\in J_1.
      \end{cases}
\]
Here $B_{86}$ has addition \texttt{000020000} and multiplication
\texttt{000012020}. The operation tables give
\[
 q_j(a+b)=q_j(a)+q_j(b),\qquad q_j(ab)=q_j(a)q_j(b).
\]
Each $q_j$ is surjective, its fibre over $0$ is $I_j$, and its
other fibres are singletons. The set $I_j$ is an ideal for both
operations, so $\ker q_j$ is exactly the Rees congruence $\rho_{I_j}$.
Consequently $A_j/\rho_{I_j}\cong B_{86}$ and
Lemma~\ref{lem:ideal-decomposition} gives
\[
 \Var(A_j)=\Var(I_j)\vee\Var(B_{86}).
\]
The finite basis conclusion follows from
Theorem~\ref{f12:thm:threshold}.

\endgroup
\subsection{Quadratic rectangles and bounded nonlinear parts}\label{sec:family13}
\begingroup
\providecommand{\eq}{}\renewcommand{\eq}{\approx}
\providecommand{\Var}{}\renewcommand{\Var}{\operatorname{Var}}
\providecommand{\SR}{}\renewcommand{\SR}{\mathsf{SR}}
\subsubsection{Two general proof devices}
Let $S$ be a finite commutative semigroup of order $n$, written
additively. Fix positive $r,q$ such that $(r+q)x\eq rx$ holds, and
put $L=r+q-1$. A linear form is a nonempty sum of variables, with each
positive coefficient reduced to $1,\ldots,L$. Some slots in a list
of forms may be formally absent; absence is not a term or a constant
of the signature.

\begin{lemma}\label{f13:lem:linear}
Given a list of $m$ linear-form slots, its variables can be compressed
to at most $n((L+1)^m-1)$ variables without changing the set of possible
evaluation vectors. The original list is recovered by substituting
sums for variables and using commutative associativity.
If $p$ variables are protected and fixed, the same assertion holds
conditionally on their values, with at most
$p+n((L+1)^m-1)$ variables in total.
\end{lemma}
\begin{proof}
Adjoin a fresh additive identity to obtain the monoid $S^1$, even
when $S$ already has an identity. For an unprotected variable, record
the coefficient vector
\[
 \alpha\in\{0,\ldots,L\}^m\setminus\{0\}.
\]
Zero coefficients are interpreted as the new identity. The map
$a\mapsto(\alpha_1a,\ldots,\alpha_ma)$ from $S$ into $(S^1)^m$
is a semigroup homomorphism. Its image $T_\alpha$ is consequently a
subsemigroup of size at most $n$. By Lemma~\ref{f3:lem:powers}, its setwise sums satisfy
$kT_\alpha=nT_\alpha$ for $k\geq n$.

The independent contribution of $k$ variables of profile $\alpha$
is exactly $kT_\alpha$. Replace every group of more than $n$ variables
by $n$ variables of the same profile. The contributions of different
profiles add independently, so the complete evaluation-vector set is
unchanged. With protected variables fixed, their contribution is a
fixed vector, and the same argument applies. There are at most
$(L+1)^m-1$ nonzero profiles.

To restore a group originally having $k>n$ variables, replace one of
its retained variables by a sum of $k-n+1$ fresh variables and replace
the other $n-1$ by individual fresh variables. Every coefficient then
distributes over that sum by commutative associativity, restoring
exactly the original profile. Different groups use disjoint variables.
Protected variables and absent slots remain unchanged.
\end{proof}

Let $h$ be a fixed unary semiring term. Consider the finite structural
identities
\begin{align}
 h(h(x))&\eq h(x),\label{f13:eq:ret}\\
 h(x\mathbin\circ y)&\eq h(x)\mathbin\circ h(y)
                   &&(\circ\in\{+,\cdot\}),\label{f13:eq:endo}\\
 x+h(y)&\eq h(x)+h(y),\label{f13:eq:plusideal}\\
 xh(y)&\eq h(x)h(y),\qquad h(x)y\eq h(x)h(y).\label{f13:eq:mulideal}
\end{align}
They hold whenever $h$ is an endomorphic retraction onto an ideal
for both operations. For example, $xh(y)$ lies in the ideal and is
therefore fixed by $h$, which proves the first multiplicative equation.
These equations also propagate the identity $T\eq h(T)$ from a
subterm to every term containing it.

\enlargethispage{1.2\baselineskip}
\begin{lemma}\label{f13:lem:lift}
Let $A$ have an endomorphic term retraction $h$ onto an FB ideal $I$.
Suppose a finite set $\Omega$ of valid identities, containing
\eqref{f13:eq:ret}--\eqref{f13:eq:mulideal}, normalizes every term either to
a member of a class $\mathcal P$ or to a term provably fixed by $h$.
Suppose a finite subset $\mathcal F\subseteq\mathcal P$ has the
following two properties, with renaming of variables allowed:
\begin{enumerate}
\item Every valid identity between members of $\mathcal P$ is
deducible, using $\Omega$ and substitution, from a valid identity
between members of $\mathcal F$.
\item Every valid identity $P\eq h(P)$ with $P\in\mathcal P$ is
similarly deducible from one with $P\in\mathcal F$.
\end{enumerate}
Then $A$ is FB.
\end{lemma}
\begin{proof}
Choose a finite basis $\Sigma$ for $I$. Adjoin to $\Omega$ the
identities $U(h(\mathbf x))\eq V(h(\mathbf x))$ for $U\eq V\in\Sigma$,
all valid identities $P\eq Q$ with $P,Q\in\mathcal F$, and all valid
$P\eq h(P)$ with $P\in\mathcal F$. This is a finite set of identities
of $A$. In any model $B$, the lifted basis ensures $h(B)\in\Var(I)$.

Normalize an arbitrary identity of $A$. If both sides are in
$\mathcal P$, property (1) proves it in $B$. If both are fixed by $h$,
their restricted identity holds in $I$, hence in $h(B)$, so the original
identity holds in $B$. Finally suppose one side is $P\in\mathcal P$
and the other side $T$ is fixed by $h$. Then $A\models P\eq h(P)$.
Property (2) gives this identity in $B$, and the restricted identity
holds in $h(B)$. Consequently
$P\eq h(P)\eq h(T)\eq T$ holds in $B$. This proves completeness.
\end{proof}

\subsubsection{Completing rectangles of quadratic monomials}
\begin{theorem}\label{f13:thm:rect}
Let $A$ be a finite semiring with a unary term defining an endomorphic
retraction $h$ onto an FB ideal $I$ for both operations. Assume
$A^3\subseteq I$ and
\begin{align}
 xy+xy&\eq xy,\label{f13:eq:prodidem}\\
 xy+zt&\eq xy+zt+xt+zy.\label{f13:eq:rectangle}
\end{align}
Then $A$ is FB. If $|A|=n$ and $(r+q)x\eq rx$, the finite sets in
Lemma~\ref{f13:lem:lift} may be constructed on an alphabet of
\[
 N=n\bigl((r+q)^6-1\bigr)
\]
variables.
\end{theorem}
\begin{proof}
Take $\SR$, \eqref{f13:eq:ret}--\eqref{f13:eq:mulideal},
\eqref{f13:eq:prodidem}--\eqref{f13:eq:rectangle}, $(r+q)x\eq rx$, and
\begin{equation}\label{f13:eq:triple}
 xyz\eq h(xyz)
\end{equation}
as $\Omega$. All are valid by the hypotheses.
Expand a term distributively into a sum of words. If any word has
length at least three, equation~\eqref{f13:eq:triple} fixes its first
three factors under $h$, and \eqref{f13:eq:mulideal} propagates this
property to the entire word. Equation~\eqref{f13:eq:plusideal} then
propagates it to the entire sum.

Otherwise every word has length one or two. Collect the length-one
words into a linear form $L_0$, possibly absent. For a nonempty sum
of quadratic monomials, list the distinct left letters as $X$ and
the distinct right letters as $Y$. If $x\in X$ and $y\in Y$, choose
existing monomials $xv$ and $uy$. Equation~\eqref{f13:eq:rectangle}
adjoins $xy$ and $uv$ without changing their sum. If these selected
monomials coincide, $xy$ was already present. Repeating for all pairs
and removing duplicates with \eqref{f13:eq:prodidem} completes the full
rectangle $X\times Y$. Distributivity then gives
\begin{equation}\label{f13:eq:rectform}
 L_0+L_1L_2,
\end{equation}
where $L_1=\sum_{x\in X}x$ and $L_2=\sum_{y\in Y}y$.
The linear part may be absent, or the product part may be absent,
but not both. The factors in a present product are nonempty. Reduce
coefficients using $(r+q)x\eq rx$.

Let $\mathcal P$ be the class of all expressions
\eqref{f13:eq:rectform} with these absence conventions and reduced
positive coefficients at most $r+q-1$. On a fixed alphabet of $N$
variables there are finitely many such expressions; let this set
be $\mathcal F$. To compress an identity between two members of
$\mathcal P$, list their at most six linear-form slots. Apply
Lemma~\ref{f13:lem:linear} without protected variables. It preserves
the entire slot-value vector set, and therefore the set of pairs of
values after applying the two fixed expression shapes. The compressed
identity is valid and lies in $\mathcal F$, up to renaming. Recovery
by sums substituted for variables restores all slots simultaneously,
and hence restores the original expressions. This proves property (1)
of Lemma~\ref{f13:lem:lift}.

For $P\eq h(P)$, compress the at most three slots of $P$. This
preserves the range of $P$, and hence preserves the assertion that
every value is fixed by $h$. The same recovery substitution gives
the original absorption identity, proving property (2). Apply
Lemma~\ref{f13:lem:lift} to obtain a complete finite basis.
\end{proof}

The two conditions \eqref{f13:eq:prodidem} and \eqref{f13:eq:rectangle} are
both used. The latter alone does not justify removing repeated
monomials. Also, $A^3\subseteq I$ need not mean that triple products
are all equal; four applications below have two-element images.

\subsubsection{Bounding the number of nonlinear monomials}
\begin{theorem}\label{f13:thm:bounded}
Let $A$ be a finite semiring with an endomorphic term retraction $h$
onto an FB ideal $I$ for both operations. Suppose there are integers
$k,d\geq2$ such that every product of $k$ elements is in $I$, and
every sum of $d$ elements of $A^2$ is in $I$. Then $A$ is FB.
If $|A|=n$ and $(r+q)x\eq rx$, a finite basis is obtained from
Lemma~\ref{f13:lem:lift} using an alphabet of size
\[
 N=2(d-1)(k-1)+n\bigl((r+q)^2-1\bigr).
\]
\end{theorem}
\begin{proof}
Take $\SR$, $(r+q)x\eq rx$, and
\eqref{f13:eq:ret}--\eqref{f13:eq:mulideal}, together with
\begin{align}
 x_1\cdots x_k&\eq h(x_1\cdots x_k),\label{f13:eq:k}\\
 x_1y_1+\cdots+x_dy_d&\eq h(x_1y_1+\cdots+x_dy_d),\label{f13:eq:d}
\end{align}
as $\Omega$. Distributive expansion produces a nonempty sum of words.
If any word has length at least $k$, use \eqref{f13:eq:k} and the ideal
identities to fix the whole sum under $h$. If there are at least $d$
words of length at least two, factor each of $d$ of them into two
nonempty factors, apply \eqref{f13:eq:d}, and again propagate the fixed
property to the whole sum. All remaining terms have the form
\begin{equation}\label{f13:eq:boundedform}
 L_0+W_1+\cdots+W_s,\qquad
 0\leq s<d,\quad 2\leq |W_i|<k,
\end{equation}
where $L_0$ is an optional reduced linear form. The whole expression
is nonempty, and words need not commute.

Let $\mathcal P$ be this class and $\mathcal F$ its finite restriction
to the stated alphabet. In an identity between two such expressions,
protect every variable in the nonlinear words. There are at most
$2(d-1)(k-1)$ protected variables. Only the two linear-form slots
contain unprotected variables. Apply Lemma~\ref{f13:lem:linear} with
the protected variables fixed. Conditional on their values, all the
nonlinear words are fixed and compression preserves the possible
pairs of linear-part values. It therefore preserves the possible
pair of full expression values. The compressed identity is valid
and belongs to $\mathcal F$, up to renaming. Recovery substitutes
sums only for unprotected variables, so it leaves every nonlinear
word unchanged and restores both linear forms. This proves (1) of
Lemma~\ref{f13:lem:lift}.

For $P\eq h(P)$, protect the nonlinear variables of $P$ and compress
its one possible linear form. This preserves the range of $P$ for
every assignment to the protected variables, and hence globally.
It preserves the fixed-point assertion and admits the same recovery,
establishing (2). Lemma~\ref{f13:lem:lift} proves the result.
\end{proof}

Neither theorem requires multiplication to be commutative. Both allow
arbitrarily many linear summands. Their finite sets of valid identities
are effectively specified by evaluation on $A$, though the resulting
sets may be extremely large. No small or minimal basis is claimed.

\subsubsection{Applications and explicit maps}
The complete operation tables and principal witnesses are given in
Table~\ref{f13:tab:targets}. A string lists the table row by row on
$\{0,1,2,3\}$; labels have no implicit numerical meaning. All image
maps are either constant with value $0$ or have image $\{0,1\}$.
The source-to-image isomorphisms are the evident maps from the defining tables; the maps for $A_{853}$ and $A_{1255}$ exchange labels
$0$ and $1$ of their two-element sources.

The one-element image is FB with basis $x\eq y$. For the other
images the tables are
\[
 \begin{array}{c|cc}
 &+&\cdot\\\hline
 E_9&\texttt{0001}&\texttt{0101}\\
 E_{10}&\texttt{0001}&\texttt{0011}
 \end{array}
\]
Thus $E_9$ has right-projection multiplication and $E_{10}$ has
left-projection multiplication. A basis is $\SR$ together with
$x+x\eq x$ and $xy\eq y$ or $xy\eq x$, respectively. To verify
completeness directly, eliminate products by the projection identity.
Terms then reduce to nonempty sums of distinct selected variables.
Two such sums are equal in the nontrivial two-element semilattice
exactly when their variable sets agree: a variable in the symmetric
difference separates them by assigning it the absorbing semilattice
value and assigning the other variables the other value. Hence the
listed identities form a complete finite basis. No external
finite-basis classification is needed for these image inputs.

\begin{table}[htbp]\centering\footnotesize
\caption{All thirteen additional FB classes.}\label{f13:tab:targets}
\setlength{\tabcolsep}{3.5pt}
\begin{tabular}{rllclll}\toprule
$j$&Addition&Multiplication&Image&$h(x)$&$h$&Method\\\midrule
254 & \texttt{0000000000020023} & \texttt{0000020000000000} & $\{0\}$ & $x^3$ & \texttt{0000} & B \\
259 & \texttt{0000000200020223} & \texttt{0000020000000000} & $\{0\}$ & $x^3$ & \texttt{0000} & B \\
261 & \texttt{0000032302020323} & \texttt{0000020000000000} & $\{0\}$ & $x^3$ & \texttt{0000} & B \\
307 & \texttt{0000010300200303} & \texttt{0000003000000000} & $\{0\}$ & $x^3$ & \texttt{0000} & R \\
308 & \texttt{0000010000230033} & \texttt{0000003000000000} & $\{0\}$ & $x^3$ & \texttt{0000} & R \\
318 & \texttt{0000012002200003} & \texttt{0000033000000000} & $\{0\}$ & $x^3$ & \texttt{0000} & R \\
321 & \texttt{0000012302230333} & \texttt{0000033000000000} & $\{0\}$ & $x^3$ & \texttt{0000} & R \\
484 & \texttt{0000012002200003} & \texttt{0000030003000000} & $\{0\}$ & $x^3$ & \texttt{0000} & R \\
487 & \texttt{0000012302230333} & \texttt{0000030003000000} & $\{0\}$ & $x^3$ & \texttt{0000} & R \\
852 & \texttt{0000010000230033} & \texttt{0000111100300000} & $E_{10}$ & $x^3$ & \texttt{0100} & R \\
853 & \texttt{0100111101230133} & \texttt{0000111100300000} & $E_{10}$ & $x^3$ & \texttt{0100} & R \\
1254 & \texttt{0000010000230033} & \texttt{0100010001300100} & $E_{9}$ & $x^3$ & \texttt{0100} & R \\
1255 & \texttt{0100111101230133} & \texttt{0100010001300100} & $E_{9}$ & $x^3$ & \texttt{0100} & R \\
\bottomrule\end{tabular}
\end{table}

Here R means Theorem~\ref{f13:thm:rect} and B means
Theorem~\ref{f13:thm:bounded}. In the R cases, addition is idempotent,
so $r=q=1$, and the theorem gives $N=4(2^6-1)=252$.
In the B cases, $3x\eq2x$, so $r=2$, $q=1$; also $k=3$, $d=2$.
The stated bound is $N=4+4(3^2-1)=36$.
All three B cases have all triple products equal to $0$ and satisfy
$xy+zt\eq0$, with $0$ represented by the unary term $x^3$.
This notation does not add a constant to the signature.

For example, $A_{307}$ has only one nonzero product, $1\cdot2=3$.
Its additive reduct satisfies $0+a=0$, and $3+3=3$.
For any quadratic monomials to have a sum different from $0$, all
their left variables must be $1$ and all their right variables $2$.
This illustrates why completing the rectangle preserves their sum.
Equation~\eqref{f13:eq:rectangle} follows directly from the fact that all nonzero quadratic products must have left value $1$ and right value $2$. The ideal image is
the singleton $\{0\}$, so Theorem~\ref{f13:thm:rect} gives finite basability.

The four examples $A_{852},A_{853},A_{1254},A_{1255}$ have proper
two-element ideal images. They were among the remaining ideal-retraction
cases, whose known FB images and known FB quotients alone did not
justify a positive conclusion. Theorem~\ref{f13:thm:rect} supplies the
missing completeness argument. We do not use an unproved closure
claim for arbitrary joins of FB varieties.

\endgroup
\section{Commutative support and coefficient forms}\label{sec:commutative-forms}
Throughout this section, use the nonlinear support monomials $m_U$
of Lemma~\ref{lem:support-calculus}: $m_{\{x\}}=x^2$ and
$m_U=\prod_{x\in U}x$ for $|U|>1$. When $x^2\eq x$ is imposed,
$m_{\{x\}}=x$. All stated word reductions follow from that lemma
and distributivity. The recovery principles and conventions of
Section~\ref{sec:preliminaries} apply throughout; catalogue transfers
are collected in Corollary~\ref{cor:catalogue-transfers}.

\subsection{Multiplicative support antichains}\label{sec:family17}
\begingroup
\providecommand{\eq}{}\renewcommand{\eq}{\approx}
\providecommand{\Var}{}\renewcommand{\Var}{\operatorname{Var}}
\providecommand{\SR}{}\renewcommand{\SR}{\mathsf{SR}}
\emph{Generators (basis size):} $A_{157}$ (12), $A_{173}$ (10), $A_{1123}$ (10).

All three have nonidempotent addition. For $A_{157}$ and $A_{173}$,
a product of at least two factors is $1$ exactly when every factor is
$1$, and is $0$ otherwise. For $A_{1123}$, such a product is $0$
exactly when every factor differs from $1$, and is $1$ otherwise.
These statements follow directly from the defining tables.

\subsubsection{Reduction of multiplicative words}
Let $\mathcal C$ consist of $\SR$ and
\begin{equation}\label{f17:eq:common}
 xy\eq yx,\qquad x^2y\eq xy,\qquad 3x\eq2x.
\end{equation}
Thus $\mathcal C$ has eight identities.

\subsubsection{An antichain basis for two generators}
\begin{theorem}\label{f17:thm:antichain}
Each of $A_{173}$ and $A_{1123}$ has the following ten-identity basis:
$\mathcal C$ together with
\begin{equation}\label{f17:eq:antichain}
 xy+xy\eq xy,\qquad x+xy\eq x.
\end{equation}
Consequently $\Var(A_{173})=\Var(A_{1123})$, although the algebras
are not isomorphic.
\end{theorem}
\begin{proof}
The defining operations satisfy the identities. By
Lemma~\ref{lem:support-calculus}(i), record linear coefficients by
$c:X\to\{0,1,2\}$ and put $L=\{x:c(x)>0\}$. Nonlinear
idempotence removes repeated monomials, and $x+xy\eq x$ deletes
every nonlinear support meeting $L$. For $U\subsetneq V$, substitute
$x=m_U$, $y=\prod_{z\in V\setminus U}z$ in that law to obtain
$m_U+m_V\eq m_U$. Thus the form is
\begin{equation}\label{f17:eq:anti-normal}
 \sum_{x\in L}c(x)x+\sum_{U\in\mathcal H}m_U,
\end{equation}
where $\mathcal H$ is an antichain of nonempty subsets of $X\setminus L$
and at least one part is present.

In $A_{173}$, $\alpha_{\{x\}}(3,0)$ gives $3$, $2$, or $0$
according as $c(x)=1,2$, or $0$. In $A_{1123}$ use
$\alpha_{\{x\}}(3,1)$: the first two outputs are again $3,2$,
whereas $c(x)=0$ gives a value in $\{0,1\}$. These tests recover $c$.
Put $(e,d)=(0,1)$ for $A_{173}$ and $(e,d)=(1,0)$ for $A_{1123}$.
In either case $e$ is an additive identity, $d+d=d$, and a nonlinear
product on $\{e,d\}$ is $d$ exactly when all factors are $d$.
For $D\subseteq X\setminus L$, evaluation at $\alpha_D(d,e)$ is
$d$ exactly when
\begin{equation}\label{f17:eq:upset}
 (\exists U\in\mathcal H)\ U\subseteq D;
\end{equation}
otherwise it is $e$. Lemma~\ref{lem:finite-set-recovery}(ii) recovers
$\mathcal H$, including the empty case. The common normal form proves
completeness and equality of varieties. The generators are not isomorphic:
the numbers of ordered pairs whose product is not the multiplicative
absorber are respectively one and nine.
\end{proof}

\subsubsection{A basis with one surviving nonlinear support}
\begin{theorem}\label{f17:thm:single}
The algebra $A_{157}$ has the following twelve-identity basis:
$\mathcal C$ together with
\begin{align}
 xy+zu&\eq v^2+v^2,\label{f17:eq:two-zero}\\
 u^2+u^2+x&\eq u^2+u^2,\label{f17:eq:zero-absorb}\\
 x+xy&\eq u^2+u^2,\label{f17:eq:incident-zero}\\
 x+yz&\eq 2x+yz.\label{f17:eq:mixed-double}
\end{align}
\end{theorem}
\begin{proof}
In $A_{157}$, both addition and multiplication absorb the element
$0$. The sum of any two elements of $\{0,1\}$ is $0$. All
nonlinear monomials take values in $\{0,1\}$, so
\eqref{f17:eq:two-zero} and \eqref{f17:eq:zero-absorb} hold.
If $xy=1$, then $x=1$ and $x+xy=1+1=0$; if $xy=0$, the sum is
again $0$. This proves \eqref{f17:eq:incident-zero}. For
\eqref{f17:eq:mixed-double}, it suffices to consider $yz=1$, since
$0$ is additively absorbing. Directly from the addition table,
$x+1=2x+1$ for $x\in\{0,1,2,3\}$. The identities in
$\mathcal C$ also follow from the tables.

Equation \eqref{f17:eq:two-zero} implies that the term $v^2+v^2$ is
independent of $v$. Denote this term by $Z$ as an abbreviation,
without introducing a constant symbol. It absorbs addition by
\eqref{f17:eq:zero-absorb}. It also absorbs multiplication:
distributivity expresses $Zx$ as the sum of two nonlinear monomials,
which is $Z$ by \eqref{f17:eq:two-zero}; multiplication is commutative.

Expand a term by Lemma~\ref{lem:support-calculus}(i). If two or more nonlinear
summands occur, including repeated copies of the same monomial,
then \eqref{f17:eq:two-zero} and \eqref{f17:eq:zero-absorb} reduce the
polynomial to $Z$. If precisely one nonlinear monomial $m_U$
occurs and its support meets the linear support, write it as $xv$
with such a linear variable $x$. Equation \eqref{f17:eq:incident-zero}
then gives $Z$, which absorbs all remaining summands. If those
supports are disjoint, \eqref{f17:eq:mixed-double} makes every
positive linear coefficient two.

\Needspace{4\baselineskip}
The three kinds of normal form are
\begin{align}
 &Z;\label{f17:eq:nf-zero}\\
 &\sum_{x\in L}c(x)x,
   \quad L\ne\varnothing,\quad c(x)\in\{1,2\};\label{f17:eq:nf-linear}\\
 &\sum_{x\in L}2x+m_U,
   \quad U\ne\varnothing,\quad L\cap U=\varnothing.\label{f17:eq:nf-mixed}
\end{align}
The set $L$ in \eqref{f17:eq:nf-mixed} may be empty.

To separate the three kinds, assign every variable the value $2$.
A pure linear form evaluates to $2$, whereas both $Z$ and every
mixed form evaluate to $0$. Each mixed form has a nonzero evaluation:
assign all variables in $U$ the value $1$ and all variables in $L$
the value $2$. The result is $1$ if $L$ is empty and $3$ otherwise.
Thus a mixed form is not equal to $Z$ as a function.

For two pure linear forms, the assignment $x=0$ with all other
variables $2$ gives $0$ exactly when $x$ belongs to the linear
support. Once the support is fixed, use $x=1$ and all other
variables $2$. A coefficient of two gives $0$; a coefficient of
one gives $1$ when $x$ is the only variable in the support and $3$
otherwise. This recovers all pure linear coefficients.

For a mixed form, a full assignment has nonzero value exactly when
\begin{equation}\label{f17:eq:rectangle}
 x=2\quad(x\in L),\qquad x=1\quad(x\in U).
\end{equation}
Indeed, failure of the second condition makes $m_U=0$ and therefore
the whole sum $0$. Under the second condition $m_U=1$.
Since $2a=2$ exactly for $a=2$, and $2a=0$ for the other three
elements, the first condition is then necessary and sufficient.
The variables outside $L\cup U$ are unrestricted. The nonzero
assignment set is consequently a nonempty Cartesian product whose
coordinate projections are $\{2\}$ on $L$, $\{1\}$ on $U$, and
$\{0,1,2,3\}$ elsewhere. Its coordinate projections uniquely
determine the disjoint sets $L$ and $U$.

The forms are separated on every finite alphabet. Lemma~\ref{lem:normal-form-completeness} proves completeness.
\end{proof}

\endgroup
\subsection{Marked and parity supports}\label{sec:family20}
\begingroup
\providecommand{\eq}{}\renewcommand{\eq}{\approx}
\providecommand{\Var}{}\renewcommand{\Var}{\operatorname{Var}}
\providecommand{\SR}{}\renewcommand{\SR}{\mathsf{SR}}
\emph{Generators (basis size):} $A_{540}$ (11), $A_{564}$ (13), $A_{569}$ (13), $A_{575}$ (13), $A_{1093}$ (13), $A_{1096}$ (11), $A_{1102}$ (11), $A_{1103}$ (13), $A_{1122}$ (13).

\subsubsection{Word reduction and common support tests}
All our bases contain $\SR$ and
\begin{equation}\label{f20:eq:words}
 xy\eq yx,\qquad x^2y\eq xy.
\end{equation}
The marked-support and parity-support bases below also contain
\begin{equation}\label{f20:eq:support}
 2xy\eq xy,\qquad xy+xyz\eq xy,\qquad x+y+xy\eq x+y.
\end{equation}
The first equation reduces nonlinear multiplicities to one. The
second gives
\begin{equation}\label{f20:eq:inclusion}
 m_U+m_V\eq m_U\quad\text{if }\varnothing\ne U\subseteq V.
\end{equation}
This is Lemma~\ref{lem:support-calculus}(iii).

If two distinct linear variables $x,y$ occur, any $m_U$ with
$\{x,y\}\subseteq U$ can be deleted. Insert $xy$ using the last
equation of \eqref{f20:eq:support}, use \eqref{f20:eq:inclusion} to delete
$m_U$, and then delete $xy$ again. This argument also covers
$U=\{x,y\}$ directly.

For the seven marked-support and parity-support generators, use labels $e,d,a,b$:
\[
 (e,d,a,b)=
 \begin{cases}
 (0,1,2,3),&j=540,564,575,\\
 (1,0,3,2),&j=1096,1102,1103,1122.
 \end{cases}
\]
The element $e$ is an additive identity and a multiplicative
absorbing element; $d+d=d$. Products of at least two factors from
$\{e,d,a\}$ equal $d$ when every factor is in $\{d,a\}$, and equal
$e$ otherwise. Also $a+d\ne a$ and $b\notin\{e,d\}$.

\begin{lemma}[Recovering supports meeting linear variables]\label{f20:lem:recover}
Let $X$ be a finite variable set, let $L\subseteq X$ be a fixed
set of linear variables, and let $\mathcal H$
be an antichain of nonempty supports such that $|U\cap L|\leq1$
for every $U\in\mathcal H$. Suppose terms under consideration
are sums of their linear parts and the monomials $m_U$.
The term function on any of the seven tables above determines
the members of $\mathcal H$ disjoint from $L$, and the members
meeting $L$ in a variable whose coefficient is one. In particular,
it determines all of $\mathcal H$ if no member meets a linear
variable of coefficient greater than one.
\end{lemma}
\begin{proof}
Use the families $\mathcal H_0,\mathcal H_x$ of
Lemma~\ref{lem:marked-antichain-recovery}. For $D\subseteq X\setminus L$,
$\alpha_D(d,e)$ gives $d$ iff $D$ contains a member of $\mathcal H_0$.
If $x$ has coefficient one and $D$ contains no such member,
$\beta_{x,D\cup\{x\}}(a,d,e)$ gives $a+d$ iff $D$ contains a
member of $\mathcal H_x$, and otherwise gives $a$.
The two outcomes differ. The cited lemma recovers the stated supports;
under the final hypothesis they exhaust $\mathcal H$.
\end{proof}

\subsubsection{Marked and parity-support bases}
\begin{theorem}\label{f20:thm:marked}\label{f20:thm:parity}
Put $J_{\rm m}=\{564,575,1103,1122\}$ and
$J_{\rm p}=\{540,1096,1102\}$. For $j$ in either group, a complete
basis of $A_j$ consists of $\SR$, \eqref{f20:eq:words},
\eqref{f20:eq:support}, and the corresponding list
\begin{align}
 J_{\rm m}:&\quad 3x\eq2x,\qquad 2x+xy\eq2x,\qquad
                    x+x^2\eq2x;\label{f20:eq:marked}\\
 J_{\rm p}:&\quad 2x\eq x^2.\label{f20:eq:parity}
\end{align}
The basis sizes are thirteen and eleven, respectively. Each group
generates a single variety. The normal forms are
\begin{equation}\label{f20:eq:markednf}
 \sum_{x\in L}c(x)x+\sum_{U\in\mathcal H}m_U,
\end{equation}
where $\mathcal H$ is an antichain, $|U\cap L|\leq1$ for
$U\in\mathcal H$, at least one part is present, and
\[
\begin{array}{ccl}
\toprule
\text{Group}&c(x)\ (x\in L)&\text{Additional restrictions}\\\midrule
J_{\rm m}&1,2&
 U\cap\{x:c(x)=2\}=\varnothing,\quad \{x\}\notin\mathcal H\ (x\in L)\\
J_{\rm p}&1&\text{None}\\\bottomrule
\end{array}
\]
\end{theorem}
\begin{proof}
For these generators a nonlinear product is $d$ precisely when all
its factors belong to $H$, and is $e$ otherwise, where
$H=\{d,a\}$ for $j=540,564,575$ and $H=\{d,a,b\}$ otherwise.
Here $e+d=d+d=d$, and sums of two members of $H$ absorb $d$.
These rules verify \eqref{f20:eq:words} and \eqref{f20:eq:support}.
In $J_{\rm m}$, the additive tables give $3r=2r$ and
$r+r^2=2r$ for every $r$, and $2r+d=2r$ for $r\in H$.
Outside $H$ all products with $r$ are $e$.
In $J_{\rm p}$ the tables give $2r=r^2$.
Thus the respective added identities hold.

Reduce words and nonlinear multiplicities first. In $J_{\rm m}$,
truncate linear coefficients at two; a nonlinear singleton at a
linear variable converts its coefficient to two and is then absorbed.
The law $2x+xy\eq2x$ deletes every support meeting a doubled linear
variable. In $J_{\rm p}$, replace pairs of linear copies by squares.
In either case, the argument after \eqref{f20:eq:inclusion} deletes
supports meeting $L$ twice, and inclusion absorption leaves an
antichain. This derives precisely the displayed restrictions.
For $J_{\rm p}$, a singleton at $x\in L$ remains: it represents
an odd coefficient at least three. Here $4x\eq2x$, whereas
$3b=b\ne2b$.

The test $\alpha_{\{x\}}(b,e)$ has value $b$ exactly when $x\in L$:
in $J_{\rm m}$, $2b=b$ and no nonlinear singleton remains at $x$;
in $J_{\rm p}$, $b+b^2=b$. Otherwise every contribution lies in
$\{e,d\}$. For $J_{\rm m}$ and $x\in L$,
$\alpha_{\{x\}}(a,e)$ distinguishes coefficients one and two by
the distinct values $a,2a$. Lemma~\ref{f20:lem:recover} recovers
$\mathcal H$ in both groups, including the permitted empty tails
in $J_{\rm p}$. Lemma~\ref{lem:normal-form-completeness} proves
completeness; the common basis within each group gives its variety equality.
\end{proof}

\subsubsection{Two generators requiring a nonlinear guard}
\begin{theorem}\label{f20:thm:guarded}
For each $j\in\{569,1093\}$, a complete thirteen-identity basis
consists of $\SR$, \eqref{f20:eq:words}, and
\begin{align}
 3x&\eq2x,&2xy&\eq xy,&xy+xyz&\eq xy,\label{f20:eq:guard-common}\\
 x+x^2&\eq x^2,&2x+yz&\eq x^2+yz,\label{f20:eq:guard-linear}\\
 x+y+xy+zu&\eq x+y+zu.&&&&\label{f20:eq:guard}
\end{align}
Consequently $\Var(A_{569})=\Var(A_{1093})$.
\end{theorem}
\begin{proof}
The identities hold and still imply \eqref{f20:eq:inclusion}.
With another nonlinear summand $p$, insert $xy$, absorb $m_U$, and
delete $xy$ by the guard law to obtain
\begin{equation}\label{f20:eq:guard-delete}
 x+y+m_U+p\eq x+y+p\quad(\{x,y\}\subseteq U,\ x\ne y).
\end{equation}
Without another summand, use $p=m_U$ to insert $xy$ and absorb $m_U$:
\begin{equation}\label{f20:eq:guard-reduce}
 x+y+m_U\eq x+y+xy.
\end{equation}
Thus a pure linear form is $\sum_{x\in L}c(x)x$ with nonempty $L$
and $c(x)\in\{1,2\}$. When a nonlinear summand is present,
$2x+yz\eq x^2+yz$ replaces doubled linear variables by squares,
and $x+x^2\eq x^2$ removes linears at nonlinear singletons.
Let $L$ be the remaining single linears. If a nonlinear support meets
$L$ at most once, retain it as a guard while deleting supports meeting
$L$ twice, then minimize by \eqref{f20:eq:inclusion}. This gives
\begin{equation}\label{f20:eq:guard-nf}
 \sum_{x\in L}x+\sum_{U\in\mathcal H}m_U,
\end{equation}
where $\mathcal H$ is nonempty, is an antichain, meets $L$ at most
once, and excludes $\{x\}$ for $x\in L$.
Otherwise \eqref{f20:eq:guard-reduce} leaves one pair product.
Any two pairs on $L$ are interchangeable by inserting one with the
other as guard and then deleting the first. Choose the first two
variables $x_0,y_0$ of $L$, giving
\begin{equation}\label{f20:eq:guard-empty}
 \sum_{x\in L}x+x_0y_0,\qquad |L|\geq2.
\end{equation}
Record this case by $\mathcal H=\varnothing$.

The uniform value $3$ in $A_{569}$, respectively $2$ in $A_{1093}$,
separates pure forms from all others: the respective pure outputs are
$3,2$ and the nonpure output is zero. For pure forms in $A_{569}$,
$\alpha_{\{x\}}(2,3)$ gives $3,2,1$ for $c(x)=0,1,2$.
In $A_{1093}$, $\alpha_{\{x\}}(0,2)$ detects occurrence, and
$\alpha_{\{x\}}(3,1)$ distinguishes coefficients one and two by $3,0$.
For nonpure forms, $\alpha_{\{x\}}(2,3)$ in $A_{569}$ and
$\alpha_{\{x\}}(3,1)$ in $A_{1093}$ give respectively $2,3$ iff
$x\in L$; otherwise the output lies in $\{0,1\}$.

For $D\subseteq X\setminus L$, $\alpha_D(1,3)$ in $A_{569}$
gives $1$ iff a member of $\mathcal H_0$ lies in $D$, and otherwise $0$.
In $A_{1093}$ use $\alpha_D(0,1)$, with respective outcomes $0,1$.
On tests avoiding $\mathcal H_0$, change a fixed $x\in L$ to $2$
in the first algebra or $3$ in the second. Containment of a tail in
$\mathcal H_x$ gives respectively $1$ or $0$; absence gives $2$ or $3$.
Lemma~\ref{lem:marked-antichain-recovery} recovers $\mathcal H$.
A retained pair guard contributes $0$ or $1$ respectively and changes
none of these tests; if all recovered families are empty,
\eqref{f20:eq:guard-empty} is fixed by $L$. This separates every case.
\end{proof}

\endgroup
\subsection{Unrestricted and weighted supports}\label{sec:family21}
\begingroup
\providecommand{\eq}{}\renewcommand{\eq}{\approx}
\providecommand{\Var}{}\renewcommand{\Var}{\operatorname{Var}}
\providecommand{\SR}{}\renewcommand{\SR}{\mathsf{SR}}
\subsubsection{Generators}
For the following generators,
$1$ is an additive identity and a multiplicative zero, and $0$ is
additively absorbing. A nonlinear product is $0$ when none of its
factors is $1$, and is $1$ otherwise.

\emph{Generators (basis size):} $A_{1094}$ (11), $A_{1108}$ (11), $A_{1121}$ (11), $A_{1127}$ (11), $A_{1095}$ (13).

\subsubsection{Common equational reductions}
Let $\mathcal B$ consist of $\SR$ and the following five identities:
\begin{align}
xy&\eq yx,&x^2y&\eq xy,\label{f21:eq:words}\\
2xy&\eq xy,&xy+xyz&\eq xy,\label{f21:eq:products}\\
x+x^2&\eq x^2.&&\label{f21:eq:square}
\end{align}
Thus $\mathcal B$ contains ten identities.

\begin{lemma}\label{f21:lem:reduction}
Modulo $\mathcal B$, every nonlinear word of support $U$ equals $m_U$.
Nonlinear summands are additively idempotent, and
$m_U+m_V\eq m_U$ whenever $U\subseteq V$. If the nonlinear singleton
$m_{\{x\}}=x^2$ is present, all linear occurrences of $x$ may be removed.
\end{lemma}
\begin{proof}
Apply Lemma~\ref{lem:support-calculus}(i),(iii), then
$2xy\eq xy$ and $x+x^2\eq x^2$ to the nonlinear and linear copies,
respectively.
\end{proof}

\subsubsection{Unrestricted antichains and two coefficient laws}
\begin{theorem}\label{f21:thm:unrestricted}
For each $j\in\{1094,1108,1121\}$, a complete eleven-identity basis
for $A_j$ is $\mathcal B$ together with
\begin{equation}\label{f21:eq:threshold}
3x\eq2x.
\end{equation}
Consequently the three algebras generate the same variety.
For $A_{1127}$, a complete eleven-identity basis is instead
$\mathcal B$ together with
\begin{equation}\label{f21:eq:period}
3x\eq x.
\end{equation}
\end{theorem}
\begin{proof}
The identities hold. Lemma~\ref{f21:lem:reduction} leaves an antichain
$\mathcal H$ of nonlinear supports and removes all linears at its
singletons. The two coefficient laws reduce each positive coefficient
to one or two, with parity determining the choice for $A_{1127}$.
The form is
\begin{equation}\label{f21:eq:normal}
\sum_{x\in X}c(x)x+\sum_{U\in\mathcal H}m_U,
\end{equation}
where $c:X\to\{0,1,2\}$, $\{x\}\in\mathcal H$ implies $c(x)=0$,
and at least one part is present. Supports may meet the linear set
without any further restriction.

The test $\alpha_D(2,1)$ gives zero iff some $U\in\mathcal H$ lies
in $D$: an active nonlinear monomial is the absorber $0$; otherwise
all contributions lie in $\{1,2\}$ and $2+2=2$. Thus
Lemma~\ref{lem:finite-set-recovery}(ii) recovers $\mathcal H$.
If $\{x\}\notin\mathcal H$, every nonlinear summand at
$\alpha_{\{x\}}(3,1)$ has value $1$. Its outputs are
\[
\begin{array}{c|ccc}
j&c(x)=0&c(x)=1&c(x)=2\\\hline
1094,1108&1&3&0\\1121,1127&1&3&2
\end{array}
\]
and distinguish all remaining coefficients. The common basis also
gives the asserted variety equality for the first three generators.
\end{proof}

\subsubsection{A basis with additive weights}
\begin{theorem}\label{f21:thm:weighted}
A complete thirteen-identity basis for $A_{1095}$ consists of
$\mathcal B$ together with
\begin{align}
3x&\eq x^2,\label{f21:eq:triple}\\
2x+y+xy&\eq2x+y,\label{f21:eq:weight21}\\
x+y+z+xyz&\eq x+y+z.\label{f21:eq:weight111}
\end{align}
\end{theorem}
\begin{proof}
The identities hold. Replace three linear copies by $x^2$ and absorb
any further copies; then apply Lemma~\ref{f21:lem:reduction}. Put
$w_c(D)=\sum_{x\in D}c(x)$. If $w_c(U)\geq3$, either $U$ contains
a doubled linear and another linear, or it contains three distinct
linears. The corresponding added law inserts $xy$ or $xyz$;
this absorbs $m_U$, and the inserted product can then be deleted.
The resulting form is \eqref{f21:eq:normal}, with
\begin{equation}\label{f21:eq:admissible}
\{x\}\in\mathcal H\Longrightarrow c(x)=0,
\qquad w_c(U)\leq2\quad(U\in\mathcal H).
\end{equation}
Deleting supports preserves the antichain, and the linear context
keeps the whole form nonempty.

For each $x$, the outputs at $\alpha_{\{x\}}(2,1)$ and
$\alpha_{\{x\}}(3,1)$ are respectively
\[
\begin{array}{c|cccc}
 &c(x)=0,\ \{x\}\notin\mathcal H&c(x)=1&c(x)=2&\{x\}\in\mathcal H\\\hline
\text{pair}&(1,1)&(2,3)&(3,0)&(0,0).
\end{array}
\]
These recover $c$ and the singleton supports. For $w_c(D)\leq2$,
$\alpha_D(2,1)$ gives zero iff $D$ contains a member of $\mathcal H$;
without an active nonlinear support it gives $1,2,3$ according as
$w_c(D)=0,1,2$. Every member of $\mathcal H$ is admissible by
\eqref{f21:eq:admissible}, so Lemma~\ref{lem:finite-set-recovery}(ii)
recovers the remaining supports.
\end{proof}

\endgroup
\subsection{Union and interval closures}\label{sec:family22}
\begingroup
\providecommand{\eq}{}\renewcommand{\eq}{\approx}
\providecommand{\Var}{}\renewcommand{\Var}{\operatorname{Var}}
\providecommand{\SR}{}\renewcommand{\SR}{\mathsf{SR}}
\emph{Generators (basis size):} $A_{620}$ (11), $A_{720}$ (11), $A_{738}$ (11), $A_{1199}$ (11), $A_{1200}$ (10), $A_{1209}$ (11), $A_{1213}$ (10), $A_{1586}$ (11), $A_{1588}$ (11), $A_{1783}$ (10), $A_{1785}$ (11).

\subsubsection{Common word reductions}
Let $\mathcal B$ consist of $\SR$ and
\begin{equation}\label{f22:eq:common}
xy\eq yx,\qquad x^2y\eq xy,\qquad 2x\eq x^2,
\qquad x+x^2\eq x^2.
\end{equation}
Thus $\mathcal B$ contains nine identities.

\begin{lemma}\label{f22:lem:words}
Modulo $\mathcal B$, each nonlinear word of support $U$ equals $m_U$.
Every term is a nonempty sum of distinct linear variables and distinct
nonlinear monomials $m_U$. A linear $x$ is absent whenever $m_{\{x\}}$
is present. Every nonlinear monomial is additively idempotent.
\end{lemma}
\begin{proof}
Lemma~\ref{lem:support-calculus}(i) reduces nonlinear words.
Replace repeated linear variables by squares using $2x\eq x^2$,
and absorb remaining copies by $x+x^2\eq x^2$.
Finally $2xy\eq(2x)y\eq x^2y\eq xy$ gives nonlinear idempotence.
\end{proof}
The term's content is its set of occurring variables. All identities
used below have the same content on both sides.

\subsubsection{Two union-closure bases}
For a union-closed family $\mathcal K$, recall
\begin{equation}\label{f22:eq:fixed}
 D\in\mathcal K\quad\Longleftrightarrow\quad
 G_{\mathcal K}(D)=D\qquad(D\ne\varnothing).
\end{equation}
\begin{theorem}\label{f22:thm:linearunion}\label{f22:thm:nonlinearunion}
A complete ten-identity basis for $A_{1200}$, respectively $A_{1213}$,
is $\mathcal B$ together with the corresponding identity
\begin{align}
 A_{1200}:&\quad x+y+xy\eq x+y;\label{f22:eq:linearunion}\\
 A_{1213}:&\quad xy+zu+xyzu\eq xy+zu.\label{f22:eq:nonlinearunion}
\end{align}
The respective normal forms are
\begin{align}
 &\sum_{x\in L}x+
   \sum_{U\in\mathcal F\setminus\{\{x\}:x\in L\}}m_U,
   \label{f22:eq:linearunionnormal}\\
 &\sum_{x\in L}x+\sum_{U\in\mathcal H}m_U,
   \label{f22:eq:nonlinearunionnormal}
\end{align}
where $\mathcal F$ is nonempty and union-closed, with
$\{x\}\in\mathcal F$ for $x\in L$, whereas $\mathcal H$ is
union-closed with $\{x\}\notin\mathcal H$ for $x\in L$.
In the second form at least one part is nonempty.
\end{theorem}
\begin{proof}
Both algebras have additive identity $1$ and absorber $0$, with
$2+2=2$ and $3+3=0$. A nonlinear product is $1$ if a factor
is $1$, $2$ if every factor is $2$, and $0$ otherwise.
These rules verify $\mathcal B$. Addition differs only at
$3+2$: its value is $0$ in $A_{1200}$ and $3$ in $A_{1213}$.
The law $u+v+uv=u+v$ holds on $\{0,1,2\}$ in both algebras
and also when a value is $3$ in $A_{1200}$. Since products lie
in $\{0,1,2\}$, this proves validity of both proposed bases.

Apply Lemma~\ref{f22:lem:words}. For $A_{1200}$ associate a
singleton to each linear summand and its support to each nonlinear
summand. Identity~\eqref{f22:eq:linearunion} inserts the product
of distinct summands, closing this family under unions.
Equal supports need no insertion, so a linear summand is never
duplicated. For $A_{1213}$, \eqref{f22:eq:nonlinearunion} closes
only the nonlinear supports under unions. In either case retain
the linear set $L$ separately, giving the stated forms.

The test $\alpha_{\{x\}}(3,1)$ has value $3$ exactly when $x\in L$.
Write $\mathcal K=\mathcal F$ or $\mathcal H$, respectively.
For $x\in D$, use $\beta_{x,D}(a,2,1)$, where
$a=0$ for $A_{1200}$, and for $A_{1213}$ take $a=3$ if
$x\in L$ and $a=0$ otherwise. A nonlinear summand has value
$0$ precisely when its support contains $x$ and lies in $D$.
In the first algebra a linear $x$ also supplies $0$, recording
its singleton in $\mathcal F$. In the second, the possible single
linear value $3$ absorbs $1,2$ but not $0$. Thus in both cases
the output is $0$ exactly when $x\in G_{\mathcal K}(D)$.
Equation~\eqref{f22:eq:fixed} recovers the family, including
$\mathcal H=\varnothing$, and Lemma~\ref{lem:normal-form-completeness}
proves completeness.
\end{proof}

\subsubsection{Interval saturation with recorded content}
Consider the identity
\begin{equation}\label{f22:eq:interval}
xy+zu+xyz\eq xy+zu.
\end{equation}
It permits partial enlargement of a support using variables of a
different nonlinear summand.

\Needspace{4\baselineskip}
\begin{lemma}\label{f22:lem:interval}
Suppose a polynomial has nonempty nonlinear support family $\mathcal K$,
and let $C=\bigcup\mathcal K$. Modulo $\mathcal B$ and
\eqref{f22:eq:interval}, its nonlinear part can be replaced by
\[
\sum_{U\in\min\mathcal K}m_U+m_C,
\]
with a duplicate $m_C$ omitted. Equivalently, one may retain every
support between a member of $\min\mathcal K$ and $C$.
\end{lemma}
\begin{proof}
Apply the interval assertion of Lemma~\ref{lem:support-calculus}(ii).
Lemma~\ref{f22:lem:words} supplies nonlinear idempotence.
\end{proof}

\begin{theorem}\label{f22:thm:disjointinterval}
For each $j\in\{720,738\}$, a complete eleven-identity basis for
$A_j$ is $\mathcal B$, \eqref{f22:eq:interval}, and
\begin{equation}\label{f22:eq:incident}
x+xy\eq x^2+xy.
\end{equation}
Consequently $\Var(A_{720})=\Var(A_{738})$.
\end{theorem}
\begin{proof}
Both defining tables satisfy the identities. After
Lemma~\ref{f22:lem:words}, the incident law replaces a linear $x$
belonging to a nonlinear support by $x^2$; a containing monomial
can be factored as $xv$, with $v=x$ for a singleton.
Thus the surviving linear set $L$ is disjoint from the nonlinear
content $C$. Lemma~\ref{f22:lem:interval} gives
\begin{equation}\label{f22:eq:disjointnormal}
 \sum_{x\in L}x+\sum_{U\in\mathcal H}m_U+m_C,
\end{equation}
where $\mathcal H$ is a nonempty antichain in $2^C\setminus\{\varnothing\}$.
Omit duplicate $m_C$; if $C=\varnothing$, only the nonempty linear part remains.

In both generators, $1$ is an additive identity, $0$ an absorber,
$2+2=2$, and multiplication on $\{1,2\}$ is minimum for $1<2$.
A nonlinear product involving $3$ or $0$ is zero.
Thus $\alpha_{\{x\}}(3,1)$ gives $0$, $3$, or $1$ according as
$x\in C$, $x\in L$, or $x$ is absent, recovering $C,L$.
For $D\subseteq C$, $\alpha_D(2,1)$ has value $2$ exactly when
$D$ contains a member of $\mathcal H$. The extra $m_C$ adds no minimal
accepted set. Lemma~\ref{lem:finite-set-recovery}(ii) completes recovery.
\end{proof}

\begin{theorem}\label{f22:thm:linearinterval}
A complete eleven-identity basis for $A_{1199}$ is $\mathcal B$
together with \eqref{f22:eq:interval} and \eqref{f22:eq:linearunion}.
\end{theorem}
\begin{proof}
The identities hold. Let $C$ be the full content and $L$ the linear
set after Lemma~\ref{f22:lem:words}. The linear-union law inserts
all pairs of distinct variables of $L$ and supplies linear variables
to nonlinear supports. Lemma~\ref{f22:lem:interval} then applies
with this full $C$, giving
\begin{equation}\label{f22:eq:linearnormal}
 \sum_{x\in L}x+\sum_{U\in\mathcal H}m_U+m_C,
\end{equation}
where $\mathcal H$ is a nonempty antichain of nonempty subsets of $C$,
$L\subseteq C$, and
\begin{equation}\label{f22:eq:linearconditions}
 \{x\}\notin\mathcal H\ (x\in L),\qquad
 \{x,y\}\in\mathcal H\ (x,y\in L,\ x\ne y).
\end{equation}
Indeed, no nonlinear singleton at $x\in L$ survives, so the inserted
pairs are minimal; every larger support meeting $L$ twice is redundant.
Omit duplicate $m_C$. The only case without nonlinear terms is the
single variable $x$, with $C=L=\{x\}$ and $\mathcal H=\varnothing$.

The addition of $A_{1199}$ is flat on $\{0,1,2\}$, with absorber $0$,
and $3+1=3$, $3+2=0$. Multiplication agrees with $A_{1200}$.
The tests $\alpha_{\{x\}}(0,2)$ and $\alpha_{\{x\}}(3,1)$
give respectively $0$ iff $x\in C$ and $3$ iff $x\in L$.

Let $\mathcal H_0=\{U\in\mathcal H:U\cap L=\varnothing\}$.
For $D\subseteq C\setminus L$, the test $\alpha_D(2,1)$ differs
from $1$ exactly when $D$ contains a member of $\mathcal H_0$.
Recover $\mathcal H_0$ by Lemma~\ref{lem:finite-set-recovery}(ii).
For each $x\in L$, put
\[
 \mathcal T_x=\{U\setminus\{x\}:U\in\mathcal H,\ U\cap L=\{x\}\}.
\]
Restrict $D\subseteq C\setminus L$ to sets containing no member of
$\mathcal H_0$, and use $\beta_{x,D\cup\{x\}}(3,2,1)$.
All monomials disjoint from $x$ or meeting another linear variable are
inactive; a monomial through $x$ is zero precisely when its tail lies
in $D$. The output is therefore $0$ if $D$ contains a member of
$\mathcal T_x$, and $3$ otherwise. Every actual tail is an allowed
test, since otherwise a member of $\mathcal H_0$ would be properly
contained in its corresponding support. Recover $\mathcal T_x$ by
the same lemma. The remaining members of $\mathcal H$ are the prescribed
pairs in \eqref{f22:eq:linearconditions}. In both test families $m_C$
adds no minimal accepted set. This also distinguishes the single-variable case.
\end{proof}

\endgroup
\subsection{Parity with recorded content}\label{sec:family23}
\begingroup
\providecommand{\eq}{}\renewcommand{\eq}{\approx}
\providecommand{\Var}{}\renewcommand{\Var}{\operatorname{Var}}
\providecommand{\SR}{}\renewcommand{\SR}{\mathsf{SR}}
\emph{Generators (basis size):} $A_{619}$ (11), $A_{736}$ (11), $A_{1597}$ (11), $A_{635}$ (10), $A_{765}$ (10).

\subsubsection{Parity reduction without loss of content}
In both directly treated generators $A_{736}$ and $A_{765}$,
$1$ is an additive identity, $0$ is an additive absorber, and
$\{1,2\}$ is an additive group of order two, with $2+2=1$.
Multiplication on $\{1,2\}$ is minimum for $1<2$. A product of
two or more factors involving $0$ or $3$ equals $0$. The tables
are identical except that $3+3=0$ in $A_{736}$ and $3+3=3$ in
$A_{765}$. 

Let $\mathcal C$ consist of $\SR$ and
\begin{align}
xy&\eq yx,&x^2y&\eq xy,\label{f23:eq:word}\\
x+xy&\eq x^2+xy,\label{f23:eq:incident}\\
2xy+2zu&\eq2xyzu.\label{f23:eq:guard}
\end{align}
Thus $\mathcal C$ has nine identities.

\begin{lemma}\label{f23:lem:content}
Assume \eqref{f23:eq:word}, \eqref{f23:eq:guard}, the semiring axioms, and
\begin{equation}\label{f23:eq:nonlinearperiod}
3xy\eq xy.
\end{equation}
For any nonempty sum of nonlinear monomials, let $C$ be the union
of their supports and let $\mathcal H$ be the set of supports
occurring an odd number of times. That sum is equationally equal to
\begin{equation}\label{f23:eq:guardnormal}
\sum_{U\in\mathcal H}m_U+2m_C.
\end{equation}

\end{lemma}
\begin{proof}
Apply the common word reduction. The law $3xy\eq xy$ retains one
copy at each odd support and permits a doubled copy there; even
positive coefficients reduce to two. The guard law combines all
doubled contributions into $2m_C$, leaving the odd family unchanged.
\end{proof}

For a finite family $\mathcal H$ of nonempty subsets of $C$, define
\[
F_{\mathcal H}(D)=\sum_{U\subseteq D}[U\in\mathcal H]
\quad\text{in }\mathbb F_2,\qquad D\subseteq C.
\]
Here brackets denote indicator values. Lemma~\ref{lem:finite-set-recovery}(i)
recovers $\mathcal H$ by
\begin{equation}\label{f23:eq:inversion}
[D\in\mathcal H]=F_{\mathcal H}(D)+
\sum_{U\subsetneq D}[U\in\mathcal H]\quad\text{in }\mathbb F_2.
\end{equation}
The empty-set coefficient is zero.

\Needspace{4\baselineskip}

\subsubsection{Two bases with a common content form}
\begin{theorem}\label{f23:thm:parity}\label{f23:thm:periodic}
A complete eleven-identity basis for $A_{736}$ is $\mathcal C$ together with
\begin{equation}\label{f23:eq:paritybasis}
2x\eq2x^2,\qquad3xy\eq xy.
\end{equation}
A complete ten-identity basis for $A_{765}$ is $\mathcal C$ together with
\begin{equation}\label{f23:eq:addperiod}
3x\eq x.
\end{equation}
\end{theorem}
\begin{proof}
The identities hold. In either basis, \eqref{f23:eq:incident} replaces
each linear occurrence in a nonlinear support by its square. In
$A_{736}$, pairs of linears first become doubled squares and hence
join the nonlinear part. In $A_{765}$, linears outside that part
retain coefficient one or two according to parity. Both bases imply
$3xy\eq xy$, so Lemma~\ref{f23:lem:content} gives
\begin{equation}\label{f23:eq:paritynormal}
\sum_{x\in L}c(x)x+\sum_{U\in\mathcal H}m_U+2m_C,
\qquad L\cap C=\varnothing,
\end{equation}
where $c(x)=1$ for $A_{736}$ and $c(x)\in\{1,2\}$ for $A_{765}$,
and $\mathcal H$ is any family of nonempty subsets of $C$.
If $C=\varnothing$, omit the nonlinear part and require $L\ne\varnothing$.

The test $\alpha_{\{x\}}(3,1)$ gives $0$, $3$, or $1$ according as
$x\in C$, $x\in L$, or $x$ is absent. For $A_{765}$ and $x\in L$,
$\alpha_{\{x\}}(2,1)$ gives $2$ or $1$ according as $c(x)=1$ or $2$.
Finally, for $D\subseteq C$, $\alpha_D(2,1)$ makes every linear and
the doubled content monomial equal to the additive identity $1$.
Its output is $2$ iff $F_{\mathcal H}(D)=1$. Subset inversion recovers
$\mathcal H$, completing separation in both cases.
\end{proof}

\endgroup
\subsection{Exceptional monomials}\label{sec:family24}
\begingroup
\providecommand{\eq}{}\renewcommand{\eq}{\approx}
\providecommand{\Var}{}\renewcommand{\Var}{\operatorname{Var}}
\providecommand{\SR}{}\renewcommand{\SR}{\mathsf{SR}}
\emph{Generators (basis size):} $A_{616}$ (11), $A_{735}$ (11), $A_{1197}$ (13), $A_{1587}$ (13).

\subsubsection{Common reductions}
The two direct bases contain $\SR$ and the word identities
\begin{equation}\label{f24:eq:word}
xy\eq yx,\qquad x^2y\eq xy.
\end{equation}
Both bases also contain
\begin{equation}\label{f24:eq:mixed}
x+yz\eq x+2yz.
\end{equation}
If $w$ is any nonlinear monomial, split $w$ into two nonempty
factors and substitute $x=w$ in \eqref{f24:eq:mixed}. This gives
\begin{equation}\label{f24:eq:threshold}
2w\eq3w.
\end{equation}
Consequently every multiplicity at least two reduces to two;
in particular, $2w+2w\eq2w$. 

\Needspace{4\baselineskip}

\subsubsection{Collapse to a doubled content monomial}
\begin{theorem}\label{f24:thm:collapse}
A complete eleven-identity basis for $A_{735}$ consists of
$\SR$, \eqref{f24:eq:word}, \eqref{f24:eq:mixed}, and
\begin{align}
2x&\eq2x^2,\label{f24:eq:double}\\
xy+zu&\eq2xyzu,\label{f24:eq:collapse}\\
x+xy&\eq x^2+xy.\label{f24:eq:incident}
\end{align}
\end{theorem}
\begin{proof}
The defining operations verify the identities. A nonlinear product is
minimum on $\{1,2\}$ and is zero if a factor is $0$ or $3$.
Addition has absorber $0$, sums two members of $\{1,2\}$ to $1$,
and satisfies $3+1=3$, $3+3=0$.

The collapse law gives $w+v\eq2wv$ for nonlinear $w,v$.
Writing $w=ab$ and using $x=2a$, $y=b$ also gives
$2w+v\eq4wv\eq2wv$ by \eqref{f24:eq:threshold}.
Induction therefore collapses any sum of at least two nonlinear
monomials to $2m_C$ on their combined support.
Replace repeated linear variables by doubled squares using
\eqref{f24:eq:double}, and replace every linear copy incident with a
nonlinear support by a square using \eqref{f24:eq:incident}.
The surviving linear set $L$ is disjoint from $C$. If a lone nonlinear
summand remains beside a linear one, \eqref{f24:eq:mixed} doubles it.
The forms are exactly
\begin{align}
 &\sum_{x\in L}x &&(L\ne\varnothing),\label{f24:eq:pure}\\
 &m_C &&(C\ne\varnothing),\label{f24:eq:sole}\\
 &\sum_{x\in L}x+2m_C
 &&(C\ne\varnothing,\ L\cap C=\varnothing),\label{f24:eq:dcontent}
\end{align}
with an absent linear part omitted. The test $\alpha_{\{x\}}(3,1)$
gives $0$, $3$, or $1$ according as $x\in C$, $x\in L$, or $x$ is
absent. The only remaining ambiguity, $m_C$ versus $2m_C$ with no
linears, is resolved by the uniform value $2$, giving outputs $2,1$.
\end{proof}

\Needspace{4\baselineskip}

\subsubsection{Antichains of doubled supports}
\begin{theorem}\label{f24:thm:antichain}
A complete thirteen-identity basis for $A_{1197}$ consists of
$\SR$, \eqref{f24:eq:word}, \eqref{f24:eq:mixed}, \eqref{f24:eq:double}, and
\begin{align}
xy+zu&\eq2xy+2zu,\label{f24:eq:pairdouble}\\
2xy+2xyz&\eq2xy,\label{f24:eq:absorb}\\
x+x^2&\eq2x^2,\label{f24:eq:linearsquare}\\
x+y+xy&\eq x+y.\label{f24:eq:linearpair}
\end{align}
\end{theorem}
\begin{proof}
The identities hold. A nonlinear product in $A_{1197}$ is $1$ if a
factor is $1$, $2$ if all factors are $2$, and zero otherwise; its double
is $1$ in the first case and zero otherwise. Adding either $0$ or $2$
to any element gives zero, while $1+1=1$, $1+3=3$, $3+3=0$.

Retain a sole nonlinear monomial $m_U$ as an exceptional form.
Otherwise \eqref{f24:eq:pairdouble} doubles nonlinear summands in pairs,
and \eqref{f24:eq:mixed} doubles one accompanied by a linear summand.
The threshold law removes excess copies. Repeated linears become
doubled squares by \eqref{f24:eq:double};
\eqref{f24:eq:linearsquare} then removes remaining linears at doubled
singletons. For $\varnothing\ne U\subsetneq V$, factoring $m_U$ into
two nonempty words and substituting the product on $V\setminus U$
in \eqref{f24:eq:absorb} gives
\begin{equation}\label{f24:eq:setabsorb}
 2m_U+2m_V\eq2m_U.
\end{equation}
If $U$ contains two distinct surviving linears $x,y$, insert $xy$
by \eqref{f24:eq:linearpair}, double it, absorb $2m_U$, then delete
each copy of $xy$. Finally minimize the doubled supports. The forms are
\begin{equation}\label{f24:eq:antichainnormal}
 \sum_{x\in L}x+\sum_{U\in\mathcal H}2m_U,
\end{equation}
where $\mathcal H$ is an antichain, at least one part is present, and
\begin{equation}\label{f24:eq:constraints}
 |U\cap L|\leq1\quad(U\in\mathcal H),\qquad
 \{x\}\notin\mathcal H\quad(x\in L).
\end{equation}

Uniform values $2,3$ give the output pair $(2,3)$ for a lone linear,
$(2,0)$ for a sole nonlinear monomial, and first output zero for every
other form. Thus exceptional forms are separated from the rest; their
supports are recovered by $\alpha_{\{x\}}(0,2)$.
For \eqref{f24:eq:antichainnormal}, $\alpha_{\{x\}}(3,1)$ gives
$3$ iff $x\in L$, by singleton exclusion.
For $D\subseteq X\setminus L$, $\alpha_D(2,1)$ gives zero iff
$D$ contains a member of $\mathcal H_0$, and otherwise gives $1$.
On tests avoiding $\mathcal H_0$,
$\beta_{x,D\cup\{x\}}(3,2,1)$ for $x\in L$ gives zero iff
$D$ contains a tail in $\mathcal H_x$, and otherwise gives $3$.
Lemma~\ref{lem:marked-antichain-recovery} recovers the entire antichain.
\end{proof}

\endgroup
\subsection{Linear coefficients and doubled support guards}\label{sec:family25}\label{sec:family26}
\begingroup
\emph{Generators (basis size):}
$T=A_{1198}$ (11), $A_{1797}$ (11), $R=A_{1227}$ (11),
$P=A_{1215}$ (10), $A_{730}$ (10), $M=A_{1244}$ (12), $A_{759}$ (12),
$V=A_{1206}$ (12), $A_{1793}$ (12), $Q=A_{1235}$ (10).

All six direct representatives have the same multiplication: a nonlinear
product is $1$ if a factor is $1$, $2$ if all factors are $2$, and $0$
otherwise. In addition, $1$ is an identity. On $\{0,2\}$, the first
four have $0$ absorbing and $2+2=0$, whereas $V,Q$ have the group of
order two with identity $0$. The remaining sums are
\[
\begin{array}{c|ccc}
 &3+0&3+2&3+3\\\hline
T&0&0&0\\R&0&0&3\\P&3&3&0\\M&3&3&3\\V&0&2&0\\Q&0&2&3
\end{array}
\]
All bases below contain $\SR$ and
\begin{equation}\label{f25:eq:word}
xy\eq yx,\qquad x^2y\eq xy.
\end{equation}

\subsubsection{Four bases with single and doubled antichains}
Let $\mathcal C$ consist of these seven identities and
\begin{equation}\label{f25:eq:triggerabsorb}
2xy+xyz\eq2xy.
\end{equation}
\begin{theorem}\label{f25:thm:triggers}\label{f26:thm:triggers}
The following are complete bases, with respectively eleven, eleven,
ten and twelve identities:
\begin{align}
T:&\quad\mathcal C,\quad 2x\eq2x^2,\quad x+x^2\eq2x^2,
       \quad x+y+xy\eq x+y;\label{f25:eq:triggersquare}\\
R:&\quad\mathcal C,\quad 3x\eq2x,\quad x+x^2\eq2x^2,
       \quad x+yz+xyz\eq x+yz;\label{f26:eq:R}\\
P:&\quad\mathcal C,\quad 2x\eq2x^2,
       \quad x+y+xy\eq x+y;\label{f26:eq:P}\\
M:&\quad\mathcal C,\quad 3x\eq2x,\quad x+x^2\eq2x,
       \quad2x+xy\eq2x,\quad x+y+xy\eq x+y.\label{f26:eq:M}
\end{align}
\end{theorem}
\begin{proof}
The defining operations verify the identities. Substituting $z=x$ in
\eqref{f25:eq:triggerabsorb} gives $3xy\eq2xy$; the support calculus gives
\begin{equation}\label{f25:eq:absorbset}
2m_U+m_W\eq2m_U\qquad(\varnothing\ne U\subseteq W).
\end{equation}
The mixed insertion in $R$ and the full insertion in $T,P,M$ insert
products whenever one summand is nonlinear; repetition doubles them.
In particular, $x+xy\eq x+2xy$ and
$m_U+m_W\eq m_U+2m_W$ for $U\subsetneq W$.

In $T,P$, pairs of linears become doubled squares; in $R,M$, linear
coefficients truncate at two. A nonlinear singleton removes its linear
variable in $T,R$, since $x+2x^2\eq3x^2\eq2x^2$. In $P$ a single
linear may coexist with its doubled singleton. In $M$ such a singleton
makes the linear coefficient two, and $2x+xy\eq2x$ then deletes all
incident nonlinear summands.

Write $L_i=\{x:c(x)=i\}$ and $L=L_1\cup L_2$. Double single nonlinear
supports meeting $L$ (delete those meeting $L_2$ in $M$) or properly
containing another single support. In $T,P,M$, supports meeting two
linears are deleted by inserting $2xy$, applying
\eqref{f25:eq:absorbset}, and deleting $2xy$. Insert doubled unions
of surviving singles and their unions with nonabsorbing linears;
minimize doubled supports and delete singles above them. The forms are
\begin{equation}\label{f25:eq:triggernormal}
\sum_{x\in L}c(x)x+\sum_{U\in\mathcal A}m_U+
\sum_{W\in\mathcal B}2m_W,
\end{equation}
where $\mathcal A,\mathcal B$ are antichains,
$U\cap L=\varnothing$ and $U\notin\uparrow\mathcal B$ for
$U\in\mathcal A$, and at least one part is present. Further conditions are
\begin{center}
\begin{tabular}{clll}\toprule
Case&$c(x)$&Conditions on $W\in\mathcal B$&$L_*$\\\midrule
$T$&1&$|W\cap L|\leq1$, $W\ne\{x\}$ for $x\in L$&$L$\\
$R$&1 or 2&$W\ne\{x\}$ for $x\in L$&$L$\\
$P$&1&$|W\cap L|\leq1$&$L$\\
$M$&1 or 2&$W\cap L_2=\varnothing$, $|W\cap L_1|\leq1$,&$L_1$\\
&&$W\ne\{x\}$ for $x\in L$&\\\bottomrule
\end{tabular}
\end{center}
and
\begin{equation}\label{f26:eq:unions}
U\cup W\in\uparrow\mathcal B\quad(U,W\in\mathcal A,\ U\ne W),
\qquad U\cup\{x\}\in\uparrow\mathcal B\quad(U\in\mathcal A,\ x\in L_*).
\end{equation}
Union insertion creates no singleton at a retained linear. Deleting a
guarded single also makes its previously inserted unions redundant,
so the conditions persist.

In all four algebras, $\alpha_{\{x\}}(3,1)$ gives $3$ iff $x\in L$:
$T,R$ exclude nonlinear singletons there, $P$ has $3+0=3$, and in $M$
the value $3$ absorbs nonlinear values. For $R,M$ and $x\in L$,
$\alpha_{\{x\}}(2,1)$ gives $2$ or zero according as $c(x)=1$ or $2$.

In $R$, assign $3$ on $D\cap L$, $2$ on $D\setminus L$, and $1$
elsewhere. In $T,P,M$, use $\alpha_D(2,1)$ with
$D\cap L_2=\varnothing$ and $|D\cap L_1|\leq1$. In either test the
output is zero exactly when $D\in\uparrow\mathcal B$: an active doubled
support gives zero; two active singles, or a single with an active
linear, force an active doubled support by \eqref{f26:eq:unions}.
Otherwise the output is nonzero. Every member of $\mathcal B$ is an
admissible test, so minimal zero tests recover it.

Finally, for $D\subseteq X\setminus L$ outside $\uparrow\mathcal B$,
$\alpha_D(2,1)$ gives $2$ iff a member of $\mathcal A$ lies in $D$,
and $1$ otherwise: at most one such member is active by
\eqref{f26:eq:unions}. Minimal accepted tests recover $\mathcal A$.
\end{proof}

\subsubsection{Two bases with parity and doubled guards}
Let $\mathcal D$ consist of $\SR$, \eqref{f25:eq:word}, and
\begin{equation}\label{f25:eq:guardabsorb}
x+2x^2\eq x^2,\qquad 2xy+2xyz\eq2xy.
\end{equation}
\begin{theorem}\label{f25:thm:parity}\label{f26:thm:parity}
A complete twelve-identity basis for $V$ is $\mathcal D$ together with
\begin{equation}\label{f25:eq:paritysquare}
2x\eq2x^2,\qquad3xy\eq xy,\qquad x+y+2xy\eq x+y.
\end{equation}
A complete ten-identity basis for $Q$ is $\mathcal D$ together with
\begin{equation}\label{f26:eq:Q}
3x\eq x.
\end{equation}
\end{theorem}
\begin{proof}
The identities hold. Both bases imply $3xy\eq xy$. Record the odd
contribution of each nonlinear support and supply its doubled guard by
$m_U\eq m_U+2m_U$. At a singleton guard, $x+2x^2\eq x^2$ removes
a linear copy and toggles the odd singleton contribution; restoring
the guard and repeating removes all linears there. Minimize guards
using \eqref{f25:eq:guardabsorb}. In $V$, pairs of linears become
doubled squares; a guard meeting two remaining linears is deleted by
inserting $2xy$, absorbing the guard, and deleting $2xy$. In $Q$,
linear coefficients reduce to one or two according to parity. The forms are
\begin{equation}\label{f25:eq:paritynormal}
\sum_{x\in L}c(x)x+\sum_{U\in\mathcal H}m_U+
\sum_{W\in\mathcal G}2m_W,
\end{equation}
where $\mathcal G$ is an antichain excluding $\{x\}$ for $x\in L$,
$\mathcal H$ is a family of nonempty supports, and at least one part
is present. The remaining conditions are
\begin{equation}\label{f25:eq:guardconditions}
\begin{array}{c|c|l|l}
 &c(x)&W\in\mathcal G&U\in\mathcal H\\\hline
V&1&|W\cap L|\leq1&U\in\uparrow\mathcal G\text{ or }|U\cap L|\geq2\\
Q&1\text{ or }2&\text{no further restriction}&U\in\uparrow\mathcal G
\end{array}
\end{equation}
Indeed, every odd support initially had a guard; a removed guard is
either replaced by a smaller one or, in $V$, supported on two linears.

The test $\alpha_{\{x\}}(3,1)$ gives $3$ iff $x\in L$; for $Q$,
$\alpha_{\{x\}}(2,1)$ distinguishes coefficients one and two by $2,0$.
At $\alpha_D(2,1)$ the indicator of output $2$ is
\[
\sum_{x\in D\cap L}c(x)+
\sum_{\substack{U\in\mathcal H\\U\subseteq D}}1
\quad\text{in }\mathbb F_2.
\]
Subtracting the known linear part and applying
Lemma~\ref{lem:finite-set-recovery}(i) recovers $\mathcal H$.

To recover $\mathcal G$, assign $3$ on $D\cap L$, $2$ on
$D\setminus L$, and $1$ elsewhere; restrict to $|D\cap L|\leq1$
for $V$, and allow every $D$ for $Q$. With no active guard the output
is $1$ if $D\cap L=\varnothing$ and $3$ otherwise. Indeed, the coverage
condition excludes an active odd support in this case. With an active
guard the output lies in $\{0,2\}$. Thus departure from the indicated
baseline tests $D\in\uparrow\mathcal G$. Every member of $\mathcal G$
is admissible, so its minimal accepted tests recover it.
\end{proof}
\endgroup

\subsection{Weighted support closure}\label{sec:family27}
\begingroup
\providecommand{\eq}{}\renewcommand{\eq}{\approx}
\providecommand{\Var}{}\renewcommand{\Var}{\operatorname{Var}}
\providecommand{\SR}{}\renewcommand{\SR}{\mathsf{SR}}
\emph{Generators (basis size):} $A_{2006}$ (10), $A_{2113}$ (10), $A_{2114}$ (9), $A_{2255}$ (9), $A_{2259}$ (10).

\subsubsection{Equational support saturation}
Let $\mathcal B$ consist of $\SR$ and
\begin{equation}\label{f27:eq:B}
xy\eq yx,\qquad x^2\eq x,\qquad
3x\eq2x,\qquad x+y+xy\eq x+y.
\end{equation}
Thus $\mathcal B$ has nine identities; positive coefficients reduce
to one or two.

A family $\mathcal F$ of sets is \emph{union-closed} if $U\cup V\in
\mathcal F$ for all $U,V\in\mathcal F$. Write $\min\mathcal F$
for its inclusion-minimal members. The last identity in
\eqref{f27:eq:B} inserts the product of two summands. If these summands
have distinct supports $U,V$, it therefore inserts $m_{U\cup V}$.
Repeated insertion supplies two copies of that monomial.

\begin{lemma}\label{f27:lem:weighted}
Modulo $\mathcal B$, every term has a form
\begin{equation}\label{f27:eq:weighted}
\sum_{U\in\mathcal F}c(U)m_U,
\end{equation}
where $\mathcal F$ is a nonempty union-closed family of nonempty
supports, $c(U)\in\{1,2\}$ for $U\in\min\mathcal F$, and
$c(U)=2$ for $U\in\mathcal F\setminus\min\mathcal F$.
\end{lemma}
\begin{proof}
Expand the term and reduce its coefficients. Insert unions until
its support family is union-closed; the process terminates because
the ambient variable set is finite. If $T\subsetneq U$ are present,
the product of their summands is $m_U$. Inserting it changes a
coefficient-one $m_U$ into a coefficient-two $m_U$. Thus every
nonminimal support may be doubled. Newly inserted unions cannot
create new minimal supports, and the original coefficient at every
minimal support is retained. All coefficients exceeding two are
reduced by $3x\eq2x$.
\end{proof}

\subsubsection{A complete basis from weighted union closure}
\begin{theorem}\label{f27:thm:union}
The nine identities $\mathcal B$ form a complete basis for $A_{2114}$.
\end{theorem}
\begin{proof}
The identities hold, and Lemma~\ref{f27:lem:weighted} supplies the forms.
Multiplication is a semilattice with absorber $1$; a product over
$\{0,2,3\}$ is $2$ or $3$ if all factors have that value and is zero
otherwise. Addition has identity $1$, absorber $0$, $2+2=2$, and $3+3=0$.
For $x\in D$, $\beta_{x,D}(0,2,1)$ makes a monomial zero exactly
when its support contains $x$ and lies in $D$; all other values are
$1,2$ and are unchanged by positive multiplicities. The output therefore
tests $x\in G_{\mathcal F}(D)$ and recovers $\mathcal F$.
For $U\in\min\mathcal F$, $\alpha_U(3,1)$ gives $3$ when $c(U)=1$
and $0$ when $c(U)=2$: only $m_U$ has value $3$.
These tests recover all coefficients not fixed by the normal form.
\end{proof}

\subsubsection{A complete basis from intervals above doubled supports}
Add the identity
\begin{equation}\label{f27:eq:I}
2x+yz+xy\eq2x+yz
\end{equation}
to $\mathcal B$.

\begin{theorem}\label{f27:thm:interval}
The ten identities $\mathcal B$ and \eqref{f27:eq:I} form a complete
basis for $A_{2113}$.
\end{theorem}
\begin{proof}
The identities hold and preserve content $C$. Start from
Lemma~\ref{f27:lem:weighted}. If $2m_U$ occurs and $z$ belongs to
another support $V$, write $m_V=zv$ (using $v=z$ for a singleton).
The added law inserts $m_{U\cup\{z\}}$; applying it twice doubles
this monomial. Retaining the sources saturates all doubled supports
upwards within $C$. A single support that acquires a doubled copy
becomes doubled by $3x\eq2x$.

Let $\mathcal D$ be the minimal doubled supports and $\mathcal A$
the remaining single supports. Both are antichains and
\begin{equation}\label{f27:eq:conditions}
 \mathcal A\cap\uparrow\mathcal D=\varnothing,\qquad
 U\cup V\in\uparrow\mathcal D
 \quad(U,V\in\mathcal A,\ U\ne V).
\end{equation}
The union condition follows from doubled union insertion.
When $\mathcal D\ne\varnothing$, the compressed form is
\begin{equation}\label{f27:eq:intervalnormal}
 \sum_{U\in\mathcal A}m_U+
 \sum_{V\in\mathcal D}2m_V+2m_C,
\end{equation}
omitting $2m_C$ if $C\in\mathcal D$. Retaining this full-content
monomial reconstructs every intermediate doubled support, justifying
compression by reversing insertions. If $\mathcal D=\varnothing$,
the union condition leaves a single support $C$ and form $m_C$.
The displayed conditions ensure that further insertions leave these data unchanged.

Multiplication agrees with $A_{2114}$. Addition has absorber $0$,
$2+2=2$, and restricts on $\{0,1,3\}$ to identity $1$ with $3+3=0$.
Thus $\alpha_{\{x\}}(0,2)$ gives zero iff $x\in C$.
Under $\alpha_E(3,1)$ a monomial gives $3$ iff its support lies in $E$.
An active doubled support forces zero; two active single supports do
also, but their union contains a member of $\mathcal D$ by
\eqref{f27:eq:conditions}. Hence the output is zero exactly when
$E\in\uparrow\mathcal D$, recovering $\mathcal D$.
On the remaining tests at most one member of $\mathcal A$ is active;
the output is $3$ iff one is active, and $1$ otherwise.
Every member of $\mathcal A$ is an allowed test, so minimal accepted
tests recover $\mathcal A$. All cases, including $\mathcal D=\varnothing$,
are separated.
\end{proof}

\endgroup
\subsection{Periodic support families}\label{sec:family28}
\begingroup
\providecommand{\eq}{}\renewcommand{\eq}{\approx}
\providecommand{\Var}{}\renewcommand{\Var}{\operatorname{Var}}
\providecommand{\SR}{}\renewcommand{\SR}{\mathsf{SR}}
\emph{Generators (basis size):} $A_{2119}$ (8), $A_{2123}$ (8), $A_{2251}$ (8), $A_{2263}$ (8).

\subsubsection{Common consequences of semilattice multiplication}
Let $\mathcal S$ consist of $\SR$ and
\begin{equation}\label{f28:eq:semilattice}
xy\eq yx,\qquad x^2\eq x.
\end{equation}
Thus $\mathcal S$ contains seven identities.

\begin{lemma}\label{f28:lem:common}
Modulo $\mathcal S$, every term is a nonempty sum of support
monomials. Moreover,
\begin{equation}\label{f28:eq:derived}
4x\eq2x,\qquad x+y+2xy\eq x+y.
\end{equation}
In particular, a doubled union may be inserted whenever two
summands with distinct supports are present.
\end{lemma}
\begin{proof}
Distributivity expands terms into sums of words, and
\eqref{f28:eq:semilattice} reduces each word to its support monomial.
For the two displayed consequences, expand squares:
\[
2x\eq(2x)^2\eq4x^2\eq4x,
\qquad
x+y\eq(x+y)^2\eq x^2+2xy+y^2\eq x+y+2xy.
\]
The product of monomials of supports $U,V$ is $m_{U\cup V}$.
Using the second identity on those two summands inserts
$2m_{U\cup V}$ and retains the original summands.
\end{proof}

The first consequence reduces a positive coefficient to $1$, $2$
or $3$: coefficient one stays one, a positive even coefficient
becomes two, and an odd coefficient at least three becomes three.
An additional identity $3x\eq x$ reduces all positive odd
coefficients to one.

The two direct representatives $R=A_{2119}$ and $P=A_{2123}$
have the same multiplication. A product containing $1$ is $1$;
a product entirely of $2$'s is $2$, one entirely of $3$'s is $3$,
and any other product over $\{0,2,3\}$ is $0$. The element $1$
is an additive identity in both algebras. In $R$, $\{0,2\}$ is
an additive group of order two with identity $0$, and
$3+3=3+0=0$. In $P$, $\{0,3\}$ is an additive group of order
two with identity $0$, while $2+2=2$ and $0+2=0$.

\subsubsection{Odd supports with an antichain of doubled supports}
For a family $\mathcal B$ of nonempty subsets of a finite variable
set $X$, put
\[
\uparrow\mathcal B=\{D\subseteq X:(\exists B\in\mathcal B)\ B\subseteq D\}.
\]
An antichain is recovered from its upward closure by taking minimal
members.

\begin{theorem}\label{f28:thm:antichain}
A complete eight-identity basis for $R=A_{2119}$ is
$\mathcal S$ together with
\begin{equation}\label{f28:eq:absorb}
2x+2xy\eq2x.
\end{equation}
\end{theorem}
\begin{proof}
The identities hold. For nonempty $U\subseteq V$,
\begin{equation}\label{f28:eq:supportabsorb}
 2m_U+2m_V\eq2m_U:
\end{equation}
for strict inclusion substitute $x=m_U$, $y=m_{V\setminus U}$
in \eqref{f28:eq:absorb}; for equality use $4m_U\eq2m_U$.
Expand using Lemma~\ref{f28:lem:common}, recording each coefficient
as an odd copy, a doubled copy, or both. The displayed law saturates
doubled supports upwards within a fixed finite $X$;
\eqref{f28:eq:derived} inserts doubled unions of distinct supports.
These operations preserve parity and stabilize after finitely many support insertions.
Compress the doubled family to its minimal antichain $\mathcal B$, obtaining
\begin{equation}\label{f28:eq:antiform}
 \sum_{U\in\mathcal O}m_U+\sum_{V\in\mathcal B}2m_V,
\end{equation}
with at least one part present. For $\mathcal A=\mathcal O\setminus\uparrow\mathcal B$,
\begin{equation}\label{f28:eq:antiunion}
 U\cup V\in\uparrow\mathcal B
 \quad(U,V\in\mathcal A,\ U\ne V).
\end{equation}
Thus $\mathcal A$ is an antichain. These conditions are sufficient:
unions involving a guarded support remain guarded, and unions of two
unguarded supports are guarded by the displayed condition. Further
saturation changes no recorded data. If $\mathcal B$ is empty, only
one odd monomial remains. Content need not be recorded separately,
since the absorption law may delete variables.

Write $p$ for the displayed form. At $\alpha_D(2,1)$, active monomials
give $2$ and inactive ones $1$.
The group $\{0,2\}$ gives
\begin{equation}\label{f28:eq:inversion}
 [p(\alpha_D(2,1))=2]=
 \sum_{\substack{U\in\mathcal O\\U\subseteq D}}1
 \quad\text{in }\mathbb F_2.
\end{equation}
With empty-set coefficient zero, subset inversion recovers $\mathcal O$.
At $\alpha_D(3,1)$, an active doubled support gives zero, absorbing
all other values. Without one, at most one unguarded odd support is active
by \eqref{f28:eq:antiunion}, so the output is $1$ or $3$.
Thus zero occurs exactly on $\uparrow\mathcal B$, recovering
$\mathcal B$ from the minimal zero tests.
\end{proof}

\subsubsection{Odd supports inside a union-closed family}
\begin{theorem}\label{f28:thm:union}
A complete eight-identity basis for $P=A_{2123}$ is
$\mathcal S$ together with
\begin{equation}\label{f28:eq:period}
3x\eq x.
\end{equation}
\end{theorem}
\begin{proof}
The identities hold. By Lemma~\ref{f28:lem:common}, expand into support
monomials and reduce coefficients by parity. The law $3x\eq x$
supplies a doubled copy at every odd support, and
\eqref{f28:eq:derived} supplies doubled unions. Closing under unions
and using $4x\eq2x$ gives
\begin{equation}\label{f28:eq:unionform}
 \sum_{U\in\mathcal O}m_U+\sum_{V\in\mathcal F}2m_V,
\end{equation}
where $\mathcal F$ is nonempty and union-closed and
$\mathcal O\subseteq\mathcal F$. The displayed coefficient three at
an odd support equals one; saturation leaves parity unchanged.

The test $\alpha_D(3,1)$ has output $3$ precisely when an odd number
of members of $\mathcal O$ lie in $D$, because $\{0,3\}$ is an
additive group of order two. Subset inversion recovers $\mathcal O$.
For $x\in D$, $\beta_{x,D}(0,2,1)$ gives zero exactly when
$x\in G_{\mathcal F}(D)$: a monomial vanishes exactly when its support
contains $x$ and lies in $D$, and addition on $\{0,1,2\}$ is idempotent
with identity $1$ and absorber $0$. These tests recover $\mathcal F$.
\end{proof}

The two varieties are incomparable. The identity $3x\eq x$
fails in $R$ at $x=3$, with left value $0$ and right value $3$.
The identity $2x+2xy\eq2x$ fails in $P$ at $x=2,y=3$,
with left value $0$ and right value $2$.

\endgroup
\subsection{Support partitions}\label{sec:family29}
\begingroup
\providecommand{\eq}{}\renewcommand{\eq}{\approx}
\providecommand{\Var}{}\renewcommand{\Var}{\operatorname{Var}}
\providecommand{\SR}{}\renewcommand{\SR}{\mathsf{SR}}
\emph{Generators (basis size):} $A_{617}$ (12), $A_{1228}$ (11), $A_{1242}$ (10), $A_{1610}$ (11), $A_{1808}$ (11).

\subsubsection{Support reduction and the idempotent representatives}
All proposed bases include
\begin{equation}\label{f29:eq:words}
xy\eq yx,\qquad x^2y\eq xy.
\end{equation}
Let $\mathcal C$ consist of $\SR$, \eqref{f29:eq:words}, and
\begin{equation}\label{f29:eq:ai}
x+x\eq x,\qquad x+x^2\eq x^2.
\end{equation}
It has nine identities. Modulo $\mathcal C$, all summands are
idempotent and a linear variable disappears if its nonlinear
singleton is present. Thus a polynomial may be written as
$\sum_{x\in L}x+\sum_{U\in\mathcal H}m_U$, where
$\{x\}\notin\mathcal H$ for $x\in L$.

The direct representatives $A_{1242}$ and $A_{1228}$ have the same
multiplication: a nonlinear product is $1$ if some factor is $1$,
is $2$ if every factor is $2$, and is $0$ otherwise. Both additions
are idempotent and have absorber $0$, with $1+3=3$. In $A_{1242}$,
$1$ is an additive identity and $2+3=3$. In $A_{1228}$,
$1+2=2+3=0$.

\subsubsection{Nonlinear union closure}
\begin{theorem}\label{f29:thm:union}
A complete ten-identity basis for $A_{1242}$ is $\mathcal C$
together with
\begin{equation}\label{f29:eq:union}
xy+zu+xyzu\eq xy+zu.
\end{equation}
\end{theorem}
\begin{proof}
The defining table satisfies the identities. The support calculus gives
\begin{equation}\label{f29:eq:unionnf}
 \sum_{x\in L}x+\sum_{U\in\mathcal H}m_U,
\end{equation}
where $\mathcal H$ is union-closed, $\{x\}\notin\mathcal H$ for
$x\in L$, and at least one part is nonempty. The two tests in the proof
of Theorem~\ref{f22:thm:nonlinearunion} recover $L$ and
$G_{\mathcal H}$ here as well: the operations of $A_{1242}$ and
$A_{1213}$ differ only at $3+3$, whereas those tests have at most
one summand with value $3$. They therefore give identical outcomes,
and \eqref{f22:eq:fixed} recovers $\mathcal H$.
\end{proof}

\subsubsection{Intervals retaining arbitrary linear intersections}
\begin{theorem}\label{f29:thm:interval}
A complete eleven-identity basis for $A_{1228}$ is $\mathcal C$
together with
\begin{equation}\label{f29:eq:interval}
xy+zu+xyz\eq xy+zu,\qquad
x+yz+xyz\eq x+yz.
\end{equation}
\end{theorem}
\begin{proof}
The identities hold and preserve full content $C$. Apply
Lemma~\ref{lem:support-calculus}(ii): the first added law supplies a
variable from a nonlinear summand and the second supplies one from a
linear summand. The interval closure therefore gives
\begin{equation}\label{f29:eq:intervalnf}
 \sum_{x\in L}x+\sum_{U\in\mathcal H}m_U+m_C,
\end{equation}
where $\mathcal H$ is a nonempty antichain of nonempty subsets of $C$,
$L\subseteq C$, and $\{x\}\notin\mathcal H$ for $x\in L$.
There is no bound on $|U\cap L|$. Omit duplicate $m_C$; with no
nonlinear term the form is the nonempty linear sum and $C=L$.

The tests $\alpha_{\{x\}}(0,2)$ and $\alpha_{\{x\}}(3,1)$ give,
respectively, zero iff $x\in C$ and $3$ iff $x\in L$.
For any $D\subseteq X$, assign $3$ on $D\cap L$, $2$ on
$D\setminus L$, and $1$ elsewhere. An active nonlinear support
(one contained in $D$) gives zero if it meets $L$ and $2$ otherwise;
inactive supports give $1$. If $D\cap L\ne\varnothing$, the linear
part is $3$, so the output is zero iff a nonlinear support is active.
If $D\cap L=\varnothing$, the output differs from $1$ iff one is active.
Thus these tests detect the supersets of $\mathcal H$; $m_C$ adds no
minimal accepted set. Lemma~\ref{lem:finite-set-recovery}(ii) recovers
$\mathcal H$, including the empty case.
\end{proof}

\subsubsection{A partition into single supports and one doubled block}
\begin{theorem}\label{f29:thm:partition}
A complete twelve-identity basis for $A_{617}$ consists of
$\SR$, \eqref{f29:eq:words}, and
\begin{gather}
3x\eq2x,\qquad 2x\eq2x^2,\qquad x+xy\eq2xy,
\label{f29:eq:linearpartition}\\
xy+xz\eq2xyz,\qquad
2xy+2zu\eq2xyzu.\label{f29:eq:merges}
\end{gather}
\end{theorem}
\begin{proof}
The identities hold. Reduce words to support monomials and coefficients
to one or two. Replace doubled linear variables by doubled squares;
if a surviving linear $x$ meets a nonlinear support, $x+xy\eq2xy$
removes it and doubles that support. For intersecting $U,V$, factor
$m_U=xu$, $m_V=xv$ at a common variable; the overlap law replaces
their sum by $2m_{U\cup V}$, including singleton cases. Extra copies
are absorbed because, for $U\subseteq W$,
$m_U+2m_W\eq3m_W\eq2m_W$.

Each connected component of the support-intersection graph therefore
collapses to a doubled union unless it was one undoubled support.
Distinct components are disjoint. The last merging law combines all
doubled components, giving
\begin{equation}\label{f29:eq:partitionnf}
 \sum_{x\in L}x+\sum_{U\in\mathcal H}m_U+2m_D,
\end{equation}
where $L,D$ and the nonempty members of $\mathcal H$ are mutually
disjoint, the doubled term is omitted if $D=\varnothing$, and at least
one part is present. The content is $C=L\cup D\cup\bigcup\mathcal H$.

A nonlinear product in $A_{617}$ is $1$ when all factors are $1$,
$2$ when all are $2$, and zero otherwise. Addition has absorber $0$,
$1+1=1$, $1+2=1+3=3$, and $2+2=2+3=3+3=0$.
Thus $\alpha_{\{x\}}(0,1)$ gives zero iff $x\in C$, while
$\alpha_{\{x\}}(3,1)$ gives $3$ on $L$, zero on nonlinear blocks,
and $1$ at absent variables, recovering $C,L$.

For nonempty $E\subseteq C\setminus L$, use $\alpha_E(2,1)$.
Meeting $D$, partially meeting a member of $\mathcal H$, or containing
two whole members forces zero. Otherwise the output lies in $\{2,3\}$
exactly when $E$ is one whole member of $\mathcal H$; all remaining
summands then give $1$. This recovers $\mathcal H$ and hence
$D=C\setminus(L\cup\bigcup\mathcal H)$.
\end{proof}

\endgroup
\subsection{Degree layers and pair saturation}\label{sec:family39}
\begingroup
\providecommand{\eq}{}\renewcommand{\eq}{\approx}
\providecommand{\SR}{}\renewcommand{\SR}{\mathsf{SR}}
\providecommand{\supp}{}\renewcommand{\supp}{\operatorname{supp}}
\subsubsection{Four additions on one multiplicative semigroup}

Put $M=\{0,a,a^2,a^3\}$, where $0$ is a multiplicative zero and
\begin{equation}\label{f39:eq:mult}
 a^i a^j=a^{\min(3,i+j)}\qquad(1\le i,j\le3).
\end{equation}
In each of the following additions, $0$ is neutral. On positive
powers define
\begin{equation}\label{f39:eq:additions}
\begin{aligned}
 a^i+_{\min}a^j&=a^{\min(i,j)},&
 a^i+_{\max}a^j&=a^{\max(i,j)},\\
 a^i+_{\mathrm P}a^j&=a^3,&
 a^i+_{\mathrm D}a^j&=a^{\min(3,i+j)}.
\end{aligned}
\end{equation}
Write $M_{\min},M_{\max},M_{\mathrm P},M_{\mathrm D}$ for the
resulting semirings, with local indices listed below. For all four rows, the map
$(0,a,a^2,a^3)\mapsto(1,2,3,0)$ is an isomorphism.
\emph{Generators (basis size):} $M_{\mathrm P}=A_{1260}$ (12), $M_{\mathrm D}=A_{1261}$ (14), $M_{\min}=A_{1268}$ (10), $M_{\max}=A_{1271}$ (12).

\subsubsection{Common word reduction and explicit bases}
Let $\Gamma$ consist of $\SR$ and
\begin{equation}\label{f39:eq:common}
 xy\eq yx,\qquad xyz\eq x^2yz.
\end{equation}
Every word of length at least three reduces to
\begin{equation}\label{f39:eq:cubic}
 c(C)=\prod_{x\in C}x^3,\qquad C=\supp(w)\ne\varnothing.
\end{equation}
Indeed, a selected occurrence can be the first factor in a
factorization into three nonempty words. The second equation in
\eqref{f39:eq:common} increases its exponent. Repetition raises each
exponent to at least three, and $x^4\eq x^3$, obtained by setting
$x=y=z$, reduces higher exponents. Thus a canonical word is a
linear variable, a quadratic word $xy$ (allowing $x=y$), or $c(C)$.
Call their degrees $1,2,3$, respectively. The word $c(C)$ may
have actual length greater than three; degree here denotes its layer.

The bases for the idempotent additions are $\Gamma$ together
with $x+x\eq x$ and the appropriate following list:
\begin{equation}\label{f39:eq:minbasis}
 \Sigma_{\min}:\qquad x+xy\eq x,\qquad x^2+xy\eq x^2;
\end{equation}
\begin{equation}\label{f39:eq:maxbasis}
\begin{aligned}
 \Sigma_{\max}:\qquad x+x^2&\eq x^2,&x^3+xy&\eq x^3,\\
 xyz+xyzu&\eq xyz,&x^2+y^2+xy&\eq x^2+y^2.
\end{aligned}
\end{equation}
Here and below $\Sigma$ denotes the full basis, including the
specified common identities. The sizes are $7+1+2=10$ and
$7+1+4=12$.

For the nonidempotent additions, use $\Gamma$, the three identities
\begin{equation}\label{f39:eq:top}
 x^3+x\eq x^3,\qquad x^3+xy\eq x^3,\qquad
 xyz+xyzu\eq xyz,
\end{equation}
and the appropriate following list:
\begin{equation}\label{f39:eq:pairbasis}
 \Sigma_{\mathrm P}:\qquad
 x+x\eq x^3,\qquad x+y\eq x+y+x^3y^3;
\end{equation}
\begin{equation}\label{f39:eq:degreebasis}
\begin{aligned}
 \Sigma_{\mathrm D}:\qquad x+x&\eq x^2,&x+x+x&\eq x^3,\\
 x+yz&\eq x+yz+x^3y^3z^3,\\
 x+y+z&\eq x+y+z+x^3y^3z^3.
\end{aligned}
\end{equation}
These lists have $7+3+2=12$ and $7+3+4=14$ identities.
Every displayed identity is valid by \eqref{f39:eq:mult} and
\eqref{f39:eq:additions}. In particular, the inserted product in a
saturation identity vanishes when any relevant factor is zero;
otherwise the sum on its left already has value $a^3$.
No identity uses more than four variables.

Fix a finite variable set $X$, and let $\mathcal W_X$ be the
finite set of canonical words above. For a polynomial
$p=\sum_{i\in I}w_i$ and nonempty $C\subseteq X$, let
$I_C=\{i:\supp(w_i)\subseteq C\}$. Repeated occurrences are
retained unless addition is idempotent. Empty parts of a sum
are omitted; neither an empty sum nor an empty word is substituted.

\subsubsection{The two idempotent additions}
\begin{theorem}\label{f39:thm:idem}
$\Sigma_{\min}$ and $\Sigma_{\max}$ are complete finite identity
bases for $M_{\min}$ and $M_{\max}$, respectively.
\end{theorem}
\begin{proof}
We first characterize absorption of a canonical target $w$,
of support $C$ and degree $d$, by a polynomial $p$.
For $M_{\min}$ the criterion is
\begin{equation}\label{f39:eq:mincriterion}
 p+w\eq p\quad\Longleftrightarrow\quad
 \text{some }i\in I_C\text{ has }\deg(w_i)\le d.
\end{equation}
Necessity follows by putting $a$ on $C$ and $0$ outside $C$.
For sufficiency, if the target vanishes there is nothing to
prove. Otherwise all variables of $C$ are nonzero. If $d=1$,
the eligible source is the same linear variable. If $d=3$,
the target has value $a^3$, absorbed by every nonzero source.
If $d=2$, the only value of the target below layer three is
$a^2$; every variable in its support must then have value $a$.
An eligible source of degree at most two has value at most
$a^2$ in exponent order. This proves \eqref{f39:eq:mincriterion}.

Every indicated absorption is derivable. A linear target is
already present. For a quadratic target $xy$, an eligible
linear source absorbs it by $x+xy\eq x$. Apart from the
target itself, the eligible quadratic sources are $x^2$ or
$y^2$, and \eqref{f39:eq:minbasis} handles these. This also covers
the case $x=y$. For a cubic target, multiply the eligible
source by a nonempty word until its support is $C$ and its
length is at least three. Factor absorption inserts that
word, and the common word reduction identifies it with $c(C)$.
For a cubic source already on $C$, additive idempotence suffices.

For $M_{\max}$ the criteria are as follows:
\begin{equation}\label{f39:eq:maxcriterion}
\begin{array}{c|l}
d&\text{Condition for }p+w\eq p\\\hline
1&I_C\ne\varnothing\\
3&\text{an eligible cubic word exists}\\
2&\text{an eligible cubic word exists, or the eligible}\\[-2pt]
 &\text{quadratic supports have union }C.
\end{array}
\end{equation}
For a linear target the only possible eligible supports are
its singleton support, and every positive power dominates that
variable. For a cubic target, evaluation at $a$ on $C$ forces
an eligible cubic word; such a word is also sufficient.
For a quadratic target, assume no eligible cubic exists.
The all-$a$ assignment requires an eligible quadratic. Changing
one target variable to $a^2$ raises the target to $a^3$; an
eligible quadratic must contain that variable. This proves
necessity of the union condition. Conversely, if the target
is $a^2$, an eligible quadratic supplies $a^2$. If the target
is $a^3$, some variable of its support is above $a$, and an
eligible quadratic containing it supplies $a^3$.

These absorptions are derivable from \eqref{f39:eq:maxbasis}.
The power-ascent equation and its product by $x$ give
$x+x^2\eq x^2$, $x^2+x^3\eq x^3$, and hence
$x+x^3\eq x^3$. A cubic source extends to any larger
cubic support by $xyz+xyzu\eq xyz$. For a quadratic target
$xy$, a cubic source can have support $\{x\}$, $\{y\}$,
or $\{x,y\}$. The first two cases use $x^3+xy\eq x^3$.
In the third case the source is $(xy)^3$, and power ascent
with the substitution $x\mapsto xy$ absorbs $xy$.
When $x=y$, use $x^2+x^3\eq x^3$.
Without a cubic source, a quadratic target is already present,
or $x\ne y$ and both squares $x^2,y^2$ occur. The last
identity in \eqref{f39:eq:maxbasis} inserts $xy$.

For either algebra, let $D(p)$ be the set of all words in
$\mathcal W_X$ absorbed by $p$. Each original word belongs
to $D(p)$ by additive idempotence. The derivations just given
therefore reduce $p$ to the sum of all words in $D(p)$.
The term function determines $D(p)$ itself. Thus two sides
of any valid identity, on a common finite variable set,
reduce to the same finite sum. This proves completeness
for arbitrary numbers of variables.
\end{proof}

\subsubsection{The two nonidempotent additions}
For a polynomial $p$ define an upward-closed family of supports
$\mathcal S(p)$ by the appropriate condition
\begin{equation}\label{f39:eq:saturation}
\begin{aligned}
C\in\mathcal S_{\mathrm P}(p)
 &\Longleftrightarrow
 \bigl(\exists i\in I_C:\deg(w_i)=3\bigr)\ \text{or}\ |I_C|\ge2,\\
C\in\mathcal S_{\mathrm D}(p)
 &\Longleftrightarrow \sum_{i\in I_C}\deg(w_i)\ge3.
\end{aligned}
\end{equation}
The cardinality and sum count occurrences, not just distinct words.

\begin{theorem}\label{f39:thm:weighted}
$\Sigma_{\mathrm P}$ and $\Sigma_{\mathrm D}$ are complete
finite identity bases for $M_{\mathrm P}$ and $M_{\mathrm D}$.
In each case $p+c(C)\eq p$ holds exactly when
$C\in\mathcal S(p)$, and every such absorption is derivable.
\end{theorem}
\begin{proof}
Put $a$ on $C$ and $0$ outside $C$. The value of $p$ is
$a^3$ exactly under the corresponding condition
\eqref{f39:eq:saturation}; since adding $a^3$ changes every
other value, this proves necessity. If the condition holds
and all variables of $C$ are nonzero, the indicated witnesses
make $p$ equal to $a^3$ at that assignment. If some variable
of $C$ is zero, $c(C)$ is zero. This proves sufficiency.

We record consequences of \eqref{f39:eq:top}. A cubic word of
support $D$ absorbs every canonical word whose support
contains $D$. For a cubic target this follows by support
extension, factoring the source into three nonempty words.
For a linear target only a singleton $D$ is possible, so
use $x^3+x\eq x^3$. For a quadratic target $xy$, singleton
supports use $x^3+xy\eq x^3$, and support $\{x,y\}$ uses
$u^3+u\eq u^3$ at $u=xy$. The repeated-variable case
uses $u=x^2$ and common word reduction. Cubic words are
additively idempotent: substitute $u=x^3$ in
$u^3+u\eq u^3$ and reduce, then substitute a word for $x$.

To derive a required cubic insertion, an eligible cubic
source extends to $C$. Otherwise, for $M_{\mathrm P}$
choose two eligible occurrences. Substitute their words
into the second equation in \eqref{f39:eq:pairbasis}; the
inserted cube has their union support, and then extends
to $C$. For $M_{\mathrm D}$, either two eligible occurrences
have total degree at least three, or three eligible linear
occurrences exist. In the former case at least one word
has length at least two. Factor that word into two nonempty
words and use the mixed equation in \eqref{f39:eq:degreebasis}.
In the latter case use its three-summand equation. Repeated
occurrences are legitimate separate summands. Again the
inserted cubic has the union support and extends to $C$.
All intermediate insertions are absorbed in the original
context, so they may be removed after inserting $c(C)$.

Insert $c(C)$ for every $C\in\mathcal S(p)$. Remove all
original cubic summands and all lower-degree occurrences
whose own support is in $\mathcal S(p)$, using the cubic
absorption consequences above. For $M_{\mathrm P}$ each
remaining word occurs once, since two copies would saturate
their own support. For $M_{\mathrm D}$ the only possible
remaining repetition is two copies of a linear variable;
replace them by its square using $x+x\eq x^2$. A repeated
quadratic or three linear copies would saturate their own
support. We have derived the finite, nonempty expression
\begin{equation}\label{f39:eq:normal}
 N(p)=\sum_{C\in\mathcal S(p)}c(C)+q(p),
\end{equation}
where $q(p)$ is a sum of distinct linear and quadratic words,
all with unsaturated own supports.

It remains to prove uniqueness from the function. For
nonempty $C\subseteq X$ evaluate at $a$ on $C$ and $0$
outside $C$, and encode $0,a,a^2,a^3$ by $0,1,2,3$.
Call the resulting integer $F_p(C)$. By
\eqref{f39:eq:saturation}, $\mathcal S(p)$ is precisely
the set on which $F_p(C)=3$.
If $C\notin\mathcal S(p)$, no subset of $C$ is saturated.
For $M_{\mathrm P}$ there is at most one eligible lower
word; for $M_{\mathrm D}$ their total degree is at most two.
Consequently, in either case
\begin{equation}\label{f39:eq:triangular}
 F_p(C)=\sum_{\substack{w\text{ in }q(p)\\\supp(w)\subseteq C}}
 \deg(w)\qquad(C\notin\mathcal S(p)),
\end{equation}
as an ordinary integer equality without truncation.
Lemma~\ref{lem:finite-set-recovery}(i) recovers the total degree
on each exact unsaturated support. On $\{x\}$, totals one and two mean
$x$ and $x^2$, respectively. On $\{x,y\}$ with $x\ne y$,
the only possible nonzero total on that exact support is
two, meaning $xy$. Larger exact supports contribute no
lower word. Thus the function determines $q(p)$ uniquely.
The normal-form criterion in Lemma~\ref{lem:normal-form-completeness}
now proves completeness.
\end{proof}

\endgroup
\section{Quadratic graphs and multiplicity forms}\label{sec:quadratic-forms}
After reducing word lengths, the nonlinear summands in these
families are represented by pairs of variables and their multiplicities.
The proofs identify the permitted absorptions and then use separating
evaluations to establish completeness.

\subsection{Short-polynomial joins}\label{sec:family14}
\begingroup
\providecommand{\eq}{}\renewcommand{\eq}{\approx}
\providecommand{\Var}{}\renewcommand{\Var}{\operatorname{Var}}
\providecommand{\SR}{}\renewcommand{\SR}{\mathsf{SR}}
\providecommand{\F}{}\renewcommand{\F}{\mathbb F_2}
\providecommand{\D}{}\renewcommand{\D}{\mathcal D}
\providecommand{\Qc}{}\renewcommand{\Qc}{Q_{\mathrm c}}
\providecommand{\Qf}{}\renewcommand{\Qf}{Q_{\mathrm f}}
\subsubsection{Generators}
Put $h(x)=x^3$ throughout this subsection.

Let $D_2$ have addition join and multiplication meet on a two-element
chain, and write $\D=\Var(D_2)$. Let $\F$ have its usual field
operations. Define $\Qc,\Qf$ on $\{0,1,2\}$ by the following tables,
encoded row by row:
\begin{center}
\begin{tabular}{cll}\toprule
Algebra&Addition&Multiplication\\\midrule
$\F$&\texttt{0110}&\texttt{0001}\\
$D_2$&\texttt{0111}&\texttt{0001}\\
$\Qc$&\texttt{000012022}&\texttt{000020000}\\
$\Qf$&\texttt{000010002}&\texttt{000020000}\\\bottomrule
\end{tabular}
\end{center}
Thus both $Q$-algebras have only one nonzero product, $1\cdot1=2$,
and every triple product is $0$. In $\Qc$ the additive order is
$1<2<0$; in $\Qf$, $1$ and $2$ are incomparable and their sum is
the top element $0$. These labels name table entries, not constant
symbols of the language. All four algebras are semirings.

We study
\[
 \Var(\F)\vee\Var(Q),\qquad \D\vee\Var(Q)
 \qquad(Q\in\{\Qc,\Qf\}).
\]

\subsubsection{Common identities and short-polynomial invariants}
Let $\mathcal H$ consist of $\SR$, commutativity of multiplication,
and the seven identities
\begin{align}
 h(h(x))&\eq h(x),\label{f14:eq:ret}\\
 h(x+y)&\eq h(x)+h(y),\label{f14:eq:addendo}\\
 h(xy)&\eq h(x)h(y),\label{f14:eq:mulendo}\\
 x+h(y)&\eq h(x)+h(y),\label{f14:eq:addideal}\\
 xh(y)&\eq h(x)h(y),\label{f14:eq:rightideal}\\
 h(x)y&\eq h(x)h(y),\label{f14:eq:leftideal}\\
 xyz&\eq h(xyz).\label{f14:eq:triple}
\end{align}
Hence $|\mathcal H|=13$. All these identities hold in each of the
four joins: on the first generator $h$ is the identity, while on
either $Q$-generator $h$ is the constant absorbing element $0$.
In a model they make $h$ an endomorphic retraction onto an ideal
for both operations. A term fixed by $h$ will be called an image term.
The ideal identities propagate this property from a subterm to the
whole term.

If a word has length at least three, \eqref{f14:eq:triple} and the ideal
identities make the entire sum an image term. Otherwise it is a
\emph{short polynomial}: a nonempty sum of variables and commutative
quadratic monomials. Repetition of summands is allowed.

For a short polynomial $P$, let $L(P)$ be the set of variables occurring
as linear monomials and $N(P)$ the set occurring in quadratic monomials.
Put $T(P)=L(P)\cup N(P)$. In the field case, let $b(P)$ be its
squarefree Boolean polynomial: reduce coefficients modulo two and
replace $x^2$ by $x$. This polynomial can be zero. In the lattice
case, let $m(P)$ be the family of inclusion-minimal variable supports
of its monomials. A linear term or a square contributes a singleton;
a product of distinct variables contributes a two-element set.

\begin{lemma}\label{f14:lem:invariants}
The following characterize identities between short polynomials.
\begin{enumerate}
\item In $\F$, equality is equivalent to equality of $b(P)$; in
$D_2$, it is equivalent to equality of $m(P)$.
\item In $\Qc$, a linear polynomial is distinguished from one with a
quadratic monomial. Among linear polynomials, equality is equivalent
to equality of $T(P)$. Among polynomials with a quadratic monomial,
it is equivalent to equality of both $T(P)$ and $N(P)$.
\item In $\Qf$, every polynomial with $L(P)\cap N(P)\ne\varnothing$
is identically $0$. Otherwise it is not identically $0$. Among these
nonzero polynomials, equality is equivalent to equality of $L(P)$
and $N(P)$.
\end{enumerate}
\end{lemma}
\begin{proof}
Over $\F$, the squarefree polynomial of degree at most two is uniquely
determined by evaluations: assignments with just one variable equal
to $1$ determine the linear coefficients, and assignments with just
two variables equal to $1$ then determine the quadratic coefficients.
For $D_2$, a polynomial takes value $1$ exactly when the set of
variables assigned $1$ contains some monomial support. The minimal
such sets are precisely $m(P)$.

In $\Qc$, a short polynomial with a quadratic monomial has value $2$
exactly when every variable in $N(P)$ is $1$ and every variable in
$T(P)\setminus N(P)$ belongs to $\{1,2\}$. In all other cases it
has value $0$. Assigning all variables $1$ distinguishes it from a
linear polynomial. Starting from that assignment, changing one variable
to $0$ detects its membership in $T(P)$, and changing it to $2$ detects
its membership in $N(P)$. Linear polynomial values similarly determine
and are determined by their supports.

In $\Qf$, if a variable occurs both linearly and quadratically, its
linear summand conflicts with any nonzero quadratic value. The entire
polynomial is therefore $0$. If $L(P),N(P)$ are disjoint and $N(P)$
is nonempty, its value is $2$ exactly when all variables in $N(P)$
are $1$ and all variables in $L(P)$ are $2$; otherwise it is $0$.
This provides a nonzero evaluation and determines both sets. More
explicitly, equality of two such functions first forces their total
supports to agree by changing a single variable to $0$ in a common
nonzero evaluation. Their prescribed values $1$ and $2$ then force
the two roles to agree. Linear polynomials are separated from pure
quadratic and mixed polynomials by the all-$1$ assignment, which
gives respectively $1$, $2$, and $0$. Their own identities are exactly
equality of their linear supports.
\end{proof}

The nonzero short polynomials in this lemma are also separated from
image terms on the corresponding $Q$-generator. This observation will
handle the potentially troublesome mixed identities.

\subsubsection{Joins with the field of two elements}
Adjoin to $\mathcal H$ the following three image identities:
\begin{align}
 h(x)^2&\eq h(x),\label{f14:eq:fbool}\\
 2h(x)&\eq2h(y),\label{f14:eq:fzero}\\
 h(x)+2h(y)&\eq h(x).\label{f14:eq:fneutral}
\end{align}
They put the image in $\Var(\F)$. Indeed, in the image the common
double is an additive identity and every element is its own additive
inverse. The image is a Boolean ring. Distributivity, coefficient
reduction modulo two, and $x^2=x$ give squarefree polynomial normal
forms. A nonzero squarefree polynomial is separated from zero on
$\F$ by induction on its number of variables, writing it as
$U+xV$. Thus these identities, together with the inherited semiring
axioms and multiplicative commutativity, form a complete basis for
the image theory in the constant-free signature.

\begin{theorem}\label{f14:thm:field}
A basis of $\Var(\F)\vee\Var(\Qc)$ consists of $\mathcal H$,
\eqref{f14:eq:fbool}--\eqref{f14:eq:fneutral}, and
\begin{align}
 3x&\eq x,\label{f14:eq:period}\\
 2xy&\eq2x^2+2y^2,\label{f14:eq:evenedge}\\
 x^2&\eq x+2x^2,\label{f14:eq:squarelinear}\\
 2x+2x^2&\eq2x^2.\label{f14:eq:evenshadow}
\end{align}
This basis has twenty identities. A basis of
$\Var(\F)\vee\Var(\Qf)$ has nineteen identities: keep
$\mathcal H$, \eqref{f14:eq:fbool}--\eqref{f14:eq:fneutral},
\eqref{f14:eq:period}--\eqref{f14:eq:evenedge}, and replace the last two
identities by
\begin{equation}\label{f14:eq:overlap}
 x+xy\eq h(x+xy).
\end{equation}
\end{theorem}
\begin{proof}
Each displayed identity holds on its two specified generators, so
soundness follows. We prove completeness by polynomial normalization.
Equation~\eqref{f14:eq:period} gives $3U\eq U$ and $4U\eq2U$ for
every term $U$. A positive coefficient $c$ of a monomial $W$ can
therefore be written as its parity contribution, either $W$ or an
omitted term, plus $2W$. A coefficient with even parity retains only
the double. Equation~\eqref{f14:eq:evenedge} replaces a doubled edge by
the doubles of its endpoint squares.

For $\Qc$, \eqref{f14:eq:squarelinear} changes the parity contribution
of a square into that of its linear variable, retaining its doubled
square. Equation~\eqref{f14:eq:evenshadow} allows a doubled linear term
to be freely adjoined when its doubled square is present. Every short
polynomial consequently reduces to
\begin{equation}\label{f14:eq:fieldchainnf}
 \sum_{W\in b(P)}W+
 \sum_{x\in T(P)}2x+
 \sum_{x\in N(P)}2x^2.
\end{equation}
The first sum lists the monomials with coefficient one in the Boolean
polynomial. To see that collection does not lose information, start
monomial by monomial. A linear monomial supplies its parity contribution
and its double; a square supplies its linear parity contribution and
its doubled square; an edge supplies its parity contribution and doubled
squares of both endpoints. Add the doubled linear terms at every
nonlinear endpoint using \eqref{f14:eq:evenshadow}. When parity terms
cancel in pairs, their doubles are already present, or split into
already present doubled squares by \eqref{f14:eq:evenedge}. Repeated
doubles merge by $4U\eq2U$. This proves the stated reduction using
only the finite listed identities.

In the flat case, if $L(P)\cap N(P)$ is nonempty, select a linear
occurrence $x$ and a quadratic occurrence containing $x$.
Multiplicative commutativity orients the latter as $xy$.
Equation~\eqref{f14:eq:overlap} makes their sum an image term, and the
ideal identities do the same for all of $P$. Otherwise the two
supports are disjoint. Keep squares as squares and reduce coefficients
as before. This gives
\begin{equation}\label{f14:eq:fieldflatnf}
 B_{L,N}(P)+\sum_{x\in L(P)}2x+\sum_{x\in N(P)}2x^2,
\end{equation}
where $B_{L,N}(P)$ is obtained from $b(P)$ by writing a singleton
monomial $x$ as $x$ when $x\in L(P)$, and as $x^2$ when
$x\in N(P)$. Quadratic monomials with distinct endpoints are unchanged.
The roles are unambiguous because $L(P)\cap N(P)$ is empty.
Every nonlinear occurrence contributes doubled endpoint squares using
\eqref{f14:eq:evenedge}; collection of parity terms and repeated doubles
is therefore justified exactly as above.

All sums in \eqref{f14:eq:fieldchainnf}--\eqref{f14:eq:fieldflatnf} are
subject to omission when empty, not to insertion of a zero constant.
The entire expression remains nonempty because the original polynomial
has nonempty total support and contributes at least one doubled term.

Consider any identity of the relevant join in a model of the proposed
basis. Normalize its sides. If both are image terms, it follows from
the image basis and endomorphism of $h$. If both are short nonimage
terms, Lemma~\ref{f14:lem:invariants} gives identical invariants on their
two generators. Their displayed normal forms therefore coincide.
There is no mixed case: a nonimage short polynomial has a nonzero
evaluation on its $Q$-generator, where every image term is zero.
Thus every identity of the join follows from the proposed finite basis.
\end{proof}

\subsubsection{Joins with the two-element distributive lattice}
For the lattice image, use
\begin{align}
 2h(x)&\eq h(x),\label{f14:eq:dadd}\\
 h(x)^2&\eq h(x),\label{f14:eq:dmul}\\
 h(x)+h(x)h(y)&\eq h(x).\label{f14:eq:dabs}
\end{align}
Together with the inherited axioms these define distributive lattices:
the second absorption law follows by distribution from the first and
multiplicative idempotence. Distributive expansion and removal of
nonminimal monomial supports give a complete basis for $\D$.
The two-element lattice separates distinct minimal-support families,
so it generates the variety.

\Needspace{4\baselineskip}
\begin{theorem}\label{f14:thm:lattice}
The variety $\D\vee\Var(\Qc)$ has the following nineteen-identity
basis: $\mathcal H$, \eqref{f14:eq:dadd}--\eqref{f14:eq:dabs}, and
\begin{align}
 2x&\eq x,\label{f14:eq:ai}\\
 x^2&\eq x+x^2,\label{f14:eq:looplinear}\\
 x+xu+yv&\eq x+xu+yv+xy.\label{f14:eq:edgewitness}
\end{align}
The variety $\D\vee\Var(\Qf)$ also has a nineteen-identity basis:
$\mathcal H$, \eqref{f14:eq:dadd}--\eqref{f14:eq:dabs}, \eqref{f14:eq:ai},
\eqref{f14:eq:overlap}, and
\begin{equation}\label{f14:eq:loopedge}
 x^2+yv\eq x^2+yv+xy.
\end{equation}
\end{theorem}
\begin{proof}
All listed identities hold on the specified generators. Expand terms
and remove repeated summands using \eqref{f14:eq:ai}. Terms with a long
word are image terms. In the flat case, the same holds for a short
polynomial with overlapping linear and nonlinear supports by
\eqref{f14:eq:overlap}. All other short polynomials have nonzero evaluations
on the $Q$-generator and cannot equal an image term.

Let $P\eq Q$ be a valid identity between two remaining short
polynomials. In an additive semilattice it suffices to adjoin every
monomial of $Q$ to $P$ without changing $P$, and conversely. Indeed
this gives $P\eq P+Q\eq Q$. Suppose first that $Q$ contains a linear
monomial $y$. Equality on $D_2$ implies that $P$ already contains
a monomial with support $\{y\}$, hence $y$ or $y^2$. In the chain
case \eqref{f14:eq:looplinear} adjoins $y$ if needed. In the flat case
Lemma~\ref{f14:lem:invariants} gives $L(P)=L(Q)$, so $y$ is already present.

Now let $xy$ be a quadratic monomial of $Q$, with $x=y$ allowed.
Equality on $D_2$ implies that $P$ contains either the same quadratic
monomial, up to commutativity, or a singleton-support monomial at one
of its endpoints. The first case needs no change. In the other case
rename endpoints so that the singleton is $x$ or $x^2$.
In the chain case, $N(P)=N(Q)$, so both $x$ and $y$ have quadratic
occurrences in $P$, say $xu$ and $yv$. If necessary adjoin the linear
term $x$ using \eqref{f14:eq:looplinear}; then \eqref{f14:eq:edgewitness}
adjoins $xy$. Coincident witness monomials may be duplicated using
additive idempotence before applying the rule.

In the flat case, $L(P)$ and $N(P)$ are disjoint and both endpoints
belong to $N(P)=N(Q)$. Thus the singleton monomial at the chosen
endpoint cannot be linear: it is $x^2$. A quadratic occurrence $yv$
exists, and \eqref{f14:eq:loopedge} adjoins $xy$. Again coincident
witnesses cause no difficulty. This proves every required adjunction
in both directions. Image/image identities follow from the image basis,
and mixed identities are excluded as above. Completeness follows.
\end{proof}

\subsubsection{A nilpotent specialization}
Let $T$ be the four-element semiring with table strings
\[
 +:\ \texttt{0000012002200003},\qquad
 \cdot:\ \texttt{0000033003000000}.
\]
Thus $0$ is absorbing for both operations, $1+2=2$, and $3$ is
incomparable with $1,2$ below the additive top $0$. Its nonzero
products are $1\cdot1=1\cdot2=2\cdot1=3$.

\begin{theorem}\label{f14:thm:T}
The variety generated by $T$ is exactly
\[
 \Var(T)=\bigl(\D\vee\Var(\Qf)\bigr)
          \cap\operatorname{Mod}(x^3\eq y^3).
\]
Consequently a complete basis for $T$ consists of the nineteen
identities in the flat case of Theorem~\ref{f14:thm:lattice}, together
with $x^3\eq y^3$.
\end{theorem}
\begin{proof}
Direct evaluation in the defining table verifies these twenty
identities. Let $B$ be any model of them. Its image $h(B)$ is a
singleton, so all image terms agree. As in the flat case above,
every term with a long word or a linear/nonlinear overlap reduces
to an image term. A remaining short polynomial has disjoint supports
$L,N$. Pure linear polynomials have exactly their semilattice support
identities. A polynomial with a quadratic monomial takes value $3$
in $T$ exactly when every variable in $L$ is assigned $3$, every
variable in $N$ is assigned $1$ or $2$, and no quadratic monomial
has both endpoints assigned $2$. A loop therefore forces its vertex
to have value $1$. In every other assignment its value is $0$.

Such a polynomial is not identically zero: assign $L$ the value $3$
and $N$ the value $1$. Its function determines $L$ and $N$. Changing
a variable to $0$ in a nonzero assignment detects total support;
assignments with linear values $3$ and nonlinear values $1$ distinguish
the two roles. For fixed $L,N$, view the quadratic monomials as an
undirected graph with possible loops on $N$. The sets of vertices
that may simultaneously be assigned $2$ are exactly the sets
containing no monomial support. Their minimal forbidden subsets are
the loop singletons and the edges not containing a loop vertex.
They determine exactly the inclusion-minimal supports of the
quadratic monomials. Explicitly, if the minimal-support families differ,
choose a minimal forbidden set in one family that contains no forbidden
set from the other; assigning precisely its vertices $2$ separates
the two functions.

Since $L\cap N$ is empty, adjoining the singleton supports from $L$
shows that these invariants are exactly $L(P),N(P),m(P)$ from
Lemma~\ref{f14:lem:invariants}. Pure linear, pure quadratic, and mixed
polynomials are distinguished by the all-$1$ assignment, whose values
are $1$, $3$, and $0$, respectively. Thus every identity between
nonimage short polynomials of $T$ is an identity of
$\D\vee\Var(\Qf)$ and follows from Theorem~\ref{f14:thm:lattice}.
No such polynomial equals an image term, because it has a nonzero
evaluation whereas all image terms of $T$ are $0$. The image/image
case was already handled, completing the proof.
\end{proof}

The equality is proved here, rather than inferred merely from
membership of $T$ in a finitely based variety. In particular,
$D_2\notin\Var(T)$: the constant-cube identity holds in $T$ and fails
in $D_2$. Thus $T$ does not generate the whole join.

\subsubsection{Applications and explicit maps}
Table~\ref{f14:tab:targets} gives the four exact join realizations and the
nilpotent specialization. All subscripts are local enumeration labels;
$A_{513}=T$. For the four joins, the table gives onto homomorphisms
$f:A_j\to K$ and $g:A_j\to Q$, where $K$ is the indicated field
or lattice generator. Each map is encoded by its four values.
The pairs $(f(a),g(a))$ are all distinct. Therefore
\[
 A_j\hookrightarrow K\times Q,
 \qquad \Var(A_j)=\Var(K)\vee\Var(Q).
\]
The forward inclusion uses the embedding, and the reverse inclusion
uses the two onto projections. Theorems~\ref{f14:thm:field} and
\ref{f14:thm:lattice} establish FB for these four classes, and
Theorem~\ref{f14:thm:T} establishes FB for $A_{513}$.

\begin{table}[htbp]\centering\small
\caption{Five applications and their homomorphisms.}\label{f14:tab:targets}
\setlength{\tabcolsep}{3pt}
\begin{tabular}{rllcll}\toprule
$j$&Addition&Multiplication&Variety&$f$&$g$\\\midrule
513&\texttt{0000012002200003}&\texttt{0000033003000000}&$\Var(T)$&---&---\\
659&\texttt{0020012022020023}&\texttt{0000030000200000}&$\Var(\F,\Qf)$&\texttt{0010}&\texttt{0102}\\
661&\texttt{0020012322020323}&\texttt{0000030000200000}&$\Var(\F,\Qc)$&\texttt{0010}&\texttt{0102}\\
662&\texttt{0020012322220323}&\texttt{0000030000200000}&$\Var(D_2,\Qc)$&\texttt{0010}&\texttt{0102}\\
1269&\texttt{0000010000230033}&\texttt{0100111101300100}&$\Var(D_2,\Qc)$&\texttt{1011}&\texttt{0012}\\
\bottomrule\end{tabular}
\end{table}

\endgroup
\subsection{Quadratic graphs and two-colourings}\label{sec:family15}
\begingroup
\providecommand{\eq}{}\renewcommand{\eq}{\approx}
\providecommand{\Var}{}\renewcommand{\Var}{\operatorname{Var}}
\providecommand{\AI}{}\renewcommand{\AI}{\mathsf{AI}}
\subsubsection{The six semirings}
Put $\AI=\SR\cup\{x+x\eq x\}$.
For the six ai-semirings below, multiplication
is commutative except in $A_{320}$ and $A_{486}$, which are opposites.

\emph{Generators (basis size):} $A_{320}$ (10), $A_{486}$ (10), $A_{495}$ (13), $A_{515}$ (10), $A_{516}$ (10), $A_{835}$ (12).

All triple products are the table element $0$. For $A_{320},A_{486},
A_{515}$ the additive order is $0<3<2<1$. For $A_{516}$ it is
$1<2<3<0$. For $A_{495},A_{835}$ it has $1<3<0$ and $2<0$, with
$2$ incomparable with $1,3$. The element $0$ absorbs multiplication
in every case.

\Needspace{4\baselineskip}

\subsubsection{Polynomial and graph conventions}
A \emph{short polynomial}
has only linear or quadratic words. Let $L$ be its linear support
and $N$ the set of variables in its quadratic words. Quadratic words
form a directed graph in the noncommutative cases and an undirected
graph, with loops allowed, in the commutative cases. Vertices in
$L\cap N$ are called marked. A square is a loop.

We use the alternative identities
\begin{align}
 xyz+u&\eq u,\label{f15:eq:bottom}\\
 xyz+u&\eq xyz.\label{f15:eq:top}
\end{align}
Under \eqref{f15:eq:bottom}, every triple word is an additive identity;
under \eqref{f15:eq:top}, every triple word is an additive absorbing
element. Commutativity of addition shows that any two triple words
are equal, by comparing their sum. Substitution shows that every
longer word is equal to the same term. Thus in either case $x^3\eq y^3$
is derivable. We refer to this common term as the zero term, but do
not add a constant symbol. Under \eqref{f15:eq:bottom} long words may
be deleted from a polynomial; if none remain, the result is the zero
term. Under \eqref{f15:eq:top} any long word makes the entire polynomial
the zero term.

\subsubsection{Triple products at the additive minimum}
\begin{theorem}\label{f15:thm:bottom}
A basis of $A_{320}$ consists of $\AI$, \eqref{f15:eq:bottom}, and
\begin{align}
 x+xy&\eq x,\label{f15:eq:absright}\\
 x+yx&\eq x,\label{f15:eq:absleft}\\
 x^2+xy&\eq x^2.\label{f15:eq:out}
\end{align}
A basis of $A_{486}$ replaces \eqref{f15:eq:out} by
\begin{equation}\label{f15:eq:in}x^2+yx\eq x^2.\end{equation}
A basis of $A_{515}$ consists of $\AI$, multiplicative commutativity,
\eqref{f15:eq:bottom}, \eqref{f15:eq:absright}, and
\begin{equation}\label{f15:eq:twoloops}
 x^2+y^2+xy\eq x^2+y^2.
\end{equation}
Each basis has ten identities.
\end{theorem}
\begin{proof}
The defining tables in Appendix~\ref{app:catalogue} verify soundness.
We give canonical forms and separating evaluations to prove completeness.

For $A_{320}$, delete long words. Keep the set $L$ of linear
monomials and delete every quadratic monomial incident with $L$ by
\eqref{f15:eq:absright}--\eqref{f15:eq:absleft}. Among the remaining vertices
keep all loops, and use \eqref{f15:eq:out} to delete every nonloop edge
whose initial vertex has a loop. No other edges are deleted. All
these are applications of the listed identities, so they give a
canonical polynomial determined by a linear set, a loop set, and a
set of directed edges with no loop at their initial vertices.

These data are determined by the term function. To detect $x\in L$,
assign $x=1$ and all other variables $0$: the value is $1$ exactly
when $x\in L$. Once $L$ is known, the same assignment at a vertex
outside $L$ has value $3$ exactly when it has a loop. To detect a
remaining directed edge $xy$ with $x\ne y$, assign $x=1$, $y=2$,
and all other variables $0$. Its initial vertex is unlooped, and a
loop at $y$ contributes $0$. The reversed edge also contributes $0$.
Thus the value is $3$ exactly when $xy$ is present. These tests
recover every part of the canonical form. The zero term is the form
with no surviving monomial. Equal term functions therefore have
identical canonical forms. The proof for $A_{486}$ reverses all words
and uses \eqref{f15:eq:in}.

For $A_{515}$, delete the edges incident with $L$, and then delete
every nonloop edge both of whose endpoints have loops, using
\eqref{f15:eq:twoloops}. Keep all remaining loops and edges. The same
one-variable assignments identify $L$ and the loop set. A surviving
edge has an unlooped endpoint. Give that endpoint value $1$, its
other endpoint value $2$, and all other variables value $0$. No loop
or other edge contributes $3$, so this detects precisely the selected
edge. Hence these normal forms too are uniquely determined by the
term function, proving completeness.
\end{proof}

\subsubsection{Independent-set constraints}
\begin{theorem}\label{f15:thm:independent}
A ten-identity basis of $A_{516}$ consists of $\AI$, multiplicative
commutativity, \eqref{f15:eq:top}, and
\begin{align}
 xy&\eq xy+x,\label{f15:eq:addendpoint}\\
 x^2+yv&\eq x^2+yv+xy.\label{f15:eq:loopwitness}
\end{align}
\end{theorem}
\begin{proof}
Soundness is checked in the table. Every polynomial with a long word
reduces to the zero term. Pure linear polynomials have exactly their
semilattice support identities. A short polynomial with a quadratic
word has value $3$ precisely when all variables of $N$ belong to
$\{1,2\}$, no quadratic word has both endpoints assigned $2$, and
all variables in $L\setminus N$ are nonzero. Otherwise its value is
$0$. In particular it is not the zero term: assigning every variable
$1$ gives $3$. That assignment also separates it from pure linear
polynomials, whose value is $1$.

Write $T=L\cup N$ and let $J$ be the loop set. The function determines
$T$ by changing individual variables from $1$ to $0$, and determines
$N$ by changing them from $1$ to $3$. Among assignments of $1,2$ to
$N$, the permitted sets of vertices assigned $2$ have as minimal
forbidden subsets exactly the singletons in $J$ and the edges between
two vertices of $N\setminus J$. Thus the function determines $J$
and the set $E$ of such edges.

Conversely these data determine the function. They also specify a
derivable canonical form. Use \eqref{f15:eq:addendpoint} to adjoin every
variable of $N$ linearly, so the linear support is $T$. For each
$x\in J$ and $y\in N$, select an occurrence $yv$ in the original
polynomial. Identity~\eqref{f15:eq:loopwitness} adjoins $xy$; if witnesses
coincide, additive idempotence duplicates them before application.
The final graph consists of the loops $J$, every edge incident with
$J$ inside $N$, and the edges $E$. No further term-function data are
needed. Equal functions consequently reduce to the same polynomial.
The zero, linear, and nonlinear cases are separated as above, completing
the proof.
\end{proof}

\Needspace{4\baselineskip}

\subsubsection{Equality constraints with marked vertices}
\begin{theorem}\label{f15:thm:equal}
A twelve-identity basis of $A_{835}$ consists of $\AI$, multiplicative
commutativity, \eqref{f15:eq:top}, and
\begin{align}
 xy&\eq xy+x^2,\label{f15:eq:endloop}\\
 xy+yz&\eq xy+yz+xz,\label{f15:eq:transitive}\\
 x+xy&\eq x+xy+y,\label{f15:eq:markprop}\\
 x+x^2+y+y^2&\eq x+x^2+y+y^2+xy.\label{f15:eq:eqbridge}
\end{align}
\end{theorem}
\begin{proof}
The table verifies soundness. Long words give the zero term. For a
short polynomial with a quadratic word, the value is $3$ exactly when
each connected component of its quadratic graph is assigned one
constant value from $\{1,2\}$, every marked component has value $1$,
and each variable of $L\setminus N$ belongs to $\{1,3\}$. Otherwise
the value is $0$. Indeed, nonzero products are precisely $1^2=2^2=3$,
and a linear summand can be added to $3$ without giving $0$ precisely
when it is $1$ or $3$.

We first derive a canonical form. By \eqref{f15:eq:endloop}, adjoin a
loop at every vertex in $N$. By \eqref{f15:eq:transitive}, complete every
connected component to a complete graph, keeping all its loops.
By \eqref{f15:eq:markprop}, a marked vertex makes every vertex of its
component linear. For marked vertices in two different components,
their loops and \eqref{f15:eq:eqbridge} adjoin a connecting edge. Repeating
this and the preceding two operations merges all marked components
into one complete, fully marked component. Unmarked components stay
separate and complete. The variables in $L\setminus N$ remain linear
and outside the graph.

This form is uniquely determined by the function. Its total support
$T=L\cup N$ and nonlinear support $N$ are detected by changing a
variable to $0$ or $3$ in the all-$1$ assignment. Fix the outside
linear variables at $1$. The sets of graph vertices that may be
assigned $2$ while all others are $1$ are exactly the unions of
unmarked components. The inclusion-minimal nonempty such sets recover
the unmarked components, and the vertices in no such set are precisely
the merged marked component. The data determine the canonical form.
Every such polynomial is nonzero under the all-$1$ assignment and
has value $3$ there, whereas a pure linear polynomial has value $1$.
Pure linear identities are support identities. These observations
also handle identities involving the zero term and finish completeness.
\end{proof}

\Needspace{4\baselineskip}

\subsubsection{Two-colourings and parity closure}
\begin{theorem}\label{f15:thm:bipartite}
A thirteen-identity basis of $A_{495}$ consists of $\AI$, multiplicative
commutativity, \eqref{f15:eq:top}, and
\begin{align}
 x^2+u&\eq x^2,\label{f15:eq:squaretop}\\
 x+y+xy&\eq x^3,\label{f15:eq:conflict}\\
 xy+yz+zt&\eq xy+yz+zt+xt,\label{f15:eq:oddpath}\\
 x+xy+yz&\eq x+xy+yz+z,\label{f15:eq:evenmark}\\
 x+xu+y+yv&\eq x+xu+y+yv+xv.\label{f15:eq:bipbridge}
\end{align}
\end{theorem}
\begin{proof}
All identities hold in the table. For a short nonlinear polynomial,
nonzero evaluation requires every edge to join a $1$ to a $2$.
Every marked vertex must be $1$, and each outside linear variable
must be $1$ or $3$. These conditions are also sufficient for value
$3$; all other evaluations have value $0$. Thus a nonzero function
requires a bipartite quadratic graph with all marked vertices of each
component on a single side of its bipartition. We prove both that
inconsistent graphs reduce to the zero term and that the consistent
ones have a derivable, unique form.

Identity~\eqref{f15:eq:oddpath} adjoins the endpoints of any walk of length
three as an edge. Induction on odd walk length shows that the endpoints
of every odd walk may be joined: shorten its first three edges to one
and continue. Repeated edges and vertices are allowed by substitution
and additive idempotence. An odd closed walk therefore produces a
loop. In the presence of a loop, \eqref{f15:eq:squaretop} and
\eqref{f15:eq:top} make the polynomial the zero term. More explicitly,
$x^2+y^3\eq x^2$ by \eqref{f15:eq:squaretop} and
$x^2+y^3\eq y^3$ by \eqref{f15:eq:top}, so $x^2\eq y^3$.
Consequently every nonbipartite graph is derivably zero. In a bipartite
graph, the same odd-walk operation completes each component to the
complete bipartite graph on its two sides.

Identity~\eqref{f15:eq:evenmark} propagates a mark along a walk of length
two. Iteration marks every vertex an even distance from an existing
mark. If a component has marks in opposite sides, the completed graph
has an edge joining two marked vertices. Identity~\eqref{f15:eq:conflict}
then produces a zero subterm, and \eqref{f15:eq:top} makes the whole
polynomial zero. This accounts for all inconsistent mark assignments.

It remains to consider consistent graphs. Each marked component has
one fully marked side, forced to be $1$, and one unmarked side, forced
to be $2$. To merge two marked components, select marked vertices
$x,y$ and neighbours $u,v$, respectively. Such neighbours exist because
every vertex in $N$ occurs in an edge and the graph has no loops.
Identity~\eqref{f15:eq:bipbridge} adjoins $xv$. Odd-walk completion and
even-mark propagation then merge the two components into one complete
bipartite component, with all its forced-$1$ vertices marked. Repeat
until all marked components have merged. The unmarked components
remain separate complete bipartite graphs, with their sides unordered.
The outside linear set is unchanged. This is the canonical form.

We verify uniqueness without a bound on the number of vertices. A
consistent graph has a nonzero evaluation, obtained by a proper
two-colouring with every marked side $1$ and every outside linear
variable $1$. Equal nonzero functions have the same total support:
in a common nonzero assignment, changing a variable to $0$ destroys
the value exactly when that variable occurs. They also have the same
nonlinear support: changing a variable to $3$ destroys the value
exactly when it belongs to $N$. Fixing the outside linear variables
at $1$, their common set of proper colourings identifies the forced-$1$
and forced-$2$ vertices. These give the two sides of the merged marked
component, if any. Among the remaining vertices, two vertices belong
to the same unmarked component exactly when their relative colour
is fixed in all permitted colourings. Within a component this relation
recovers its two sides, up to interchange. Between distinct unmarked
components either one may be flipped independently, so their relative
colour is not fixed. The allowed colourings therefore determine the
whole canonical form.

All consistently coloured nonlinear polynomials have some nonzero
evaluation, separating them from the zero term. Their value at the
all-$1$ assignment is $0$, whereas any pure linear polynomial has value
$1$. Pure linear polynomials have only their support identities. All
possible identities are consequently covered by the listed reductions
and separation arguments, proving completeness.
\end{proof}

\endgroup
\subsection{Weighted quadratic graphs}\label{sec:family16}
\begingroup
\emph{Generators (basis size):} $A_{241}$ (11), $A_{264}$ (10),
$A_{867}$ (13), $A_{870}$ (11).

The four multiplications have the common description
$ab=q$ if $a,b\in H$, and $ab=0$ otherwise, where
\[
\begin{array}{c|cc}
\toprule j&H&q\\\midrule
241,264&\{1\}&2\\
867,870&\{1,2\}&3\\\bottomrule
\end{array}
\]
In every case $0$ is an additive identity and a multiplicative zero.
Since $q\notin H$, all triple products are $0$. Let $\mathcal C$
consist of $\SR$, multiplicative commutativity, and
\begin{equation}\label{f16:eq:bottom}
 xyz+u\eq u.
\end{equation}
As in Section~\ref{sec:family15}, long words reduce to a common
zero term; after their deletion a polynomial has only linear and
quadratic monomials. A zero term represents the case with none left.

\subsubsection{Three bases and their normal forms}
Use the following identities:
\begin{align}
 3x&\eq2x,\label{f16:eq:threshold}\\
 2xy&\eq xy,\label{f16:eq:prodidem}\\
 x+xy&\eq x,\label{f16:eq:absorb}\\
 x^2+xy&\eq x^2,\label{f16:eq:loop}\\
 2x+xy&\eq2x,\label{f16:eq:doubleabsorb}\\
 x+y+xy&\eq x+y,\label{f16:eq:pairabsorb}\\
 x+x^2&\eq2x,\label{f16:eq:loopdouble}\\
 3x&\eq x,\label{f16:eq:period}\\
 2xy&\eq x^3.\label{f16:eq:doublezero}
\end{align}
\begin{theorem}\label{f16:thm:threshold}\label{f16:thm:marked}\label{f16:thm:parity}
Adjoining the indicated identities to $\mathcal C$ gives complete bases:
\[
\begin{array}{ccl}
\toprule \text{Case}&\text{Generators}&\text{Added identities}\\\midrule
{\rm T}&A_{241},A_{870}&
 \eqref{f16:eq:threshold}\text{--}\eqref{f16:eq:loop}\\
{\rm M}&A_{867}&
 \eqref{f16:eq:threshold},\eqref{f16:eq:prodidem},
 \eqref{f16:eq:loop}\text{--}\eqref{f16:eq:loopdouble}\\
{\rm P}&A_{264}&
 \eqref{f16:eq:absorb},\eqref{f16:eq:period},\eqref{f16:eq:doublezero}\\\bottomrule
\end{array}
\]
Their sizes are eleven, thirteen and ten, respectively, and
$\Var(A_{241})=\Var(A_{870})$. Put
$L=\{x:c(x)>0\}$ and $S=\{x:c(x)=1\}$ for $c:X\to\{0,1,2\}$.
The normal forms are
\begin{equation}\label{f16:eq:normal}
 \sum_{x\in L}c(x)x+\sum_{x\in J}x^2+
 \sum_{\{x,y\}\in E}xy,
\end{equation}
where $J\subseteq X\setminus L$ and $E$ is a simple undirected
edge set with the following allowed endpoints:
\[
\begin{array}{cl}
\toprule \text{Case}&\text{Allowed edges}\\\midrule
{\rm T}&\text{Both endpoints outside }L\cup J\\
{\rm M}&\text{Both outside }L\cup J,\text{ or one in }S
                      \text{ and one outside }L\cup J\\
{\rm P}&\text{Both endpoints outside }L\\\bottomrule
\end{array}
\]
If all three parts are absent, use the zero term.
\end{theorem}
\begin{proof}
The multiplication description proves $\mathcal C$. In T and M,
$q+q=q$, so \eqref{f16:eq:prodidem} holds; a nonzero $xy$
requires $x\in H$, whence $x^2=q$, proving \eqref{f16:eq:loop}.
In T and P, $r+q=r$ for $r\in H$, proving
\eqref{f16:eq:absorb}. In M, $2r+q=2r$ and
$r+s+q=r+s$ for $r,s\in H$, while $r+r^2=2r$ for every
element. These prove \eqref{f16:eq:doubleabsorb}--\eqref{f16:eq:loopdouble}.
Finally, addition gives $3r=2r$ in T and M, and $3r=r$,
$q+q=0$ in P. Hence all listed identities are valid.

Remove long words and reduce positive linear coefficients to one or
two, using parity in P. Quadratic coefficients become one in T and
M, and reduce modulo two in P. In T and P,
\eqref{f16:eq:absorb} removes every quadratic monomial incident
with $L$, selecting one linear copy when its coefficient is two.
In M, \eqref{f16:eq:loopdouble} converts a linear variable with
a loop into a doubled linear variable. Then
\eqref{f16:eq:doubleabsorb} removes all edges and loops at doubled
variables, and \eqref{f16:eq:pairabsorb} removes edges joining two
linear variables. In T and M, \eqref{f16:eq:loop} removes all
nonloop edges incident with the remaining loop set $J$.
These steps give exactly \eqref{f16:eq:normal}; no loop absorption
is used in P.

We recover its parameters by assignments with all unmentioned
variables equal to $0$. In T and M, setting $x=1$ gives $1$
for $x\in L$, $q$ for $x\in J$, and $0$ otherwise.
For $x\in L$, the following assignments distinguish its coefficient:
\[
\begin{array}{c|ccc}
\toprule j&x&c(x)=1&c(x)=2\\\midrule
241&3&3&2\\
867,870&2&2&1\\
264&1&1&3\\\bottomrule
\end{array}
\]
In P first set $x=3$: its output is $3$ exactly when $x\in L$.
For $x\notin L$, setting $x=1$ gives $2$ exactly when $x\in J$.
Thus $c$ and $J$ are determined in every case.

For T and M, an edge outside $L\cup J$ is detected by assigning
both endpoints $1$: its presence gives $q$, and its absence gives $0$.
For a marked edge in M, assign $x=2$ at its endpoint in $S$ and
$y=1$ at its endpoint outside $L\cup J$. The output is $1$
if the edge occurs and $2$ otherwise, since $2+3=1$.
For P, assign distinct $x,y\notin L$ the value $1$. The output
is $2$ precisely when the number of present loops at $x,y$
and the edge $xy$ is odd. Known loop membership therefore
determines the edge. These tests separate every normal form,
including the zero term, on arbitrary finite $X$.
Lemma~\ref{lem:normal-form-completeness} completes the proof.
\end{proof}

The generators $A_{241}$ and $A_{870}$ are not isomorphic:
they have respectively one and four ordered pairs with nonzero product.
In $A_{867}$, $x=2,y=1$ separates $x$ from $x+xy$, so the marked
edges retained in M cannot be discarded.
\endgroup

\subsection{Quadratic graphs and cubic constraints}\label{sec:family30}
\begingroup
\providecommand{\eq}{}\renewcommand{\eq}{\approx}
\providecommand{\SR}{}\renewcommand{\SR}{\mathsf{SR}}
\providecommand{\Var}{}\renewcommand{\Var}{\operatorname{Var}}

Throughout, $\Sigma$ consists of $\SR$ and
\begin{equation}\label{f30:eq:common}
x+x\eq x,\qquad xy\eq yx,\qquad xyzu\eq v^4.
\end{equation}
Thus $\Sigma$ has eight identities. The expression $\omega=x^4$
denotes a constant term function, not a signature constant.
Indeed, \eqref{f30:eq:common} makes all fourth powers equal.
It also gives $\omega t=t^5\eq\omega$, using the four factors
$t,t,t,t^2$. All words of length at least four therefore equal
$\omega$. Distributivity reduces every other term to a sum of
linear, quadratic and cubic words, with repeated summands deleted.

The following indices refer to Appendix~\ref{app:catalogue}.
All three have multiplication
\[
i\cdot j=\begin{cases}i+j,&i,j\ne0\text{ and }i+j\leq3,\\
0,&\text{otherwise}.
\end{cases}
\]
Here the addition on the right is ordinary integer addition.
Thus $0$ is a multiplicative zero and the powers of $1$ are
$1,2,3,0$. Element labels are not constants of the language.
\emph{Generators (basis size):} $A_{502}$ (17), $A_{503}$ (15), $A_{504}$ (13).

Give the labels $1,2,3,0$ respective weights $1,2,3,4$.
Addition in $A_{504}$ takes the maximum weight; addition in
$A_{503}$ takes the minimum weight. In $A_{502}$, addition is
idempotent, $0$ absorbs everything, $1+3=3$, and $1+2=2+3=0$.
Basis sizes include $\SR$ and are not asserted minimal.
An absent sum is omitted, not interpreted as a constant.

\subsubsection{Maximum weights and quadratic graphs}
\begin{theorem}\label{f30:thm:max}
A complete thirteen-identity basis for $A_{504}$ is $\Sigma$ and
\begin{gather}
x^4+y\eq x^4,\qquad x+xy\eq xy,\label{f30:eq:maxabs}\\
xyz\eq x^3+y^3+z^3,\label{f30:eq:maxcubes}\\
x^2+yz\eq x^2+yz+xy+xz,\qquad
x^3+y^2\eq x^3+y^3.\label{f30:eq:maxgraph}
\end{gather}
\end{theorem}
\begin{proof}
The identities hold by the weight description, or by checking the
finite tables. A term containing $\omega$ reduces to $\omega$.
Otherwise replace every cubic word by the sum of the cubes of its
variables. If any cube occurs, the last identity in
\eqref{f30:eq:maxgraph} replaces all squares by cubes. In the absence
of cubes, retain the square variables. Let $d=3$ in the first
case, and $d=2$ in the second; write $G$ for the resulting cube
or square variable set. Let $Q$ be the set of all other quadratic
variables, $E$ the set of quadratic edges with both endpoints in
$Q$, and $L$ the remaining linear variables. Absorption deletes
every linear variable occurring nonlinearly. A pure linear sum is
a separate form with $d=1$.

If $G\ne\varnothing$, choose its least variable $g$. The first
identity in \eqref{f30:eq:maxgraph} inserts $gx$ for every quadratic
variable $x$. For a cubic guard, first insert its square using
$g^2+g^3\eq g^3$, which follows from absorption. Edges between
two guard variables are redundant: substituting $z=y$ in that
identity gives $x^2+y^2+xy\eq x^2+y^2$; cubes absorb the
corresponding squares. Edges from any guard to $Q$ are equivalent
to the inserted edges from $g$: with both guard squares present,
either choice reconstructs the other by the same insertion law.
Consequently the normal form is
\begin{equation}\label{f30:eq:maxnf}
\sum_{x\in L}x+\sum_{g\in G}g^d+
\sum_{x\in Q}gx+\sum_{\{x,y\}\in E}xy,
\end{equation}
where the third sum is used only when $G\ne\varnothing$ and
its $g$ is the chosen guard. The sets $L,G,Q$ are pairwise
disjoint. When $G=\varnothing$, $d=2$ and $E$ has vertex set
$Q$ with no isolated vertices. When $d=3$, $G$ is nonempty.
For $d=2$ at least one quadratic summand is present.

We recover the form on any finite ambient variable set. The value
at the all-$1$ assignment is $d$; it is $0$ for $\omega$.
For $d=2$, changing just $x$ to $2$ gives $0$ precisely for
$x\in G$, and gives $3$ precisely for $x\in Q$. For $d=3$,
changing just $x$ to $2$ gives $0$ precisely for $x\in G$;
changing just $x$ to $3$ gives $0$ precisely for
$x\in G\cup Q$. In either case, changing just $x$ to $0$
recovers the full content, hence $L$. For distinct $x,y\in Q$,
changing precisely these variables to $2$ gives $0$ exactly when
$\{x,y\}\in E$. All other quadratic summands have weight at
most three in this test. Pure linear sums are recovered by their
content. Thus every datum in \eqref{f30:eq:maxnf} is determined by
the term function. Equational reduction and separation prove
completeness for arbitrarily many variables.
\end{proof}

\subsubsection{Minimum weights and cubic antichains}
For $\varnothing\ne U$ with $|U|\leq3$, let $c_U$ be a
cubic word with support $U$: use $x^3$ for a singleton, $x^2y$
for a pair with $x<y$, and the squarefree product for a triple.
\begin{theorem}\label{f30:thm:min}
A complete fifteen-identity basis for $A_{503}$ is $\Sigma$ and
\begin{gather}
x^4+y\eq y,\qquad x+xy\eq x,\qquad x^2y\eq xy^2,
\label{f30:eq:minbasic}\\
x^2+xyz\eq x^2,\qquad x^3+xyz\eq x^3,\qquad
x^2y+xyz\eq x^2y,\label{f30:eq:mincubes}\\
x^2+y^3+xy\eq x^2+y^3.\label{f30:eq:minedge}
\end{gather}
\end{theorem}
\begin{proof}
The minimum-weight description verifies the identities.
Delete $\omega$ summands, retaining $\omega$ alone if no other
summand remains. Absorption deletes every nonlinear word
containing a retained linear variable. Write $L$ for these
linear variables, $S$ for the square variables, and $E$ for the
simple quadratic edges. The cubic support identity in
\eqref{f30:eq:minbasic} identifies the two multiplicities on a
two-variable cubic support.

The identities in \eqref{f30:eq:mincubes} remove a cubic support
meeting $S$ and remove any cubic support properly containing
another cubic support. Absorption removes cubic supports
containing a quadratic edge. Let $\mathcal H$ be the remaining
antichain of nonempty supports of size at most three. It is
disjoint from $L\cup S$, and none of its members contains an
edge of $E$. Put
$T=\{x:\{x\}\in\mathcal H\}$.
Equation~\eqref{f30:eq:minedge} deletes edges between $S$ and
$S\cup T$. When the second endpoint lies in $S$, insert its
cube using absorption before applying that equation. Such a
deleted edge could not have removed a surviving cubic support,
because it has an endpoint in $S$. These reductions therefore
give the form
\begin{equation}\label{f30:eq:minnf}
\sum_{x\in L}x+\sum_{x\in S}x^2+
\sum_{\{x,y\}\in E}xy+\sum_{U\in\mathcal H}c_U.
\end{equation}
Here $L\cap S=\varnothing$, edges avoid $L$, and no edge joins
$S$ to $S\cup T$. At least one summand is present.

Fix any finite variable set. Assign $x=1$ and all other variables
$0$. The output is $1$ exactly for $x\in L$, $2$ exactly for
$x\in S$, and $3$ exactly for $x\in T$. Thus these sets are
recovered. For two variables outside $L\cup S$, set them both
to $1$ and the rest to $0$. Output $2$ detects their edge.
For $x\in S$ and $y\notin L\cup S\cup T$, instead set
$x=2,y=1$ and the rest to $0$. Output $3$ detects the edge
$xy$: the square at $x$ is now $0$, and there is no singleton
cubic at $y$. All other possible edges incident with $S$ are
excluded in the normal form. This recovers $E$.

Finally, for $D\subseteq X$, assign $D$ the value $1$ and its
complement $0$. Call $D$ admissible if it misses $L\cup S$
and contains no edge of $E$. On admissible sets, the output is
$3$ exactly when $D$ contains a member of $\mathcal H$.
Every member of $\mathcal H$, and each of its subsets, is
admissible. Hence the minimal admissible sets producing $3$
recover $\mathcal H$. The term $\omega$ is distinguished
because every nonzero form takes a nonzero value at all $1$.
All normal-form data are recovered, proving completeness.
\end{proof}

\subsubsection{Bipartite components and prescribed values}
\begin{theorem}\label{f30:thm:fork}
A complete seventeen-identity basis for $A_{502}$ is $\Sigma$ and
\begin{gather}
x^4+y\eq x^4,\qquad x+xy\eq x^3+xy,\qquad
xyz\eq x^3+y^3+z^3,\label{f30:eq:forkbasic}\\
x^2+yz\eq x^2+y^2+z^2,\label{f30:eq:forksquare}\\
xy+yz+zu\eq xy+yz+zu+xu,\label{f30:eq:forkpath}\\
x^3+xy+yz\eq x^3+xy+yz+z^3,\label{f30:eq:forkanchor}\\
x^3+z^3+xy\eq x^3+z^3+xy+zy,\label{f30:eq:forkchange}\\
x^2+y^3\eq z^4,\qquad x^3+y^3+xy\eq z^4.
\label{f30:eq:forkzero}
\end{gather}
\end{theorem}
\begin{proof}
The identities hold in the defining tables. Replace cubic words
by cubes of all their variables. Regard a square as a loop and
a squarefree quadratic word as an undirected edge. Replace every
linear variable incident with a quadratic term by its cube using
\eqref{f30:eq:forkbasic}. A linear variable already cubed disappears:
substituting $y=x^2$ into its middle identity yields
$x+x^3\eq x^3$. Thus remaining linear variables avoid all
quadratic and cubic variables.

If a square occurs, \eqref{f30:eq:forksquare} turns all quadratic
terms into squares on their full variable set $Q$. If any cube
also occurs, \eqref{f30:eq:forkzero} gives $\omega$. Otherwise the
form is
\begin{equation}\label{f30:eq:squarenf}
\sum_{x\in L}x+\sum_{x\in Q}x^2,
\qquad Q\ne\varnothing,\quad L\cap Q=\varnothing.
\end{equation}

Suppose next that there are no squares. Equation~\eqref{f30:eq:forkpath}
inserts the edge joining the endpoints of any three-edge walk.
Iteration inserts the endpoints of every odd-length walk, by
induction reducing its length by two. In a connected nonbipartite
component, an odd closed walk consequently inserts a loop, and
the preceding square reduction applies to all quadratic terms.
In a bipartite component with parts $A,B$, odd walks connect
exactly opposite parts. The same insertions therefore complete
the component to the graph with edge set $A\times B$.
Every part is nonempty; parts and components are ordered
canonically to remove the immaterial choice of their names.

If no cube occurs, the resulting form is
\begin{equation}\label{f30:eq:graphnf}
\sum_{x\in L}x+\sum_{(A,B)\in\mathcal P}
\sum_{x\in A,\,y\in B}xy,
\end{equation}
where $\mathcal P$ is a nonempty family of components with
mutually disjoint vertex sets, disjoint also from $L$.

If cubes occur, write $C$ for their variable set. Within a
bipartite component, \eqref{f30:eq:forkanchor} propagates a cube
along every even-length path. Thus a cube at one vertex inserts
cubes on its whole part. Cubes in both parts give $\omega$
by the second identity in \eqref{f30:eq:forkzero}. Otherwise put
the cubed part into $C$ and its opposite part into $D$.
Retain in $\mathcal P$ only the components containing no cube.
Equation~\eqref{f30:eq:forkchange} permits changing the cubed
endpoint of any edge to any other cubed variable. Hence the
edges of all components meeting $C$ can be replaced by $gd$
for $d\in D$, where $g$ is the least member of $C$.
The replaced edges are reconstructed by the same identity,
which justifies their deletion. The form is
\begin{equation}\label{f30:eq:cubicnf}
\sum_{x\in L}x+\sum_{x\in C}x^3+
\sum_{y\in D}gy+
\sum_{(A,B)\in\mathcal P}\sum_{x\in A,\,y\in B}xy.
\end{equation}
Here $C\ne\varnothing$; $L,C,D$ and the component vertex
sets are pairwise disjoint. Both $D$ and $\mathcal P$ may
be empty. The remaining possibilities are $\omega$ and a
nonempty pure linear sum. This proves existence of all forms
by equational transformations using the proposed basis.

We separate them over an arbitrary finite ambient set $X$.
A pure linear sum takes the value $1$ at all $1$; none of the
other forms ever takes value $1$. Such sums are distinguished
by setting a chosen variable to $0$ and all others to $1$.
The term $\omega$ is identically $0$, whereas each other form
admits a nonzero assignment.

Forms \eqref{f30:eq:squarenf} and \eqref{f30:eq:graphnf} admit output
$2$. Its inverse image is exactly the set of assignments
\[
Q\mapsto1,\qquad L\mapsto2,
\]
with all absent variables arbitrary, where $Q$ is the full
quadratic vertex set. Coordinate projections of this inverse
image recover $Q$ and $L$. The square form never takes value
$3$. The graph form takes value $3$ exactly when every
component is properly colored by $1,2$, with its two parts
receiving opposite values, and every variable of $L$ lies in
$\{1,3\}$. The choice of orientation is independent in
different components. Two vertices of $Q$ receive different
values in every such assignment exactly when they are in
opposite parts of one component. This recovers its complete
bipartite edge relation, hence $\mathcal P$.

Form \eqref{f30:eq:cubicnf} takes only values $0,3$. Its inverse
image of $3$ consists exactly of assignments satisfying
\[
C\mapsto1,\quad D\mapsto2,\quad L\mapsto\{1,3\},
\]
with each component of $\mathcal P$ properly colored by
$1,2$ and absent variables arbitrary. This set is nonempty.
Its coordinate projections are respectively $\{1\}$,
$\{2\}$, $\{1,3\}$, $\{1,2\}$, and the whole algebra.
They recover $C,D,L$, the unanchored vertex set, and the absent
variables. On the unanchored vertices, the pairs that always
receive different values recover the complete bipartite
components as before. All forms are therefore uniquely
determined by their functions. This proves the asserted
arbitrary-variable completeness.
\end{proof}

\endgroup
\subsection{Directed stars and rectangles}\label{sec:family31}
\begingroup
\providecommand{\eq}{}\renewcommand{\eq}{\approx}
\providecommand{\SR}{}\renewcommand{\SR}{\mathsf{SR}}
\providecommand{\Var}{}\renewcommand{\Var}{\operatorname{Var}}

\emph{Generators (basis size):} $A_{277}$ (13), $A_{281}$ (14), $A_{289}$ (11), $A_{293}$ (12), $A_{297}$ (10), $A_{400}$ (11), $A_{413}$ (13), $A_{417}$ (14), $A_{425}$ (11), $A_{429}$ (12), $A_{433}$ (10), $A_{990}$ (10), $A_{1343}$ (11), $A_{1627}$ (10).

\subsubsection{Directed stars with two treatments of linear variables}
Let $\mathcal S$ consist of $\SR$, additive idempotence, and
\begin{gather}
xyz\eq x^2+yz,\qquad xy\eq x^2+xy,\label{f31:eq:starword}\\
x+yx\eq x+y^2,\qquad
x^2+yz\eq x^2+y^2+xz.\label{f31:eq:startransport}
\end{gather}
This is a set of ten identities.
\begin{theorem}\label{f31:thm:stars}
The set $\mathcal S$ is a complete basis for $A_{297}$.
The set $\mathcal S$ together with
\begin{equation}\label{f31:eq:linearremove}
x+x^2\eq x^2
\end{equation}
is a complete eleven-identity basis for $A_{400}$.
\end{theorem}
\begin{proof}
Both assertions of validity follow from the tables. For a word
$x_1\cdots x_n$ of length $n\geq2$, repeated use of
\eqref{f31:eq:starword} gives
\[
x_1\cdots x_n\eq x_1^2+\cdots+x_{n-2}^2+x_{n-1}x_n.
\]
For $n=2$ only its last term is used. Insert the square of
each quadratic source using the second identity in
\eqref{f31:eq:starword}. Let $P$ be this nonempty source set
and $L$ the linear variable set. Under the extra identity
\eqref{f31:eq:linearremove}, remove $L\cap P$; without it,
retain those linear variables.

If a terminal variable belongs to $L$, the first identity in
\eqref{f31:eq:startransport} removes its incident incoming edge
while retaining the source square. If a terminal belongs to
$P$, the second identity, applied with its square as guard,
replaces that edge by the two endpoint squares. Let $T$ be
the remaining terminal set; it avoids $L\cup P$.

Choose the least $g\in P$. With the squares of both sources
present, the second identity in \eqref{f31:eq:startransport}
permits transporting any edge to source $g$, and also
reconstructing its original source. Thus all terms containing
a nonlinear word reduce to
\begin{equation}\label{f31:eq:starnf}
\sum_{x\in L}x+\sum_{p\in P}p^2+\sum_{t\in T}gt,
\qquad P\ne\varnothing,\quad T\cap(L\cup P)=\varnothing.
\end{equation}
For $A_{400}$ also require $L\cap P=\varnothing$.
A term with no nonlinear word is a pure nonempty linear sum.

In both generators a nonlinear product can be nonzero only
if every nonterminal factor is $1$. In $A_{297}$ its
terminal must be $1$ or $2$, and the value is that terminal.
In $A_{400}$ the terminal can also be $3$, which gives
value $2$. In $A_{297}$ addition is the maximum in the chain
$1<2<0<3$; in $A_{400}$ it is the maximum in $1<2<3<0$.

Fix a finite ambient variable set $X$. For $A_{297}$,
assigning $x=3$ and all other variables $1$ gives $3$
exactly when $x\in L$, since a linear $3$ absorbs even
zero-valued products. Assigning just $x=2$ gives $0$
exactly when $x\in P$. Having recovered $L,P$, for
$x\notin L\cup P$ the first test gives $0$ exactly
when $x\in T$. For $A_{400}$, assigning just $x=2$
gives $0$ exactly for $x\in P$, and assigning just
$x=3$ gives $3$ exactly for $x\in L$. Outside these
sets, assigning $x=2$ gives $2$ exactly for $x\in T$.
These tests also distinguish the case $P=\varnothing$,
which is precisely a pure linear form and has $T=\varnothing$.
All data in \eqref{f31:eq:starnf} are therefore recovered for
arbitrary finite $X$. This proves both completeness claims.
\end{proof}

\subsubsection{Directed rectangles and additive multiplicities}
Let $\mathcal R$ be the following five identities:
\begin{gather}
xyz\eq xz+yz,\label{f31:eq:rectword}\\
xy+zu\eq xy+zu+xu+zy,\label{f31:eq:rectcomplete}\\
x^2+yz\eq x^2+y^2+z^2,\label{f31:eq:rectsquare}\\
x+xy\eq x^2+y^2,\qquad x+yx\eq yx.
\label{f31:eq:rectlinear}
\end{gather}
Define
\begin{align*}
\mathcal R_I&=\SR\cup\mathcal R\cup\{x+x\eq x\},\\
\mathcal R_T&=\SR\cup\mathcal R\cup
\{3x\eq2x,\ xy+xy\eq xy,\ 2x+yz\eq yx+yz\}.
\end{align*}
They have eleven and thirteen identities, respectively.
Finally, write
\begin{equation}\label{f31:eq:rectcollapse}
x+y^2\eq x^2+y^2.
\end{equation}
\begin{theorem}\label{f31:thm:rectangles}
The sets $\mathcal R_I$ and $\mathcal R_T$ are complete bases
for $A_{289}$ and $A_{277}$, respectively. Adding
\eqref{f31:eq:rectcollapse} to these sets gives complete bases
for $A_{293}$ and $A_{281}$, respectively. The four basis
sizes are eleven, thirteen, twelve and fourteen in this order.
\end{theorem}
\begin{proof}
The tables verify the proposed identities. All four generators
have multiplication $1\cdot1=1$, $1\cdot2=2$, and all other
products zero. Consequently a nonlinear word has value $1$
when all its variables are $1$, value $2$ when its nonterminal
factors are all $1$ and its terminal is $2$, and value $0$
otherwise.

For clarity, the four additions admit one description. Zero
absorbs everything; $1,2$ are idempotent and $1+2=0$.
Let $e=1$ for $A_{289},A_{277}$ and $e=2$ for $A_{293},A_{281}$.
Then $e+3=3$ and $(3-e)+3=0$. Finally $3+3=3$ in the
idempotent cases and $3+3=0$ in the other two cases.

First reduce each long word to the sum of the directed edges
from its nonterminal letters to its final letter, using
\eqref{f31:eq:rectword} inductively. Nonlinear repetitions are
deleted. Linear coefficients are one in the idempotent cases
and one or two in the threshold cases. When a nonlinear
edge $yz$ is present, the identity $2x+yz\eq yx+yz$
replaces every doubled linear variable by an incoming edge.
Thus all remaining linear variables have coefficient one
whenever nonlinear terms occur.

Regard the quadratic terms as directed edges, with a square
as a loop. Equation~\eqref{f31:eq:rectcomplete} completes their
edge set to $P\times T$, where $P,T$ are the sets of sources
and terminals. A common vertex of $P,T$ therefore produces
a loop. A linear variable in $P$ also produces a square by
the first identity in \eqref{f31:eq:rectlinear}. Once a square
is present, \eqref{f31:eq:rectsquare} replaces all quadratic
terms by squares on their full vertex set $Q$.
Linear variables in $Q$ disappear, since the second identity
in \eqref{f31:eq:rectlinear} gives $x+x^2\eq x^2$.
With the extra identity \eqref{f31:eq:rectcollapse}, the remaining
linear variables are also converted into squares.

The resulting square forms are
\begin{equation}\label{f31:eq:rectsquarenf}
\sum_{x\in L}x+\sum_{q\in Q}q^2,
\qquad Q\ne\varnothing,\quad L\cap Q=\varnothing,
\end{equation}
with the further restriction $L=\varnothing$ in the two
cases satisfying \eqref{f31:eq:rectcollapse}. If no square is
created, $P,T$ are nonempty and disjoint, and linear variables
in $T$ disappear by \eqref{f31:eq:rectlinear}. The form is
\begin{equation}\label{f31:eq:rectnf}
\sum_{x\in L}x+\sum_{p\in P,\,t\in T}pt,
\qquad L,P,T\text{ pairwise disjoint},\quad P,T\ne\varnothing.
\end{equation}
The remaining possibility is a pure nonempty linear sum
$\sum c_xx$, with $c_x\in\{0,1\}$ or $c_x\in\{0,1,2\}$
as appropriate. These reductions use only the claimed basis.

We now separate all forms over any finite ambient variable set.
Square forms never take value $2$. Pure linear and rectangle
forms do take that value. The inverse image of $2$ for a
pure linear form consists exactly of assignments giving
every present variable value $2$. For a rectangle form it
consists exactly of assignments
\begin{equation}\label{f31:eq:recttwo}
P\mapsto1,\qquad T\cup L\mapsto2,
\end{equation}
with absent variables arbitrary. The nonempty set of
coordinates forced to $1$ distinguishes rectangle forms
from pure linear forms. Thus the three types are distinguished.

Within pure linear forms, assigning one variable $0$ and
all others $1$ recovers content. In a threshold case, assign
a chosen $x$ the value $3$ and all others $e$. The outputs
for coefficient one, coefficient two and absence are $3,0,e$,
respectively. Thus all pure linear coefficients are recovered.

When $e=1$, assigning $x=3$ and all others $1$ gives $3$
exactly for a retained linear variable, in both
\eqref{f31:eq:rectsquarenf} and \eqref{f31:eq:rectnf}.
For a square form, the analogous test with $x=0$ recovers
the union $L\cup Q$, hence $Q$. For a rectangle form,
\eqref{f31:eq:recttwo} recovers $P$ and $T\cup L$, hence $T$.

When $e=2$, a square form has no linear part and has value
$1$ exactly when every variable in $Q$ is $1$; this recovers
$Q$. For a rectangle form, \eqref{f31:eq:recttwo} again recovers
$P$ and $U=T\cup L$. Assign $P$ the value $1$, every
variable in $U$ the value $2$, except a chosen $u\in U$
which receives $3$. The output is $3$ if $u\in L$ and
$0$ if $u\in T$. This works also when $3+3=0$, since
there is exactly one linear contribution equal to $3$.
It recovers $L,T$. All normal-form data are determined,
and arbitrary-variable completeness follows.
\end{proof}

\endgroup
\subsection{Directed support constraints}\label{sec:family32}
\begingroup
\providecommand{\eq}{}\renewcommand{\eq}{\approx}
\providecommand{\SR}{}\renewcommand{\SR}{\mathsf{SR}}
\providecommand{\Var}{}\renewcommand{\Var}{\operatorname{Var}}
\subsubsection{Scope, generators and common word identities}

\emph{Generators (basis size):} $A_{387}$ (13), $A_{390}$ (16), $A_{395}$ (11), $A_{398}$ (14), $A_{1330}$ (13), $A_{1333}$ (16), $A_{1338}$ (11), $A_{1341}$ (14), $A_{1348}$ (17), $A_{1356}$ (15).

The first four generators have multiplication
\[
 1\cdot1=1,\qquad 1\cdot2=1\cdot3=2,
 \qquad ab=0\text{ for all other pairs }(a,b).
\]
A nonlinear word has value \(1\) when all its letters have value
\(1\), value \(2\) when all nonfinal letters have value \(1\)
and the final letter has value \(2\) or \(3\), and value \(0\)
otherwise. The four generators all satisfy
\begin{equation}\label{f32:eq:word}
 xyz\eq yxz,\qquad x^2y\eq xy,\qquad xy^2\eq yx^2.
\end{equation}
Write \(\mathcal W=\SR+\eqref{f32:eq:word}\).

Fix a finite variable set \(X\) and an ordering of it. For nonempty
\(P\subseteq X\), let \(w_P\) be the ordered product on \(P\).
For \(t\notin P\), put \(e_{P,t}=w_Pt\). For nonempty
\(Q\subseteq X\), put \(s_Q=w_Qq\), where \(q=\min Q\).
The notation \(s_Q\) represents a nonlinear square-support word,
including \(s_{\{q\}}=q^2\).
In displayed polynomials, zero coefficients and empty parts of a sum
mean omitted summands. Every represented term is nonempty; no empty
sum or product is introduced into the language.

\begin{lemma}\label{f32:lem:words}
Modulo \(\mathcal W\), every nonlinear word has exactly one of the
following prescribed representatives: \(e_{P,t}\), when its final
letter \(t\) is absent from its prefix support \(P\); or \(s_Q\),
when its final letter occurs earlier and \(Q\) is its full support.
Furthermore \(s_Q\eq w_Q^2\).
\end{lemma}
\begin{proof}
The first law of \eqref{f32:eq:word}, applied in a context, exchanges
any adjacent nonfinal letters. The second removes repeated prefix
letters. If the final letter \(t\) already belongs to the prefix,
the word can be written with suffix \(it^2\) for any other prefix
letter \(i\). The third law replaces this suffix by \(ti^2\).
Thus the repeated final letter can be chosen arbitrarily from the
support. Choose its minimum and sort the remaining prefix. Applying
the same argument to \(w_Q^2\) proves the last assertion.
The case of a one-letter support follows from \(x^3\eq x^2\),
an instance of the second law. No empty word is used.
\end{proof}

\subsubsection{Two bases from directed antichains}
For \(A_{395}\), addition is maximum in the chain
\(0<3<2<1\). The addition of \(A_{387}\) is the same except
that \(3+3=2\). In both, sums of nonlinear word values are
computed by maximum. Define
\begin{equation}\label{f32:eq:anti}
 x+xy\eq x,\qquad xy+zxy\eq xy.
\end{equation}

\begin{theorem}\label{f32:thm:anti}
A complete basis for \(A_{395}\) is
\(\mathcal W+\eqref{f32:eq:anti}+\{2x\eq x\}\), of size eleven.
A complete basis for \(A_{387}\) is
\(\mathcal W+\eqref{f32:eq:anti}\) together with
\begin{equation}\label{f32:eq:threshold}
 3x\eq2x,\qquad 2xy\eq xy,\qquad 2x+yx\eq2x,
\end{equation}
of size thirteen.
\end{theorem}
\begin{proof}
The tables verify the proposed identities. Distributivity and
Lemma~\ref{f32:lem:words} reduce every term to linear summands and the
monomials \(s_Q,e_{P,t}\). Linear coefficients are truncated at
one for \(A_{395}\), and at two for \(A_{387}\). Nonlinear
multiplicities are one. Write these linear coefficients as \(c_x\).

For a nonlinear monomial \(m\), denote its support by \(C(m)\),
so \(C(s_Q)=Q\) and \(C(e_{P,t})=P\cup\{t\}\).
The pointwise order on nonlinear monomial functions is
\begin{align}
 s_Q\preceq n&\quad\Longleftrightarrow\quad C(n)\subseteq Q,
                                                        \label{f32:eq:domsquare}\\
 e_{P,t}\preceq n&\quad\Longleftrightarrow\quad
 C(n)\subseteq P\ \text{or}\
 [n=e_{R,t}\text{ with }R\subseteq P].                  \label{f32:eq:domedge}
\end{align}
Sufficiency follows from the stated word evaluations. For necessity
in \eqref{f32:eq:domsquare}, assign \(1\) on \(Q\) and \(0\)
elsewhere. For necessity in \eqref{f32:eq:domedge}, assign \(1\)
on \(P\), \(3\) to \(t\), and \(0\) elsewhere. The possible
nonlinear summands with value at least \(2\) give precisely the
alternatives displayed. In particular these relations distinguish
all nonlinear monomial functions and define a partial order.

Every absorption prescribed by this order is derivable. If
\(C(n)\subseteq P\), extend the word \(n\) on the right to
have prefix support \(P\) and final letter \(t\), then use
\(x+xy\eq x\) and Lemma~\ref{f32:lem:words}. The same construction,
ending in a repeated letter, deals with
\(C(n)\subseteq Q\) for a square-support word. The equal-support
case \(n=s_Q\) uses nonlinear idempotence. If
\(n=e_{R,t}\) and \(R\subsetneq P\), prepend the product on
\(P\setminus R\) and use \(xy+zxy\eq xy\); equality uses
idempotence. These substitutions have nonempty factors whenever
the corresponding extension is required.

A linear summand \(x\) absorbs \(s_Q\) when \(x\in Q\),
and absorbs \(e_{P,t}\) when \(x\in P\), by the first law of
\eqref{f32:eq:anti} and word rearrangement. In \(A_{387}\), a doubled
linear summand \(2t\) also absorbs \(e_{P,t}\), by the last law
of \eqref{f32:eq:threshold}. Delete all monomials absorbed by these
linear summands. Among the remaining monomials, retain only the
maximal elements for \(\preceq\). We obtain
\begin{equation}\label{f32:eq:antinormal}
 \sum_{x\in X}c_xx+\sum_{m\in H}m,
\end{equation}
where \(H\) is an antichain and no member is absorbed by a displayed
linear summand. At least one summand is present. The deletions are
derivable; a deleted monomial cannot become necessary later because
absorption is transitive.

We prove uniqueness for arbitrary finite \(X\). Assign \(x=3\)
and all other variables \(0\). Every nonlinear word vanishes,
and the value is \(0,3,2\) according as \(c_x=0,1,2\);
the last case occurs only for \(A_{387}\). Thus the function
determines all linear coefficients.

Suppose two canonical polynomials with these same coefficients
define the same function. Take \(s_Q\) in the first antichain and
use the witness for \eqref{f32:eq:domsquare}. The common value is \(1\).
Some summand on the second side must be \(1\). A linear variable
in \(Q\) would have deleted \(s_Q\) on the first side, so a
nonlinear summand \(n\) has \(C(n)\subseteq Q\).

Next take \(e_{P,t}\) in the first antichain and use the witness
for \eqref{f32:eq:domedge}. The common value is at least \(2\) in
the chain order. Linear variables in \(P\), or a doubled \(t\),
would have deleted this edge. All remaining linear contributions
are zero, except possibly a single copy of \(t=3\). Hence some
nonlinear summand \(n\) on the second side has value at least
\(2\), and \(e_{P,t}\preceq n\). This argument also covers
the threshold addition: there are not two surviving linear
contributions of value \(3\) under this witness.

Each monomial on either side is therefore below one on the other.
Two antichains in a partial order with this property coincide:
if \(m\preceq n\preceq m'\), where \(m,m'\) belong to one
antichain, then \(m=m'=n\). The canonical forms are identical.
Reduction and uniqueness prove the two completeness assertions.
\end{proof}

\subsubsection{Two bases from combined prefix constraints}
In \(A_{398}\) and \(A_{390}\), \(0\) absorbs addition,
\(1+1=1\), and \(1+2=1+3=0\). On \(\{2,3\}\), addition
is respectively maximum with \(2<3\), or the group operation
with identity \(2\) and \(3+3=2\). Define the five laws
\begin{align}
 xy+zu&\eq xzy+xzu,                                  \label{f32:eq:union}\\
 x^2+yz&\eq(xyz)^2,                                  \label{f32:eq:squareunion}\\
 x+y^2&\eq(xy)^2,                                    \label{f32:eq:squarelinear}\\
 x+xy&\eq(xy)^2,                                     \label{f32:eq:linearprefix}\\
 x+yz&\eq x+yz+yx.                                   \label{f32:eq:linearterminal}
\end{align}

\Needspace{4\baselineskip}
\begin{theorem}\label{f32:thm:constraints}
A complete basis for \(A_{398}\) consists of \(\mathcal W\),
\eqref{f32:eq:union}--\eqref{f32:eq:linearterminal}, and \(2x\eq x\).
It has fourteen identities. Replacing \(2x\eq x\) by
\begin{equation}\label{f32:eq:parity}
 3x\eq x,\qquad 2xy\eq xy,\qquad 2x+yz\eq yz+yx
\end{equation}
gives a complete sixteen-identity basis for \(A_{390}\).
\end{theorem}
\begin{proof}
The tables verify every identity. Expand terms and apply
Lemma~\ref{f32:lem:words}. Linear coefficients lie in \(\{0,1\}\)
or \(\{0,1,2\}\), respectively; in the latter case positive
coefficients reduce according to parity. Nonlinear multiplicities
are one. A term with no nonlinear word is a pure nonempty linear sum.

First suppose that a square-support summand occurs. Express it as
\(w_Q^2\). Law \eqref{f32:eq:squareunion}, with \(x=w_Q\),
merges it with any other nonlinear word \(yz\), giving a square
on the union of the two supports. Law \eqref{f32:eq:squarelinear}
then incorporates every linear summand, one copy at a time. The
result is \(s_C\), where \(C\) is the full content of the term.

Otherwise each nonlinear summand is \(w_Pt\) with \(t\notin P\).
Law \eqref{f32:eq:union} merges two summands by taking the union of
their prefix supports. More generally, substitute a sum of terminals
for \(y\) in that law and distribute. Induction on the number of
nonlinear summands gives
\[
 \sum_{t\in T}w_Pt,
\]
where \(P\) is the union of all original prefix supports and
\(T\) the set of original terminals. If \(P\cap T\ne\varnothing\),
a square-support summand is present, so the preceding reduction
applies. If a linear variable \(x\) belongs to \(P\), rearrange
one nonlinear word to begin with \(x\), apply
\eqref{f32:eq:linearprefix}, and again finish with the square reduction.
This also works when the linear coefficient is two: use one copy
first and then absorb the other.

It remains that the prefix support \(P\) is disjoint from every
terminal and linear variable. Both \(P\) and \(T\) are nonempty.
For a single linear summand \(x\), law
\eqref{f32:eq:linearterminal}, with \(y=w_P\) and any existing
terminal substituted for \(z\), adjoins \(w_Px\).
In the parity case, the last law of \eqref{f32:eq:parity} replaces a
doubled linear variable by this terminal monomial while retaining
the existing nonlinear summand. Consequently the normal form is
\begin{equation}\label{f32:eq:proper}
 \sum_{x\in L}x+\sum_{t\in D}w_Pt,
 \qquad \varnothing\ne P,\quad \varnothing\ne D,\quad
 P\cap D=\varnothing,\quad L\subseteq D.
\end{equation}
Here \(D\) is the union of the original terminals and linear
variables. The set \(L\) consists of all original linear variables
in the idempotent case, or just those with odd positive coefficients
in the parity case. This gives three types: pure linear sums,
square-support words, and forms \eqref{f32:eq:proper}.

We now separate these forms on an arbitrary finite \(X\).
Assigning every variable \(3\) gives a value in \(\{2,3\}\)
on a pure linear sum, but \(0\) on either nonlinear type.
For pure linear sums, support is recovered by assigning one chosen
variable \(0\) and all others \(1\): the value is \(0\)
exactly when that variable is present. In the parity case, assign
one present variable \(3\) and all others \(2\); the value
is \(3\) or \(2\) according as its coefficient is one or two.

A square-support word \(s_C\) takes value \(1\) exactly when
every variable of \(C\) is \(1\), and otherwise takes value
\(0\). It has no assignments with value in \(\{2,3\}\).
A proper form has value \(1\) exactly when all variables in
\(P\cup D\) are \(1\). It has value in \(\{2,3\}\)
exactly on the nonempty set of assignments
\begin{equation}\label{f32:eq:box}
 a(P)=\{1\},\qquad a(D)\subseteq\{2,3\},
\end{equation}
with all other coordinates unrestricted. Indeed, each nonlinear
summand must be nonzero; mixed values \(1\) and \(2\) in the
sum give \(0\). Linear variables belong to \(D\), so they
obey the same restriction.

The inverse image of \(\{2,3\}\) therefore distinguishes the
two nonlinear types. Its coordinate projections recover \(P\)
and \(D\): the allowed values are respectively \(\{1\}\),
\(\{2,3\}\), and \(\{0,1,2,3\}\) outside their union.
To recover \(L\), fix \(P\) at \(1\), fix all of \(D\)
at \(2\), and change a chosen \(d\in D\) to \(3\).
The result is \(3\) exactly when \(d\in L\), for both
additions. Finally, the inverse image of \(1\) recovers the
support \(C\) of a square form by the same single-zero tests
used above. Thus distinct canonical forms define distinct functions.
Every valid identity has the same canonical form on both sides
and is derivable from the proposed finite basis.
\end{proof}

\subsubsection{A completed multiplicative sector}

The last two listed generators have commutative
multiplication. They are term equivalent to two of our direct
generators, as follows.
\begin{theorem}\label{f32:thm:terms}
The semirings \(A_{1348}\) and \(A_{1356}\) have complete bases
of sizes seventeen and fifteen, respectively. With \(A,B\)
equal to \((A_{1348},A_{387})\) or \((A_{1356},A_{395})\),
retain addition and define
\begin{equation}\label{f32:eq:termdefs}
 x\circ y=(x\cdot x)\cdot y\quad\text{in }A,
 \qquad x\cdot y=x\circ y+y\circ x\quad\text{in }B.
\end{equation}
These operations recover exactly the tables of \(B\) and \(A\),
respectively, on the same element labels.
\end{theorem}
\begin{proof}
The operation tables verify \eqref{f32:eq:termdefs}. Let $\sigma$
interpret multiplication of $B$ in $A$ by $u^2v$. With $\Sigma_B$
from Theorem~\ref{f32:thm:anti}, Proposition~\ref{prop:mutual-interpretations}
gives the complete basis
\begin{equation}\label{f32:eq:translatedbasis}
\SR\ \cup\ \sigma(\Sigma_B)\ \cup\ \{xy\eq x^2y+y^2x\}.
\end{equation}
Counting the two common additive axioms once gives
$5+(|\Sigma_B|-2)+1$, namely seventeen or fifteen identities.
\end{proof}

\endgroup
\subsection{Cliques, matchings and coefficients}\label{sec:family34}
\begingroup
\providecommand{\eq}{}\renewcommand{\eq}{\approx}
\providecommand{\SR}{}\renewcommand{\SR}{\mathsf{SR}}
\providecommand{\Var}{}\renewcommand{\Var}{\operatorname{Var}}
\providecommand{\St}{}\renewcommand{\St}{\operatorname{St}}
\emph{Generators (basis group):} $A_{439}$ (III), $A_{443}$ (III), $A_{461}$ (I), $A_{653}$ (VI), $A_{654}$ (VI), $A_{1350}$ (III), $A_{1353}$ (IV), $A_{1361}$ (II), $A_{1470}$ (V), $A_{1772}$ (I).

All six bases include $\SR$ and
\begin{equation}\label{f34:eq:quadratic}
 xy\eq yx,\qquad xyz\eq xy+xz+yz.
\end{equation}
For Groups I and II add
\begin{equation}\label{f34:eq:idem}
 x+x\eq x,\qquad x^2+y+xy\eq x^2+y,\qquad
 x^2+yz+xy\eq x^2+yz.
\end{equation}
Group II has the additional identity
\begin{equation}\label{f34:eq:collapse}
 x+x^2\eq x^2.
\end{equation}
Thus Groups I and II have ten and eleven identities, respectively.

For Groups III and IV instead add to $\SR$ and
\eqref{f34:eq:quadratic} the identities
\begin{equation}\label{f34:eq:weighted}
 2xy\eq xy,\qquad x^2+2y+xy\eq x^2+2y.
\end{equation}
Here $kt$ denotes a sum of $k$ copies of the term $t$.
Group III further includes
\begin{equation}\label{f34:eq:threshold}
 3x\eq2x,\qquad x+xy\eq xy,
\end{equation}
giving eleven identities. Group IV further includes
\begin{equation}\label{f34:eq:parity}
 3x\eq x,\qquad 2x+xy\eq xy,\qquad x+x^2\eq x^2,
\end{equation}
giving twelve identities. Group V instead adds to $\SR$ and \eqref{f34:eq:quadratic}
\begin{equation}\label{f34:eq:clique}
 x+x\eq x,\qquad x+xy\eq xy,\qquad x^2+y^2+xy\eq x^2+y^2,
\end{equation}
giving ten identities. Group VI adds to $\SR$ and \eqref{f34:eq:quadratic}
\begin{equation}\label{f34:eq:matching-coeff}
 3x\eq2x,\qquad 2x\eq2x^2,
\end{equation}
\begin{equation}\label{f34:eq:matching-edges}
 xy+xz\eq2x^2+2y^2+2z^2,\qquad x+xy\eq2x^2+2y^2,
\end{equation}
giving eleven identities. Denote the six displayed bases by
$\Sigma_{\mathrm I},\ldots,\Sigma_{\mathrm{VI}}$.
Direct evaluation of the tables verifies the relevant identities
on each generator.
The proofs below establish their completeness without a variable bound.

In Groups I--V, quadratic monomials are additively idempotent.
Distributivity first expands terms into nonempty sums of words.
For a word of length at least three, apply
$xyz\eq xy+xz+yz$ with $x,y$ the first two letters and $z$ the
remaining nonempty word. Each resulting word is shorter. Induction
therefore expresses every nonlinear word as a sum of quadratic
monomials in all six groups. In Groups I--V these are the products of every pair of
distinct positions of the original word, with duplicate quadratic
monomials removed. Repeated letters can produce squares.
In Group VI repetitions must be retained and are handled by its
coefficient identities. No empty product or empty sum is introduced.

\subsubsection{Idempotent addition and marked graphs}
Fix any finite variable set $X$. For $Q\subseteq T\subseteq X$ put
\[
 \St_Q(T)=\{\{x,y\}\subseteq T:x\ne y,\ \{x,y\}\cap Q\ne\varnothing\}.
\]
This is the set of edges meeting a marked vertex.
For a graph $F$, write $v(F)$ for its incident vertices and
$e(F)=\sum_{\{x,y\}\in F}xy$.
In all displayed polynomials, an empty part is omitted. Only
nonempty resulting terms are used.

\begin{theorem}\label{f34:thm:idem}
The basis $\Sigma_{\mathrm I}$ is complete for both $A_{461}$ and
$A_{1772}$. The basis $\Sigma_{\mathrm{II}}$ is complete for $A_{1361}$.
In particular $\Var(A_{461})=\Var(A_{1772})$.
\end{theorem}
\begin{proof}
After quadratic reduction let $L$ be the set of linear variables,
$C$ the set of all occurring variables, and $Q$ the variables whose
squares occur. Let $E$ be the set of nonsquare edges having neither
endpoint in $Q$.
The two absorption identities in \eqref{f34:eq:idem} adjoin every
edge in $\St_Q(C)$: for a marked vertex $x$, another variable
$y\in C$ occurs either as a linear term $y$ or in a quadratic
term $yz$, allowing $z=y$ for a square. Thus the polynomial is
equivalent to
\begin{equation}\label{f34:eq:nf-idem}
 N(L,C,Q,E)=\sum_{x\in L}x+\sum_{x\in Q}x^2
             +e(E)+e(\St_Q(C)).
\end{equation}
Here $C\ne\varnothing$, $L,Q\subseteq C$,
$E\subseteq\binom{C\setminus Q}{2}$, and the displayed terms
must together have content $C$. For Group II,
\eqref{f34:eq:collapse} removes the linear terms in $Q$, so require
$L\cap Q=\varnothing$. These operations preserve $C,Q,E$ and
the allowed remaining $L$. Every parameter tuple satisfying these
conditions already gives a normal form.

We recover the parameters from the term function using explicit
assignments. A perturbation of a baseline $b$ means that all
variables except the specified ones receive $b$. Use
\[
 (b,z)=(1,0)\text{ for }A_{461},A_{1361},\qquad
 (b,z)=(0,1)\text{ for }A_{1772}.
\]
At the baseline the polynomial has value $b$. Changing one variable
$x$ to $z$ gives value $z$ exactly when $x\in C$, recovering $C$.
Changing $x\in C$ to $2$ gives value $z$ exactly when $x\in Q$,
and otherwise gives $2$, recovering $Q$.
For distinct $x,y\in C\setminus Q$, change both to $2$.
The value is $z$ exactly when $\{x,y\}\in E$; otherwise it is
$2$. These assignments recover $E$.

It remains to recover $L$. In $A_{461}$ a perturbation $x=3$
has value $3$ exactly when $x\in L$: a nonlinear occurrence of
$x$ gives $0$, which is absorbed by a linear value $3$.
In $A_{1772}$ the same perturbation gives $3$ exactly when
$x\in L$, and otherwise $0$. In $A_{1361}$, for
$x\in C\setminus Q$, the perturbation gives $3$ exactly when
$x\in L$, and otherwise $2$. Linear terms at marked vertices
have already been removed in this last case.
All these statements follow from the defining tables and remain
valid when some polynomial parts are absent.

The parameters are therefore separated, proving completeness by
Lemma~\ref{lem:normal-form-completeness}. The common basis in Group I
also gives equality of the two generated varieties.
\end{proof}

\subsubsection{Nonidempotent addition and coefficients}
The graph representation persists, but linear multiplicities cannot
be discarded. The two coefficient laws lead to different normal forms.

\begin{theorem}\label{f34:thm:weighted}
The basis $\Sigma_{\mathrm{III}}$ is complete for each of
$A_{439},A_{443},A_{1350}$, and these three algebras generate the
same variety. The basis $\Sigma_{\mathrm{IV}}$ is complete for $A_{1353}$.
\end{theorem}
\begin{proof}
For Group III, reduce each positive linear coefficient to $1$ or $2$
using $3x\eq2x$. Let $D$ consist of all variables occurring either
in a quadratic monomial or as a linear term of coefficient two.
Every quadratic monomial absorbs its linear endpoints by
$x+xy\eq xy$, so it also absorbs their doubles. Adjoin these
doubles. A coefficient-one linear term whose variable is now in
$D$ can be removed; equivalently, its coefficient together with
the adjoined double reduces from three to two. Let $L$ be the
remaining coefficient-one variables. Then $L\cap D=\varnothing$.
As before let $Q$ record squares and let $E$ record edges not
meeting $Q$. The identity $x^2+2y+xy\eq x^2+2y$ adjoins every
edge in $\St_Q(D)$. We obtain the normal form
\begin{equation}\label{f34:eq:nf-threshold}
 N_3(L,D,Q,E)=\sum_{x\in L}x+\sum_{x\in D}2x
              +\sum_{x\in Q}x^2+e(E)+e(\St_Q(D)),
\end{equation}
where $Q\subseteq D$, $E\subseteq\binom{D\setminus Q}{2}$,
and $L\cup D\ne\varnothing$. Isolated vertices of $D$ are
permitted because their doubled linear terms are retained. A marked
vertex is never joined to $L$ by this construction.

To separate these forms, use baseline $1$ for all three generators.
Changing $x$ to $0$ gives $0$ exactly when $x\in C=L\cup D$.
Changing $x\in C$ to $2$ gives $0$ exactly when $x\in Q$.
Changing distinct $x,y\in C\setminus Q$ to $2$ gives $0$
exactly when $\{x,y\}\in E$; otherwise it gives $2$.
Finally, changing $x\in C$ to $3$ gives $3$ exactly when
$x\in L$. For a vertex of $D$, the doubled linear term has
value $0$ in $A_{439}$ and $A_{443}$, and value $2$ in
$A_{1350}$; the remaining terms cannot change that value to $3$.
Thus the function recovers $C,Q,E,L$, and then $D=C\setminus L$.
The common complete basis follows as in Theorem~\ref{f34:thm:idem}.

For Group IV, the law $3x\eq x$ reduces a positive linear
coefficient to one when odd and two when even. In particular
$4x\eq2x$, so a doubled linear term can be adjoined to any
existing linear term. The law $2x+xy\eq xy$ adjoins a double
for every endpoint of a quadratic term. If $C$ is the total content,
we can therefore retain $2x$ for every $x\in C$.
The remaining odd linear terms at square vertices disappear by
$x+x^2\eq x^2$. Let $O\subseteq C\setminus Q$ be the other
odd linear variables. Adjoining the edges meeting squares gives
\begin{equation}\label{f34:eq:nf-parity}
 N_4(C,O,Q,E)=\sum_{x\in C}2x+\sum_{x\in O}x
              +\sum_{x\in Q}x^2+e(E)+e(\St_Q(C)),
\end{equation}
where $C\ne\varnothing$, $Q\subseteq C$, and
$E\subseteq\binom{C\setminus Q}{2}$. The displayed three
copies at a vertex of $O$ are equivalent to one copy, but retaining
the double explicitly makes the content available for edge insertion.

Use baseline $1$ in $A_{1353}$. The perturbations by $0$ and
$2$ just used recover $C$, $Q$, and $E$. For
$x\in C\setminus Q$, change $x$ to $3$. Quadratic occurrences
of $x$ have value $2$, its doubled linear term has value $2$,
and all terms omitting it have value $1$. As $2+2=2$ and
$2+3=3$, the total is $3$ exactly when $x\in O$, and is
$2$ otherwise. This recovers $O$ and proves separation.
The normal-form derivations and separation again prove completeness
for every finite variable set. Empty components are omitted throughout;
each normal form contains at least one actual monomial.
\end{proof}

\subsubsection{Cliques and matching blocks}
\begin{theorem}\label{f34:thm:extra}
The basis $\Sigma_{\mathrm V}$ is complete for $A_{1470}$.
The basis $\Sigma_{\mathrm{VI}}$ is complete for both $A_{653}$
and $A_{654}$, which consequently generate the same variety.
\end{theorem}
\begin{proof}
For Group V, quadratic reduction and $x+xy\eq xy$ remove any
linear term occurring in a quadratic monomial. Let $C$ be the total
content and $Q$ the set of square vertices. The last identity in
\eqref{f34:eq:clique} completes the clique on $Q$. Retain the set
$E$ of all other edges; each such edge has at least one endpoint
outside $Q$. The normal form is
\begin{equation}\label{f34:eq:nf-clique}
 N_5(C,Q,E)=\sum_{x\in C\setminus(Q\cup v(E))}x
             +\sum_{x\in Q}x^2+e(E)+e\left(\binom Q2\right),
\end{equation}
where $C\ne\varnothing$, $Q\subseteq C$, and
$E\subseteq\binom C2\setminus\binom Q2$.
Use baseline $1$ in $A_{1470}$. Changing a variable $x$ to $0$
gives $0$ precisely when $x\in C$. For $x\in C$, changing it
to $3$ gives $0$ precisely when $x\in Q$, and otherwise gives
$3$. For an unordered pair with an endpoint outside $Q$, orient
it as $(x,y)$ with $x\notin Q$, and put $x=3,y=2$.
The value is $0$ precisely when $\{x,y\}\in E$, and otherwise
is $3$. This recovers $C,Q,E$ and proves completeness.

For Group VI, reduce all positive coefficients to one or two by
$3x\eq2x$. A doubled linear term becomes $2x^2$ by
\eqref{f34:eq:matching-coeff}; call such a vertex forced. The first
identity in \eqref{f34:eq:matching-edges} changes two quadratic terms
sharing a vertex into doubled squares of their endpoints. It also
handles a repeated quadratic term by setting $z=y$. A linear term
incident to a quadratic term forces both endpoints by the second
identity. A forced vertex propagates across any incident edge,
because the shared-edge identity with $z=x$, followed by adding
$x^2$ and applying the coefficient law, gives
\[
 2x^2+xy\eq2x^2+2y^2.
\]
Thus each affected graph component reduces to doubled squares of
all its vertices. What remains is a family $\mathcal H$ of
mutually disjoint blocks of size one or two. A singleton block
$\{x\}$ represents $x^2$, and a two-element block $\{x,y\}$
represents $xy$. Each occurs once. These blocks are disjoint from
the forced set $F$ and from the set $L$ of remaining linear variables.
We obtain
\begin{equation}\label{f34:eq:nf-matching}
 N_6(L,F,\mathcal H)=\sum_{x\in L}x+\sum_{x\in F}2x^2
                       +\sum_{H\in\mathcal H}m_H.
\end{equation}
Here $L,F$, and the blocks in $\mathcal H$ are pairwise disjoint,
and their union $C$ is nonempty. The reduction is finite: a
collision of edges or an incident linear term forces its endpoints,
and propagation removes incident edges while enlarging $F$.

Use baseline $2$ in either generator. Changing $x$ to $0$ gives
$0$ precisely when $x\in C$. Changing $x\in C$ to $3$ gives
$3$ precisely when $x\in L$, and otherwise gives $0$.
For any nonempty $T\subseteq C\setminus L$, put the variables
of $T$ at $1$, retaining baseline $2$ elsewhere. The value is
$3$ precisely when $T$ is one block of $\mathcal H$.
Indeed, meeting $F$ produces a zero doubled square; meeting only
part of a two-element block produces a zero product; and including
two complete blocks produces $3+3=0$. Exactly one complete block
produces $3$, with all other terms having the additive neutral value
$2$. The element $0$ absorbs all sums in both generators.
Tests with $|T|=1,2$ therefore recover every block, and
$F=C\setminus(L\cup\bigcup\mathcal H)$ is then determined.
The complete normal forms are separated for arbitrary $X$, proving
both assertions by equational reduction.
\end{proof}

\endgroup
\subsection{Component partitions and Boolean parity}\label{sec:family35}
\begingroup
\providecommand{\eq}{}\renewcommand{\eq}{\approx}
\providecommand{\SR}{}\renewcommand{\SR}{\mathsf{SR}}
\providecommand{\Var}{}\renewcommand{\Var}{\operatorname{Var}}
\emph{Generators (basis size):} $A_{663}$ (13), $A_{1000}$ (15), $A_{1441}$ (12), $A_{1910}$ (15), $A_{1944}$ (11), $A_{1948}$ (12).

For a fixed finite variable set $X$ and $\varnothing\ne C\subseteq X$,
$K\subseteq C$, write
\[
 m(C,K)=\prod_{x\in C}x^{\epsilon_x},\qquad
 \epsilon_x=2\ (x\in K),\quad \epsilon_x=1\ (x\in C\setminus K).
\]
Fix an ordering of $X$ to order the factors. Multiplication is
commutative in the four direct generators. Empty parts of displayed
polynomials are omitted; every whole term is nonempty. We never
substitute an empty word for a variable.

\subsubsection{Cubic constraints and component partitions}
Let $\Sigma_C$ consist of $\SR$ and the following eight identities:
\begin{equation}\label{f35:eq:component-common}
 xy\eq yx,\qquad x+x\eq x,
\end{equation}
\begin{equation}\label{f35:eq:component-cubes}
 xyz\eq x^3+y^3+z^3,\qquad x^3y^3\eq x^3+y^3,
\end{equation}
\begin{equation}\label{f35:eq:component-edges}
 xy+x^2\eq xy,\qquad xy+yz+xz\eq xy+yz,
\end{equation}
\begin{equation}\label{f35:eq:component-forcing}
 x^3+xy\eq x^3+y^3,\qquad x+xy\eq x^3+y^3.
\end{equation}
These thirteen identities are valid in $A_{663}$ by its tables.

\begin{theorem}\label{f35:thm:component}
The set $\Sigma_C$ is a complete identity basis for $A_{663}$.
\end{theorem}
\begin{proof}
Distributivity expands a term into a sum of words. By commutativity
and the cubic product law, the cube of any word equals the sum of
the cubes of its letters. Apply the triple law with the first two
letters and the remaining nonempty word: every word of length at
least three becomes a sum of unary cubes. Additive idempotence
removes repetitions. Thus only linear terms, squares, edges $xy$
with $x\ne y$, and cubes remain.

Regard the quadratic part as an undirected graph, allowing loops
for squares. An edge absorbs the square at each endpoint, so these
loops may be inserted. The second identity in
\eqref{f35:eq:component-edges} inserts the third edge of any two-edge
path. Consequently every connected component can be completed to
the graph containing every pair and every loop on its vertices.
A component of one vertex is represented by its square.

Call a vertex forced if its cube occurs. The first identity in
\eqref{f35:eq:component-forcing} replaces an edge touching a forced
vertex by cubes at both endpoints, so forcing propagates throughout
its component. A linear variable incident to a quadratic term
starts the same propagation by the second identity. Cubes also
absorb their own linear terms: set $y=x^2$ in $x+xy\eq x^3+y^3$
and use $x^6\eq x^3$, which follows from the cubic product law.
A cube absorbs its own square by setting $y=x$ in the other forcing
identity. These operations reduce the term to
\begin{equation}\label{f35:eq:component-normal}
 N(L,F,\mathcal H)=\sum_{x\in L}x+\sum_{x\in F}x^3
       +\sum_{H\in\mathcal H}\ \sum_{\substack{x,y\in H\\x\le y}}xy.
\end{equation}
The sets $L,F$ and the nonempty blocks $H\in\mathcal H$ are
pairwise disjoint. Their union $C$ is nonempty. Unforced quadratic
components have become the blocks of $\mathcal H$, while all other
quadratic components have been absorbed into $F$.

Use assignments with baseline value $2$. Changing a single variable
$x$ to $0$ gives value $0$ exactly when $x\in C$. Changing
$x\in C$ to $3$ gives $3$ exactly when $x\in L$, and otherwise
$0$. Thus the function determines $C$ and $L$.
For $\varnothing\ne T\subseteq C\setminus L$, set exactly the
variables of $T$ to $1$, keeping all others at $2$. A block term
has value $3$ if all its variables are $1$, value $2$ if all are
$2$, and value $0$ otherwise. A cube is $0$ at $1$ and $2$ at $2$.
Since $0$ absorbs sums and $3+3=3$, the polynomial has value $3$
precisely when $T$ is a nonempty union of whole blocks of
$\mathcal H$. Its inclusion-minimal such subsets are exactly
the blocks. Recovering them also recovers
$F=C\setminus(L\cup\bigcup\mathcal H)$.
These tests apply to arbitrarily large blocks. Distinct normal forms
therefore define distinct functions. Every valid identity reduces
to the same normal form on its two sides, proving completeness for
arbitrary $X$.
\end{proof}

\subsubsection{Doubled parts and Boolean parity}
Put $U=A_{1944}$, and let $A$ be either $A_{1441}$ or $A_{1948}$
when discussing the second basis. Both bases contain $\SR$ and
\begin{equation}\label{f35:eq:period-common}
 xy\eq yx,\qquad x^3\eq x^2,\qquad 3x\eq x.
\end{equation}
Here $kt$ denotes the sum of $k$ copies of a term $t$.
For $U$ add
\begin{equation}\label{f35:eq:upper-raise}
 x+x^2\eq2x,\qquad 2x+2xy\eq2x,
\end{equation}
\begin{equation}\label{f35:eq:upper-meet}
 2xy^2+2x^2y+2xy\eq2xy^2+2x^2y.
\end{equation}
This gives an eleven-identity basis $\Sigma_U$.
For the other two generators instead add
\begin{equation}\label{f35:eq:lower-lower}
 x+2x^2\eq x^2,
\end{equation}
\begin{equation}\label{f35:eq:lower-triangle}
 2xyz\eq2xy+2xz+2yz,
\end{equation}
\begin{equation}\label{f35:eq:lower-edges}
 2x+2xy\eq2xy,\qquad 2x^2+2y+2xy\eq2x^2+2y.
\end{equation}
This gives a twelve-identity basis $\Sigma_D$.
Every displayed identity is verified on the respective tables.

The law $3x\eq x$ implies $4x\eq2x$. The map $\delta(a)=2a$
therefore preserves both operations: in particular
$(2a)(2b)=4ab=2ab$. Its image is the subsemiring $T=\{0,1,3\}$,
which it fixes pointwise. In all three generators, $T$ has
multiplicative identity $1$, zero $0$, and $3^2=0$.
Addition on $T$ is maximum in the order $0<3<1$ for $U$, and
maximum in $1<3<0$ for the other two generators.
Identities between doubled polynomials can thus be analyzed on
these three-element images: a doubled polynomial evaluated at an
assignment equals its value at the assignment followed by $\delta$.

\begin{theorem}\label{f35:thm:parity}
The basis $\Sigma_U$ is complete for $A_{1944}$.
The basis $\Sigma_D$ is complete for both $A_{1441}$ and $A_{1948}$;
these two algebras generate the same variety.
\end{theorem}
\begin{proof}
Commutativity and $x^3\eq x^2$ reduce words to $m(C,K)$.
The additive period law reduces each positive coefficient to one
or two. We first prove a finite-basis normal-form theorem for the
doubled part in each case, then treat the odd coefficients.

\emph{The doubled part for $U$.}
Let $d=\sum_{i\in I}2m(C_i,K_i)$ be nonempty, and set
$I_C=\{i:C_i\subseteq C\}$. Direct evaluation on $T$, ordered
as $0<3<1$, shows that
\begin{equation}\label{f35:eq:upper-criterion}
 d+2m(C,K)\eq d\quad\Longleftrightarrow\quad
 I_C\ne\varnothing\ \hbox{and}\ \bigcap_{i\in I_C}K_i\subseteq K.
\end{equation}
For necessity put $1$ on $C$ and $0$ outside $C$, proving that
$I_C$ is nonempty. For $t\in C\setminus K$, instead set $t=3$.
The proposed absorbed term has value $3$, so some eligible monomial
must omit $t$ or contain it unsquared. This proves the intersection
condition. Conversely a nonzero target value is either $1$, when
all of $C$ is at $1$, or $3$, when exactly one unsquared variable
is at $3$ and all others are at $1$. An eligible monomial supplied
by the respective condition has value at least the target value.
This proves the criterion on $T$ and hence for doubled terms in $U$.

Every absorption in \eqref{f35:eq:upper-criterion} has a derivation
from $\Sigma_U$. Doubling $x+x^2\eq2x$ gives
$2x+2x^2\eq2x$. Thus squared-variable sets can be enlarged in
an absorbed doubled monomial. The law $2x+2xy\eq2x$ extends an
eligible support $C_i$ to $C$ by multiplying by the product of
the missing variables, retaining its squared set $K_i$.
If no variable is missing, no multiplication is needed.

For two monomials on the same support $C$, we can intersect their
squared sets. Factor their common exponent-one and exponent-two
parts into a word $g$. Let $p$ and $q$ be the squarefree products
of the variables with opposing exponents. If both opposing groups
are nonempty, their two remaining factors have the form $pq^2$
and $p^2q$. Apply \eqref{f35:eq:upper-meet}, with a multiplicative
context $g$ when it is nonempty, to adjoin $2gpq$. This has squared
set equal to the intersection. If either opposing group is empty,
the desired monomial is already one of the two sources. No empty
word substitution occurs. Iterating intersections and then enlarging
the squared set proves every absorption in the criterion.

Let $S(d)$ contain all pairs $(C,K)$ satisfying this criterion.
It is nonempty and finite for fixed $X$. Adjoining all these terms
gives a canonical doubled polynomial
\[
 D_U(S)=\sum_{(C,K)\in S}2m(C,K).
\]
Two doubled polynomials define the same function exactly when they
have the same absorbed set $S$: equality implies identical absorption
tests, and their common saturated sum gives the converse and its
equational derivation.

\emph{The doubled part for $A_{1441}$ and $A_{1948}$.}
Use \eqref{f35:eq:lower-triangle} repeatedly to reduce every doubled
word of length at least three to doubled quadratic monomials.
Repeated doubled terms disappear by $4t\eq2t$. Insert $2x$ for
every occurring variable by the first law in
\eqref{f35:eq:lower-edges}, including the case $y=x$.
Let $C$ be the total content, $Q$ the square vertices, and $E$ the
edges whose endpoints are outside $Q$. The second edge law removes
every edge meeting $Q$, since a doubled linear term is available
at its other endpoint. Thus the doubled part has canonical form
\begin{equation}\label{f35:eq:lower-normal}
 D_D(C,Q,E)=\sum_{x\in C}2x+\sum_{x\in Q}2x^2
                      +\sum_{\{x,y\}\in E}2xy.
\end{equation}
Here $C\ne\varnothing$, $Q\subseteq C$ and
$E\subseteq\binom{C\setminus Q}{2}$.
On $T$, use baseline $1$. A perturbation $x=0$ detects $x\in C$.
For $x\in C$, the perturbation $x=3$ gives $0$ exactly when
$x\in Q$. For distinct $x,y\in C\setminus Q$, putting both at
$3$ gives $0$ exactly when $\{x,y\}\in E$; otherwise it gives
$3$. These tests recover $C,Q,E$. They prove completeness of the
doubled reduction for both generators, whose doubled images are
identical.

\emph{Decomposition and separation of the odd part.}
In the upper case, $x+x^2\eq2x$ and $3x\eq x$ give
$x\eq x^2+2x$. Multiplication by a word context and repeated
exponent raising therefore yield
\begin{equation}\label{f35:eq:upper-decomposition}
 m(C,K)\eq m(C,C)+2m(C,K).
\end{equation}
The intermediate doubled terms are absorbed by $2m(C,K)$ using
the doubled exponent-raising law already proved.
In the other case, $x^2\eq x+2x^2$ instead gives
\begin{equation}\label{f35:eq:lower-decomposition}
 m(C,K)\eq m(C,\varnothing)+2m(C,K).
\end{equation}
Here doubling the same identity gives $2x+2x^2\eq2x^2$, so the
intermediate doubled terms are again absorbed, now in the reverse
exponent direction. These arguments also cover one-letter words
without a multiplicative context.

For any nonempty polynomial $f$, these decompositions give
\begin{equation}\label{f35:eq:full-normal}
 f\eq D+\sum_{C\in\mathcal O}r_C,
\end{equation}
where $D$ is the canonical doubled form of $2f$,
$r_C=m(C,C)$ in the upper case and $r_C=m(C,\varnothing)$ in
the other case. The set $\mathcal O$ records, modulo two, the
number of original summands with support $C$.
For each such support, $D$ absorbs $2r_C$ by the doubled exponent
law in the appropriate direction. Even copies of $r_C$ can
therefore be absorbed into $D$, leaving one copy exactly for
$C\in\mathcal O$. Also $2D\eq D$, and doubling
\eqref{f35:eq:full-normal} gives $D$ again.

If two normal forms define the same function, doubling their
functions recovers the same $D$, whose parameters were separated
above. To recover $\mathcal O$, use the two-element subring
$\{1,2\}$ for $A_{1944}$ and $A_{1948}$, or $\{0,2\}$ for
$A_{1441}$. In each case this is the field of two elements, with
identity $2$ and zero $\epsilon$ equal to $1$ or $0$, respectively.
On these assignments the nonempty doubled part has constant value
$\epsilon$. Put exactly the variables of $Y\subseteq X$ at $2$
and all others at $\epsilon$. The resulting parity is
\[
 b(Y)=\sum_{C\subseteq Y}\mathbf1_{C\in\mathcal O}\pmod2.
\]
Lemma~\ref{lem:finite-set-recovery}(i) recovers $\mathcal O$.
Both components are therefore determined by the function, proving
completeness by the criterion in Lemma~\ref{lem:normal-form-completeness}.
The same doubled parameters, odd representatives, and proof for
$A_{1441}$ and $A_{1948}$ give their common complete basis.
\end{proof}

\subsubsection{Two transfers by mutual term definitions}
On $A=A_{1000}$ or $A=A_{1910}$, with original product written
$\cdot$, define
\begin{equation}\label{f35:eq:forward}
 x\circ y=2(x\cdot y)+y\cdot x.
\end{equation}
With unchanged addition, this gives exactly $U=A_{1944}$ on the
same labels. In $U$ the original products are recovered by
\begin{equation}\label{f35:eq:backward}
 x\cdot y=(x\circ x)\circ y\quad(A=A_{1000}),\qquad
 x\cdot y=(x\circ y)\circ y\quad(A=A_{1910}).
\end{equation}
Both claims are checked on all sixteen ordered pairs.

\begin{theorem}\label{f35:thm:transfer}
Both $A_{1000}$ and $A_{1910}$ have complete bases of fifteen
identities. Let $\sigma$ recursively translate the $U$ product by
\eqref{f35:eq:forward}, and let $\tau$ translate the $A$ product by
its definition in \eqref{f35:eq:backward}. A complete basis is
\[
 \SR\ \cup\ \sigma(\Sigma_U)\ \cup\
 \{x\cdot y\eq\sigma\tau(x\cdot y)\}.
\]
The two duplicated additive commutativity and associativity axioms
are counted once.
\end{theorem}
\begin{proof}
The table definitions verify all displayed identities. The count is
$5+(11-2)+1=15$. The interpretations \eqref{f35:eq:forward}--\eqref{f35:eq:backward}
and the recovery identity satisfy
Proposition~\ref{prop:mutual-interpretations}. The complete source
basis from Theorem~\ref{f35:thm:parity} therefore gives completeness.
\end{proof}

\endgroup
\subsection{Multiplicity chains and forks}\label{sec:family36}
\begingroup
\providecommand{\eq}{}\renewcommand{\eq}{\approx}
\providecommand{\SR}{}\renewcommand{\SR}{\mathsf{SR}}
\providecommand{\Var}{}\renewcommand{\Var}{\operatorname{Var}}
\subsubsection{Generators}
Write $F=A_{1940}$ and $H=A_{1943}$.

\emph{Generators (basis size):} $A_{1940}$ (260), $A_{1943}$ (260).

In both generators, $0$ is a zero for multiplication and an identity
for addition, $2$ is the multiplicative identity, $1^2=1$,
$1\cdot3=3$, and $3^2=0$. Addition satisfies
$1+a=1$ for every $a$, $2+2=1$, and $3+3=3$.
The sole difference is $2+3=1$ in $F$ and $2+3=2$ in $H$.
The map $a\mapsto2a$ fixes $0,1,3$ and sends $2$ to $1$.
Its image has addition given by maximum in $0<3<1$.

Fix any finite set $X$ of variables. For $\varnothing\ne C\subseteq X$
and $K\subseteq C$, put
\[
 m(C,K)=\prod_{x\in C}x^{\epsilon_x},\qquad
 \epsilon_x=2\ (x\in K),\quad \epsilon_x=1\ (x\in C\setminus K).
\]

\subsubsection{Two explicitly finite bases}
Both bases contain $\SR$ and the following six identities:
\begin{equation}\label{f36:eq:common}
 xy\eq yx,\qquad x^3\eq x^2,\qquad 3x\eq2x,
\end{equation}
\begin{equation}\label{f36:eq:extend}
 2x+2xy\eq2x,\qquad 2x+2x^2\eq2x,
\end{equation}
\begin{equation}\label{f36:eq:meet}
 2xy^2+2x^2y+2xy\eq2xy^2+2x^2y.
\end{equation}
Here $kt$ denotes the sum of $k$ copies of a term $t$.

The remaining equations form a finite list, indexed by subsets of
the eight-element set
\[
 E=\{0,1,2\}^2\setminus\{(0,0)\}.
\]
For each $P\subseteq E$ for which some first coordinate and some
second coordinate are nonzero, introduce distinct variables
$z_{ab}$, $(a,b)\in P$, and define the nonempty words
\[
 u_P=\prod_{(a,b)\in P,\ a>0}z_{ab}^{a},\qquad
 v_P=\prod_{(a,b)\in P,\ b>0}z_{ab}^{b}.
\]
Define, for $(a,b)\in E$,
\begin{equation}\label{f36:eq:profile}
 f_F(a,b)=\begin{cases}1&1\in\{a,b\},\\2&\text{otherwise},\end{cases}
 \qquad
 f_H(a,b)=\begin{cases}2&a=b=2,\\1&\text{otherwise}.\end{cases}
\end{equation}
For $A\in\{F,H\}$ put
$w_{A,P}=\prod_{(a,b)\in P}z_{ab}^{f_A(a,b)}$ and include
\begin{equation}\label{f36:eq:pair}
 u_P+v_P\eq u_P+v_P+2w_{A,P}.
\end{equation}
Call the resulting basis $\Sigma_A$. There are
$255-3-3=249$ admissible subsets $P$, so the list has
$5+6+249=260$ equations. Each equation uses at most eight variables.

\subsubsection{The absorption criterion}
Distributivity and \eqref{f36:eq:common} expand every term into a
nonempty polynomial
\[
 p=\sum_{i\in I}m(C_i,K_i),
\]
in which each distinct monomial occurs once or twice. We regard
$I$ as an occurrence index set, so that two copies have distinct
indices. For a nonempty $C\subseteq X$, set
$I_C=\{i\in I:C_i\subseteq C\}$ and $U_i=C_i\setminus K_i$.

\begin{lemma}\label{f36:lem:absorption}
For $K\subseteq C$, the identity $p+2m(C,K)\eq p$ holds in $H$
if and only if
\begin{equation}\label{f36:eq:chain-criterion}
 |I_C|\ge2,\qquad \bigcap_{i\in I_C}K_i\subseteq K.
\end{equation}
It holds in $F$ if and only if $|I_C|\ge2$ and, for every
$t\in C\setminus K$, at least one of the following holds:
\begin{equation}\label{f36:eq:fork-criterion}
 t\in U_i\text{ for some }i\in I_C,
 \qquad
 |\{i\in I_C:t\notin C_i\}|\ge2.
\end{equation}
\end{lemma}
\begin{proof}
Assign $2$ to every variable in $C$ and $0$ outside $C$.
Exactly the occurrences indexed by $I_C$ have value $2$; all
others have value $0$. The added doubled monomial has value $1$.
Thus its absorption requires at least two eligible occurrences.

Now assign $3$ to $t\in C\setminus K$, $2$ to the other variables
of $C$, and $0$ outside $C$. The added doubled monomial has value
$3$. An eligible occurrence has value $3$, $0$, or $2$ according
as $t\in U_i$, $t\in K_i$, or $t\notin C_i$.
In $H$, absorption fails precisely when every eligible occurrence
has value $0$. This gives \eqref{f36:eq:chain-criterion}.
In $F$, absorption fails when there is no value $3$ and there
are fewer than two values $2$. This gives \eqref{f36:eq:fork-criterion}.

For sufficiency consider any assignment. If $2m(C,K)$ has value
$0$, it is absorbed. If it has value $1$, all variables of $C$
lie in $\{1,2\}$. Each eligible occurrence then has value $1$
or $2$, and any two such values sum to $1$.
Finally, if the doubled monomial has value $3$, exactly one variable
$t\in C\setminus K$ has value $3$ and all other variables of $C$
lie in $\{1,2\}$. An eligible occurrence with $t\in U_i$ has
value $3$. One omitting $t$ has value $1$ or $2$.
In $H$ either kind suffices to absorb $3$. In $F$ a value $3$
or two occurrences omitting $t$ suffice. The respective criterion
supplies these occurrences. Absorption by a subsum implies
absorption by the whole sum, by associativity and commutativity.
\end{proof}

The six common identities are valid by the tables. In particular,
the three doubled identities can be checked on $\{0,1,3\}$,
where addition is maximum in $0<3<1$.
For a pair of occurrences $i,j$, put $C_{ij}=C_i\cup C_j$ and
\begin{equation}\label{f36:eq:pair-sets}
 K^H_{ij}=K_i\cap K_j,\qquad
 K^F_{ij}=(K_i\cup K_j)\setminus(U_i\cup U_j).
\end{equation}
Lemma~\ref{f36:lem:absorption} shows that the sum of these two
occurrences absorbs $2m(C_{ij},K^A_{ij})$. For $H$, this is
exactly its intersection condition. For $F$, every variable of
$C_{ij}\setminus K^F_{ij}$ is unsquared in an occurrence.
Equations \eqref{f36:eq:pair} are precisely these valid pair insertions
when the variables are grouped by their two exponent profiles.
This proves the soundness of every equation of $\Sigma_A$.

\subsubsection{Derivations and completeness}
\begin{lemma}\label{f36:lem:derive}
Every absorption in Lemma~\ref{f36:lem:absorption} is derivable
from the corresponding basis $\Sigma_A$.
\end{lemma}
\begin{proof}
First, any two occurrences $m(C_i,K_i),m(C_j,K_j)$ permit insertion
of $2m(C_{ij},K^A_{ij})$. Group the variables of their union
according to the exponent pairs $(a,b)\in E$. For each nonempty
group substitute its squarefree product for $z_{ab}$ in
\eqref{f36:eq:pair}. Every substituted word is nonempty. Commutativity
then gives the desired insertion, in the context of the other
summands of $p$. Distinct indices may represent equal monomials;
the same argument still applies to their two occurrences.

Three further operations are derivable for doubled monomials.
The first law in \eqref{f36:eq:extend} extends a support by any
nonempty product of missing variables. The second law increases
an exponent from one to two, using the product of the other
factors as a multiplicative context when that product is nonempty.
Neither operation removes the source summand.

Two doubled monomials on the same support can also insert the
monomial with the intersection of their squared-variable sets.
Factor their common parts into a word $g$. Let $a$ and $b$ be
the squarefree products of the variables with opposing exponents.
If both opposing groups are nonempty, the two monomials have
the form $gab^2$ and $ga^2b$. Apply \eqref{f36:eq:meet}, using the
context $g$ only when nonempty, to insert $2gab$.
If one opposing group is empty, one squared set is contained
in the other and the desired intersection is already present.

Suppose the criterion holds for $(C,K)$. Insert the pair monomial
for every pair of distinct indices from $I_C$ and extend its
support to $C$. Such a pair exists since $|I_C|\ge2$.
Intersect all their squared sets to obtain a doubled monomial
on $C$ with squared set
\[
 J_A=\bigcap_{\{i,j\}\subseteq I_C}K^A_{ij}.
\]
For $H$, this is $\bigcap_{i\in I_C}K_i$.
For $F$, a variable $t\in C$ belongs to $J_F$ exactly when
no eligible occurrence contains it unsquared and at most one
eligible occurrence omits it. Indeed an unsquared occurrence
gives a pair excluding $t$, and two omissions give such a pair
as well. Conversely, in the absence of those two possibilities,
every pair contains a squared occurrence of $t$ and no unsquared
one. The respective criterion therefore says $J_A\subseteq K$.
Increasing exponents inserts $2m(C,K)$.

Only finitely many pair insertions and intersections are needed
for any given $p,C,K$. All were made in the additive context
of $p$ and previously inserted terms. Each previous term was
already derivably absorbed by $p$, so it can be removed again.
The result is a derivation of $p+2m(C,K)\eq p$ itself.
\end{proof}

\begin{theorem}\label{f36:thm:main}
The set $\Sigma_F$ is a complete finite identity basis for
$A_{1940}$, and $\Sigma_H$ is a complete finite identity basis
for $A_{1943}$.
\end{theorem}
\begin{proof}
For a polynomial $p$ on $X$, let $D_A(p)$ be the set of pairs
$(C,K)$ for which $p+2m(C,K)\eq p$ is valid in $A$.
By Lemma~\ref{f36:lem:derive}, all these absorptions are derivable.
Notice also that $p+2m\eq p$ implies $p+m\eq p$ using
$3m\eq2m$: indeed
$p+m\eq p+3m\eq p+2m\eq p$.

Call an occurrence essential if its support $C$ contains exactly
one occurrence support of $p$, counted with multiplicity, namely
itself. Such an occurrence appears once. Equivalently, its
support is minimal among the occurrence supports and has
multiplicity one. Denote the set of its pairs $(C,K)$ by $E(p)$.
An original occurrence belongs to $D_A(p)$ if and only if it is
not essential. Necessity follows from $|I_C|\ge2$.
For sufficiency, the occurrence itself verifies the intersection
condition for $H$ and supplies an unsquared occurrence for every
$t\in C\setminus K$ in $F$.

Insert every doubled monomial indexed by $D_A(p)$. The set is
finite, containing at most $3^{|X|}-1$ elements. Every original
nonessential occurrence is then removable since its doubled copy
is present and $3m\eq2m$. Thus the basis derives the canonical
expression
\begin{equation}\label{f36:eq:normal}
 N_A(p)=\sum_{(C,K)\in D_A(p)}2m(C,K)
       +\sum_{(C,K)\in E(p)}m(C,K).
\end{equation}
It is nonempty: if $p$ has only one occurrence it is essential;
if it has at least two, a pair insertion makes $D_A(p)$ nonempty.

It remains to show that the function of $p$ determines both sets.
It determines $D_A(p)$ by the definition of absorption.
For $\varnothing\ne C\subseteq X$, evaluate $p$ by putting $2$
on $C$ and $0$ outside $C$. Its value is $2$ exactly when
$|I_C|=1$. The inclusion-minimal sets $C$ giving value $2$
are precisely the supports of essential occurrences. To see this,
if the unique eligible occurrence has support $B\subseteq C$,
the assignment on $B$ also gives $2$; minimality forces $B=C$.
Conversely a single occurrence on a minimal support is the only
eligible occurrence there and no proper subset gives $2$.

For such an essential support $C$, change one variable $t\in C$
from $2$ to $3$, keeping $2$ on $C\setminus\{t\}$ and $0$
outside $C$. Every other occurrence vanishes. The value is $0$
when $t\in K$ and $3$ when $t\notin K$. Hence the function
recovers $K$ and therefore all of $E(p)$.

If $p\eq q$ is any valid identity, take $X$ to contain all
variables of both sides. Their functions give the same $D_A$
and $E$, so both sides reduce to exactly the same expression
\eqref{f36:eq:normal}. This proves completeness for arbitrary
finite variable sets, and hence for all identities.
\end{proof}

\begin{corollary}\label{f36:cor:proper}
$\Var(A_{1943})$ is a proper subvariety of $\Var(A_{1940})$.
\end{corollary}
\begin{proof}
The common laws hold in both generators. Each pair monomial
for $F$ has a squared set containing the corresponding squared
set for $H$, so its insertion is valid in $H$ by
Lemma~\ref{f36:lem:absorption}. Thus $H$ satisfies $\Sigma_F$.
The profile identity
\[
 x^2+y\eq x^2+y+2xy
\]
holds in $H$. In $F$, the assignment $x=3,y=2$ gives $2$ on
the left and $1$ on the right. The containment is therefore strict.
\end{proof}

\endgroup
\subsection{Linear layers and truncated coefficients}\label{sec:family37}
\begingroup
\providecommand{\eq}{}\renewcommand{\eq}{\approx}
\providecommand{\SR}{}\renewcommand{\SR}{\mathsf{SR}}
\providecommand{\Var}{}\renewcommand{\Var}{\operatorname{Var}}
\emph{Generators (basis size):} $A_{980}$ (265), $A_{983}$ (265), $A_{1614}$ (265), $A_{1619}$ (265), $A_{1758}$ (261), $A_{1763}$ (261), $A_{1881}$ (261), $A_{1882}$ (71).

For a finite variable set $X$, write
$m(C,K)=\prod_{x\in C}x^{\epsilon_x}$, where
$\varnothing\ne C\subseteq X$, $K\subseteq C$, and
$\epsilon_x=2$ on $K$ and $1$ on $C\setminus K$.
Put $s_C=m(C,\varnothing)$. Fixed orderings make sums and
products unambiguous. Empty summands are omitted; no whole term
is empty and no empty word is substituted into an identity.

\subsubsection{The two retained bases}
Let $F=A_{1940}$ and $H=A_{1943}$, and use the complete
260-identity bases $\Sigma_F,\Sigma_H$ of
Theorem~\ref{f36:thm:main}, with the profiles in
\eqref{f36:eq:profile} and pair identities \eqref{f36:eq:pair}.

We will also use a consequence of the pair equations. Given $m(C_i,K_i)$ and $m(C_j,K_j)$, the $H$ equations insert
\begin{equation}\label{f37:eq:retained-pair-consequence}
 2m(C_i\cup C_j,K_i\cap K_j).
\end{equation}
To obtain it, group variables by their two exponents and substitute
each nonempty group product into the appropriate equation
\eqref{f36:eq:pair}. This derivation never uses an empty word.

\subsubsection{Making all products additively idempotent}
For $B\in\{F,H\}$ define $B^d$ on the same set, retaining addition
and replacing multiplication by
\begin{equation}\label{f37:eq:double-product}
 a\circ b=2(ab).
\end{equation}
The relabelling $g=(1,0,3,2)$, meaning $g(0)=1$, $g(1)=0$,
$g(2)=3$, $g(3)=2$, gives isomorphisms
$A_{1758}\cong F^d$ and $A_{1763}\cong H^d$.
These statements follow directly from the operation tables.

\begin{theorem}
Complete bases for $A_{1758}$ and $A_{1763}$ are, respectively,
\begin{equation}\label{f37:eq:linear-bases}
 \Sigma_F\cup\{2xy\eq xy\},\qquad
 \Sigma_H\cup\{2xy\eq xy\}.
\end{equation}
Each listed basis has 261 equations.
\end{theorem}
\begin{proof}
The law $3x\eq2x$ makes every coefficient at least two equal
to two. Consequently a word of length at least two evaluated
in $B^d$ equals twice its value in $B$; a one-letter word is
unchanged. In particular, doubled words have the same values
in the two algebras.

The semiring axioms and \eqref{f36:eq:common} hold in
$B^d$ by distributivity and the coefficient law. All doubled
identities \eqref{f36:eq:extend}--\eqref{f36:eq:meet}
have their original values. For a pair equation, its left side
in $B^d$ is its left side in $B$ with a further copy of each
nonlinear source word added. The left side in $B$ absorbs the
specified doubled word, so adding these further copies preserves
that absorption. The doubled target has the same value in both
algebras. Thus $B^d$ satisfies $\Sigma_B$, as well as $2xy\eq xy$.

Expand any identity $p\eq q$ valid in $B^d$ into polynomials.
Let $p^\sharp,q^\sharp$ double every nonlinear monomial, leaving
the linear occurrences unchanged. The added product identity
derives $p\eq p^\sharp$ and $q\eq q^\sharp$.
The values of $p$ in $B^d$ are exactly those of $p^\sharp$ in
$B$, and likewise for $q$. Hence $p^\sharp\eq q^\sharp$ is
valid in $B$ and follows from its complete basis $\Sigma_B$.
Combining the three derivations proves $p\eq q$ from the
stated basis. The argument has no bound on the number of variables.
\end{proof}

\subsubsection{Support saturation at threshold two}
Relabel $A_{1881}$ by $g=(1,2,0,3)$ and call the resulting
algebra $Q$. In $Q$, $0$ is the additive identity and a
multiplicative zero; $2$ is the multiplicative identity.
The sum of any two nonzero elements is $1$. Products of nonzero
elements remain nonzero, $1a=1$ for nonzero $a$, and $3^2=1$.

\begin{theorem}
A complete 261-equation basis for $A_{1881}$ is
\begin{equation}\label{f37:eq:flat-basis}
 \Sigma_H\cup\{2x^2\eq2x\}.
\end{equation}
\end{theorem}
\begin{proof}
The equations are valid in $Q$. For the pair equations, if their
target is nonzero both source words are nonzero and their sum
is $1$, which absorbs every element. The other equations follow
immediately from the stated tables.

Normalize words to $m(C,K)$ and positive coefficients to one
or two. Index monomial occurrences of a polynomial $p$ by $I$,
and set $I_C=\{i:C_i\subseteq C\}$. Then
\begin{equation}\label{f37:eq:flat-criterion}
 p+2m(C,K)\eq p\text{ in }Q
 \quad\Longleftrightarrow\quad |I_C|\ge2.
\end{equation}
For necessity put $2$ on $C$ and $0$ elsewhere: the doubled
target is $1$, while zero or one eligible occurrence gives $0$
or $2$. For sufficiency, a nonzero target means every variable
of $C$ is nonzero; two eligible occurrences then sum to $1$.

All these absorptions are derivable. Choose two eligible
occurrences and insert \eqref{f37:eq:retained-pair-consequence}.
The added square identity removes or inserts squared exponents
in doubled words. The support-extension law in
\eqref{f36:eq:extend} then gives the desired doubled monomial
on $C$. Previously inserted absorbed terms can be removed again.

Let $D$ contain all pairs $(C,K)$ with $|I_C|\ge2$.
Let $E$ contain the original occurrences for which $|I_{C_i}|=1$.
Insert every doubled monomial indexed by $D$. Every original
occurrence outside $E$ is removed using $3m\eq2m$, leaving
\[
 \sum_{(C,K)\in D}2m(C,K)+\sum_{(C,K)\in E}m(C,K).
\]
The function determines $D$ by absorption. Its evaluations on
$\{0,2\}$ recover the supports in $E$ as the inclusion-minimal
sets $C$ giving value $2$. On such a support, set one variable
to $3$. The resulting value is $1$ if it is squared and $3$
if it is unsquared, recovering $K$. Thus the canonical expression
is determined by the function, proving completeness for all $X$.
\end{proof}

\subsubsection{Natural numbers truncated at three}
Let $T_3=\{0,1,2,3\}$ with
\begin{equation}\label{f37:eq:truncated-operations}
 a\oplus b=\min(3,a+b),\qquad a\odot b=\min(3,ab).
\end{equation}
The relabelling $(3,1,0,2)$ is an isomorphism
$A_{1882}\cong T_3$. Thus the label $3$ here represents all
natural numbers at least three, while the signature still names
no constants.

Define $\Gamma$ to contain $\SR$ and
\begin{equation}\label{f37:eq:triple-common}
 xy\eq yx,\qquad x^3\eq x^2,\qquad4x\eq3x,
\end{equation}
\begin{equation}\label{f37:eq:triple-closure}
 3x+3xy\eq3x,\qquad3x^2\eq3x,
\end{equation}
\begin{equation}\label{f37:eq:triple-insertion}
 x+y+z+3xyz\eq x+y+z.
\end{equation}
There are 60 more equations. For a subset
$P\subseteq\{0,1,2\}\times\{1,2\}$ containing a pair with
nonzero first coordinate, define
\[
 u_P=\prod_{a>0}z_{ab}^{a},\qquad
 v_P=\prod_{(a,b)\in P}z_{ab}^{b},\qquad
 s_P=\prod_{(a,b)\in P}z_{ab}.
\]
Include
\begin{equation}\label{f37:eq:nested-profile}
 u_P+v_P\eq u_P+s_P.
\end{equation}
There are $63-3=60$ such subsets, so $\Gamma$ has 71 equations,
each using at most six variables.

\begin{theorem}
$\Gamma$ is a complete finite identity basis for $T_3$ and
hence for $A_{1882}$.
\end{theorem}
\begin{proof}
The first eleven equations follow from truncated arithmetic.
For \eqref{f37:eq:nested-profile}, $u_P$ has support contained in
that of $v_P$. If $s_P=0$, both $v_P$ and $s_P$ vanish.
If $s_P=1$, all its variables are $1$, and $v_P=1$.
If $s_P\ge2$, all its variables are positive, $u_P\ge1$ and
$v_P\ge s_P$, so both sides are $3$. This proves soundness.

Normalize a polynomial to words $m(C,K)$, with each positive
coefficient at most three, and count occurrences with multiplicity.
Write $n_C=|I_C|$ as above. The fundamental criterion is
\begin{equation}\label{f37:eq:triple-criterion}
 p+3m(C,K)\eq p\text{ in }T_3
 \quad\Longleftrightarrow\quad n_C\ge3.
\end{equation}
The assignment equal to $1$ on $C$ and $0$ outside it gives
necessity. For sufficiency, whenever the target is positive,
all variables of $C$ are positive and three eligible occurrences
already sum to $3$.

The criterion has an equational derivation from $\Gamma$.
Choose three eligible occurrences $u,v,w$. Equation
\eqref{f37:eq:triple-insertion} inserts $3uvw$ in their additive
context. The cube law reduces exponents, and $3x^2\eq3x$
removes squared exponents inside this tripled word. Extend its
support to $C$ using \eqref{f37:eq:triple-closure}, then introduce
any required squares using the same square identity in reverse.
This inserts $3m(C,K)$, and all auxiliary absorbed terms can
be removed. Moreover, $p+3m\eq p$ implies $p+m\eq p$ by
$4m\eq3m$.

There is a second useful derivation. If one occurrence $u$ has
support contained in the support $C$ of another occurrence $v$,
then
\begin{equation}\label{f37:eq:nested-consequence}
 u+v\eq u+s_C
\end{equation}
follows by grouping variables by their two exponents and applying
\eqref{f37:eq:nested-profile}. Every group substituted into a variable
is nonempty. Distinct occurrences can represent the same word.
This permits the removal of squares from any occurrence having
another eligible occurrence on its support.

Let $D=\{C:n_C\ge3\}$, an upward-closed family of nonempty
subsets of $X$. For $C\notin D$, let $c_C$ be the number of
original occurrences with support exactly $C$. Such a count is
at most two. When $c_C>0$, put $K_C$ equal to the squared set
of the unique eligible occurrence if $n_C=1$, and put
$K_C=\varnothing$ if $n_C=2$. A canonical expression is
\begin{equation}\label{f37:eq:triple-normal}
 N(p)=\sum_{C\in D}3s_C+
       \sum_{\substack{C\notin D\\c_C>0}}c_Cm(C,K_C).
\end{equation}
To derive it, insert $3s_C$ for every $C\in D$ and remove each
original occurrence supported in $D$. This removal is legitimate
even for a squared occurrence, since $3s_C\eq3m(C,K)$ and
$4m\eq3m$. The remaining occurrences have $n_C\le2$.
When $n_C=2$, apply \eqref{f37:eq:nested-consequence} to make each
occurrence on $C$ squarefree. A needed witness has not been
removed: a removed support contained in $C$ would imply $C\in D$.
These operations leave the support counts unchanged and give
\eqref{f37:eq:triple-normal}. The expression is nonempty because
any polynomial has an occurrence, and any three occurrences
give a member of $D$.

Finally the function recovers the canonical data. On the assignment
equal to $1$ on $C$ and $0$ outside, its value is
$r_C=\min(3,n_C)$. Hence $D=\{C:r_C=3\}$.
For $C\notin D$, every nonempty subset of $C$ also lies outside
$D$. Lemma~\ref{lem:finite-set-recovery}(i), over the integers,
therefore recovers the exact coefficients by
\begin{equation}\label{f37:eq:coefficient-recovery}
 c_C=r_C-\sum_{\varnothing\ne B\subsetneq C}c_B.
\end{equation}
When $n_C=r_C=1$ and $c_C>0$, set one variable of $C$ to $2$,
the others to $1$, and all outside variables to $0$. There is
only one eligible occurrence. Its value is $3$ when that variable
is squared and $2$ when it is unsquared. This recovers $K_C$.
All other nonzero coefficients in \eqref{f37:eq:triple-normal} use
squarefree representatives. Equal functions therefore give the
same canonical expression on any common finite variable set,
and every valid identity follows from $\Gamma$.
\end{proof}

\subsubsection{Four mutual term definitions}
On each of $A_{980},A_{983},A_{1614},A_{1619}$ define
\begin{equation}\label{f37:eq:symmetric-product}
 x\circ y=xy+yx.
\end{equation}
Addition is retained. The resulting algebra is isomorphic to
the source in the following table. The map $h$ is from source
labels to the corresponding derived target labels.
\begin{center}
\begin{tabular}{rrll}\toprule
Target&Source&$h$&Original product in the source language\\\midrule
980&1758&$(1,0,3,2)$&$(x\circ x)\circ y$\\
983&1763&$(1,0,3,2)$&$(x\circ x)\circ y$\\
1614&1758&$(0,1,2,3)$&$(x\circ y)\circ y$\\
1619&1763&$(0,1,2,3)$&$(x\circ y)\circ y$\\\bottomrule
\end{tabular}
\end{center}
The displayed operation definitions and recoveries follow from the
tables and specify both term interpretations.

\begin{theorem}
The four targets are finitely based. Each has a complete listed
basis of 265 equations.
\end{theorem}
\begin{proof}
Let $\sigma$ replace a source product by the term $xy+yx$, and
let $\tau$ replace a target product by the appropriate source
term in the table. Both leave addition unchanged. A basis is
\begin{equation}\label{f37:eq:translated-basis}
 \SR\ \cup\ \sigma(\Sigma)\ \cup\
 \{xy\eq\sigma\tau(xy)\},
\end{equation}
where $\Sigma$ is the corresponding 261-equation source basis.
List the two duplicated additive axioms only once, obtaining
$5+(261-2)+1=265$ equations.
The interpretations and recovery identity satisfy
Proposition~\ref{prop:mutual-interpretations}, so the displayed basis
is complete after the indicated isomorphic relabelling.
\end{proof}

\endgroup
\section{Leading letters, prefixes and support transfer}\label{sec:leading-forms}
These families require the position of a letter as well as its
support. We isolate the relevant leading-letter or prefix information
before giving the basis and its completeness proof. Opposite algebras
are treated by the reversal result of Lemma~\ref{lem:opposite}.

\subsection{Leading letters and affine parity}\label{sec:family40}
\begingroup
\providecommand{\eq}{}\renewcommand{\eq}{\approx}
\providecommand{\SR}{}\renewcommand{\SR}{\mathsf{SR}}
\providecommand{\supp}{}\renewcommand{\supp}{\operatorname{supp}}
\providecommand{\Aff}{}\renewcommand{\Aff}{\operatorname{Aff}}
\providecommand{\F}{}\renewcommand{\F}{\mathbb F_2}
\emph{Generators (basis size):} $A_{596}$ (12), $A_{598}$ (13), $A_{603}$ (11), $A_{605}$ (12), $A_{851}$ (11), $A_{1253}$ (11).

For 851 use the carrier $B=\{0,a,b,c\}$ with addition equal
to maximum in $0<c<b<a$ and multiplication
\begin{equation}\label{f40:eq:guardproduct}
 0y=0,\qquad aa=b,\qquad xy=c
 \quad\text{if }x\ne0\text{ and }(x,y)\ne(a,a).
\end{equation}
The map $(0,a,b,c)\mapsto(1,2,3,0)$ identifies $B$ with 851.
Here $0$ is neutral for addition and is a left multiplicative
zero, but $a0=c$. The semiring 1253 is exactly the opposite
of 851 on the original labels.

For the other four generators use $\{z,e,g,d\}$, identified
with catalogue labels $(0,1,2,3)$. Multiplication is commutative,
with
\begin{equation}\label{f40:eq:groupproduct}
 e^2=e,\qquad eg=g,\qquad g^2=e,\qquad
 zx=dx=z\quad(x\in\{z,e,g,d\}).
\end{equation}
Addition has absorbing element $z$, satisfies
$e+e=e$, $g+g=g$, and $e+g=z$, and is specified further by
\begin{equation}\label{f40:eq:groupaddition}
 d+h=d,\qquad d+k=z\ (k\in\{e,g\}\setminus\{h\}),
 \qquad d+d=\begin{cases}z&\varepsilon=0,\\d&\varepsilon=1.\end{cases}
\end{equation}
Write $F_{h,\varepsilon}$ for this algebra. Then $A_{596}=F_{e,0}$,
$A_{598}=F_{g,0}$, $A_{603}=F_{e,1}$ and $A_{605}=F_{g,1}$.
Every nonlinear word has values in the flat core $\{z,e,g\}$.

\subsubsection{An explicit leading-letter basis}
Let $\Sigma_B$ consist of $\SR$ and the following six identities:
\begin{equation}\label{f40:eq:guardbasis}
\begin{aligned}
 xyz&\eq x^3,&x+x&\eq x,\\
 x+xy&\eq x,&x^2+xy&\eq x^2,\\
 x^3+y^2&\eq x^3+y^2+xy,\\
 x^3+yx&\eq x^3+yx+xy.
\end{aligned}
\end{equation}
These identities are valid in $B$. A quadratic word $xy$
has value $0$ when $x=0$, value $b$ when $x=y=a$, and
value $c$ otherwise. A word of length at least three has
value $0$ when its first variable is $0$, and $c$ otherwise.
These descriptions verify the two insertion identities as
well as the remaining laws.

Fix a finite set $X$ of variables. The first identity in
\eqref{f40:eq:guardbasis} reduces every word of length at least
three to the cube of its first letter, by splitting the word
after its first and second letters. Thus the canonical words are
\begin{equation}\label{f40:eq:guardwords}
 \mathcal W_X=X\ \cup\ \{xy:x,y\in X\}\ \cup\ \{x^3:x\in X\}.
\end{equation}
For a polynomial $p$ in these words, let $L(p)$ be its linear
variables and $H(p)$ the set of first letters of its summands.
For $C\subseteq X$, let $W_C(p)$ mean that $p$ contains a
linear or quadratic word with support contained in $C$.
Repetitions may be removed by additive idempotence.

\subsubsection{Absorption and completeness for the leading-letter case}
\begin{theorem}\label{f40:thm:guard}
$\Sigma_B$ is a complete finite identity basis for 851.
For canonical targets, absorption is characterized by
\begin{equation}\label{f40:eq:guardcriterion}
\begin{aligned}
 p+x\eq p&\ \Longleftrightarrow\ x\in L(p),\\
 p+x^3\eq p&\ \Longleftrightarrow\ x\in H(p),\\
 p+xy\eq p&\ \Longleftrightarrow\ x\in H(p)\text{ and }W_{\{x,y\}}(p).
\end{aligned}
\end{equation}
Every indicated absorption is derivable from $\Sigma_B$.
\end{theorem}
\begin{proof}
For the first criterion, put $x=a$ and every other variable
at $0$. A value $a$ can only come from a linear summand $x$.
Its presence is also sufficient. For the second, put $x=c$
and all other variables at $0$. The value is $c$ precisely
when a summand starts with $x$. Every such summand dominates
$x^3$ at every assignment.

For a quadratic target $xy$, the same test with $x=c$ forces
$x\in H(p)$. Put all variables in $\{x,y\}$ at $a$ and
the others at $0$. The target is $b$, and a value at least
$b$ requires the stated low-degree witness. Conversely,
the first-letter witness dominates the target whenever its
value is $0$ or $c$. If the target is $b$, both $x,y$ are
$a$, and the low-degree witness has value $a$ or $b$.
This proves the semantic criteria, including $x=y$.

To derive them, first insert $x^3$ whenever a summand starts
with $x$. For a linear source use $x+xy\eq x$ at $y=x^2$.
For a quadratic source $xv$, factor absorption gives
$xv+xvu\eq xv$, and the word law replaces $xvu$ by $x^3$.
A cubic source already is $x^3$.

For the target $xy$ with $x\ne y$, the possible low-degree
witnesses are $x,y,x^2,y^2,xy,yx$. The witnesses $x,x^2,xy$
are handled by factor absorption, square spreading, and
idempotence. The witness $y^2$ combines with the inserted
$x^3$ using the fifth equation in \eqref{f40:eq:guardbasis}.
A linear $y$ first inserts $y^2$, since $y+y^2\eq y$.
The last equation handles the witness $yx$. When $x=y$,
only $x$ or $x^2$ is needed. All intermediate insertions
are absorbed in the context of $p$ and can be removed after
inserting the target.

Let $D(p)\subseteq\mathcal W_X$ consist of all absorbed
canonical words. Each original summand belongs to $D(p)$,
and every other member can be inserted by these derivations.
Thus $p$ reduces to the sum of all words in $D(p)$ in a fixed
order. Its term function determines $D(p)$ by absorption.
The criterion in Lemma~\ref{lem:normal-form-completeness} proves
completeness.
\end{proof}

\subsubsection{Four bases from affine parity spaces}
For $F_{h,\varepsilon}$ use $\SR$ and
\begin{equation}\label{f40:eq:affinecommon}
 xy\eq yx,\qquad xy\eq x^3y,\qquad x+xy\eq x^3+xy.
\end{equation}
Add the two core identities
\begin{equation}\label{f40:eq:corelaws}
\begin{aligned}
 xy+zu&\eq xy(zu)^2+(xy)^2zu,\\
 xy+zu+vw&\eq xy+zu+vw+xyzuvw.
\end{aligned}
\end{equation}
If $\varepsilon=0$, include
\begin{equation}\label{f40:eq:nonidem}
 xy+xy\eq xy,\qquad x+x\eq x^3;
\end{equation}
if $\varepsilon=1$, include just $x+x\eq x$.
Finally, if $h=g$, include
\begin{equation}\label{f40:eq:anchor}
 x+y^2\eq x^3+y^2.
\end{equation}
Denote this full list by $\Sigma_{h,\varepsilon}$. The sizes
for 596, 598, 603, 605 are respectively 12, 13, 11, 12.
At most six variables occur in any listed equation.

The laws are sound. Products take values in the flat core.
For core elements $r,s$, the sums $r+s$ and $rs^2+r^2s$
agree; if a sum of three core elements is nonzero, they are
equal group elements and their product equals that same element.
This proves \eqref{f40:eq:corelaws}. The incident projection in
\eqref{f40:eq:affinecommon} changes a variable only at $d$; in
that case the incident product is $z$, which absorbs the sum.
For $h=g$, the square in \eqref{f40:eq:anchor} is either $z$ or
$e$. Since $d+e=z$, projection is again valid. The remaining
laws follow directly from the tables.

For nonempty $U\subseteq X$ and $q\in\F^U$, let $m(U,q)$
have exponent two where $q_x=0$ and three where $q_x=1$.
Every nonlinear word reduces to such a word: $xy\eq x^3y$
raises an exponent one in a nonlinear context, and its
instance $y=x$ gives $x^4\eq x^2$ to reduce higher powers.

\begin{theorem}\label{f40:thm:affine}
$\Sigma_{h,\varepsilon}$ is a complete finite identity basis
for each of the four algebras $F_{h,\varepsilon}$.
Every term has a canonical form
\begin{equation}\label{f40:eq:affinenormal}
 N(A,U,Q)=\sum_{x\in A}x+\sum_{q\in Q}m(U,q),
\end{equation}
where $A\cap U=\varnothing$, and either $U=\varnothing$
and $Q=\varnothing$, or $U\ne\varnothing$ and $Q$ is a
nonempty affine subspace of $\F^U$. The expression is nonempty.
For $h=g$, require additionally that $0\in Q$ implies
$A=\varnothing$. Its parameters are determined by its term function.
\end{theorem}
\begin{proof}
We first derive a form of the indicated kind. Expand into
words. In the idempotent cases remove repetitions. In the
other cases, two occurrences of a linear $x$ become $x^3$;
further copies are absorbed by the incident projection and
idempotence of nonlinear words. Remove repeated nonlinear
words in all cases.

Let $U$ be the union of supports of the nonlinear words.
If a remaining linear $x$ lies in $U$, some nonlinear word
has a factorization $xv$. The incident projection replaces
that occurrence of $x$ by $x^3$ in the context of $xv$.
Thus the remaining bare linear variables form a set $A$
disjoint from $U$.

All nonlinear words can be given the same support $U$.
Indeed, the first core identity replaces two nonlinear words
$r,s$ by $rs^2,r^2s$; both new words have the union of their
supports and the same respective parity vectors as before.
Keep occurrences separate during this step. Successively
pairing a fixed occurrence with the others gives it support
$U$; a second pass gives every occurrence that support.
If there is just one occurrence its support already is $U$.
Normalize the powers to obtain words $m(U,q)$.

The second core identity inserts the product of any three
nonlinear summands, factoring each into two nonempty words.
After power reduction its parity vector is their sum in $\F^U$.
By Lemma~\ref{lem:affine-recovery}, iteration inserts exactly the
affine span of the original nonempty set of vectors.
Nonlinear idempotence removes repetitions.

For $h=g$, suppose $0\in Q$ and $A\ne\varnothing$.
The summand $m(U,0)$ is the square of the nonempty product
of the variables in $U$. Equation \eqref{f40:eq:anchor} therefore
projects every bare $x\in A$ to $x^3$. Repeat support
equalization and affine closure with the enlarged nonlinear
support $U\cup A$. The zero parity vector remains present,
and the new bare set is empty. This proves derivability
of a canonical form satisfying the extra condition.

We now recover its parameters from its function. Its content
$C=A\cup U$ is determined by putting one selected variable
at $z$ and all others at $e$. The result is $z$ exactly when
the selected variable belongs to $C$, and is $e$ otherwise.

A variable $x$ belongs to $A$ exactly when some assignment
with $x=d$ and all other variables in $\{e,g\}$ gives value
$d$. If $x\in U$, a nonlinear summand vanishes to $z$;
if $x\notin C$, no summand can produce $d$. Conversely,
for $h=e$ put every other variable at $e$. For $h=g$, put
the other bare variables at $g$. If $Q$ is empty this suffices.
Otherwise $A\ne\varnothing$ implies $0\notin Q$.
Lemma~\ref{lem:affine-recovery} supplies a functional equal to one
on $Q$. Assign its corresponding $e,g$ values to $U$. All nonlinear summands then equal $g$, and the whole
sum is $d$. There is exactly one $d$ summand in this test,
so the argument works for both values of $\varepsilon$.
This recovers $A$, and then $U=C\setminus A$.

Restrict next to group assignments. Encode $e,g$ by $0,1$
and embed all parity vectors in $\F^X$. Let $e_x$ denote
the coordinate vector at $x$, and put
\begin{equation}\label{f40:eq:globalspan}
 H=\Aff\bigl(Q\cup\{e_x:x\in A\}\bigr).
\end{equation}
This affine space is nonempty. If $H=v_0+L$, the polynomial
is nonzero precisely on $L^\perp$, and at such an input $r$
its group value has bit $v_0\cdot r$. Lemma~\ref{lem:affine-recovery}
therefore recovers $H$ by
\begin{equation}\label{f40:eq:dualrecovery}
 H=\{v\in\F^X:v\cdot r=\text{the output bit at }r
               \text{ for every }r\in L^\perp\}.
\end{equation}
Finally,
\begin{equation}\label{f40:eq:intersection}
 Q=H\cap\F^U,
\end{equation}
where $\F^U$ is embedded by zero coordinates outside $U$.
Indeed, an affine combination with zero coordinates on $A$
must have coefficient zero at each independent $e_x$.
The remaining coefficients sum to one and use only $Q$.
If $Q$ is empty, no such affine combination exists, so the
intersection is empty as required.

All parameters of \eqref{f40:eq:affinenormal} are recovered, so
Lemma~\ref{lem:normal-form-completeness} proves completeness.
\end{proof}

The extra projection cannot be applied to the $h=e$ cases.
For example, $x=d,y=e$ gives $x+y^2=d$ and $x^3+y^2=z$.
Thus merely satisfying a smaller collection of identities
would not justify transferring its complete-basis conclusion.

\endgroup
\subsection{Conditional and unconditional tail unions}\label{sec:family42}
\begingroup
\providecommand{\eq}{}\renewcommand{\eq}{\approx}
\providecommand{\SR}{}\renewcommand{\SR}{\mathsf{SR}}
\providecommand{\Var}{}\renewcommand{\Var}{\operatorname{Var}}
\subsubsection{Generators}
Each basis for the following generators
has eleven identities, each involving at most three variables.

\emph{Generators (basis size):} $A_{689}$ (11), $A_{933}$ (11), $A_{1167}$ (11), $A_{1293}$ (11), $A_{1303}$ (11), $A_{1487}$ (11), $A_{1573}$ (11), $A_{1724}$ (11).

Denote 1303 by $A$ and 1167 by $B$. For both, write the
catalogue elements $(0,1,2,3)$ as $(z,k,e,b)$. Addition has
absorbing element $z$, is idempotent, and satisfies
\begin{equation}\label{f42:eq:addition}
 k+b=b,\qquad e+k=e+b=z.
\end{equation}
For $A$, multiplication is
\begin{equation}\label{f42:eq:Aproduct}
 ey=y,\qquad xy=\begin{cases}k&y=k,\\z&y\ne k\end{cases}
 \quad(x\ne e).
\end{equation}
The multiplication of $B$ differs only in $eb=z$ instead of
$eb=b$. Thus a nonlinear word in $B$ is $k$ if its last
letter is $k$, is $e$ if its whole content is assigned $e$,
and is $z$ otherwise.

Denote 1573 by $C$ and 1724 by $D$. In their catalogue
labels, addition is maximum on $2<1<3<0$ for $C$ and on
$2<3<0<1$ for $D$. A nonlinear word in $C$ has value $0$
if its last letter is $0$ or $3$, value $2$ if its whole
content is $2$, and value $1$ otherwise. A nonlinear word
in $D$ has value $1$ if any variable in its content is $1$;
otherwise its value is $0$ if its last letter is $0$ or $3$,
and $2$ otherwise. These statements follow directly from
the multiplication tables and induction on word length.

\subsubsection{Two bases and common reductions}
Let $\Gamma$ consist of $\SR$ and
\begin{equation}\label{f42:eq:common}
\begin{aligned}
 x^2y&\eq xy,&xyz&\eq yxz,\\
 x+x&\eq x,&x+yx&\eq yx,\\
 xy+zy&\eq xzy.
\end{aligned}
\end{equation}
Define $\Sigma=\Gamma\cup\{\text{\eqref{f42:eq:conditional}}\}$
and $\Delta=\Gamma\cup\{\text{\eqref{f42:eq:unconditional}}\}$, where
\begin{equation}\label{f42:eq:conditional}
 xy+x+z\eq zxy+x+z
\end{equation}
and
\begin{equation}\label{f42:eq:unconditional}
 xy+z\eq zxy+z.
\end{equation}
Both lists contain eleven identities. The word evaluations
above verify that $A\models\Sigma$ and $B,C,D\models\Delta$. For example,
in $A$ the identity $xy+zy\eq xzy$ follows by testing the
last factor: if it is $k$, both sides are $k$; if it is
$e$ or $b$, both sides retain that value exactly when
both preceding factors are assigned $e$, and have the
absorbing value otherwise. An absorbing last factor
gives that absorbing value on both sides.

For clarity, soundness of \eqref{f42:eq:conditional} in $A$
can also be checked by its first variable. If $x=z$ as an
element, its summand absorbs the sum. If $x=k$ or $b$,
the sum can avoid the absorbing element only when $y=k$,
in which case both displayed products end in $k$. If $x=e$,
a nonabsorbing sum requires $y=e$ and the remaining summand
to be $e$ as well. Both sides are then $e$.
Adding the summand $x$ to \eqref{f42:eq:unconditional} derives
\eqref{f42:eq:conditional}. The converse is false in $A$:
at $(x,y,z)=(e,b,k)$ the two sides of
\eqref{f42:eq:unconditional} are $b$ and the absorbing element $z$.

Fix a finite nonempty variable set $X$. A standard word is
\begin{equation}\label{f42:eq:word}
 w(P,y)=\Bigl(\prod_{x\in P}x\Bigr)y,\qquad P\subseteq X,
 \quad y\in X.
\end{equation}
The prefix uses each variable once in a fixed order. The
convention $w(\varnothing,y)=y$ denotes a variable, not
an empty product. The last letter may belong to $P$.
The first two laws of \eqref{f42:eq:common} permute prefix
letters and delete repeated prefix letters, reducing every
nonempty word to this form.

\begin{lemma}\label{f42:lem:merge}
Modulo $\Gamma$, summands with a common last variable merge by
\begin{equation}\label{f42:eq:merge}
 w(P,y)+w(Q,y)\eq w(P\cup Q,y).
\end{equation}
Consequently every term reduces to a sum with exactly one
standard word for each of its last variables.
\end{lemma}
\begin{proof}
When both prefixes are nonempty, substitute their products
for $x,z$ in $xy+zy\eq xzy$ and reduce the new prefix.
When one is empty use $x+yx\eq yx$; when both are empty
use additive idempotence. Expand a term by distributivity
and merge repeatedly.
\end{proof}

\begin{lemma}[Prefix saturation]\label{lem:prefix-saturation}
Assume the word laws $x^2y\eq xy$, $xyz\eq yxz$,
the merging law $xy+zy\eq xzy$, and additive idempotence of
nonlinear words. Let $r=w(P,y)$ be a nonlinear summand of a
polynomial, and let $\mathcal R$ be a collection of its summands
containing $r$. If, in a fixed additive context which is retained,
the available identities permit
$r+s\eq sr+s$ for each $s\in\mathcal R$, also after enlarging
the prefix of $r$, then $r$ may be replaced by
\[
 w\left(P\cup\bigcup_{s\in\mathcal R}\supp(s),y\right)
\]
while all other summands remain. In particular,
\eqref{f42:eq:unconditional} allows $\mathcal R$ to contain
all summands; the law $xy+zu\eq zuxy+zu$ allows all nonlinear summands.
\end{lemma}
\begin{proof}
Retain a copy of $r$ by nonlinear idempotence, replace the other
copy in $r+s$ by $sr$, and merge the two words ending in $y$.
Their prefixes unite, adjoining $\supp(s)$. For $s=r$, use
three copies of the nonlinear $r$; after merging, its last
variable also belongs to the prefix. Repeat over $\mathcal R$,
retaining the original contexts until they have been used.
Only nonlinear words are duplicated. The word laws remove prefix
repetitions, and no variable outside the indicated union is introduced.
\end{proof}

\subsubsection{Conditional spreading}
For a merged sum, write $T$ for its set of last variables
and $C_0$ for its whole content. If its prefix $P_y$ meets
$T$, call that last variable marked. Put $J=\{y:P_y\cap T
\ne\varnothing\}$. For unmarked last variables write
$U_y=P_y\subseteq C_0\setminus T$.

\begin{theorem}\label{f42:thm:conditional}
$\Sigma$ is a complete finite identity basis for $A$.
Every term has a normal form
\begin{equation}\label{f42:eq:conditionalnormal}
 N(C_0,T,J,U)=\sum_{y\in T\setminus J}w(U_y,y)
                   +\sum_{y\in J}w(C_0,y),
\end{equation}
where $\varnothing\ne T\subseteq C_0\subseteq X$,
$J\subseteq T$, and $U_y\subseteq C_0\setminus T$.
If $J=\varnothing$, require
$C_0=T\cup\bigcup_{y\in T}U_y$. The parameters are
uniquely determined by the term function.
\end{theorem}
\begin{proof}
Merge summands using Lemma~\ref{f42:lem:merge}. If $P_y$ contains
$x\in T$, a word ending in $x$ permits insertion of the auxiliary
linear $x$ by $x+vx\eq vx$; a linear witness already supplies it.
Permute the target prefix to write the word as $xv$. In this context
\eqref{f42:eq:conditional} permits $xv+r\eq rxv+r$ for every
summand $r$, and still does so after its prefix grows.
Lemma~\ref{lem:prefix-saturation} therefore expands each marked
prefix to $C_0$. Remove the auxiliary $x$ with its original witness.
Content and last variables are unchanged; unmarked prefixes remain
disjoint from $T$. This gives \eqref{f42:eq:conditionalnormal},
including the stated content condition when $J=\varnothing$.

We recover the parameters by actual assignments in $A$.
Assign a selected variable $v$ the element $z$ and every
other variable $e$. The output is $z$ exactly when
$v\in C_0$, and is $e$ otherwise. Next assign $v=z$ and
all others $k$. The output is $z$ exactly when $v\in T$,
and is $k$ otherwise. These tests recover $C_0$ and $T$.

For $y\in T$, set $y=b$, set other members of $T$ to $k$,
and set variables outside $T$ to $e$. Every other summand
ends in $k$, hence equals $k$. The summand ending in $y$
is $z$ when marked, since its prefix contains a last
variable not assigned $e$; otherwise it is $b$. The total
output is therefore $z$ exactly when $y\in J$, and is
$b$ otherwise. For unmarked $y$, change a selected
$v\notin T$ from $e$ to $z$ in this test. Its output
changes from $b$ to $z$ precisely when $v\in U_y$.
All parameters have been recovered from the term function.

Completeness follows from Lemma~\ref{lem:normal-form-completeness}.
\end{proof}

\subsubsection{Unconditional spreading and a common variety}
\begin{theorem}\label{f42:thm:unconditional}
$\Delta$ is a complete finite identity basis for each of
$B,C,D$. Every term has a normal form
\begin{equation}\label{f42:eq:unconditionalnormal}
 M(C_0,T,L)=\sum_{y\in L}y+\sum_{y\in T\setminus L}w(C_0,y),
\end{equation}
where $\varnothing\ne T\subseteq C_0\subseteq X$ and
$L\subseteq T$. If $L=T$, require $C_0=T$.
These parameters are determined by the term function in
each generator. In particular, $\Var(B)=\Var(C)=\Var(D)$.
\end{theorem}
\begin{proof}
After Lemma~\ref{f42:lem:merge}, apply
Lemma~\ref{lem:prefix-saturation} with all summands as contexts.
Every nonlinear prefix becomes $C_0$; the remaining linears form
$L$. If every summand is linear then $C_0=T=L$, proving
\eqref{f42:eq:unconditionalnormal}.

For $B$ use the same $C_0$ and $T$ tests as in the preceding
proof. For $y\in T$, assign $y=b$, the other members of
$T$ to $k$, and variables outside $T$ to $e$. Nonlinear
words ending in $y$ are $z$, whereas a linear $y$ is $b$.
All other summands end in $k$. Thus the output is $b$
exactly when $y\in L$, and is $z$ otherwise.

For $C$ and $D$ use catalogue labels. Assign a selected
variable $v$ the value $1$ and all others $2$; the output
is $1$ exactly when $v\in C_0$, and is $2$ otherwise.
Assign $v=0$ and all others $2$; the output is $0$ exactly
when $v\in T$. In $C$, a variable occurring only in a
prefix can instead produce $1$, which is distinct from
$0$; in $D$ it leaves the value $2$. Finally, for
$y\in T$, assign $y=3$ and all others $2$. The output
is $3$ when $y\in L$ and $0$ when $y\notin L$.
Indeed, a nonlinear word ending in $3$ is $0$ in both
algebras. If $y$ is linear, other terms are at most $1$
in the chain for $C$, and are $2$ in $D$, so their sum
with that linear $3$ remains $3$.

The parameters $C_0,T,L$ are recovered, so
Lemma~\ref{lem:normal-form-completeness} proves completeness.
The common basis $\Delta$ gives equality of the three varieties.
\end{proof}

\endgroup
\subsection{Doubling and square profiles}\label{sec:family43}
\begingroup
\emph{Generators (basis size):} $A_{672}$ (11), $A_{918}$ (11),
$A_{1150}$ (11), $A_{1280}$ (11), $A_{1558}$ (11), $A_{1711}$ (11).

The direct generators $A_{1150},A_{1558},A_{1711}$ have the
same respective multiplications as $A_{1167},A_{1573},A_{1724}$
in Section~\ref{sec:family42}. In each pair the sole change to
addition is $3+3=0$ in place of $3+3=3$. Thus their nonlinear
word evaluations are unchanged and never take the value $3$.
Let $\Theta$ consist of $\SR$ and
\begin{equation}\label{f43:eq:basis}
\begin{aligned}
 x^2y&\eq xy,&xyz&\eq yxz,\\
 x+x&\eq x^2,&x+yx&\eq yx,\\
 xy+zy&\eq xzy,&xy+z&\eq zxy+z.
\end{aligned}
\end{equation}
These eleven identities involve at most three variables.
Except for doubling, they hold by the same word evaluations:
every sum in an added identity has at most one linear summand,
so the changed additive pair is never used. For doubling, both
$2a$ and $a^2$ equal $0$ at $a=3$ and equal $a$ otherwise.

\begin{lemma}\label{f43:lem:reductions}
The following identities follow from $\Theta$:
\begin{equation}\label{f43:eq:derived}
 xy+xy\eq xy,\qquad x+x+x\eq x+x,\qquad x+x^2\eq x^2.
\end{equation}
Words reduce to the standard words \eqref{f42:eq:word}, and
summands with a common last variable merge by
\begin{equation}\label{f43:eq:merge}
\begin{aligned}
 w(P,y)+w(Q,y)&\eq w(P\cup Q,y)&&\text{if }P\cup Q\ne\varnothing,\\
 y+y&\eq w(\{y\},y).&&
\end{aligned}
\end{equation}
\end{lemma}
\begin{proof}
Set $z=x$ in $xy+zy\eq xzy$ and use $x^2y\eq xy$ to obtain
nonlinear idempotence. Substituting $y=x$ in $x+yx\eq yx$
gives square absorption, and
$3x\eq x+x^2\eq x^2\eq2x$ gives the threshold law.
The word reduction and the first merging rule follow as in
Lemma~\ref{f42:lem:merge}; two empty prefixes instead use doubling.
\end{proof}

\begin{theorem}\label{f43:thm:complete}
$\Theta$ is a complete basis for each direct generator, with normal forms
\begin{equation}\label{f43:eq:normal}
 N(U,T,L)=\sum_{y\in L}y+\sum_{y\in T\setminus L}w(U,y),
\end{equation}
where
\begin{equation}\label{f43:eq:parameters}
 \varnothing\ne T\subseteq U\subseteq X,\qquad L\subseteq T,
 \qquad L=T\ \Longrightarrow\ U=T.
\end{equation}
Each term function determines its triple uniquely. Consequently
$\Var(A_{1150})=\Var(A_{1558})=\Var(A_{1711})$.
\end{theorem}
\begin{proof}
Expand with multiplicities and apply Lemma~\ref{f43:lem:reductions}
until one summand remains for each last variable. Repeated linear
copies become a square. Lemma~\ref{lem:prefix-saturation}, using all
summands as contexts, expands each nonlinear prefix to the full content
$U$; linear contexts are retained once. The surviving linears form
$L$, giving the stated restrictions, including $U=T$ when all
summands are linear.

The recovery tests for $B,C,D$ in Theorem~\ref{f42:thm:unconditional}
apply to the corresponding direct generators here. Each test assigns
at most one variable the catalogue value $3$, and that variable occurs
at most once as a linear summand in \eqref{f43:eq:normal}.
All nonlinear values differ from $3$, and adding two values different
from $3$ never gives $3$ in these algebras. Thus every tested sum can
be evaluated without using $3+3$: add the possible single $3$ last.
Its result is unchanged by the modified additive square.
Those tests therefore recover $U,T,L$ here as well.
Lemma~\ref{lem:normal-form-completeness} proves completeness and
the common basis gives the variety equalities.
\end{proof}
\endgroup

\subsection{Prefix erasure in sums}\label{sec:family44}
\begingroup
\providecommand{\eq}{}\renewcommand{\eq}{\approx}
\providecommand{\SR}{}\renewcommand{\SR}{\mathsf{SR}}
\providecommand{\Var}{}\renewcommand{\Var}{\operatorname{Var}}
\providecommand{\supp}{}\renewcommand{\supp}{\operatorname{supp}}
\subsubsection{Signature, background and tables}
Write $R=A_{1149}$ and $S=A_{1557}$. Their opposites are identified in
Table~\ref{tab:opposite-transfers}. We determine the additive contexts
in which a product's prefix can be erased.

\emph{Generators (basis size):} $A_{671}$ (10), $A_{917}$ (10), $A_{1149}$ (10), $A_{1557}$ (10).

For $R$, write $(0,1,2,3)=(z,k,e,b)$. Addition has
$k+k=k$ and $k+b=b+k=b$, with every other sum equal to $z$.
A nonlinear word has value $k$ if its last variable is $k$,
value $e$ if every variable is $e$, and value $z$ otherwise.
For $S$, addition agrees with maximum on $2<1<3<0$, except
that $2+2=1$ and $3+3=0$. A nonlinear word is $0$ if its last
variable is $0$ or $3$, is $2$ if every variable is $2$, and
is $1$ otherwise. These word descriptions follow by induction
from the tables. In particular, no nonlinear word takes value $3$.

\subsubsection{Erasure and a complete finite basis}
Let $\Xi$ consist of $\SR$ and the following five identities:
\begin{equation}\label{f44:eq:basis}
\begin{aligned}
 x^2y&\eq xy,&xyz&\eq yxz,\\
 xyx&\eq yx,&x+x+x&\eq x+x,\\
 xy+z&\eq y+y+z.
\end{aligned}
\end{equation}
The last identity erases the prefix only in the presence of
an additive context. The preceding word evaluations and the
addition tables verify all ten identities in $R,S$.
Each identity uses at most three variables.

Fix a nonempty finite set $X=\{x_1,\ldots,x_n\}$ with a linear
order. For a nonempty subset $U$ and $y\in U$, define
\begin{equation}\label{f44:eq:word}
 W(U,y)=\Bigl(\prod_{x\in U}x\Bigr)y,\qquad
 \varnothing\ne U\subseteq X,\quad y\in U.
\end{equation}
The prefix uses the fixed order. It contains $y$, so this word
has length at least two; no empty product occurs. For a vector
of coefficients define
\begin{equation}\label{f44:eq:additive}
 A(\alpha)=\sum_{i=1}^n\alpha_i x_i,\qquad
 \alpha\in\{0,1,2\}^n\setminus\{\boldsymbol0\}.
\end{equation}
A coefficient two means two additive copies, and a zero
coefficient means omission. In particular $A(e_i)=x_i$.
There is no nullary zero term.

\begin{lemma}\label{f44:lem:reduction}
Every term is provably equal modulo $\Xi$ to a form
\eqref{f44:eq:word} or \eqref{f44:eq:additive}.
\end{lemma}
\begin{proof}
The first two word laws in \eqref{f44:eq:basis} permute adjacent
prefix variables and remove repetitions from the prefix of
any nonempty word. If a word $vy$ is nonlinear, $v$ is nonempty.
The third word law gives $yvy\eq vy$, inserting the last
variable in the prefix when needed. Sorting and deleting
repetitions now yields $W(U,y)$, where $U$ is the whole support.

Expand an arbitrary term by distributivity into a nonempty
sum of nonempty words, retaining their multiplicities. A single
word is either a variable or the preceding nonlinear form.
If at least two word occurrences are present, select a
nonlinear occurrence $vy$ and let $t$ be the nonempty sum of
the remaining occurrences. Substitute $x=v,z=t$ into the
erasure identity to replace $vy+t$ by $y+y+t$. This step
eliminates one nonlinear occurrence and leaves at least two
linear or nonlinear occurrences. Repetition eliminates every
nonlinear occurrence. The additive threshold identity then
reduces all positive coefficients to one or two. No step
duplicates a nonlinear word or assumes additive idempotence.
\end{proof}

\begin{theorem}\label{f44:thm:complete}
In each of $R,S$, distinct forms in \eqref{f44:eq:word} and
\eqref{f44:eq:additive} have distinct term functions, including
forms of different types. Hence $\Xi$ is a complete identity
basis for both algebras, and $\Var(R)=\Var(S)$.
\end{theorem}
\begin{proof}
We first separate additive forms in either algebra. Give
a selected variable $x_i$ value $0$ and every other variable
value $1$. The output is $0$ precisely when $\alpha_i>0$,
and is $1$ otherwise: $0$ is additively absorbing, whereas
$1+1=1$. Thus the function recovers the support of $\alpha$.
For a supported variable, instead give it value $3$ and
all others value $1$. The output is $3$ for coefficient
one and $0$ for coefficient two, since $3+1=3$ and
$3+3=0$. This recovers every coefficient, also when only
one variable occurs in the expression.

For a nonlinear form in $R$, set a selected variable $v$
to $0$ and all others to $2$. The output is $0$ precisely
when $v\in U$, and is $2$ otherwise. For a nonlinear form
in $S$, set $v=1$ and all others to $2$; the output is
$1$ precisely when $v\in U$, and is $2$ otherwise. In
either algebra, set $v=0$ and all others to $1$. The
output is $0$ precisely when $v=y$, and is $1$ otherwise.
The nonlinear function therefore determines both $U$ and $y$.

Finally, put every variable at $2$. A nonlinear form gives
$2$. An additive form with at least two occurrences gives
$0$ in $R$ or $1$ in $S$, including arbitrarily many
occurrences by associativity and the threshold identity.
The only remaining additive forms are individual variables.
Putting every variable at $3$ gives $3$ for an individual
variable and $0$ for every nonlinear word. This separates
the two types of forms.

Together with Lemma~\ref{f44:lem:reduction}, this separation proves
completeness by Lemma~\ref{lem:normal-form-completeness}. The common
basis $\Xi$ gives equality of the generated varieties.
\end{proof}

The erasure law fails in $A_{1299}$ at $(x,y,z)=(2,3,1)$:
$xy+z$ and $y+y+z$ have values $3$ and $0$. It fails in
$A_{1300}$ at $(2,2,1)$ with the same values. Their conditional
erasure laws are treated in Theorem~\ref{f45:thm:complete}.

For the paired opposites, reversal gives $yx+z\eq y+y+z$,
equivalently $xy+z\eq x+x+z$.

\endgroup
\subsection{Tail conflict erasure}\label{sec:family45}
\begingroup
\providecommand{\eq}{}\renewcommand{\eq}{\approx}
\providecommand{\SR}{}\renewcommand{\SR}{\mathsf{SR}}
\providecommand{\Var}{}\renewcommand{\Var}{\operatorname{Var}}
\subsubsection{Signature and the two direct generators}
Write $R=A_{1299}$ and $S=A_{1300}$; see
Table~\ref{tab:opposite-transfers} for their opposites. The normal forms
retain a prefix precisely when the erasure conditions below fail.

\emph{Generators (basis size):} $A_{1299}$ (14), $A_{1300}$ (12), $A_{1481}$ (14), $A_{1482}$ (12).

The direct generators have the same multiplication. A word
with last variable $y$ evaluates to $1$ whenever $y=1$.
Otherwise, it evaluates to the value of $y$ if all its
prefix variables have value $2$, and to $0$ if that condition
fails. This includes a linear word by taking an empty prefix
condition. The description follows by induction from the tables.

In $R$, the only nonzero binary sums are $1+1=1$ and
$1+3=3+1=3$. In $S$ they are $1+1=1$ and
$1+2=2+1=1+3=3+1=3$. In both algebras every sum of at
least two elements is nonzero precisely when either all
summands are $1$, or exactly one is exceptional and all
others are $1$. The exceptional value must be $3$ in $R$
and can be $2$ or $3$ in $S$. The outputs in these cases
are $1$ and $3$, respectively. All other sums are $0$.

\subsubsection{Two finite lists and their reductions}
Let $\Gamma$ consist of $\SR$ and seven identities:
\begin{equation}\label{f45:eq:common}
\begin{aligned}
 x^2y&\eq xy,&xyz&\eq yxz,\\
 x+x+x&\eq x+x,&x+yx&\eq x+x,\\
 xy+zy&\eq y+y,&xy+x&\eq y+y+x,\\
 xy+zx&\eq y+y+zx.
\end{aligned}
\end{equation}
Thus $\Gamma$ has twelve identities. Let $\Delta$ be
$\Gamma$ together with
\begin{equation}\label{f45:eq:self}
 x^2+y\eq x+x+y,\qquad xy^2+z\eq y+y+z.
\end{equation}
All fourteen identities hold in $R$, and all twelve in
$\Gamma$ hold in $S$. The preceding word evaluations
verify these assertions. Each identity has at most three variables.

Fix a nonempty finite ordered variable set $X$. The word laws
reduce every word to
\begin{equation}\label{f45:eq:word}
 w(P,y)=\Bigl(\prod_{x\in P}x\Bigr)y,\qquad
 P\subseteq X,\quad y\in X.
\end{equation}
The prefix uses each of its variables once in the fixed order.
The last variable may also belong to $P$. When $P$ is empty,
$w(P,y)$ means $y$, not an empty product.

Besides single words, consider the expressions
\begin{equation}\label{f45:eq:sum}
 N(T,D,(P_y))=\sum_{y\in D}(y+y)
             +\sum_{y\in T\setminus D}w(P_y,y),
\end{equation}
where
\begin{equation}\label{f45:eq:parameters}
 \varnothing\ne T\subseteq X,\quad D\subseteq T,\qquad
 |T|\ge2\ \text{or}\ (|T|=1\ \text{and}\ D=T).
\end{equation}
For $\Delta$ require $P_y\subseteq X\setminus T$.
For $\Gamma$ require $P_y\subseteq (X\setminus T)\cup\{y\}$.
The restrictions apply only for $y\notin D$. Empty sub-sums
are omitted. Every expression in \eqref{f45:eq:sum} has at least
two word occurrences. No global content parameter is required:
variables whose prefixes are erased may disappear entirely.

\begin{lemma}\label{f45:lem:reduction}
Modulo $\Delta$, respectively $\Gamma$, every term derives a
single word \eqref{f45:eq:word} or an expression \eqref{f45:eq:sum}
with the corresponding prefix restrictions.
\end{lemma}
\begin{proof}
The first two identities in \eqref{f45:eq:common} permute adjacent
prefix variables and delete their repetitions. They prove the
single-word reduction without changing its last variable.
Expand a general term by distributivity, keeping multiplicities.
If just one word occurrence results, this completes reduction.

Otherwise, group occurrences by their last variable. Two
linear occurrences already give $y+y$. A linear and a
nonlinear occurrence give $y+y$ by $x+yx\eq x+x$, with
the prefix substituted for the second variable. Two nonlinear
occurrences give $y+y$ by $xy+zy\eq y+y$, substituting
their nonempty prefixes. Further occurrences are absorbed:
$(y+y)+vy\eq y+(y+vy)\eq y+(y+y)\eq y+y$, and the
linear case uses the threshold identity. Each group therefore
is either $y+y$ or a single $w(P_y,y)$. Let $T$ be the
last-variable set and $D$ the doubled groups.

Suppose that $x\ne y$ lies in $P_y\cap T$. Reorder the
prefix so the target word is $xq$, where $q$ is a nonempty
word ending in $y$. Select a witness occurrence ending in
$x$. If it is the linear word $x$, use $xq+x\eq q+q+x$.
If it is $vx$ with $v$ nonempty, use
$xq+vx\eq q+q+vx$. These are substitutions into the last
two identities of \eqref{f45:eq:common}. The two copies of $q$
then reduce to $y+y$ by the preceding group rule. Thus the
target becomes doubled and loses its prefix. The witness
is retained; if it came from a doubled group, its other
copy stays in the surrounding sum. In particular $T$ stays
unchanged throughout. Repeating eliminates all occurrences
of other last variables in surviving prefixes.

For $\Gamma$ this is the required restriction. For $\Delta$,
if $y\in P_y$, the target can be ordered as $y^2$, or
$vy^2$ with $v$ nonempty. The rest of the sum is nonempty.
Apply the respective identity of \eqref{f45:eq:self} to replace
this target by $y+y$. This removes the remaining self
occurrences from surviving prefixes. Each erasure decreases
the number of groups represented by a single word, so the
process terminates. Grouping and erasure keep at least two
occurrences, giving exactly condition \eqref{f45:eq:parameters}.
No nonlinear term is assumed additively idempotent.
\end{proof}

\subsubsection{Separation, completeness and a strict inclusion}
\begin{theorem}\label{f45:thm:complete}
$\Delta$ is a complete identity basis for $R$, and $\Gamma$
is a complete identity basis for $S$. In each algebra the
term function uniquely determines the associated normal form,
for every finite variable set.
\end{theorem}
\begin{proof}
We give all assignments in the direct catalogue labels.
When every variable is $2$, a single word has value $2$,
whereas every sum form has value $0$, since it has at least
two occurrences, each valued $2$. This separates the types.

For a single word, set a selected variable $v$ to $0$ and
all others to $1$. The result is $0$ precisely when $v=y$,
and is $1$ otherwise. Having recovered the last variable,
recover each other variable's prefix membership by setting
that variable to $0$ and all others to $2$: output $0$
means membership, output $2$ means absence. Finally set
$y=3$ and all other variables to $2$. The result is $0$
if $y\in P$, and $3$ otherwise. Hence $(P,y)$ is recovered,
including an empty prefix.

For a sum form in either algebra, give a selected variable
$v$ value $0$ and every other variable value $1$. The
output is $0$ precisely when $v\in T$, and otherwise is
$1$, since a word ending in $1$ always evaluates to $1$.
This recovers $T$ even when some prefixes have lost variables.

In $R$, for $y\in T$ set $y=3$, all other members of $T$
to $1$, and all variables outside $T$ to $2$. If $y\in D$
its two occurrences give $0$. Otherwise its prefix lies
outside $T$, its word gives $3$, and all other occurrences
give $1$. The output is then $3$. Thus this test recovers
$D$. For $y\notin D$, change an outside variable $v$ from
$2$ to $0$. The output changes from $3$ to $0$ exactly
when $v\in P_y$. Other words still end in $1$. This recovers
every surviving prefix.

In $S$, instead set $y=2$, other members of $T$ to $1$,
and outside variables to $2$. A doubled group gives $0$;
a surviving word gives $2$, and, since $|T|\ge2$ in a
surviving case, it is added to at least one $1$. The
output is $3$. Hence this detects $D$. For a surviving
group, change an outside variable to $0$ to recover its
prefix membership just as above. To detect the remaining
possibility $y\in P_y$, change $y$ from $2$ to $3$.
The output is $0$ when the prefix contains $y$, and $3$
otherwise. In the singleton-tail case the parameters
already force the doubled form $y+y$.

All normal-form parameters are uniquely determined. Together with
Lemma~\ref{f45:lem:reduction}, this proves completeness by
Lemma~\ref{lem:normal-form-completeness}.
\end{proof}

\begin{corollary}\label{f45:cor:inclusion}
The generated varieties satisfy
\begin{equation}\label{f45:eq:inclusion}
 \Var(R)\subsetneq\Var(S).
\end{equation}
\end{corollary}
\begin{proof}
The complete list $\Delta$ contains $\Gamma$, so its
models form a subvariety of the models of $\Gamma$.
In $S$, set $(x,y)=(2,1)$ in $x^2+y\eq x+x+y$.
The sides have values $3$ and $0$, respectively. Thus
$S$ is outside the smaller variety, proving strictness.
\end{proof}

\endgroup
\subsection{Incident linear conversion}\label{sec:family46}
\begingroup
\providecommand{\eq}{}\renewcommand{\eq}{\approx}
\providecommand{\SR}{}\renewcommand{\SR}{\mathsf{SR}}
\providecommand{\Var}{}\renewcommand{\Var}{\operatorname{Var}}
\subsubsection{Generators}
The ten algebras listed below have nonidempotent
addition, since $3+3=0\ne3$. The reductions convert a linear summand
to a square when its variable occurs in a nonlinear support, and
enlarge supports within the nonlinear part.

\emph{Generators (basis size):} $A_{678}$ (13), $A_{679}$ (13), $A_{775}$ (13), $A_{778}$ (13), $A_{803}$ (13), $A_{806}$ (13), $A_{916}$ (13), $A_{1156}$ (13), $A_{1157}$ (13), $A_{1556}$ (13).

Call 775, 778, 1156, 1157 and 1556 the direct generators.
Here a linear summand is a single variable, and a nonlinear
word means a word of length at least two.
For 775 and 778, a nonlinear word has value equal to its
last variable when every variable in its support belongs
to $\{1,2\}$, and value $0$ otherwise. For 1156 and 1157,
it has value $1$ if its last variable is $1$, value $2$
if all its variables are $2$, and value $0$ otherwise.
For 1556, it has value $0$ if the last variable is $0$ or
$3$, value $2$ if every variable is $2$, and value $1$
otherwise. These descriptions follow by induction from
the multiplication tables. A nonlinear word never has value $3$.

\subsubsection{Conversion and normalization}
Let $\Sigma$ consist of $\SR$ and eight identities:
\begin{equation}\label{f46:eq:basis}
\begin{aligned}
 x^2y&\eq xy,&xyz&\eq yxz,\\
 xyx&\eq yx,&x+x&\eq x^2,\\
 x+yx&\eq yx,&xy+zy&\eq xzy,\\
 xy+x&\eq xy+x^2,&xy+zu&\eq zuxy+zu.
\end{aligned}
\end{equation}
The preceding word evaluations and the addition tables in
Appendix~\ref{app:catalogue} verify all thirteen identities in each direct generator. The last identity uses four variables; all others use
at most three. We make no assertion that four variables or
thirteen identities are necessary.

\begin{lemma}\label{f46:lem:reduction}
The basis $\Sigma$ implies
\begin{equation}\label{f46:eq:derived}
 xy+xy\eq xy,\qquad x+x^2\eq x^2,\qquad x+x+x\eq x+x.
\end{equation}
Fix a finite nonempty ordered variable set $X$. Every term
derives a normal form
\begin{equation}\label{f46:eq:normal}
 N(U,T,L)=\sum_{x\in L}x+\sum_{y\in T}W(U,y),\qquad
 W(U,y)=\Bigl(\prod_{x\in U}x\Bigr)y,
\end{equation}
with parameters
\begin{equation}\label{f46:eq:parameters}
 T\subseteq U\subseteq X,\quad L\subseteq X\setminus U,
 \quad(T=\varnothing\ \Longleftrightarrow\ U=\varnothing),
 \quad U\cup L\ne\varnothing.
\end{equation}
Every displayed $W(U,y)$ is nonlinear; its prefix contains
each variable of $U$ once in the fixed order, including $y$.

\end{lemma}
\begin{proof}
Set $z=x$ in $xy+zy\eq xzy$ and use $x^2y\eq xy$
to derive idempotence of nonlinear words. Set $y=x$ in
$x+yx\eq yx$ to derive $x+x^2\eq x^2$. Together with
doubling into a square this gives the threshold identity.

The first two word laws permute adjacent prefix variables
and delete repetitions. In a nonlinear word $vy$, the
nonempty prefix $v$ allows $yvy\eq vy$ to be used backwards,
inserting $y$ in the prefix. Sorting and removing repetitions
therefore reduce every nonlinear word to $W(V,y)$ with
$V$ its support.

Expand a general term by distributivity, preserving
multiplicities. Repeated linear copies become a square,
and further copies are absorbed. Nonlinear occurrences
are idempotent. Terms ending in the same variable merge:
nonlinear prefixes unite by $xy+zy\eq xzy$, while a
linear term with that last variable is absorbed by
$x+yx\eq yx$.

Suppose a remaining linear $x$ occurs in the support of
a nonlinear word. If it is that word's last variable,
it is absorbed by the preceding rule. Otherwise put $x$
first in its prefix, writing the word as $xq$ with $q$
nonempty, and apply $xq+x\eq xq+x^2$. Thus the linear
summand becomes a square. Repeat until all remaining
linear variables are outside every nonlinear support.
Each step removes a remaining linear summand and introduces
no new variable into those supports, so the process terminates.

Let $U$ be the union of nonlinear supports and $T$ their last-variable
set. The last law of \eqref{f46:eq:basis} permits
Lemma~\ref{lem:prefix-saturation} with precisely the nonlinear summands
as contexts. It expands every nonlinear prefix to $U$ without
changing its last variable. The surviving linears form
$L\subseteq X\setminus U$. If no nonlinear summand remains,
$U=T=\varnothing$ and $L\ne\varnothing$, giving the stated restrictions.
\end{proof}

\subsubsection{Complete bases}
\begin{theorem}\label{f46:thm:complete}
$\Sigma$ is a complete identity basis for each of 775,
778, 1156, 1157 and 1556. Each term function uniquely
determines its triple in \eqref{f46:eq:parameters}; hence these
five generators generate the same variety $\mathcal V$.
\end{theorem}
\begin{proof}
We recover the parameters using the original table labels.
First assign a selected variable $v$ value $3$ and every
other variable the background value $b$ in the following table.
The output columns describe all possible memberships:
\begin{center}
\begin{tabular}{rccccc}\toprule
Generator&$b$&$v\notin U\cup L$&$v\in L$&$v\in T$&$v\in U\setminus T$\\\midrule
775, 778&1&1&3&0&0\\
1156, 1157&2&2&3&0&0\\
1556&2&2&3&0&1\\\bottomrule
\end{tabular}
\end{center}
The preceding word descriptions and the addition tables
verify every column. For instance, in 775 or 778, a
variable $3$ in $U$ makes every nonlinear summand $0$,
which absorbs the remaining linear values $1$. A variable
in $L$ instead supplies a single $3$, added only to $1$s.
For 1156 and 1157 use $2$ as background and the same
argument. In 1556, a prefix-only occurrence gives $1$
in every nonlinear summand; a last-variable occurrence
gives at least one $0$, which is additively absorbing.
The outside linear values are $2$ in both cases.

This test recovers $U$ and $L$ in every generator, and
also $T$ in 1556. For the other four generators use a
second test: assign $v=a$ and all other variables $b$.
The output equals $a$ precisely when $v\in T\cup L$,
and equals $b$ otherwise, with the following choices:
\begin{center}
\begin{tabular}{rcc}\toprule
Generator&$a$&$b$\\\midrule
775&2&1\\778&1&2\\1156&0&1\\1157&1&0\\\bottomrule
\end{tabular}
\end{center}
For 775, the last variable $2$ is detected against $1$s;
for 778 the absorbing choice among $\{1,2\}$ is instead $1$.
For 1156, a last variable $0$ is detected against $1$s;
for 1157 the absorbing choice among $\{0,1\}$ is instead $1$.
Prefix-only occurrences do not change the other last values
in these assignments. Since $L$ has already been recovered
and is disjoint from $T$, this recovers $T$.

Distinct allowed triples give distinct term functions in every
direct generator. Lemmas~\ref{f46:lem:reduction} and
\ref{lem:normal-form-completeness} prove completeness; the common
basis $\Sigma$ gives equality of the five varieties.
\end{proof}

\endgroup
\subsection{Minimal-support antichains}\label{sec:family47}
\begingroup
\providecommand{\eq}{}\renewcommand{\eq}{\approx}
\providecommand{\SR}{}\renewcommand{\SR}{\mathsf{SR}}
\providecommand{\Var}{}\renewcommand{\Var}{\operatorname{Var}}
\providecommand{\Min}{}\renewcommand{\Min}{\operatorname{Min}}
\providecommand{\Sing}{}\renewcommand{\Sing}{\operatorname{Sing}}
\providecommand{\ac}{}\renewcommand{\ac}{\operatorname{ac}}
\subsubsection{Generators}
The four generators listed below satisfy
$3+3=0\ne3$. Their polynomial normal forms are determined by the
minimal nonlinear supports and the linear part.

\emph{Generators (basis size):} $A_{676}$ (11), $A_{1154}$ (11), $A_{1635}$ (11), $A_{1736}$ (11).

In 1154, a nonlinear word has value $1$ if its last variable
is $1$, value $2$ if all its variables are $2$, and $0$ otherwise.
On the nonlinear values, addition is maximum in the chain
$1<0<2$; moreover $1$ is an additive identity and $3+3=0$.
In 1736, a nonlinear word has value $2$ if its last variable is
$2$, value $1$ if its last variable is $1$, and otherwise value
$0$ when every variable in its support is in $\{0,3\}$, or $1$
when some such variable is in $\{1,2\}$. Addition on nonlinear
values is maximum in $2<1<0$, and $2$ is an additive identity.
These descriptions follow by induction from the defining tables.

\subsubsection{Minimal-support normal forms}
Let $\Gamma$ consist of $\SR$ and six identities:
\begin{equation}\label{f47:eq:basis}
\begin{aligned}
 x^2y&\eq xy,&xyz&\eq yxz,&xyx&\eq yx,\\
 x+x&\eq x^2,&x+yx&\eq x+x,&x+yz&\eq x+yz+xz.
\end{aligned}
\end{equation}
All eleven identities hold in 1154 and 1736 by the preceding
word evaluations and the addition tables.
Every identity uses at most three variables. Neither the number
of identities nor the number of variables is claimed minimal.

Fix a finite nonempty ordered variable set $X$. Write $c(w)$ for
the support of a word $w$. The notation $\Min\mathcal B$ means
the family of inclusion-minimal members of a nonempty finite
family $\mathcal B$ of nonempty sets. For an antichain $\mathcal A$,
put $\Sing\mathcal A=\{x:\{x\}\in\mathcal A\}$.

\begin{lemma}\label{f47:lem:normal}
The basis $\Gamma$ implies
\begin{equation}\label{f47:eq:derived}
\begin{aligned}
 (xy)^2&\eq xy,&xy+xy&\eq xy,&x+x^2&\eq x^2,\\
 x+x+x&\eq x+x,&xy+zxy&\eq xy.
\end{aligned}
\end{equation}
Every nonlinear word reduces to
\begin{equation}\label{f47:eq:word}
 W(U,y)=\Bigl(\prod_{x\in U}x\Bigr)y,\qquad y\in U\subseteq X,
\end{equation}
where the product uses the fixed order. Every term reduces to a
normal form indexed by parameters
\begin{equation}\label{f47:eq:parameters}
\begin{gathered}
 \varnothing\ne\mathcal A\subseteq 2^X\setminus\{\varnothing\}
 \text{ is an antichain},\qquad D\subseteq X,\\
 U\in\mathcal A,\ |U|>1\ \Longrightarrow\ U\cap D\ne\varnothing,
 \qquad L=\Sing\mathcal A\setminus D.
\end{gathered}
\end{equation}
The corresponding form is
\begin{equation}\label{f47:eq:normal}
 N(\mathcal A,D)=\sum_{x\in L}x+
 \sum_{y\in D}\sum_{U\in\mathcal A}W(U\cup\{y\},y).
\end{equation}
Empty sub-sums are omitted; the expression as a whole is nonempty.
\end{lemma}
\begin{proof}
Twice applying $xyx\eq yx$ gives $xyxy\eq yxy\eq xy$.
Doubling $xy$ into its square proves its additive idempotence.
Setting $y=x$ in $x+yx\eq x+x$ gives $x+x^2\eq x^2$;
the threshold law follows. Finally, substitute a nonlinear
word $w$ for $x$ in $x+yx\eq x+x$: this gives
$w+vw\eq w+w\eq w$. This includes the last derived identity.

The prefix laws permute adjacent prefix variables and remove
their repetitions. For a nonlinear word $vy$, use
$yvy\eq vy$ backwards to insert its last variable in the prefix.
This gives the standard word \eqref{f47:eq:word} with $U=c(vy)$.
If $y\in U\subseteq V$, the larger word $W(V,y)$ is equivalent
to a left multiple of $W(U,y)$, unless $U=V$. Thus
$W(U,y)+W(V,y)\eq W(U,y)$ in either case.

Expand a term by distributivity, retaining multiplicities.
Repeated copies of a linear variable become its square; further
copies are absorbed. If a linear $x$ and any nonlinear word
ending in $x$ coexist, $x+vx\eq x^2$ replaces them by a square.
Any further word ending in $x$ is absorbed by that square.
After these steps, let $L$ be the surviving linear variables
and $D$ the nonlinear last variables, so $L\cap D=\varnothing$.
Let $\mathcal A$ be the minimal supports among all summands,
including the singleton supports of linears.

For each $U\in\mathcal A$, retain an original summand $u_U$
with that support. For each $y\in D$, choose an existing
nonlinear summand $vy$ ending in $y$. Substitution in the
sharing identity gives
$u_U+vy\eq u_U+vy+u_Uy$. Hence we may first add all words
$W(U\cup\{y\},y)$ for $U\in\mathcal A$ and $y\in D$.
Every original nonlinear word $W(V,y)$ has some $U\in\mathcal A$
with $U\subseteq V$. Its newly added companion has support
$U\cup\{y\}\subseteq V$ and the same last variable, so it
absorbs the original word. Deleting these redundant summands
leaves exactly \eqref{f47:eq:normal}.

Every nonsingleton minimal support comes from a nonlinear
summand and therefore meets $D$. Every surviving linear gives
a singleton in $\mathcal A$. Conversely, a singleton support
$\{x\}$ is represented either by the surviving linear $x$ or
by a nonlinear word ending in $x$; these two cases are disjoint.
Consequently $L=\Sing\mathcal A\setminus D$. If $D$ is empty,
all members of $\mathcal A$ are singletons and $L$ is nonempty.
If $D$ is nonempty, the double sum is nonempty. No empty product
or nullary zero is used.
\end{proof}

\subsubsection{Complete bases}
\begin{theorem}\label{f47:thm:complete}
$\Gamma$ is a complete identity basis for 1154 and for 1736.
Every allowed pair $(\mathcal A,D)$ is realized, and distinct
pairs induce distinct term functions in either generator.
In particular, $\Var(A_{1154})=\Var(A_{1736})=\mathcal W$.
\end{theorem}
\begin{proof}
For an allowed pair, the minimal supports of \eqref{f47:eq:normal}
are exactly $\mathcal A$. A singleton in $L$ is explicitly
present. Every other $U\in\mathcal A$ meets $D$: choose
$y\in U\cap D$ to obtain the summand $W(U,y)$. Every remaining
displayed support contains a member of $\mathcal A$.
Also the nonlinear last set is exactly $D$ and the linear
set is exactly $L$. Thus every allowed pair is realized.

In 1154, assign a selected variable $v$ value $3$ and every
other variable value $1$. The result is $3$ if $v\in L$,
$0$ if $v\in D$, and $1$ otherwise. All other summands end
in $1$, so prefix occurrences of $v$ do not interfere.
This recovers $D$ and $L$.

For 1736, assign $v=1$ and all other variables $2$.
The result is $1$ exactly when $v\in D\cup L$ and is $2$
otherwise. Next assign $v=3$ and the others $2$. The result
is $3$ exactly when $v\in L$: nonlinear words never take
value $3$, while all other last variables have value $2$.
When $v\in D$, the result is $0$ or $1$; when $v$ is in
neither set, the result is $2$. These tests recover both sets.

To recover $\mathcal A$, choose any subset $Q\subseteq X$.
In 1154 assign $2$ on $Q$ and $1$ outside $Q$. A summand
takes value $2$ exactly when its support is contained in $Q$;
$2$ is additively absorbing. In 1736 assign $0$ on $Q$ and
$2$ outside $Q$. A summand takes value $0$ exactly when its
support is contained in $Q$; now $0$ is additively absorbing.
In either case the designated value is obtained exactly
when some $U\in\mathcal A$ is contained in $Q$. The minimal
such $Q$ therefore recover $\mathcal A$.

The tests separate all parameters for every finite nonempty $X$.
Lemma~\ref{lem:normal-form-completeness} proves completeness of $\Gamma$
in both generators, and hence equality of their varieties.
\end{proof}

\endgroup
\section{Boolean and capped coefficient constructions}\label{sec:coefficient-forms}
The next constructions use finite coefficient systems together
with support data. Their normal forms retain precisely the distinctions
visible in the addition tables; each recovery argument applies to an
arbitrary finite alphabet.

\subsection{Boolean-ring joins and affine lifts}\label{sec:family49}
\begingroup
\providecommand{\eq}{}\renewcommand{\eq}{\approx}
\providecommand{\SR}{}\renewcommand{\SR}{\mathsf{SR}}
\providecommand{\Var}{}\renewcommand{\Var}{\operatorname{Var}}
\providecommand{\Aff}{}\renewcommand{\Aff}{\operatorname{Aff}}
\providecommand{\supp}{}\renewcommand{\supp}{\operatorname{supp}}
\providecommand{\F}{}\renewcommand{\F}{\mathbb F_2}
\emph{Generators (basis size):} $A_{1324}$ (12), $A_{1326}$ (12), $A_{1667}$ (22), $A_{1899}$ (22).

For a word $w$, write $c(w)$ for its nonempty support.
If $X$ is a fixed finite variable set, its parity vector
$q(w)\in\F^X$ records the exponents modulo two.

\subsubsection{A finite basis construction for Boolean-ring joins}
Write $R$ for the two-element field in the signature
$(+,\cdot)$, and $\mathcal B=\Var(R)$.
Word expansion and $a^2=a$ in $R$ show that an identity
$P\eq Q$ holds in $R$ precisely when, for every nonempty
support $D$, the numbers of summands having support $D$ in
$P$ and $Q$ have the same parity. Necessity follows from
linear independence of the squarefree monomial functions
on $\{0,1\}^X$: evaluating on the characteristic vectors of
subsets of $X$ and inducting on subset size recovers every
coefficient. Sufficiency follows by direct evaluation.

\begin{lemma}\label{f49:lem:double}
In any semiring satisfying $3x\eq x$, the map $h(x)=2x$
is an endomorphic retraction onto an ai-subsemiring.
Every polynomial $P$ is equivalent to $2P+L_P$, where
$L_P$ consists of the word occurrences left with odd
multiplicity after pairing equal words. If no such
occurrences remain, the expression is $2P$.
\end{lemma}
\begin{proof}
The identity implies $4x\eq2x$, $x+2x\eq x$, and
$2x+2x\eq2x$. Thus $h$ is idempotent and its image is
additively idempotent. It preserves addition, and
$h(x)h(y)=4xy\eq2xy=h(xy)$.
Apply $P\eq P+2P$, expand, and use the idempotence of the
doubled background to remove every paired copy. This proves
the polynomial assertion without using an empty term.
\end{proof}

For positive integers $r,s$ put $\ell=r+s-1$.
Let $\mathcal T_\ell$ be the following finite set.
For each nonempty subset
$J\subseteq\{1,\ldots,\ell\}^2$, introduce one variable
$z_{ab}$ for each $(a,b)\in J$, and put
\begin{equation}\label{f49:eq:profiles}
 u_J=\prod_{(a,b)\in J}z_{ab}^{a},\qquad
 v_J=\prod_{(a,b)\in J}z_{ab}^{b}.
\end{equation}
If $u_J=v_J$, omit the resulting trivial identity.
Otherwise include
\begin{equation}\label{f49:eq:transfer}
 2u_J+2v_J+u_J\eq2u_J+2v_J+v_J.
\end{equation}
There are at most $2^{\ell^2}-1$ such identities,
each using at most $\ell^2$ variables.

\begin{lemma}\label{f49:lem:profiles}
Commutativity of multiplication, $x^{r+s}\eq x^r$, and
$\mathcal T_\ell$ imply
\begin{equation}\label{f49:eq:alltransfer}
 2u+2v+u\eq2u+2v+v
\end{equation}
for all nonempty words $u,v$ with $c(u)=c(v)$.
\end{lemma}
\begin{proof}
First reduce every positive exponent to
$\{1,\ldots,\ell\}$. Partition the common support according
to the pair $(a,b)$ of its two exponents.
For each nonempty block substitute its product of distinct
variables for $z_{ab}$ in \eqref{f49:eq:transfer}.
Commutativity recovers the two reduced words exactly.
If all pairs are diagonal, the words already agree.
The blocks are nonempty, so no multiplicative identity
has been used in the substitution.
\end{proof}

\begin{theorem}\label{f49:thm:join}
Let $\mathcal K$ be a commutative ai-semiring variety
with a finite identity basis $\Sigma$, and suppose
$x^{r+s}\eq x^r$ holds in $\mathcal K$ for positive $r,s$.
Then $\mathcal K\vee\mathcal B$ is finitely based.
A complete basis consists of $\SR$ and
\begin{equation}\label{f49:eq:joinbasis}
 3x\eq x,\quad xy\eq yx,\quad x^{r+s}\eq x^r,\quad
 \mathcal T_\ell,\quad
 \{P(2\boldsymbol x)\eq Q(2\boldsymbol x):
                        (P\eq Q)\in\Sigma\}.
\end{equation}
In particular the hypothesis on the power identity holds
whenever $\mathcal K$ is generated by a finite commutative
ai-semiring.
\end{theorem}
\begin{proof}
All laws in \eqref{f49:eq:joinbasis} are valid in both factors.
For \eqref{f49:eq:transfer}, addition is idempotent in
$\mathcal K$, so either side is $u_J+v_J$.
In $R$, the doubled terms vanish and words with the same
support agree. The other laws are immediate; lifted
identities are evaluated at zero in $R$.

For completeness let $A$ satisfy the proposed basis.
Lemma~\ref{f49:lem:double} gives a retraction $h$, and the
lifted identities give $h(A)\in\mathcal K$.
Consequently any identity $P\eq Q$ of the join implies
$2P\eq2Q$ in $A$. Equationally this follows by lifting
a $\Sigma$-derivation and using the homomorphism laws for $h$.

Expand $P,Q$ into words, reduce exponents, and retain the
common doubled background
\[
 b=2P+2Q\eq2P\eq2Q.
\]
The background is additively idempotent and absorbs
$2w$ for every word occurrence $w$ in either polynomial.
Lemma~\ref{f49:lem:double} therefore gives $P\eq b+L_P$
and $Q\eq b+L_Q$.

For each support $D$ occurring in either expansion choose
one representative word $v_D$ from the two expansions.
In the presence of $b$, equation \eqref{f49:eq:alltransfer}
replaces every residual word $u$ with support $D$ by $v_D$:
insert the already absorbed $2u+2v_D$, use the identity,
and remove those doubled terms again.
Pairs of $v_D$ then disappear into $b$.
Only the parity of the number of occurrences of support
$D$ remains. Those parities agree for $P,Q$ because the
identity holds in $R$. Both sides have thus been reduced
to the same expression. This proves completeness for all
finite sets of variables.

For a finite generating algebra, the sequence of unary
multiplicative powers eventually repeats, giving the
required positive $r,s$.
\end{proof}

The transfer laws are essential to this argument.
Containment in a finitely based variety by itself is not
an FB criterion. The theorem is a construction for the
specified join; it does not assert hereditary finite basability.

\subsubsection{The flat group semiring and two applications}
Let $F=\{e,g,z\}$, where $\{e,g\}$ is the group of order
two under multiplication and $z$ is multiplicatively
absorbing. Addition is flat:
$e+e=e$, $g+g=g$, $e+g=z$, and $z+x=z$.
Its local representative is $B_{80}$, with operation strings
\[
 F_+=\texttt{011111112},\qquad
 F_\cdot=\texttt{012111210}
\]
in the order $(e,z,g)$.
The following proof supplies its finite basis directly.

\begin{lemma}\label{f49:lem:flat}
A complete basis of $F$ consists of $\SR$ and
\begin{equation}\label{f49:eq:flatbasis}
\begin{aligned}
 x+x&\eq x,&xy&\eq yx,&x^3&\eq x,\\
 x+y&\eq xy^2+x^2y,&
 x+y+z&\eq x+y+z+xyz.
\end{aligned}
\end{equation}
The free algebra on a finite nonempty set $X$ consists of
pairs $(C,H)$, where $\varnothing\ne C\subseteq X$
and $H$ is a nonempty affine subspace of $\F^C$.
\end{lemma}
\begin{proof}
Soundness follows immediately from flat addition and the
group of exponent two. Normalize a polynomial into words.
The fourth law in \eqref{f49:eq:flatbasis} replaces a pair
$u,v$ by $uv^2,u^2v$, giving both words the union of their
supports without changing their respective parity vectors.
For several word occurrences, first pair a fixed occurrence
successively with all others to give it the full support
$C$, and then make a second pass to give all occurrences
that support. With one occurrence no pass is needed.
The law $x^3\eq x$ reduces the positive exponents to
one and two.

The last law inserts the sum modulo two of any three parity
vectors. By Lemma~\ref{lem:affine-recovery}, the resulting canonical
expression contains all monomials with support $C$ and parity in
the affine span $H$.

Recover $C$ by assigning a selected variable $z$ and all
others $e$: the output is $z$ exactly when the selected
variable is in $C$. On group assignments encode $e,g$
as $0,1$. Write $H=q_0+L$. The polynomial is nonzero
exactly on $L^\perp$, where its value has bit $q_0\cdot t$.
Lemma~\ref{lem:affine-recovery} recovers
\begin{equation}\label{f49:eq:dual}
 H=\{q\in\F^C:q\cdot t=q_0\cdot t
                         \text{ for all }t\in L^\perp\}.
\end{equation}
Thus $(C,H)$ is determined, proving completeness by
Lemma~\ref{lem:normal-form-completeness}.
\end{proof}

For $\ell=2$, the profile set consists of
$(1,1),(1,2),(2,1),(2,2)$.
The twelve nontrivial subsets are recorded compactly
in Table~\ref{f49:tab:transfers}. A string such as
$11,12,22$ means the three indicated profiles, assigned
in that order to $x,y,z$ in \eqref{f49:eq:profiles}.
\begin{table}[htbp]\centering
\caption{The twelve instances of \eqref{f49:eq:transfer}
used in the Boolean join.}\label{f49:tab:transfers}
\begin{tabular}{llll}\toprule
$12$&$11,12$&$21$&$11,21$\\
$12,21$&$11,12,21$&$12,22$&$11,12,22$\\
$21,22$&$11,21,22$&$12,21,22$&$11,12,21,22$\\
\bottomrule\end{tabular}
\end{table}

\begin{theorem}\label{f49:thm:applications}
The algebras $A_{1667}$ and $A_{1899}$ both generate
$\Var(F)\vee\mathcal B$. A common complete basis
has 22 identities: $\SR$, the twelve transfers in
Table~\ref{f49:tab:transfers}, and the five identities
\begin{equation}\label{f49:eq:22}
\begin{aligned}
 3x&\eq x,&xy&\eq yx,&x^3&\eq x,\\
 2x+2y&\eq (2x)(2y)^2+(2x)^2(2y),\\
 2x+2y+2z&\eq 2x+2y+2z+(2x)(2y)(2z).
\end{aligned}
\end{equation}
\end{theorem}
\begin{proof}
Apply Theorem~\ref{f49:thm:join} with $r=1,s=2$ and
the basis in Lemma~\ref{f49:lem:flat}.
The lifted semiring laws are consequences of $\SR$ and
Lemma~\ref{f49:lem:double}; the lifted idempotence,
commutativity, and cube laws also follow from the three
unlifted identities in the first line of \eqref{f49:eq:22}.
This leaves exactly the two lifted laws in its last lines.
All 22 identities hold in both displayed four-element
tables by direct substitution.
The theorem therefore gives the forward containment
in the claimed join.

For the reverse containment, $F$ embeds in $A_{1667}$
by $(e,z,g)\mapsto(0,1,3)$, and $R$ embeds by
$(0,1)\mapsto(1,2)$.
In $A_{1899}$ the corresponding embeddings are
$(e,z,g)\mapsto(0,2,3)$ and $(0,1)\mapsto(0,1)$.
Direct checking of both operations proves these claims.
Both factors therefore belong to each generated variety.
\end{proof}

\subsubsection{Two lifts of the flat group semiring}
Use the carrier $\{e,z,g,d\}$ with labels $(0,1,2,3)$.
Let $\pi$ fix $F$ and send $d$ to $g$, and define
\[
 ab=\pi(a)\pi(b)\quad(a,b\in F\cup\{d\}).
\]
Thus every product of length at least two lies in $F$.
Addition extends that of $F$ by
\[
 d+z=z,\qquad d+e=z,\qquad d+g=d,\qquad
 d+d=\begin{cases}g&\text{in }E_{\mathrm c},\\
                  d&\text{in }E_{\mathrm i}.
       \end{cases}
\]
These are exactly $A_{1324}$ and $A_{1326}$ respectively.
The map $\pi(a)=a^3$ is a homomorphism onto $F$ in both.

Let $\Gamma_\epsilon$ consist of $\SR$ and the seven
identities below:
\begin{equation}\label{f49:eq:12}
\begin{aligned}
 xy&\eq yx,&xy&\eq x^3y,\\
 xy+zt&\eq xy(zt)^2+(xy)^2zt,\\
 x+y+z&\eq x+y+z+xyz,\\
 x+y^2&\eq x^3+y^2,&x+x^3&\eq x,\\
 x+x&\eq x^3 &&(\epsilon=\mathrm c),\\[-2pt]
 x+x&\eq x &&(\epsilon=\mathrm i).
\end{aligned}
\end{equation}
Only one of the last two alternatives is included.

For $\varnothing\ne C\subseteq X$ and $q\in\F^C$,
let $m(C,q)$ be the monomial with exponent two at
coordinates with bit zero and exponent three at coordinates
with bit one. Write $e_x$ for the coordinate vector at $x$.
Define parameters
\begin{equation}\label{f49:eq:liftparameters}
\begin{gathered}
 \varnothing\ne C\subseteq X,\quad
 H\subseteq\F^C\text{ a nonempty affine subspace},\\
 A\subseteq\{x\in C:e_x\in H\},\qquad
 0\in H\Longrightarrow A=\varnothing .
\end{gathered}
\end{equation}
Their proposed normal form is
\begin{equation}\label{f49:eq:liftnormal}
 N(C,H,A)=\sum_{x\in A}x+\sum_{q\in H}m(C,q).
\end{equation}

\begin{theorem}\label{f49:thm:lifts}
$\Gamma_{\mathrm c}$ and $\Gamma_{\mathrm i}$ are
complete twelve-identity bases for $E_{\mathrm c}$ and
$E_{\mathrm i}$ respectively. The forms
\eqref{f49:eq:liftnormal}, subject to \eqref{f49:eq:liftparameters},
give their free algebras on every finite nonempty $X$.
\end{theorem}
\begin{proof}
First check soundness. The products lie in $F$, giving
the period law and the law for two nonlinear summands.
For the ternary insertion, a sum not equal to $z$ is
either made entirely of $e$'s, or entirely of elements
of $\{g,d\}$. The product is respectively $e$ or $g$,
which is absorbed by the sum. In the remaining case
the sum is $z$ and absorbs every product.
A square is $e$ or $z$, so $x+y^2$ is unchanged by
replacing $d$ by $g$ in its first argument.
Cube absorption and the two alternative doubling laws
follow from $d+g=d$ and the displayed values of $d+d$.
This checks all seven additional laws.

We prove normalization by equational deductions.
The law $xy\eq x^3y$ gives $x^4\eq x^2$ by putting
$y=x$. It changes any exponent one in a nonlinear
word to three and reduces larger exponents to two or
three. In particular $(uv)^3\eq uv$. Substitution
in the doubling law shows that every nonlinear word
is additively idempotent, in either algebra.

Expand the term as a polynomial. In $E_{\mathrm i}$
retain at most one occurrence of each bare linear
variable. In $E_{\mathrm c}$ replace pairs of linear
occurrences $x+x$ by $x^3$ and use $x+x^3\eq x$
when a further bare occurrence is present.
In both cases let $A$ record the remaining linear
variables. Insert $x^3$ beside every remaining bare
$x$ using cube absorption. There is now at least
one nonlinear summand, and the union $C$ of their
supports is the content of the original polynomial.
Their parity vectors include $e_x$ for each $x\in A$.

Use the law for two nonlinear summands in
\eqref{f49:eq:12}, factoring each chosen summand into
two nonempty words. The two-pass procedure from
Lemma~\ref{f49:lem:flat} gives all nonlinear summands
the full support $C$, preserving their parities.
The ternary insertion, followed by power reduction,
then adds exactly the affine closure $H$ of these
vectors. Nonlinear idempotence removes repetitions.
This derives \eqref{f49:eq:liftnormal} except possibly
for the final constraint in \eqref{f49:eq:liftparameters}.

If $0\in H$, the summand $m(C,0)$ is a square
of a nonempty word. The square projection in
\eqref{f49:eq:12} replaces every bare $x\in A$ by
$x^3$. Equalizing its support again and using
$e_x\in H$ makes it redundant. Thus all bare
linears can be deleted, giving $A=\varnothing$.
No other parameter is changed.

We now separate the resulting forms.
Putting one variable at $z$ and all others at $e$
recovers $C$, as in Lemma~\ref{f49:lem:flat}.
On group assignments the bare parity vectors $e_x$
already lie in $H$, so \eqref{f49:eq:dual} recovers $H$.
If $0\in H$, the parameter $A$ is prescribed to be empty.
Otherwise Lemma~\ref{lem:affine-recovery} gives a functional $t$
equal to one on $H$.

For a coordinate $x$ with $e_x\in H$, this forces
$t_x=1$. Use the group assignment specified by $t$
on the other coordinates, and put $x=d$.
Every nonlinear summand has value $g$.
Every bare summand other than $x$ also has value $g$.
The total is $d$ exactly when $x\in A$, and is $g$
otherwise. This remains true for $E_{\mathrm c}$
because there is at most one $d$ summand in this test.
Coordinates with $e_x\notin H$ cannot belong to $A$.
Thus the function determines $A$ as well.

All allowed parameters are realized by the displayed
expression. Its nonlinear summands have support $C$,
their affine span is exactly $H$, and the displayed
bare set satisfies the constraints. Completeness follows from
Lemma~\ref{lem:normal-form-completeness}.
\end{proof}

The two varieties in Theorem~\ref{f49:thm:lifts} are distinct:
only $E_{\mathrm i}$ satisfies additive idempotence.
Their free algebras nevertheless have the same cardinalities.
This also shows why comparing free-algebra sizes alone
would not prove equality of varieties.

\endgroup
\subsection{Capped addition}\label{sec:family50}
\begingroup
\providecommand{\eq}{}\renewcommand{\eq}{\approx}
\providecommand{\SR}{}\renewcommand{\SR}{\mathsf{SR}}
\providecommand{\Var}{}\renewcommand{\Var}{\operatorname{Var}}
\providecommand{\SL}{}\renewcommand{\SL}{\mathsf{SL}}
\providecommand{\supp}{}\renewcommand{\supp}{\operatorname{supp}}
\providecommand{\FF}{}\renewcommand{\FF}{\mathbb F_2}
\providecommand{\zz}{}\renewcommand{\zz}{\mathbf 0}
\providecommand{\PP}{}\renewcommand{\PP}{\mathcal P}
\providecommand{\SSS}{}\renewcommand{\SSS}{\mathcal S}

\subsubsection{The saturation construction}
Let $S$ be a nonempty semigroup. Adjoin two distinct new elements
$\zz$ and $\omega$, and put
\[
C(S)=S\mathbin{\dot\cup}\{\zz,\omega\}.
\]
Addition and multiplication are defined as follows:
\begin{align}
\zz+a&=a+\zz=a,&
a+b&=\omega\quad(a,b\ne\zz),\label{f50:eq:addition}\\
\zz a&=a\zz=\zz,&
\omega a&=a\omega=\omega\quad(a\ne\zz),\label{f50:eq:multiplication}
\end{align}
with multiplication on $S$ unchanged. The symbols $\zz,\omega$
name elements in this description only. They are not nullary
operations or constants in any basis below.
An element that was a zero inside $S$, if one exists, is distinct
from the newly adjoined $\zz$.

\begin{lemma}\label{f50:lem:construction}
The algebra $C(S)$ is a semiring. Its multiplicative reduct
$C(S)_\cdot$ satisfies a word identity $u\eq v$ precisely when
$c(u)=c(v)$ and $S\models u\eq v$. Consequently,
\[
\Var(C(S)_\cdot)=\Var(S)\vee\SL,
\]
where $\SL$ denotes the semigroup variety of semilattices.
\end{lemma}
\begin{proof}
A sum is $\zz$ if all its summands are $\zz$, is its sole nonzero
summand if there is exactly one, and is $\omega$ if at least two
are nonzero. This proves additive associativity and commutativity.
A product with a $\zz$ factor is $\zz$; otherwise a product with
an $\omega$ factor is $\omega$; all remaining products are those
in $S$. This proves multiplicative associativity.
Multiplication by a nonzero element preserves the distinction
between zero and nonzero. Thus multiplying a sum with two nonzero
summands gives $\omega$ on both sides of either distributive law.
The cases with zero or one nonzero summand are immediate.

Every identity of the reduct holds in $S$. If a variable occurs
on only one side, assign it $\zz$ and assign all other variables
a fixed member of $S$. One side is $\zz$ and the other belongs
to $S$, a contradiction. Conversely, when supports agree, any
assignment with a used variable equal to $\zz$ gives $\zz$ on
both sides. If no such value occurs but a used variable is
$\omega$, both sides give $\omega$. Otherwise the identity is
evaluated inside $S$. This proves the characterization.
Semilattices satisfy precisely the word identities with equal
supports. The identity theory of the reduct is therefore the
intersection of those of $S$ and $\SL$, which proves the last
assertion.
\end{proof}

Let $\mathcal E$ consist of the following five laws:
\begin{equation}\label{f50:eq:extras}
\begin{aligned}
3x&\eq2x,&2x^2&\eq2x,&2xy&\eq2yx,\\
x+y+xy&\eq x+y,&x+y+xyz&\eq x+y.
\end{aligned}
\end{equation}
Every law holds in $C(S)$. The first three depend only on which
variables are zero. For either absorption law, the added product
is zero if $x$ or $y$ is zero, while $x+y=\omega$ if both are
nonzero.

\begin{lemma}\label{f50:lem:doubled}
From $\SR$ and the second and third laws of
\eqref{f50:eq:extras} one can deduce
\[
2u\eq2v\qquad\text{whenever }c(u)=c(v).
\]
\end{lemma}
\begin{proof}
Inside a doubled word, the law $2xy\eq2yx$ permits any adjacent
transposition. Indeed, multiplication on either side by a context
word, followed by distributivity, gives
$2pxyq\eq2pyxq$; absent left or right contexts are simply omitted.
Sort the letters in a fixed order, then remove adjacent
repetitions using $2x^2\eq2x$ in the same way. The result is
twice the fixed squarefree word on the support.
\end{proof}

\begin{theorem}\label{f50:thm:transfer}
For any nonempty semigroup $S$, the semiring $C(S)$ is finitely
based if and only if its multiplicative reduct is finitely based.
More explicitly, if $\Theta$ is a finite basis of
$C(S)_\cdot$ relative to associativity, then
\begin{equation}\label{f50:eq:generalbasis}
\SR\ \cup\ \Theta\ \cup\ \mathcal E
\end{equation}
is a finite basis of $C(S)$.
\end{theorem}
\begin{proof}
We first prove sufficiency, including completeness for arbitrarily
many variables. Soundness has already been established.
Fix a finite nonempty ambient variable set $X$ and expand a
term $P$ into its multiset of word occurrences. For
$\varnothing\ne Y\subseteq X$ set
\begin{equation}\label{f50:eq:profile}
b_P(Y)=\min\bigl(2,\#\{w\text{ occurring in }P:c(w)\subseteq Y\}\bigr).
\end{equation}
Define the saturated family
$\SSS_P=\{Y:b_P(Y)=2\}$, which is upward closed.

For every $Y\in\SSS_P$, choose two distinct occurrences $u,v$
with $c(u)\cup c(v)\subseteq Y$. If their union is $Y$,
the first absorption law inserts $uv$. Otherwise let $z$
be a nonempty word containing exactly the missing variables;
the second law inserts $uvz$. Repeat the insertion to obtain
two copies. Lemma~\ref{f50:lem:doubled} changes them into $2m_Y$,
where $m_Y$ is the fixed squarefree word with support $Y$.
All original occurrences are retained during this process.
Any support at which the inserted word is active contains $Y$
and was already saturated. Thus these insertions do not change
the profile.

After inserting $2m_Y$ for every saturated $Y$, remove every
original occurrence $w$ whose support $U=c(w)$ is saturated:
\[
2m_U+w\eq2w+w\eq2w\eq2m_U.
\]
The required doubled marker remains available after each removal,
even when several identical occurrences are removed.
Every remaining original word has an unsaturated support, and
there is exactly one occurrence at that support.
These supports form an antichain $\PP_P$. If one were contained
in another, the latter would have two active original occurrences
and would be saturated. Similarly, for distinct $U,V\in\PP_P$,
their union belongs to $\SSS_P$.

For each $U\in\PP_P$, replace its remaining word by a fixed
representative $w_{U,\lambda_U}$ of its identity class on $S$,
among words with support exactly $U$.
Lemma~\ref{f50:lem:construction} and the completeness of $\Theta$
justify this replacement equationally. We have derived
\begin{equation}\label{f50:eq:normal}
N(\SSS,\PP,\lambda)=
\sum_{U\in\PP}w_{U,\lambda_U}
+\sum_{Y\in\SSS}2m_Y.
\end{equation}
The existence of representatives is sufficient here; when $S$
is finite there are only finitely many word functions for each
fixed support.

We now show that a valid identity has the same normal form on
both sides. Fix any $s\in S$ and assign all variables of $Y$
the value $s$, with all others equal to $\zz$. The three values
of $b_P(Y)$ are distinguished: the output is $\zz$ for zero,
lies in $S$ for one, and is $\omega$ for two.
Consequently a term function determines $b_P$, and hence
$\SSS_P$. The members of $\PP_P$ are precisely the
inclusion-minimal positive supports whose profile value is one,
so they too are determined.

For $U\in\PP_P$, assign arbitrary elements of $S$ to the
variables in $U$ and $\zz$ to all others. Exactly one word
occurrence is active. Thus equality of term functions implies
that the two residual words with support $U$ define the same
word function on $S$. They consequently have the same
representative. Both sides of every valid identity therefore
reduce to \eqref{f50:eq:normal}. This proves completeness of
\eqref{f50:eq:generalbasis}.

For necessity, suppose $\Gamma$ is a finite semiring identity
basis of $C(S)$. Assign every variable an arbitrary fixed
$s\in S$. A term with no addition has its value in $S$.
A term containing addition has value $\omega$: in a polynomial
expansion it has at least two nonzero word occurrences.
Hence no valid identity has addition on exactly one side.

Let $\Gamma_\cdot$ be the finite subset of identities in
$\Gamma$ with no addition on either side. Consider an
equational deduction, expressed as successive replacements
by substitution instances in contexts, between two words.
Starting with a word, no identity containing addition can
be applied, because its matching side contains an addition
symbol. An applicable identity from $\Gamma_\cdot$ has
word substitutions for all variables on its matching side.
Its two supports agree by Lemma~\ref{f50:lem:construction}, so
its other side also becomes a word. The context is likewise
multiplicative. Inductively, the deduction uses only
$\Gamma_\cdot$, with word substitutions and multiplicative
contexts. Every identity of the reduct is thus a consequence
of this finite set. This proves necessity.
\end{proof}

The theorem concerns the particular construction $C(S)$.
Its hypothesis is finite basability of
$\Var(S)\vee\SL$, equivalently of $C(S)_\cdot$.
It does not replace that hypothesis by finite basability of
$S$ alone, and it does not assert hereditary finite basability.

\subsubsection{Three four-element applications}
Let $G=\{e,g\}$ be the cyclic group with $g^2=e$.
Let $R=\{a,b\}$ have multiplication $xy=y$, and let
$L=\{a,b\}$ have multiplication $xy=x$.
The local representatives of $C(G),C(R),C(L)$ are listed below.
Each string gives the entries of a table row by row on
$\{0,1,2,3\}$. In all three representatives $\zz$ has
label $1$, $\omega$ has label $2$, and the two core
elements have labels $0,3$. For $G$ the label $0$ is $e$.
\emph{Generators (basis size):} $A_{1841}$ (10), $A_{2302}$ (10), $A_{2324}$ (10).

\begin{lemma}\label{f50:lem:wordbases}
Relative to associativity, the respective multiplicative
reducts have the following complete bases:
\begin{align}
C(G)_\cdot:\quad &xy\eq yx,\qquad x^3\eq x;
                                                   \label{f50:eq:groupwords}\\
C(R)_\cdot:\quad &x^2\eq x,\qquad xyz\eq yxz;
                                                   \label{f50:eq:rightwords}\\
C(L)_\cdot:\quad &x^2\eq x,\qquad xyz\eq xzy.
                                                   \label{f50:eq:leftwords}
\end{align}
For a fixed nonempty support $U$, the residual invariants
are respectively a parity vector in $\FF^U$, the last
variable in $U$, and the first variable in $U$.
\end{lemma}
\begin{proof}
For $G$, a word function is determined by the parity of each
exponent. Different parity vectors are separated by assigning
$g$ to a differing variable and $e$ to the others.
By Lemma~\ref{f50:lem:construction}, the reduct additionally
records the support. Commutativity and $x^3\eq x$ reduce
each positive exponent to one or two, which gives exactly
these invariants.

On $R$ the value of a word is the value of its last variable.
The law $xyz\eq yxz$ allows adjacent transpositions anywhere
before the final position, since there is a nonempty suffix.
Move earlier occurrences of the final variable next to its
last occurrence and remove them by idempotence. Sort the
remaining prefix and remove its repetitions in the same way.
The result lists each other support variable once in fixed
order, followed by the distinguished final variable.
The support is separated by $\zz$ assignments, and distinct
final variables are separated by the two distinct values of
$R$. This proves the right-zero claim. Reversing words
proves the left-zero claim, with a fixed first variable.
\end{proof}

\begin{theorem}\label{f50:thm:three}
Each listed generator has a complete
basis of ten identities. All three bases contain $\SR$
and the two absorption identities
\begin{equation}\label{f50:eq:common}
x+y+xy\eq x+y,\qquad x+y+xyz\eq x+y.
\end{equation}
The three additional identities are, respectively,
\begin{align}
1841:\quad &xy\eq yx,\qquad x^3\eq x,\qquad
                         2x\eq2x^2; \label{f50:eq:b1841}\\
2302:\quad &x^2\eq x,\qquad xyz\eq yxz,\qquad
                         2xy\eq2yx; \label{f50:eq:b2302}\\
2324:\quad &x^2\eq x,\qquad xyz\eq xzy,\qquad
                         2xy\eq2yx. \label{f50:eq:b2324}
\end{align}
No minimality assertion is made about these basis sizes.
\end{theorem}
\begin{proof}
For $1841$, put $x=y=z$ in the second absorption
identity to obtain $2x+x^3\eq2x$. The cube identity
gives $3x\eq2x$, and doubled commutativity follows
from multiplicative commutativity.
For either band case, put $y=x$ in the first absorption
identity. The resulting $2x+x^2\eq2x$ gives
$3x\eq2x$, while $2x^2\eq2x$ follows directly from
idempotence.
Thus each proposed ten-identity basis entails every
law in \eqref{f50:eq:generalbasis} for the corresponding
word basis of Lemma~\ref{f50:lem:wordbases}.
All its laws hold in the specified generator.
Theorem~\ref{f50:thm:transfer} proves completeness.
\end{proof}

\endgroup
\subsection{Periodic boundaries and single occurrences}\label{sec:family51}
\begingroup
\providecommand{\eq}{}\renewcommand{\eq}{\approx}
\providecommand{\SR}{}\renewcommand{\SR}{\mathsf{SR}}
\providecommand{\Var}{}\renewcommand{\Var}{\operatorname{Var}}
\providecommand{\FF}{}\renewcommand{\FF}{\mathbb F_2}
\providecommand{\A}{}\renewcommand{\A}{\mathcal A}
\providecommand{\K}{}\renewcommand{\K}{\mathcal K}
\providecommand{\E}{}\renewcommand{\E}{\mathcal E}
\providecommand{\R}{}\renewcommand{\R}{\mathcal R}
\providecommand{\Sc}{}\renewcommand{\Sc}{S_{\mathrm c}(abc)}
\providecommand{\stirling}{}\renewcommand{\stirling}[2]{\genfrac{\{}{\}}{0pt}{}{#1}{#2}}
\emph{Generators (basis status):} $A_{452}$ (NFB), $A_{471}$ (NFB), $A_{1387}$ (NFB), $A_{1461}$ (FB (12)), $A_{1465}$ (FB (14)).

Write $P=A_{1461}$ and $Q=A_{1465}$. For these two algebras,
use $(0,e,d,a)$ for labels $(0,1,2,3)$. Their common multiplication
has zero $0$, identity $e$, and
\[
d^2=d,\qquad da=ad=0,\qquad a^2=0.
\]
In $P$, addition has zero $0$ and satisfies
\[
e+e=d,\quad e+d=e,\quad d+d=d,\quad
a+a=0,\quad a+e=e,\quad a+d=d.
\]
In $Q$, the element $0$ is additively absorbing, and
\[
e+e=e+d=d+d=d,\qquad a+a=0,\qquad a+e=a+d=a.
\]
These are descriptions of elements, not named constants in the signature.

For a nonempty set $D$, let $m_D$ be the product of its
variables in a fixed order, and for $x\in D$ put
\[
v_{D,x}=x\prod_{y\in D\setminus\{x\}}y^2,\qquad
\delta_{D,x}=v_{D,x}+m_D^2.
\]
The product following $x$ is omitted when $D=\{x\}$.

\subsubsection{A periodic basis and support boundaries}
Let $\Gamma_P$ consist of $\SR$ and the following seven identities:
\begin{equation}\label{f51:eq:fb12}
\begin{aligned}
xy&\eq yx,&x^3&\eq x^2,\\
3x&\eq x,&2x&\eq2x^2,\\
2x+2xy&\eq2x,\\
xy&\eq x^2y+xy^2+x^2y^2,\\
2y+xy&\eq2y+x^2y.
\end{aligned}
\end{equation}

Fix a finite nonempty variable set $X$. If $\A$ is a nonempty
antichain of nonempty subsets of $X$, define
\begin{equation}\label{f51:eq:boundary}
\begin{aligned}
\K(\A)&=\{D\subseteq X:(\exists U\in\A)\ U\subseteq D\},\\
\partial\K(\A)&=\{(D,x):D\in\K(\A),\ x\in D,\
                                  D\setminus\{x\}\notin\K(\A)\}.
\end{aligned}
\end{equation}
Thus $\partial\K(\A)$ records the edges entering the upward-closed
family when a single variable is adjoined.
Choose arbitrary subsets $\R\subseteq\K(\A)$ and
$\E\subseteq\partial\K(\A)$. The corresponding normal form is
\begin{equation}\label{f51:eq:periodicnormal}
N_P(\A,\R,\E)=
\underbrace{\sum_{U\in\A}2m_U^2}_{B_\A}
+\sum_{D\in\R}m_D^2+\sum_{(D,x)\in\E}\delta_{D,x}.
\end{equation}
The nonempty background $B_\A$ is always present.

\begin{theorem}\label{f51:thm:periodic}
$\Gamma_P$ is a complete twelve-identity basis of $P$.
The parameters $(\A,\R,\E)$ above are in bijection with the
elements of the free algebra of $\Var(P)$ on $X$.
\end{theorem}
\begin{proof}
Soundness can be checked directly from the displayed operations.
The doubled elements are $0,d$ and form a two-element distributive
lattice. The last two laws of \eqref{f51:eq:fb12} follow by considering
$0$, the two elements $e,d$, and the square-zero element $a$.
We prove completeness equationally.

First, all doubled words with the same support are equivalent.
Commutativity orders the letters, and $2x\eq2x^2$, used in
multiplicative contexts with distributivity, makes every positive
exponent two; $x^3\eq x^2$ reduces larger exponents.
The law $3x\eq x$ implies $4x\eq2x$, so doubled terms are
additively idempotent. If $U\subseteq D$, doubled absorption gives
\begin{equation}\label{f51:eq:backgroundabsorb}
2m_U^2+2m_D^2\eq2m_U^2.
\end{equation}
For strict inclusion, factor $m_D^2=m_U^2m_{D\setminus U}^2$
and use the fifth law; equality uses doubled idempotence.

Expand a term $t$ as a sum of word occurrences and let $\A$
be the minimal members of their support family.
The preceding reductions give $2t\eq B_\A$. Since
$t\eq3t=2t+t$, we may retain $B_\A$ while reducing the
remaining words. It absorbs $2w$ for every word with support
in $\K(\A)$.

Normalize the positive exponents of a word $w$ to one or two.
Let $D=c(w)$ and let $L(w)$ be its variables with exponent one.
Repeatedly splitting a product by the sixth law of
\eqref{f51:eq:fb12} reduces it, in the retained background, to
\begin{equation}\label{f51:eq:wordexpansion}
B_\A+w\eq B_\A+m_D^2+\sum_{x\in L(w)}\delta_{D,x}.
\end{equation}
For clarity, the induction separates an exponent-one variable
from the remaining factors. In $uv\eq u^2v+uv^2+u^2v^2$,
the first two terms retain the respective exponent-one variables
of $v$ and $u$, and the third is a square word. Repeated square
or marked terms occur in pairs and are absorbed by $B_\A$.
The case with no exponent-one variables is a square word.
For one such variable, the right side of
\eqref{f51:eq:wordexpansion} is $v_{D,x}+2m_D^2$, equivalent
to $v_{D,x}$ because its doubled part is absorbed.
This proves the induction, including these boundary cases.

All coefficients can now be reduced modulo two in the retained
background: $B_\A+2w\eq B_\A$ for the words in question.
If $(D,x)$ is not a boundary pair, then
$D\setminus\{x\}\in\K(\A)$, and this set is nonempty.
Put $y=m_{D\setminus\{x\}}^2$ in the last law of
\eqref{f51:eq:fb12}. Since $B_\A$ absorbs $2y$, this replaces
$v_{D,x}$ by $m_D^2$ in the background. Consequently
$\delta_{D,x}$ disappears there.
Only the boundary pairs remain, giving
\eqref{f51:eq:periodicnormal}. The square support set $\R$
is the parity of the original word-support multiplicities.
The set $\E$ is the parity of the exponent-one occurrences
at eligible pairs. This is a reduction by finitely many identities
for every finite $X$.

We next separate the parameters by evaluations.
The map $h(z)=2z$ has image $\{0,d\}$. Assign the variables
of $Y\subseteq X$ the value $d$ and all others $0$.
Then $h(t)$ is $d$ precisely when $Y$ contains a member
of $\A$. Thus the function of $t$ determines $\K(\A)$
and its minimal antichain $\A$.

The set $\{0,e,d\}$ is a subsemiring of $P$, and
\[
\chi(0)=\chi(d)=0,\qquad\chi(e)=1
\]
defines a homomorphism onto the two-element field.
On assignments from $\{0,e\}$, the background and every
$\delta_{D,x}$ have $\chi$-value zero. Hence the
$\chi$-value of \eqref{f51:eq:periodicnormal} is the
squarefree Boolean polynomial
\[
\sum_{D\in\R}\prod_{x\in D}z_x\quad\text{over }\FF.
\]
Evaluation at $Y$ sums the coefficients on subsets of $Y$;
Lemma~\ref{lem:finite-set-recovery}(i) therefore recovers $\R$.

Finally fix $x\in X$. Assign $x=a$, assign variables of
$Y\subseteq X\setminus\{x\}$ the value $e$, and assign
all others $0$. Restrict to $Y\notin\K(\A)$.
Then the background and all square terms vanish.
Marked terms indexed by a variable other than $x$ also
vanish: either they contain $x^2$, or their support would
have to lie in $Y\in\K(\A)$ to survive.
The output is therefore $a$ or $0$ according to
\begin{equation}\label{f51:eq:boundaryseparation}
\sum_{\substack{(D,x)\in\E\\D\setminus\{x\}\subseteq Y}}1
\pmod2.
\end{equation}
The allowed $Y$ form a downward-closed family, and each
$D\setminus\{x\}$ appearing in this sum belongs to that
family by the boundary condition.
Lemma~\ref{lem:finite-set-recovery}(i), applied to this family,
recovers every coefficient in \eqref{f51:eq:boundaryseparation}.
Repeating for all $x$ determines $\E$.

Every permitted parameter set is realized by its displayed
normal form. Its background has minimal supports $\A$;
its square coefficients are $\R$; and each marked pair
contributes exactly its chosen boundary coefficient.
The separation argument and equational normalization prove
the assertions.
\end{proof}

\subsubsection{A basis using variables that occur once}
Let $\Gamma_Q$ consist of $\SR$ and the following nine laws:
\begin{equation}\label{f51:eq:fb14}
\begin{aligned}
xy&\eq yx,&x^3&\eq x^2,\\
3x&\eq2x,&2x&\eq2x^2,\\
x+xy&\eq x^2+x^2y,\\
xy+xz&\eq x^2y+x^2z,\\
x^2+y&\eq2x^2+y,\\
x^2y+z&\eq x^2y+z+x^2,\\
x^2+y^2&\eq2x^2y^2.
\end{aligned}
\end{equation}

There are two types of normal forms on a finite nonempty $X$.
The first is a single word with nonempty support $C\subseteq X$
and a chosen set $L\subseteq C$:
\begin{equation}\label{f51:eq:singleword}
W(C,L)=m_Lm_{C\setminus L}^2.
\end{equation}
The exponents are one on $L$ and two elsewhere.
For the second, choose a nonempty $C\subseteq X$, a subset
$R\subseteq C$, a partition $\pi$ of $L=C\setminus R$,
and a context $D_U\subseteq R$ for each block $U\in\pi$.
If $R=\varnothing$, require $|\pi|\ge2$.
If $R\ne\varnothing$, the partition may be empty when
$L=\varnothing$. Put
\begin{equation}\label{f51:eq:occurrencenormal}
N_Q(C,R,\pi,D)=2m_R^2+\sum_{U\in\pi}m_Um_{D_U}^2,
\end{equation}
omitting the doubled background when $R=\varnothing$.

\begin{theorem}\label{f51:thm:occurrence}
$\Gamma_Q$ is a complete fourteen-identity basis of $Q$.
The forms \eqref{f51:eq:singleword} and
\eqref{f51:eq:occurrencenormal}, with the stated constraints,
are distinct representatives of all elements of the free
algebra of $\Var(Q)$ on $X$.
\end{theorem}
\begin{proof}
All nine laws hold in $Q$. The first four follow immediately
from the tables. For the two promotion laws, $x^2=x$ unless
$x=a$, and with $x=a$ both sides of either identity are
zero: two $a$ summands sum to zero, and any zero summand
is absorbing. The remaining three laws follow because a
nonzero square is $e$ or $d$. A sum with at least two
nonzero summands has value $d$, $a$, or $0$, and adding
an available nonzero square preserves it in the stated cases.
If that square is zero, an existing summand is already zero.

A polynomial with one word reduces to \eqref{f51:eq:singleword}
by commutativity and the power law. Now take a polynomial
with at least two word occurrences. Let $C$ be its support,
let $L$ contain the variables occurring exactly once in the
entire polynomial, and put $R=C\setminus L$.

We first make every occurrence of every variable in $R$
squared. An exponent already at least two is reduced to two.
When $x$ occurs in different words, the two promotion laws
apply after factoring out $x$. They cover the cases where
one word is bare $x$ and where both have nonempty remaining
factors. If both words are $x$, use $2x\eq2x^2$.
Use a word already containing $x^2$ as a partner to promote
further occurrences; the power law reduces any exponent three
that this creates. No support or variable of $L$ changes.

Suppose $R\ne\varnothing$. For each $r\in R$ choose an
original word containing $r^2$. If it has a nonempty remaining
factor, the eighth law of \eqref{f51:eq:fb14}, with another
original word as its other summand, inserts $r^2$.
Apply it twice. If the chosen word is just $r^2$, use
the seventh law twice instead. Original occurrences are
retained throughout, so the required other summand exists.
We have inserted two copies of $r^2$ for every $r\in R$.
The last law, followed by $3x\eq2x$, combines these
markers into $2m_R^2$.

Every original word containing no variable of $L$ is now
$m_D^2$ with nonempty $D\subseteq R$. It is absorbed by
the background: apply the last law to one copy of $m_R^2$
and $m_D^2$, reduce powers using $D\subseteq R$, and
then reduce three copies of $m_R^2$ to two.
Every remaining original word contains a nonempty block
$U$ of variables of $L$ and has the form $m_Um_{D_U}^2$
for some $D_U\subseteq R$.
The blocks partition $L$, since each of its variables
originally occurred exactly once.
If $R=\varnothing$, all original words are disjoint
squarefree words and already have this form; the partition
has at least two blocks. Thus the two cases give exactly
\eqref{f51:eq:occurrencenormal}.

To prove separation, assigning all variables $e$ distinguishes
the types: a single word gives $e$, and every form of the
second type gives $d$. Assigning one variable $0$ and all
others $e$ gives zero precisely when that variable lies
in $C$, so $C$ is determined.
Next assign $x=a$ and all other variables $e$.
For $x\in C$ the result is $a$ precisely when $x\in L$,
and is zero otherwise. This determines $L$ and hence $R$.

For a form of the second type and $x\in L$, assign
$x=a$, another variable $y\ne x$ the value $d$,
and all others $e$. The result is zero precisely when
$y$ belongs to the support of the unique word containing
$x$; otherwise it is $a$. Indeed, $ad=0$, whereas a
square or another word with value $d$ is additively
absorbed by the sole $a$ summand.
These tests recover that support as
$U_x\cup D_{U_x}$, and intersection with $L$ and $R$
recovers the block $U_x$ and its context. All parameters
are therefore determined.

Every allowed form realizes its stated parameters: the
background repeats exactly the variables of $R$, and
the partition blocks contain each variable of $L$ once.
Together with the reductions above, this proves the
complete basis and free-algebra assertions.
\end{proof}

\endgroup
\section{Terminal, initial and doubled profiles}\label{sec:terminal-forms}
We now record the information at the ends of words and the effect
of doubling their coefficients. For each family, the proposed identities
first reduce a polynomial to these data, and the subsequent evaluations
show that no further identifications are valid.

\subsection{Terminal support profiles}\label{sec:family52}
\begingroup
\providecommand{\eq}{}\renewcommand{\eq}{\approx}
\providecommand{\SR}{}\renewcommand{\SR}{\mathsf{SR}}
\providecommand{\Var}{}\renewcommand{\Var}{\operatorname{Var}}
\providecommand{\FF}{}\renewcommand{\FF}{\mathbb F_2}
\providecommand{\A}{}\renewcommand{\A}{\mathcal A}
\providecommand{\K}{}\renewcommand{\K}{\mathcal K}
\providecommand{\E}{}\renewcommand{\E}{\mathcal E}
\providecommand{\R}{}\renewcommand{\R}{\mathcal R}

\emph{Generators (basis size):} $P=A_{819}$ (12), $P^{\mathrm{op}}=A_{1447}$ (12), $Q=A_{822}$ (15), $Q^{\mathrm{op}}=A_{1450}$ (15), $I=A_{827}$ (13), $I^{\mathrm{op}}=A_{1455}$ (13).

Write $(0,e,d,a)$ for labels $(0,1,2,3)$. In $P,Q,I$, the element $0$
annihilates multiplication, and the other products are
\[
e^2=e,\quad ed=de=d^2=d,\quad ea=a,\quad
ae=ad=da=a^2=0.
\]
Thus $e$ is a left identity but not a right identity. In $P$, addition
has identity $0$ and satisfies
\[
e+e=d,\quad e+d=e,\quad d+d=d,\quad
a+a=0,\quad a+e=e,\quad a+d=d.
\]
In $Q$, the element $0$ is additively absorbing and
\[
e+e=e+d=d+d=d,\qquad a+a=0,\qquad a+e=a+d=a.
\]
In $I$, addition is the maximum operation on the chain
$e<d<a<0$. 

Fix a finite nonempty set $X$ of variables and an order on it. For a
nonempty $D\subseteq X$, let
\[
b_D=\prod_{x\in D}x^2,
\qquad v_{D,x}=\left(\prod_{y\in D\setminus\{x\}}y^2\right)x
\quad(x\in D).
\]
The prefix is omitted when $D=\{x\}$, so $v_{\{x\},x}=x$.

\subsubsection{A common word lemma}
All three bases will contain $\SR$ and
\begin{equation}\label{f52:eq:word}
xyz\eq yxz,\qquad x^2y\eq xy,\qquad xy^2\eq yx^2.
\end{equation}
The first identity permutes a prefix while retaining a nonempty suffix.
The second permits squaring or unsquaring a prefix variable.

\begin{lemma}\label{f52:lem:word}
Modulo \eqref{f52:eq:word} and multiplicative associativity, a word $w$
with support $D$ reduces to $v_{D,x}$ if its final variable $x$
occurs exactly once in $w$, and to $b_D$ otherwise.
Furthermore, products of square words are determined by their support:
$b_Ub_V\eq b_{U\cup V}$.
\end{lemma}
\begin{proof}
Putting $y=x$ in the second identity gives $x^3\eq x^2$,
and hence every power of exponent at least two reduces to $x^2$.
The first identity, inside contexts, allows all permutations of the
letters preceding the final letter. Group these prefix letters by
variable and use the second identity in both directions to retain
exactly two copies of each distinct prefix variable.
If the final variable does not occur in the prefix, this gives $v_{D,x}$.
If it does, move its prefix copies next to the final copy and reduce
the resulting cube to a square. All factors are now squares.
Finally
\[
x^2y^2\eq xy^2\eq yx^2\eq y^2x^2
\]
permutes square factors. Repeated factors collapse by the power law.
This gives $b_D$ and the product assertion.
\end{proof}

Every identity in \eqref{f52:eq:word} holds in the three multiplication
tables. The first two follow because every product with $a$ in a
nonfinal position is zero. In the last identity the two sides both
vanish if $x$ or $y$ is $0$ or $a$; on $\{e,d\}$ multiplication is
commutative and idempotent.

For later separation arguments, the values of the standard words are
useful. A square word $b_D$ gives $e$ when every variable of $D$ is $e$,
gives $d$ when their values are in $\{e,d\}$ with at least one $d$,
and gives $0$ otherwise. The word $v_{D,x}$ agrees with $b_D$ except
that it gives $a$ when $x=a$ and all variables of $D\setminus\{x\}$
are $e$. These descriptions also show that the standard words are
distinct functions. There are
\begin{equation}\label{f52:eq:wordcount}
M_n=2^n-1+n2^{n-1}
\end{equation}
such functions on $n$ variables.

\subsubsection{Periodic addition and support boundaries}
Let $\Gamma_P$ consist of $\SR$, \eqref{f52:eq:word}, and
\begin{equation}\label{f52:eq:periodic}
3x\eq x,\qquad 2x\eq2x^2,\qquad
2x+2xy\eq2x,\qquad 2y+yx\eq2y+yx^2.
\end{equation}
For a nonempty antichain $\A$ of nonempty subsets of $X$, define
\begin{equation}\label{f52:eq:boundary}
\begin{aligned}
\K(\A)&=\{D\subseteq X:(\exists U\in\A)\ U\subseteq D\},\\
\partial\K(\A)&=\{(D,x):D\in\K(\A),\ x\in D,
                              \ D\setminus\{x\}\notin\K(\A)\}.
\end{aligned}
\end{equation}
Choose any $\R\subseteq\K(\A)$ and $\E\subseteq\partial\K(\A)$.
Put $\delta_{D,x}=v_{D,x}+b_D$ and set
\begin{equation}\label{f52:eq:pnormal}
N_P(\A,\R,\E)=
\underbrace{\sum_{U\in\A}2b_U}_{B_\A}
+\sum_{D\in\R}b_D+\sum_{(D,x)\in\E}\delta_{D,x}.
\end{equation}

\begin{theorem}\label{f52:thm:periodic}
$\Gamma_P$ is a complete twelve-identity basis for $P$.
The forms \eqref{f52:eq:pnormal}, with the stated parameters, are distinct
representatives of the elements of the free algebra of $\Var(P)$ on $X$.
\end{theorem}
\begin{proof}
The first two laws of \eqref{f52:eq:periodic} follow from the addition
table and the fact that only $a$ changes on squaring. For doubled
absorption, if $x\in\{0,a\}$, both $2x$ and $2xy$ are $0$;
if $x\in\{e,d\}$, then $2x=d$ and $2xy\in\{0,d\}$.
For the last law, only $x=a$ needs examination. The only case where
$ya\ne0$ is $y=e$, and then $2y=d$ absorbs $a$. Thus the basis is sound.

Expand a term $t$ into a nonempty sum of word occurrences. Let $\A$
be the minimal supports of those occurrences. By distributivity,
$2x\eq2x^2$ can be used at any position of a doubled word. Lemma
\ref{f52:lem:word} then gives $2w\eq2b_{c(w)}$.
Since $3x\eq x$ implies $4x\eq2x$, doubled terms are additively
idempotent. If $U\subseteq D$, the third law gives
$2b_U+2b_D\eq2b_U$, factoring out $b_U$ for strict inclusion
and using doubled idempotence for equality.
Consequently $2t\eq B_\A$, and $t\eq3t\eq B_\A+t$.
This background absorbs $2w$ for every word whose support belongs
to $\K(\A)$.

Normalize every retained word by Lemma~\ref{f52:lem:word}. A square word
contributes $b_D$. A marked word can be written in the background as
$b_D+\delta_{D,x}$: the extra two copies of $b_D$ are absorbed.
All these coefficients can be reduced modulo two. For a nonboundary
pair, $D\setminus\{x\}\in\K(\A)$ is nonempty. Substituting
$y=b_{D\setminus\{x\}}$ in the last law of \eqref{f52:eq:periodic},
and inserting or deleting $2y$ using the background, gives
\[
B_\A+v_{D,x}\eq B_\A+b_D.
\]
Hence $\delta_{D,x}$ disappears. This proves reduction to
\eqref{f52:eq:pnormal} for every finite $X$. Specifically, $\R$ records
the parity of the original word-support multiplicities, and $\E$
records the parity of the words with each eligible unique terminal
variable.

We prove uniqueness by evaluations. Assign $d$ on a chosen subset
$Y\subseteq X$ and $0$ outside it. The output is $d$ exactly when
$Y$ contains some member of $\A$. This recovers $\K(\A)$ and
its minimal antichain $\A$.
The subsemiring $\{0,e,d\}$ admits a homomorphism to $\FF$ defined
by $\chi(0)=\chi(d)=0$ and $\chi(e)=1$. On assignments from
$\{0,e\}$, the background and every $\delta_{D,x}$ have zero
$\chi$-value, leaving the Boolean polynomial
\[
\sum_{D\in\R}\prod_{x\in D}z_x.
\]
Evaluation on $Y$ sums the coefficients on subsets of $Y$;
Lemma~\ref{lem:finite-set-recovery}(i) recovers $\R$.

Fix $x\in X$, assign $x=a$, assign $e$ on
$Y\subseteq X\setminus\{x\}$, and $0$ elsewhere, with
$Y\notin\K(\A)$. The background and square terms vanish.
A marked term distinguished by another variable either contains $x$
in its squared prefix, or would require its support to be contained
in $Y$ to survive; the latter is impossible because that support
belongs to $\K(\A)$. The remaining output is $a$ or $0$ according
to the parity of
\begin{equation}\label{f52:eq:boundarytest}
\{(D,x)\in\E:D\setminus\{x\}\subseteq Y\}.
\end{equation}
The permitted $Y$ form a downward-closed family. Every support
$D\setminus\{x\}$ in \eqref{f52:eq:boundarytest} belongs to this family.
Lemma~\ref{lem:finite-set-recovery}(i) therefore recovers every
coefficient, including the empty-subset coefficient. Repeating for
all $x$ recovers $\E$.

Every allowed parameter set is realized by its displayed term: the
background fixes its minimal supports, square terms fix $\R$, and
each $\delta_{D,x}$ changes exactly its selected boundary coefficient.
Completeness follows from Lemma~\ref{lem:normal-form-completeness}.
\end{proof}

\subsubsection{Nonperiodic addition and terminal contexts}
Let $\Gamma_Q$ consist of $\SR$, \eqref{f52:eq:word}, and
\begin{equation}\label{f52:eq:threshold}
\begin{aligned}
3x&\eq2x,&2x&\eq2x^2,\\
x+yx&\eq x^2+yx^2,&yx+zx&\eq yx^2+zx^2,\\
x^2+y&\eq2x^2+y,\\
x^2y+z&\eq x^2y+z+x^2,\\
x^2+y^2&\eq2x^2y^2.
\end{aligned}
\end{equation}
Let $\Gamma_I$ consist of $\SR$, \eqref{f52:eq:word}, and
\begin{equation}\label{f52:eq:idempotent}
\begin{aligned}
x+x&\eq x,&xy&\eq x^2+xy,\\
x^2+y^2&\eq x^2y^2,&yx+zx&\eq yzx,\\
x+yx&\eq yx.
\end{aligned}
\end{equation}
For disjoint $R,L\subseteq X$ and contexts $D_x\subseteq R$ $(x\in L)$,
define
\begin{equation}\label{f52:eq:qnormal}
N_Q(R,L,D)=2b_R+\sum_{x\in L}v_{D_x\cup\{x\},x}
\qquad(R\ne\varnothing\text{ or }|L|\ge2),
\end{equation}
and
\begin{equation}\label{f52:eq:inormal}
N_I(R,L,D)=b_R+\sum_{x\in L}v_{D_x\cup\{x\},x}
\qquad(R\cup L\ne\varnothing).
\end{equation}
Omit the background when $R=\varnothing$. For $Q$, also retain the
single-word forms $b_C$ and $v_{C,x}$ from Lemma~\ref{f52:lem:word}.

\begin{theorem}\label{f52:thm:threshold}\label{f52:thm:idempotent}
The systems $\Gamma_Q$ and $\Gamma_I$ are complete bases of sizes
fifteen and thirteen for $Q$ and $I$, respectively. The single-word
forms together with \eqref{f52:eq:qnormal} represent the free algebra
of $\Var(Q)$ on $X$ without repetition; \eqref{f52:eq:inormal} does
the same for $\Var(I)$.
\end{theorem}
\begin{proof}
\emph{Soundness.} In $Q$, the cap and doubling laws follow from
addition and squaring. In either promotion law, only $x=a$ changes
on squaring; both sides are then $0$. A square is $0$, $e$, or $d$:
the first absorbs addition, and the latter two have identical sums
with any other element. This proves duplication and exposure (a
nonzero $x^2y+z$ is $d$ or $a$). The last law follows on
$\{0,e,d\}$.
In $I$, prefix exposure holds because $xy\ge x^2$ in the additive
chain, and square merging holds on $\{0,e,d\}$. For tail merging,
$x=0$ gives $0$ on both sides; $x=a$ gives $a$ exactly when
$y=z=e$, and $0$ otherwise. For $x=e$, multiplication by $e$
sends $a$ to $0$ and fixes the other elements, so the maximum
of $ye,ze$ is $(yz)e$. For $x=d$, both sides are $d$ when
$y,z\in\{e,d\}$ and $0$ otherwise. The same four cases give
$x+yx=yx$. This verifies both bases.

\emph{Reduction for $Q$.} Expand a term into word occurrences.
A single word reduces by Lemma~\ref{f52:lem:word}. Otherwise let $C$
be the total support, let $L$ consist of variables occurring exactly
once in the polynomial, in final position, and put $R=C\setminus L$.
Normalize the words. If $r$ occurs in a prefix, permute it first
and expose $r^2$ from $r^2u$ with $u$ nonempty, using another original word
as $z$; apply exposure twice. If the chosen word is $r^2$ itself,
use duplication with another original word instead. If $r\in R$
has no prefix occurrence, it ends at least two word occurrences.
Promote a pair by the third or fourth law of
\eqref{f52:eq:threshold}, using doubling when both are bare $r$,
and then expose the markers as above. The two promoted occurrences
still provide the required other summand. Thus two markers $r^2$
can be inserted for every $r\in R$.

These markers promote all remaining terminal occurrences of $r$:
for a nonempty prefix $u$, the fourth law with $z=r$ gives
\begin{equation}\label{f52:eq:markerpromotion}
ur+r^2\eq ur^2+r^3\eq ur^2+r^2,
\end{equation}
and the third law with $y=r$ gives $r+r^2\eq2r^2$.
Each step retains a marker. Every resulting word with terminal in
$R$ is now a square word; the others are the unique marked words
indexed by $L$, with contexts $D_x\subseteq R$.
The last law, the square-product rule, and $3x\eq2x$ combine the
markers into $2b_R$. This absorbs each square word $b_D$ with
$\varnothing\ne D\subseteq R$: merge $b_D$ with one copy of
$b_R$, then cap the resulting three copies at two. This gives
\eqref{f52:eq:qnormal}. If $R=\varnothing$, the original words
were distinct bare variables, so $|L|\ge2$.

\emph{Reduction for $I$.} After word normalization, let $R$ be the
union of the square-word supports and marked prefixes, let $C$ be
the total support, and put $L=C\setminus R$. Prefix exposure,
with idempotence for a bare square, inserts all markers $r^2$
while retaining the original summands. Square merging combines
them into $b_R$ and gives $b_R+b_D\eq b_R$ for every nonempty
$D\subseteq R$.
For a terminal $r\in R$, insert a marker from the background.
Tail merging gives $ur+r^2\eq ur^2$ for a nonempty prefix $u$;
bare absorption gives $r+r^2\eq r^2$. The resulting square word
is absorbed by the retained background, as are the original square
words. For $x\in L$, merge nonempty prefixes by $ux+vx\eq uvx$
and absorb a bare $x$ if a prefixed word is present. Idempotence
removes repeated bare copies. Exactly one marked word remains
for each $x$, with $D_x$ the union of its prefixes. This is
\eqref{f52:eq:inormal}; if $R=\varnothing$, all words were bare.

\emph{Separation.} The same tests recover the parameters for both
algebras. Assign one variable $0$ and all others $e$; the output
is $0$ exactly for variables in the total support $C$. For $x\in C$,
assign $x=a$ and all others $e$: the output is $a$ for $x\in L$
and $0$ for $x\in R$. For $x\in L$, $y\in R$, assign
$x=a$, $y=d$ and all others $e$. The value is $0$ exactly when
$y\in D_x$, and otherwise $a$, which absorbs the other nonzero
summands. These tests recover $R,L$ and all $D_x$.
For $Q$, the all-$e$ assignment also distinguishes a single word
(value $e$) from a sum form (value $d$); the words themselves
were separated in Lemma~\ref{f52:lem:word}. Each allowed form
realizes its stated data, so Lemma~\ref{lem:normal-form-completeness}
proves both assertions.
\end{proof}

\endgroup
\subsection{Initial support profiles}\label{sec:family53}
\begingroup
\providecommand{\eq}{}\renewcommand{\eq}{\approx}
\providecommand{\SR}{}\renewcommand{\SR}{\mathsf{SR}}
\providecommand{\Var}{}\renewcommand{\Var}{\operatorname{Var}}
\providecommand{\FF}{}\renewcommand{\FF}{\mathbb F_2}
\providecommand{\K}{}\renewcommand{\K}{\mathcal K}
\providecommand{\E}{}\renewcommand{\E}{\mathcal E}
\providecommand{\R}{}\renewcommand{\R}{\mathcal R}

\emph{Generators (basis size):} $P=A_{1491}$ (11), $P'=A_{1652}$ (11), $Q=A_{1493}$ (13), $Q'=A_{1653}$ (13), $I=A_{1498}$ (14), $I'=A_{1658}$ (14).

For the source algebras $P,Q,I$, identify the labels $(0,1,2,3)$ with
$(0,e,d,a)$. Define a map $h$ on their common carrier by
\[
h(0)=h(a)=0,\qquad h(e)=h(d)=d.
\]
Their multiplication is exactly
\begin{equation}\label{f53:eq:mu}
uv=\begin{cases}u,&v=e,\\h(u),&v\ne e.\end{cases}
\end{equation}
In particular, $e$ is a right identity, while $0$ and $d$ are
left-zero elements. The element $0$ is not a two-sided multiplicative
zero: $e0=d$. 

The additions of $P,Q,I$ are those in Section~\ref{sec:family52},
under the same labels $(0,e,d,a)$.
For $P$ and $Q$, one has $h(u)=u+u$.

\subsubsection{Word reduction and notation}
Fix a finite nonempty variable set $X$ with a chosen order. For
$x\in X$ and $D\subseteq X\setminus\{x\}$, write
\[
u_{x,D}=x\prod_{y\in D}y,\qquad
s_{x,D}=x^2\prod_{y\in D}y.
\]
For $D=\varnothing$, these mean $x$ and $x^2$, respectively.

All three bases contain the word identities
\begin{equation}\label{f53:eq:word}
xyz\eq xzy,\qquad xy^2\eq xy.
\end{equation}
\begin{lemma}\label{f53:lem:word}
Modulo \eqref{f53:eq:word} and associativity, a word $w$ with first
variable $x$ reduces to $u_{x,c(w)\setminus\{x\}}$ if $x$ occurs
once in $w$, and to $s_{x,c(w)\setminus\{x\}}$ otherwise.
These standard words are distinct functions in each of $P,Q,I$.
There are $M_n=n2^n$ of them on $n\ge1$ variables.
\end{lemma}
\begin{proof}
The first law permutes the letters after the initial letter, using
contexts with nonempty prefixes. Group the suffix letters and use
the second law to reduce each positive suffix exponent to one.
If $x$ belongs to the suffix, place it first in that suffix to obtain
$x^2$. The substitution $y=x$ also gives $x^3\eq x^2$.
This proves the reductions. Both laws hold in \eqref{f53:eq:mu}, because
a product has its initial value when every suffix letter is $e$;
otherwise it takes the value of $h$ at that initial value.

For distinction, assign a selected variable $0$ and all others $e$.
The result is $0$ precisely when it is the initial variable; assigning
$0$ to a variable in the suffix but not initially gives $d$ instead.
With the initial variable $a$ and
all others $e$, the result is $a$ for $u_{x,D}$ and $0$ for $s_{x,D}$.
With all variables $e$ except a selected noninitial variable set to
$d$, the result is $d$ precisely when that variable lies in $D$.
These tests recover the initial variable, the type, and the context.
The count is $n\cdot2\cdot2^{n-1}=n2^n$.
\end{proof}

The values of $s_{x,D}$ are $0$ when $x\in\{0,a\}$, $d$ when
$x=d$, and, when $x=e$, either $e$ if all variables of $D$ are $e$
or $d$ otherwise. The function $u_{x,D}$ differs only at assignments
with $x=a$ and all variables of $D$ equal to $e$, when its value is $a$.

\subsubsection{Periodic addition and initial supports}
Let $\Gamma_P$ consist of $\SR$, \eqref{f53:eq:word}, and
\begin{equation}\label{f53:eq:periodic}
\begin{aligned}
3x&\eq x,&2xy&\eq2x,\\
2x+2y+x^2y&\eq2x+2y+y^2x,\\
2y+xy&\eq2y+x^2y.
\end{aligned}
\end{equation}
For a nonempty $H\subseteq X$, define
\[
\K_H=\{C\subseteq X:C\cap H\ne\varnothing\},\qquad
\E_H=\{(C,x):C\subseteq X,\ C\cap H=\{x\}\}.
\]
For $C\in\K_H$, let $j_H(C)$ be its least member in $H$, and put
\[
q_{H,C}=s_{j_H(C),C\setminus\{j_H(C)\}},\qquad
\delta_{C,x}=u_{x,C\setminus\{x\}}+s_{x,C\setminus\{x\}}.
\]
Choose arbitrary $\R\subseteq\K_H$ and $\E\subseteq\E_H$. The
normal form is
\begin{equation}\label{f53:eq:pnormal}
N_P(H,\R,\E)=
\underbrace{\sum_{x\in H}2x}_{B_H}
+\sum_{C\in\R}q_{H,C}+\sum_{(C,x)\in\E}\delta_{C,x}.
\end{equation}

\begin{theorem}\label{f53:thm:periodic}
$\Gamma_P$ is a complete eleven-identity basis for $P$. The forms
\eqref{f53:eq:pnormal} are distinct representatives of all elements of
the free algebra of $\Var(P)$ on $X$.
\end{theorem}
\begin{proof}
Soundness of $3x\eq x$ follows from the addition table, and
$2xy\eq2x$ follows from $h(xy)=h(x)$. For square transfer, if both
$x,y\in\{0,a\}$, all relevant terms vanish. Otherwise the doubled
background is $d$. Both square products are $e$ precisely when
$x=y=e$; in every other case they belong to $\{0,d\}$ and their
addition to the background gives $d$. For projection, only $x=a$
can change on squaring. The only case with $xy\ne0$ is $y=e$,
whose doubled value $d$ absorbs $a$.

Expand a term $t$ and let $H$ be the set of initial variables of its
word occurrences. Distributivity and $2xy\eq2x$ imply $2w\eq2x$
for every word beginning in $x$. Since $3x\eq x$ implies
$4x\eq2x$, one obtains $2t\eq B_H$ and $t\eq B_H+t$.
The background absorbs two copies of every word beginning in $H$.

Normalize words by Lemma~\ref{f53:lem:word}. A word $u_{x,D}$ can be
replaced, in the background, by $s_{x,D}+\delta_{D\cup\{x\},x}$,
because the added two copies of $s_{x,D}$ are absorbed.
Suppose $y\in H\cap D$. Order the suffix with $y$ first and apply
the last law of \eqref{f53:eq:periodic}. If there is a remaining suffix
word $v$, multiply that identity on the right by $v$; its doubled
summand remains $2y$ by the second law. Thus
\[
B_H+u_{x,D}\eq B_H+s_{x,D}.
\]
Consequently a delta term disappears whenever its support contains
an initial variable other than its own. The only possible surviving
deltas are those in $\E_H$.

If $x,y\in H\cap C$, multiply the square-transfer law by the word
on $C\setminus\{x,y\}$, or use it without a context if this set is
empty. Again doubled summands reduce to $2x$ and $2y$. Suffix
permutation gives
\[
B_H+s_{x,C\setminus\{x\}}\eq
B_H+s_{y,C\setminus\{y\}}.
\]
Hence every square word of support $C$ can be replaced by $q_{H,C}$.
Reduce all remaining coefficients modulo two in the background.
This yields \eqref{f53:eq:pnormal}. The set $\R$ is the parity of the
original support multiplicities. The set $\E$ is the parity of
the uniquely initial words at eligible pairs $(C,x)$.

For separation, assign $d$ to one variable and $0$ to all others.
The output is $d$ precisely when that variable belongs to $H$.
Indeed, words with other initial variables vanish, and the doubled
summand $2d=d$ ensures survival for an initial variable in $H$.
This recovers $H$.
The map $\chi:P\to\FF$ defined by $\chi(e)=1$ and
$\chi(0)=\chi(d)=\chi(a)=0$ is a homomorphism, as follows
directly from the two tables. On assignments from $\{0,e\}$,
the $\chi$-value of the normal form is
\[
\sum_{C\in\R}\prod_{x\in C}z_x.
\]
Lemma~\ref{lem:finite-set-recovery}(i) recovers $\R$ from these
subset sums.

Fix $x\in H$. Assign $x=a$, set every other variable of $H$ to $0$,
assign $e$ on $Y\subseteq X\setminus H$, and set the remaining
variables to $0$. The background and all square words vanish.
Deltas with other initial variables vanish as well. The output is
$a$ or $0$ according to the parity of
\[
\{(C,x)\in\E:C\setminus\{x\}\subseteq Y\}.
\]
Lemma~\ref{lem:finite-set-recovery}(i), on $2^{X\setminus H}$,
recovers these coefficients, including the empty support; varying
$x$ recovers $\E$.
Each displayed form realizes its parameters: its doubled background
fixes $H$, the square summands fix support parity, and each delta
affects exactly its indicated initial coefficient. Reduction and
separation prove completeness for every finite $X$.
\end{proof}

\subsubsection{Nonidempotent addition and erasure of square contexts}
Let $\Gamma_Q$ consist of $\SR$, \eqref{f53:eq:word}, and
\begin{equation}\label{f53:eq:threshold}
\begin{aligned}
3x&\eq2x,&2xy&\eq2x,\\
x+xy&\eq x^2+x^2y,&xy+xz&\eq x^2y+x^2z,\\
x^2y+z&\eq x^2+z,&x^2+y&\eq2x^2+y.
\end{aligned}
\end{equation}
One type of normal form is a single standard word from
Lemma~\ref{f53:lem:word}. For the sum type, choose disjoint sets
$R,L\subseteq X$, with $R\ne\varnothing$ or $|L|\ge2$, and
arbitrary contexts $D_x\subseteq X\setminus\{x\}$ for $x\in L$.
Set
\begin{equation}\label{f53:eq:qnormal}
N_Q(R,L,D)=\sum_{r\in R}2r^2+\sum_{x\in L}u_{x,D_x}.
\end{equation}
Context variables need not be initial variables, and may be initial
variables of other summands.

\begin{theorem}\label{f53:thm:threshold}
$\Gamma_Q$ is a complete thirteen-identity basis for $Q$. Its
standard words and forms \eqref{f53:eq:qnormal} give distinct
representatives of the free algebra of $\Var(Q)$ on $X$.
\end{theorem}
\begin{proof}
The cap law holds in the additive table; doubling a product gives
$h(xy)=h(x)$. For either promotion law, only $x=a$ needs attention.
Every left-hand summand is then $0$ or $a$, and their sum is $0$,
as is the right side. For square-context erasure, $x\in\{0,a\}$
makes both sides zero. Otherwise $x^2$ and $x^2y$ differ, if at
all, only between $e$ and $d$, which are indistinguishable after
addition of any other element in $Q$. The same observation proves
duplication of a square in a sum. Thus the basis is sound.

A polynomial with one word reduces by Lemma~\ref{f53:lem:word}.
Take a polynomial with at least two word occurrences. Let $H$ be
its initial-variable set. Let $L$ comprise those $x\in H$ occurring
initially in exactly one word and nowhere else in that word; put
$R=H\setminus L$. Occurrences of $x$ in words with a different
initial variable do not exclude $x$ from $L$.

Normalize each word. For $r\in R$, if a word already has repeated
initial variable $r$, it is $r^2v$ or $r^2$. In the former case,
erase $v$ using the fifth law of \eqref{f53:eq:threshold} and another
retained word as $z$. In the latter case no erasure is needed.
If no such word exists, at least two words begin in $r$.
The two promotion laws turn a chosen pair into square-initial words;
when both are bare $r$, use $2r\eq2r^2$, obtained by setting $y=r$
in the doubling law. Then erase their contexts.

The available $r^2$ promotes every further word beginning in $r$.
For a nonempty context $v$, the fourth law with $z=r$ gives
$rv+r^2\eq r^2v+r^3\eq r^2v+r^2$; erase $v$ afterward.
For a bare $r$, the third law with $y=r$ gives
$r+r^2\eq2r^2$. These operations always retain at least two
word occurrences in the full polynomial, so the erasure step has
the required other summand. After treating all $r\in R$, only
bare squares with these initials and the unique words $u_{x,D_x}$
for $x\in L$ remain. Duplication of a square and $3x\eq2x$
make the coefficient of each $r^2$ exactly two. This is
\eqref{f53:eq:qnormal}. If $R$ is empty, the initial variables of
the original summands were distinct, and there are at least two.

For separation, the all-$e$ assignment gives $e$ on a word and
$d$ on a sum form. For a sum form, assigning a selected variable
$0$ and all others $e$ gives $0$ exactly when the selected variable
lies in $H=R\cup L$; a purely contextual $0$ can only change an
$e$ summand to $d$, leaving the sum $d$. This recovers $H$.
Assigning $x=a$ and all others $e$ gives $a$ for $x\in L$ and
$0$ for $x\in R$. For $x\in L$ and $y\ne x$, assign $x=a$,
$y=d$ and all others $e$. The sum is $0$ exactly when $y\in D_x$;
otherwise the sole $a$ summand absorbs the other $e$ or $d$ values.
Thus every context is recovered. The word forms were separated in
Lemma~\ref{f53:lem:word}, and every permitted sum form realizes the
specified data. The normal-form reduction proves completeness.
\end{proof}

\subsubsection{Idempotent addition and transfer of square contexts}
Let $\Gamma_I$ consist of $\SR$, \eqref{f53:eq:word}, and
\begin{equation}\label{f53:eq:idempotent}
\begin{aligned}
x+x&\eq x,&xy+xz&\eq xyz,&x+xy&\eq xy,\\
x^2y+z^2&\eq x^2+z^2y,\\
x^2y+y&\eq x^2+y,\\
x^2y+zy&\eq x^2+zy,\\
x^2y+yz&\eq x^2+yz.
\end{aligned}
\end{equation}
Choose disjoint $R,L\subseteq X$ with $H=R\cup L\ne\varnothing$,
and any $D_x\subseteq X\setminus\{x\}$ for $x\in L$. Put
\[
C_0=H\cup\bigcup_{x\in L}D_x.
\]
If $R=\varnothing$, set $K=\varnothing$. If $R\ne\varnothing$,
choose any $K\subseteq X\setminus C_0$ and let $r_0$ be the least
member of $R$. The normal form is
\begin{equation}\label{f53:eq:inormal}
N_I(R,L,D,K)=\sum_{r\in R}r^2+
\sum_{x\in L}u_{x,D_x}+s_{r_0,K},
\end{equation}
where the final summand is included only when $K\ne\varnothing$.
This convention also omits $r_0$ when $R=\varnothing$.

\begin{theorem}\label{f53:thm:idempotent}
$\Gamma_I$ is a complete fourteen-identity basis for $I$.
The forms \eqref{f53:eq:inormal} are distinct representatives of the
free algebra of $\Var(I)$ on $X$.
\end{theorem}
\begin{proof}
We first check soundness. Multiplication \eqref{f53:eq:mu} makes
$xy+xz\eq xyz$ valid: for $x=e$ both sides are $e$ precisely
when $y=z=e$, and otherwise $d$; for $x=a$ they are $a$ precisely
in that same case and otherwise $0$; for $x\in\{0,d\}$ they
equal $x$. Also $xy$ is at least $x$ in the additive chain,
giving bare absorption. For square-context transfer, if either
$x$ or $z$ is $0$ or $a$, a square summand makes both sides $0$.
If $x,z\in\{e,d\}$ and at least one is $d$, both sides are $d$;
if $x=z=e$, the two sides coincide by symmetry.
For each of the last three laws, $y=e$ gives $x^2y=x^2$.
If $y\ne e$, the companion $y$, $zy$, or $yz$ is at least $d$.
It therefore makes the possible change from $e$ to $d$ in
$x^2y$ irrelevant; if $x^2=0$, both sides already vanish.
This verifies all identities.

Expand a term and group its word occurrences by initial variable.
Repeatedly applying $xy+xz\eq xyz$ merges all nonempty contexts
of a fixed initial variable by union, using Lemma~\ref{f53:lem:word}
to remove repeated suffix letters. A bare initial variable is
absorbed by any nonempty context; repeated bare copies collapse
by idempotence. There is now one word for each initial variable.
Let $R$ consist of the initials that also occur in their merged
suffix, and let $L$ be the other initials. The words indexed by
$L$ have the form $u_{x,D_x}$ and will be retained throughout.
The others have the form $s_{r,E_r}$.

From each $s_{r,E_r}$ insert a bare $r^2$ using bare absorption
in reverse, or duplicate an existing $r^2$ by idempotence.
Using $xy+xz\eq xyz$ in reverse with $x=r^2$, split a nonempty
square context into individual contributions $r^2y$ for $y\in E_r$.
Thus all bare markers $r^2$ are available, and all remaining square
contexts occur as single-letter contributions.

If $R\ne\varnothing$, use the fourth law to move each contribution
to the chosen initial $r_0$:
\[
r^2y+r_0^2\eq r^2+r_0^2y.
\]
Idempotence lets the bare markers be retained while doing this.
We can now erase every $r_0^2y$ for which $y\in C_0$.
If $y\in R$, a retained $y^2$ serves as companion in the last law
with $z=y$. If $y\in L$, use the fifth law when its word is bare
$y$, or the last law with its suffix as $z$ when it is not bare.
If $y$ is in the context of another retained live word, suffix
permutation puts it last, making that word $vy$ with $v$ nonempty;
the sixth law applies with $z=v$. These cases cover $C_0$ and keep
all the required companions. A contribution with $y=r_0$ also
reduces directly by $r_0^3\eq r_0^2$.

Only the distinct contributions $r_0^2y$ with
$y\in K\subseteq X\setminus C_0$ remain. Merge them back into
$s_{r_0,K}$ when $K$ is nonempty. The resulting expression is
\eqref{f53:eq:inormal}. If $R=\varnothing$, no square word or square
context was present, so the construction already consists solely
of the live words, with $K=\varnothing$.

The single-$0$, single-$a$, and paired $a,d$ assignments in
Theorem~\ref{f53:thm:threshold} recover $H,L$ and the contexts
$D_x$ here as well: a contextual $0$ gives $d$, not $0$, and
a live $a$ survives exactly when its context contains only $e$.
Thus $R$ and $C_0$ are also known. For $y\in X\setminus C_0$,
assigning $y=d$ and all others $e$ gives $d$ precisely when
$y\in K$, and $e$ otherwise. This recovers $K$. Each allowed
form realizes its parameters, completing the proof.
\end{proof}

\endgroup
\subsection{Doubled flat and chain backgrounds}\label{sec:family54}
\begingroup
\providecommand{\eq}{}\renewcommand{\eq}{\approx}
\providecommand{\SR}{}\renewcommand{\SR}{\mathsf{SR}}
\providecommand{\Var}{}\renewcommand{\Var}{\operatorname{Var}}
\providecommand{\FF}{}\renewcommand{\FF}{\mathbb F_2}
\providecommand{\R}{}\renewcommand{\R}{\mathcal R}

\emph{Generators:} $P=A_{1421}$, $\widetilde P=A_{1904}$, $P^{\mathrm{op}}=A_{709}$, $\widetilde P^{\mathrm{op}}=A_{994}$, $Q=A_{1423}$, $\widetilde Q=A_{1914}$, $Q^{\mathrm{op}}=A_{711}$, $\widetilde Q^{\mathrm{op}}=A_{1004}$.

Write $(0,e,f,a)$ for labels $(0,1,2,3)$. In the four source algebras
$P,\widetilde P,Q,\widetilde Q$, the element $0$ is multiplicatively
absorbing, $e^2=e$, $f^2=f$, $ae=a$, and $ea=fa=a^2=0$.
For $P,Q$, the additional products are $ef=fe=af=0$.
For $\widetilde P,\widetilde Q$, they are $ef=fe=e$ and $af=a$.
Every source contains the subsemiring $K=\{0,e,a\}$, whose multiplication
is $uv=u$ when $v=e$ and $uv=0$ otherwise. For $P,\widetilde P$,
addition on $K$ is idempotent with $e+a=0$ and absorbing element $0$;
write $K_{\mathrm f}$ for this semiring. For $Q,\widetilde Q$,
addition on $K$ is the maximum operation on $e<a<0$; write
$K_{\mathrm c}$.

There is also a copy of $\FF$, considered only in the signature
$(+,\cdot)$: it is $\{0,f\}$ in $P,Q$ and $\{e,f\}$ in the two
tilded algebras. In each copy $f$ corresponds to $1$ and the other
element corresponds to $0$. 

\subsubsection{Words and doubled terms}
All source bases contain $\SR$ and
\begin{equation}\label{f54:eq:word}
xyz\eq xzy,\qquad xy^2\eq xy,\qquad x^2y\eq y^2x.
\end{equation}
Fix a finite nonempty variable set $X$ and an order on it. For nonempty
$D\subseteq X$ and $x\notin T$, define
\[
b_D=\prod_{y\in D}y^2,\qquad
u_{x,T}=x\prod_{y\in T}y.
\]
For $T=\varnothing$, the latter word is $x$.

\begin{lemma}\label{f54:lem:words}
The identities \eqref{f54:eq:word} reduce a word $w$ to $u_{x,c(w)-\{x\}}$
if its initial variable $x$ does not occur again, and to $b_{c(w)}$
otherwise. In particular $w^2\eq b_{c(w)}$ and
$b_Db_E\eq b_{D\cup E}$. The standard words are distinct functions
in every source algebra. Their number on $n$ variables is
\begin{equation}\label{f54:eq:wordcount}
M_n=2^n-1+n2^{n-1}.
\end{equation}
\end{lemma}
\begin{proof}
The identities \eqref{f54:eq:word} are the reversals of
\eqref{f52:eq:word}. Reverse the derivation in
Lemma~\ref{f52:lem:word}, and use $xy^2\eq xy$ to unsquare
the marked suffix. This gives the reductions and square-product rule.

All laws hold in the displayed multiplication tables: suffix products
are insensitive to order and repetition, and a squared initial value
removes the exceptional behavior of $a$. Different supports are
separated on the field copy. For a fixed support $D$, assign $a$ to
$x\in D$ and $e$ to its other variables. Precisely the marked word
$u_{x,D-\{x\}}$ gives $a$; all other standard words with support $D$
give $0$. This proves distinctness and the count.
\end{proof}

Every source satisfies $3x\eq x$. By Lemma~\ref{f49:lem:double},
$E(t)=2t$ is an endomorphism with additively idempotent image,
here exactly $K$, and $t\eq E(t)+t$. We call a canonical
representative of $E(t)$ the doubled background of $t$.

\subsubsection{The flat additive background}
Let $\Gamma_{\mathrm f}$ consist of $\SR$, \eqref{f54:eq:word}, and
\begin{equation}\label{f54:eq:flat}
\begin{aligned}
3x&\eq x,\\
2xy+2z&\eq2xy+2zy,\\
2x^2+y&\eq2x^2+y^2.
\end{aligned}
\end{equation}
For a polynomial $t$, expand it as a nonempty sum of word occurrences.
Let $C$ be its content, $H$ the set of initial variables, and $T$ the
set of variables occurring in noninitial positions. Thus $C=H\cup T$.
Repeated word occurrences are included in this definition. There are
two background types:
\begin{equation}\label{f54:eq:flatbg}
B(C,0)=2b_C,\qquad
B(C,H)=\sum_{x\in H}2u_{x,C-H}\quad
(\varnothing\ne H\subseteq C).
\end{equation}
The symbol $0$ in the first parameter list is a type marker, not a
constant in a term.

\begin{lemma}\label{f54:lem:flatbg}
Modulo $\Gamma_{\mathrm f}$, the doubled term $E(t)$ reduces to
$B(C,H)$ if $H\cap T=\varnothing$, and to $B(C,0)$ otherwise.
For $B(C,0)$, every $2b_D$ with $\varnothing\ne D\subseteq C$ is
absorbed additively. For $B(C,H)$ with $H\ne\varnothing$, every
$2u_{x,D}$ with $x\in H$ and $D\subseteq C-H$ is absorbed.
\end{lemma}
\begin{proof}
The middle law of \eqref{f54:eq:flat} propagates any existing suffix letter
to another doubled word. Indeed, expose a suffix letter $y$ at the end
of a source word, writing it as $vy$ with $v$ nonempty. The law replaces
$2vy+2z$ by $2vy+2zy$ while retaining the source word. By
Lemma~\ref{f54:lem:words}, the source can expose each of its suffix letters.
Propagating every such letter to every doubled word produces the common
suffix set $T$. If $H\cap T=\varnothing$, the result is exactly
$B(C,H)$ after removing duplicate doubled words.

If $H\cap T\ne\varnothing$, one resulting word has its initial
variable in its suffix and hence is a square word $b_D$. Since
$b_D^2\eq b_D$, the last law of \eqref{f54:eq:flat}, with $x=b_D$,
allows every other doubled word to be squared in its presence.
Every variable of a square word can be exposed in a suffix position.
The same propagation now gives the full support $C$ to each square
word. All reduce to $b_C$, so their doubled sum is $2b_C$.

Finally add one of the doubled words specified in the statement to the
corresponding background. In the first case the content remains $C$
and a square word remains present; in the second case the initial set
remains $H$ and the suffix set remains $C-H$. The reductions just proved
give the original background again, proving absorption.
\end{proof}

For the unmarked type $H=0$, choose any family
$\R\subseteq\{D:\varnothing\ne D\subseteq C\}$ and set
\begin{equation}\label{f54:eq:flatdead}
N_{\mathrm f}(C,0,\R)=2b_C+\sum_{D\in\R}b_D.
\end{equation}
For $\varnothing\ne H\subseteq C$, choose any family
\[
\R\subseteq\{D\subseteq C:|D\cap H|=1\}.
\]
Writing $x_D$ for the unique element of $D\cap H$, set
\begin{equation}\label{f54:eq:flatlive}
N_{\mathrm f}(C,H,\R)=B(C,H)+
\sum_{D\in\R}u_{x_D,D-\{x_D\}}.
\end{equation}
Here and below $C$ is always nonempty, while $\R$ may be empty.

\begin{theorem}\label{f54:thm:flat}
The eleven identities $\Gamma_{\mathrm f}$ form a complete basis for
both $P$ and $\widetilde P$. The forms \eqref{f54:eq:flatdead} and
\eqref{f54:eq:flatlive} are distinct representatives of their free-algebra
elements on every finite variable set $X$.
\end{theorem}
\begin{proof}
For soundness, $3x\eq x$ follows from the addition tables. The doubled
image is $K_{\mathrm f}$, where $uv+z=uv+zv$: if $v=e$ both sides
are $u+z$, and otherwise $uv=zv=0$ and $0$ is additively absorbing.
This proves the middle law of \eqref{f54:eq:flat}. In the last law,
$y$ differs from $y^2$ only at $y=a$, while $2x^2\in\{0,e\}$;
either element has the same sum with $a$ as with $0$. The word laws
were checked in Lemma~\ref{f54:lem:words}. Thus the basis is sound in
both algebras.

Normalize $E(t)$ by Lemma~\ref{f54:lem:flatbg}, and use $t\eq E(t)+t$.
In the unmarked case, the last law of \eqref{f54:eq:flat} squares every
remaining word in the presence of $2b_C$, giving a square word with
the same support. In the marked case, every original word has its
initial variable in $H$ and all suffix variables in $C-H$, so it is
already of the permitted marked type after word reduction.
The absorption conclusions of the lemma reduce each coefficient modulo
two. Thus both cases give the displayed forms, with $\R$ equal to
the family of supports occurring an odd number of times in the original
polynomial. This is an equational reduction for arbitrary finite $X$.

We prove uniqueness by evaluations in the two displayed subsemirings.
First assign $0$ to one variable $x$ and $e$ to all others. The value
is $0$ exactly when $x\in C$, since $0$ is absorbing for both
operations on $K_{\mathrm f}$. This recovers $C$. Next, for every
nonempty $Y\subseteq C$, assign $a$ on $Y$ and $e$ elsewhere.
An unmarked form always gives $0$. A form marked by $H$ gives $a$
exactly when $Y=H$: each word then has initial value $a$ and suffix
values $e$. For any other nonempty $Y$, either a suffix contains $a$
and produces $0$, or initial values include both $e$ and $a$, whose
sum is $0$. The background forces that result. Hence the type and $H$
are recovered.

On the field copy, the background is identically zero and a standard
word of support $D$ has value $\prod_{x\in D}z_x$. The value of the
normal form is therefore the Boolean polynomial
\begin{equation}\label{f54:eq:boolean}
\sum_{D\in\R}\prod_{x\in D}z_x\quad\text{over }\FF.
\end{equation}
Evaluation at the characteristic vector of $Y$ sums the coefficients
on subsets of $Y$. Lemma~\ref{lem:finite-set-recovery}(i) recovers
$\R$, completing separation in both algebras and hence the proof.
\end{proof}

\subsubsection{The chain additive background}
Let $\Gamma_{\mathrm c}$ consist of $\SR$, \eqref{f54:eq:word}, and
\begin{equation}\label{f54:eq:chain}
3x\eq x,\qquad
2xy\eq2x+2y^2,\qquad
2x^2+x\eq2x^2+x^2.
\end{equation}
For a polynomial with initial set $H$ and suffix set $T$, put
$C=H\cup T$ and $L=H-T$. Thus $L$ consists of the variables whose
occurrences are all initial; there is no restriction on how many words
they initiate. Define
\begin{equation}\label{f54:eq:chainbg}
D(C,L)=\sum_{x\in L}2x+\sum_{y\in C-L}2y^2
\quad(\varnothing\ne C\subseteq X,\ L\subseteq C).
\end{equation}

\begin{lemma}\label{f54:lem:chainbg}
Modulo $\Gamma_{\mathrm c}$, $E(t)\eq D(C,L)$. This background
absorbs $2b_U$ for every nonempty $U\subseteq C-L$, and
$2u_{x,U}$ for $x\in L$, $U\subseteq C-L$. In its presence, every
word occurrence of $t$ whose initial variable belongs to $C-L$ can
be replaced by its square.
\end{lemma}
\begin{proof}
Apply the middle law of \eqref{f54:eq:chain} repeatedly to a doubled word:
it gives twice the initial variable, together with twice the square of
each suffix variable. Substituting $y=x$ in this law gives
$2x^2\eq2x+2x^2$. Doubled terms are additively idempotent, so a
variable occurring in a suffix is retained only in squared form. The
sum over all words is precisely $D(C,L)$. The same expansion shows
that each specified doubled word is a sum of components already in
the background, proving absorption.

If a word $w$ of $t$ starts in $C-L$, all its variables belong to $C-L$,
since every other variable is in a suffix position. Hence
$w^2\eq b_{c(w)}$ and $2w^2$ is absorbed by $D(C,L)$. Inserting
this absorbed term, applying the last law of \eqref{f54:eq:chain} with
$x=w$, and deleting it again gives
\[
D(C,L)+w\eq D(C,L)+2w^2+w
\eq D(C,L)+2w^2+w^2\eq D(C,L)+w^2.\qedhere
\]
\end{proof}

Choose an arbitrary family
\begin{equation}\label{f54:eq:chainallowed}
\R\subseteq\{U:\varnothing\ne U\subseteq C,\ |U\cap L|\le1\}.
\end{equation}
For its members define $v_U=b_U$ if $U\cap L=\varnothing$ and
$v_U=u_{x,U-\{x\}}$ if $U\cap L=\{x\}$. The normal form is
\begin{equation}\label{f54:eq:chainnormal}
N_{\mathrm c}(C,L,\R)=D(C,L)+\sum_{U\in\R}v_U.
\end{equation}

\begin{theorem}\label{f54:thm:chain}
The eleven identities $\Gamma_{\mathrm c}$ form a complete basis for
both $Q$ and $\widetilde Q$. The forms \eqref{f54:eq:chainnormal}, with
the stated parameter restrictions, are distinct representatives of
their free-algebra elements on every finite variable set $X$.
\end{theorem}
\begin{proof}
For soundness, the doubled image is $K_{\mathrm c}$. In this image
$uv=u+v^2$: if $v=e$ both sides equal $u$, and otherwise both equal
$0$, the greatest element in the additive order. Thus
$2xy\eq2x+2y^2$. In the last law, the only nonidempotent
multiplicative value of $x$ is $a$, and then $2x^2=0$ absorbs both
$a$ and $a^2$. All remaining laws were already checked.

Use $t\eq E(t)+t$ and Lemma~\ref{f54:lem:chainbg}. A word starting
outside $L$ becomes a square word on a subset of $C-L$. A word
starting at $x\in L$ has no further occurrence of $x$, and all suffix
variables are in $C-L$, so it becomes $u_{x,U}$ for a permitted $U$.
Absorbed doubled words reduce every coefficient modulo two. This gives
\eqref{f54:eq:chainnormal}, with $\R$ again the odd-support family.

For uniqueness, the single-$0$ assignments used in
Theorem~\ref{f54:thm:flat}, now in $K_{\mathrm c}$, recover $C$.
For each $x\in C$, assign $a$ to $x$ and $e$ to every other
variable. The background and full normal form have value $a$ if
$x\in L$, and $0$ if $x\notin L$. In the first case every added
word has value $e$ or $a$; in the second case the background already
has the absorbing value $0$. These tests recover $L$. Finally the
field copy gives \eqref{f54:eq:boolean}, which uniquely recovers $\R$.
Thus distinct normal forms are distinct functions in each algebra,
and the equational reduction is complete.
\end{proof}

\subsubsection{Generated varieties}
\begin{corollary}\label{f54:cor:join}
In the constant-free signature,
\[
\begin{aligned}
\Var(P)&=\Var(\widetilde P)=\Var(K_{\mathrm f})\vee\Var(\FF),\\
\Var(Q)&=\Var(\widetilde Q)=\Var(K_{\mathrm c})\vee\Var(\FF).
\end{aligned}
\]
\end{corollary}
\begin{proof}
The two indicated subsemirings separate the normal forms in the
preceding theorems. Apply Lemma~\ref{lem:normal-form-completeness}.
\end{proof}

\endgroup
\subsection{Support-antichain saturation}\label{sec:family55}
\begingroup
\providecommand{\eq}{}\renewcommand{\eq}{\approx}
\providecommand{\SR}{}\renewcommand{\SR}{\mathsf{SR}}
\providecommand{\Var}{}\renewcommand{\Var}{\operatorname{Var}}
\providecommand{\A}{}\renewcommand{\A}{\mathcal A}
\providecommand{\BB}{}\renewcommand{\BB}{\mathbb B}
\providecommand{\Ant}{}\renewcommand{\Ant}{\operatorname{Ant}}

\emph{Generators:} $P=A_{1424}$, $Q=A_{1916}$, $R=A_{1870}$, $P^{\mathrm{op}}=A_{712}$, $Q^{\mathrm{op}}=A_{1006}$, $R^{\mathrm{op}}=A_{1318}$.

For each of the source algebras $P,Q,R$, write $(0,e,f,a)$ for
labels $(0,1,2,3)$. They share the subsemiring $K=\{0,e,a\}$ with
addition the maximum operation on $e<a<0$ and multiplication
\begin{equation}\label{f55:eq:k}
uv=\begin{cases}u,&v=e,\\0,&v\ne e.\end{cases}
\end{equation}

Let $\BB$ denote the two-element Boolean semiring, with addition
maximum and multiplication minimum on $\{0,1\}$. All three sources
contain a copy of $\BB$. The embeddings are given in
Table~\ref{f55:tab:small}; an ordered carrier lists the images in the
order $(0,e,a)$ or $(0,1)$, respectively.
\begin{table}[htbp]\centering
\caption{Subsemiring embeddings in local carrier labels.}
\label{f55:tab:small}
\begin{tabular}{rcc}\toprule
Source&Carrier of $K$&Carrier of $\BB$\\\midrule
1424&$[0,1,3]$&$[0,2]$\\
1916&$[0,1,3]$&$[1,2]$\\
1870&$[0,1,3]$&$[2,0]$\\
\bottomrule\end{tabular}
\end{table}

In $R$, label $2$ is multiplicatively absorbing and $0\cdot2=2$.
Thus the element called $0$ in its copy of $K$ is not a zero for the
whole multiplication. All algebraic arguments respect this distinction.

\subsubsection{The basis and its word reductions}
Let $\Gamma$ consist of $\SR$ and the following eight identities:
\begin{gather}
xyz\eq xzy,\qquad xy^2\eq xy,\qquad x^2y\eq y^2x,\label{f55:eq:word}\\
x+x\eq x,\label{f55:eq:idem}\\
x+yx\eq x^2+yx,\label{f55:eq:projectbare}\\
xy+zx\eq x^2y+zx,\label{f55:eq:projectword}\\
x+yz\eq x+yz+xz,\label{f55:eq:extend}\\
xy+z^2\eq xy+z^2+xz^2.\label{f55:eq:transfer}
\end{gather}
The three identities on line \eqref{f55:eq:word} are counted separately.
Direct substitution in the three source operation tables verifies
every identity. For the two projection laws, the only carrier value changed by
squaring is $a$; the other summand makes that change invisible.
The remaining laws can also be read directly from the defining tables.

Fix a finite nonempty variable set $X$ and an order on it. Let $c(w)$
denote the support of a word. For nonempty $D\subseteq X$, and for
$x\notin U$, put
\[
b_D=\prod_{y\in D}y^2,\qquad u_{x,U}=x\prod_{y\in U}y.
\]
When $U=\varnothing$, the latter word is $x$.

\begin{lemma}\label{f55:lem:words}
Modulo \eqref{f55:eq:word}, a word with initial variable $x$ reduces to
$u_{x,c(w)-\{x\}}$ if $x$ occurs only once in that word, and to
$b_{c(w)}$ otherwise. Moreover $b_Db_E\eq b_{D\cup E}$.
The standard words are distinct functions in each of $P,Q,R$, and
there are $2^n-1+n2^{n-1}$ of them on $n$ variables.
\end{lemma}
\begin{proof}
The reductions follow from Lemma~\ref{f54:lem:words}, with the same
identities. Its separation argument also applies: use $\BB$ in
place of the field to distinguish supports, and the same $a,e$
assignments in $K$ to distinguish marks. The count is unchanged.
\end{proof}

\subsubsection{Support saturation}
For a polynomial $t$, let $H(t)$ be its set of initial variables and
$T(t)$ the union of all variables appearing in noninitial positions.
Define
\begin{equation}\label{f55:eq:parameters}
\begin{aligned}
C(t)&=H(t)\cup T(t),\qquad L(t)=H(t)-T(t),\\
\A(t)&=\min\{c(w):w\text{ is a word occurrence of }t\}.
\end{aligned}
\end{equation}
Here $\min$ means inclusion-minimal members. A variable in $L(t)$
may initiate several word occurrences; it must have no noninitial
occurrence anywhere in $t$. Notice that $C(t)-L(t)=T(t)$.

An admissible triple $(C,L,\A)$ consists of a nonempty $C\subseteq X$,
a subset $L\subseteq C$, and a nonempty antichain $\A$ satisfying
\begin{align}
&\varnothing\ne A\subseteq C,\quad |A\cap L|\le1
          &&(A\in\A),\label{f55:eq:admissible1}\\
&\text{for each }x\in L\text{ there exists }A\in\A
           \text{ with }A\subseteq(C-L)\cup\{x\}.
&&\label{f55:eq:admissible2}
\end{align}
For such a triple define
\[
\Omega(C,L,\A)=\{D:\varnothing\ne D\subseteq C,
             \ |D\cap L|\le1,\ (\exists A\in\A)\ A\subseteq D\}.
\]
For $D$ in this family set
\[
v_D=\begin{cases}
b_D,&D\cap L=\varnothing,\\
u_{x,D-\{x\}},&D\cap L=\{x\},
\end{cases}
\qquad
N(C,L,\A)=\sum_{D\in\Omega(C,L,\A)}v_D.
\]
This is a saturated normal form: every permitted extension of a minimal
support is included once.

\begin{lemma}\label{f55:lem:project}
Every polynomial $t$ is equivalent modulo $\Gamma$ to a sum of
standard words $v_D$ with the same $C(t)$, $L(t)$, and $\A(t)$.
Every word whose initial variable is outside $L(t)$ is squared;
the other words keep their unique initial variables in $L(t)$.
\end{lemma}
\begin{proof}
Put $L=L(t)$ and $T=T(t)$. A word starting in $L$ has no repeated
initial and all its suffix variables are in $T$. Lemma~\ref{f55:lem:words}
puts it in the required form. A word starting at $x\in T$ either
already repeats its initial, in which case it reduces to a square word,
or has a separate donor word with $x$ in a suffix position. Expose
that suffix occurrence at the end of the donor, writing it as $zx$
with a nonempty prefix $z$.

If the word to be projected is the bare variable $x$, use
\eqref{f55:eq:projectbare}. Otherwise write it as $xy$ with $y$ nonempty
and use \eqref{f55:eq:projectword}. The donor is retained while the
initial becomes squared. The resulting word reduces to $b_D$ with
unchanged support. Every letter in this word was already in $T$, so
the projection introduces no new noninitial variable. Conversely,
square-word normalization preserves all its noninitial support,
so the union $T$ is unchanged. Donors therefore remain available
throughout the process. If necessary, idempotence duplicates a donor
before using it. The support family is unchanged, and variables in $L$
remain exactly those occurring only initially. This proves the lemma.
\end{proof}

\begin{lemma}\label{f55:lem:saturate}
For every polynomial $t$, its parameters \eqref{f55:eq:parameters} form
an admissible triple and
\begin{equation}\label{f55:eq:normalization}
t\eq N(C(t),L(t),\A(t))
\end{equation}
is derivable from $\Gamma$. Conversely, the parameters of every
displayed normal form are its defining admissible triple.
\end{lemma}
\begin{proof}
Write $C=C(t)$, $L=L(t)$, $T=C-L$, and $\A=\A(t)$.
No word contains two variables of $L$, since only its initial position
can contain such a variable. This proves \eqref{f55:eq:admissible1}.
Every $x\in L$ initiates a word supported in $T\cup\{x\}$;
choosing a minimal support below that word proves
\eqref{f55:eq:admissible2}.

First apply Lemma~\ref{f55:lem:project}. For each $z\in T$ there is
a donor word with $z$ in a suffix position. Expose it as $yz$.
Substituting any currently present word $w$ for $x$ in
\eqref{f55:eq:extend} adjoins $wz$ while retaining both existing words.
After word reduction this extends its support by $z$. Thus any word
can be extended by any subset of $T$.

Consider a required support $D\in\Omega(C,L,\A)$ and choose
$A\in\A$ with $A\subseteq D$. A standard word on $A$ is present
after projection. If $A$ and $D$ contain the same member of $L$, or
both avoid $L$, suffix extension produces $v_D$. The only other case
is $A\subseteq T$ and $D\cap L=\{x\}$. A live word initiated
by $x$ is present. If the bare word $x$ is present, extend that word
by $D-\{x\}$. Otherwise choose such a word $xy$ with nonempty
suffix $y$. In the presence of $b_A$, put $z=b_A$ in
\eqref{f55:eq:transfer}. Since $b_A^2\eq b_A$, this adjoins $xb_A$,
which reduces to $u_{x,A}$. Extend its suffix to $D-\{x\}$.

All inserted words have permitted supports containing an original
minimal support. Every original projected word is already among the
required summands. There are finitely many required supports, so these
steps adjoin them all, and idempotence deletes duplicate occurrences.
This proves \eqref{f55:eq:normalization} for arbitrary finite $X$.

Conversely, the minimal supports of $N$ are exactly $\A$. Every
$x\in L$ appears initially by \eqref{f55:eq:admissible2}, and never
noninitially. Every $y\in C-L$ appears in a suffix: extend any
member of $\A$ by $y$; its canonical word has $y$ in a suffix,
either because it is a square word or because its initial belongs
to $L$. Thus the total content is $C$ and the variables occurring
only initially are precisely $L$.
\end{proof}

For example, on three distinct variables the parameters of $xy+z^2$
are $C=\{x,y,z\}$, $L=\{x\}$, and $\A=\{\{x,y\},\{z\}\}$.
The saturation proof gives the explicit consequence
\[
xy+z^2\eq xy+xz+xyz+z^2+y^2z^2.
\]
The set $L$ cannot be discarded: $x+xy$ and $x+xy+x^2$ have
the same content and minimal support family, but assigning $x=a$,
$y=e$ in $K$ gives $a$ and $0$, respectively.

\subsubsection{Completeness and generated varieties}
\begin{theorem}\label{f55:thm:basis}
The thirteen identities $\Gamma$ form a complete basis for each of
$P,Q,R$. The admissible triples $(C,L,\A)$ are in bijection with
the elements of the free algebra of their common generated variety
on any finite nonempty variable set $X$.
\end{theorem}
\begin{proof}
Soundness follows from the displayed operation tables, and
Lemma~\ref{f55:lem:saturate} gives equational normalization. It remains
to separate the parameters of normal forms, using only the two
subsemirings specified in Table~\ref{f55:tab:small}.

Here $K=K_{\mathrm c}$, so the single-$0$ and single-$a$ tests
of Theorem~\ref{f54:thm:chain} recover $C$ and $L$: for $x\in C$,
the latter test gives $a$ exactly when $x$ occurs only initially,
and $0$ otherwise. On $\BB$, a subset $Y$ gives value $1$
exactly when it contains a word support. Its minimal accepted
subsets recover $\A$ by Lemma~\ref{lem:finite-set-recovery}(ii).
Lemma~\ref{lem:normal-form-completeness} now proves completeness.
\end{proof}

\begin{corollary}\label{f55:cor:join}
In the two-binary-operation signature,
\[
\Var(P)=\Var(Q)=\Var(R)=\Var(K)\vee\Var(\BB).
\]
\end{corollary}
\begin{proof}
The subsemirings $K,\BB$ separate the normal forms in
Theorem~\ref{f55:thm:basis}; apply Lemma~\ref{lem:normal-form-completeness}.
\end{proof}

\endgroup
\subsection{Boundary and saturated supports}\label{sec:family57}
\begingroup
\providecommand{\eq}{}\renewcommand{\eq}{\approx}
\providecommand{\SR}{}\renewcommand{\SR}{\mathsf{SR}}
\providecommand{\Var}{}\renewcommand{\Var}{\operatorname{Var}}
\providecommand{\Sat}{}\renewcommand{\Sat}{\mathcal S}
\providecommand{\Res}{}\renewcommand{\Res}{\mathcal R}
\providecommand{\Bd}{}\renewcommand{\Bd}{\mathcal B}

The first source algebra is $A=\{\bot,e,a,o\}$. The element $\bot$
is an additive identity and a two-sided multiplicative zero. Any two
elements of $K=\{e,a,o\}$ have sum $o$. Multiplication on $K$ is
\begin{equation}\label{f57:eq:a}
uv=\begin{cases}u,&v=e,\\o,&v\ne e.\end{cases}
\end{equation}
The second source algebra is $B=\{0,t,e,a\}$. Here $0$ is an
additive identity and a two-sided multiplicative zero, while $t$
absorbs every element under addition. The remaining sums are
\begin{equation}\label{f57:eq:badd}
e+e=t,\qquad e+a=e,\qquad a+a=a.
\end{equation}
On $\{t,e,a\}$ the multiplication table is
\begin{equation}\label{f57:eq:bmul}
\begin{array}{c|ccc}
\cdot&t&e&a\\\hline
t&t&t&0\\
e&t&e&0\\
a&a&a&0
\end{array}.
\end{equation}
Thus $e$ is a right identity in both sources, but need not be a
two-sided identity. These definitions satisfy $\SR$, as may be checked
directly from the four-element tables.

The local representatives are listed below.
For $A$, the labels of $(\bot,e,a,o)$ are $(2,1,3,0)$; for $B$,
those of $(0,t,e,a)$ are $(0,1,2,3)$.

\emph{Generators (basis size):} $A=A_{1866}$ (13), $A^{\rm op}=A_{1310}$ (13), $B=A_{1909}$ (15), $B^{\rm op}=A_{999}$ (15).

All four algebras have nonidempotent addition and noncommutative multiplication.

\subsubsection{Words and support profiles}
For $\varnothing\ne D\subseteq X$ and $x\in D$, set
\begin{equation}\label{f57:eq:words}
b_D=\prod_{y\in D}y^2,\qquad
u_{x,D}=x\prod_{y\in D\setminus\{x\}}y,
\end{equation}
The second expression is $x$ when $D=\{x\}$. Call $b_D$ a square
word and $u_{x,D}$ a marked word.

\begin{lemma}\label{f57:lem:words}
The three laws
\begin{equation}\label{f57:eq:wordlaws}
xyz\eq xzy,\qquad xy^2\eq xy,\qquad x^2y\eq y^2x
\end{equation}
form, relative to associativity, a basis for the word identities of
both $A$ and $B$. A word $w$ reduces to $u_{x,c(w)}$ if its initial
variable $x$ occurs only once, and to $b_{c(w)}$ otherwise. In
particular,
\begin{equation}\label{f57:eq:squares}
w^2\eq b_{c(w)},\qquad b_Db_E\eq b_{D\cup E}.
\end{equation}
There are $2^{|X|}-1+|X|2^{|X|-1}$ distinct standard word functions.
\end{lemma}
\begin{proof}
The equational reductions and square-product formulas follow from
Lemma~\ref{f54:lem:words}, whose three word identities are identical
to those used here. These identities hold in the defining tables.

Different supports are distinguished by assigning the multiplicative
zero to a variable present on only one side and $e$ to the other
variables. On a fixed support $D$, assign $a$ to $x$ and $e$ to
the other variables of $D$. The word $u_{x,D}$ has value $a$.
Every other standard word on $D$ has value $o$ in $A$ or $0$ in
$B$. Thus all standard words are distinct as functions. The reduction
and these separating assignments prove completeness and the count.
\end{proof}

For a polynomial $p$, define its capped support profile on all
$Y\subseteq X$, including $Y=\varnothing$, by
\begin{equation}\label{f57:eq:profile}
b_p(Y)=\min\bigl(2,\#\{\text{word occurrences }w\text{ in }p:
c(w)\subseteq Y\}\bigr).
\end{equation}
Let $\Sat_p=\{Y:b_p(Y)=2\}$. This family is upward closed and
does not contain the empty set. A residual support is the support
$D$ of an occurrence for which $b_p(D)=1$. There is exactly one
occurrence at such a support, and we record its standard word in
$\Res_p$. Equivalently, residual supports are the inclusion-minimal
positive supports whose profile value is one. They form an antichain.

The profile is determined by $\Sat_p$ and the residual supports:
outside $\Sat_p$, its value is one precisely when a residual support
is contained in $Y$. Indeed, if $b_p(Y)=1$, the sole occurrence
supported inside $Y$ has a residual support. Two different residual
supports have saturated union.

\begin{lemma}\label{f57:lem:insert}
Assume \eqref{f57:eq:wordlaws} and the absorption identities
\begin{equation}\label{f57:eq:absorb}
x+y+xy\eq x+y,\qquad x+y+xyz\eq x+y.
\end{equation}
One can insert $2b_Y$ into $p$ for every $Y\in\Sat_p$, retaining
all original occurrences. Here $2w$ means $w+w$.
\end{lemma}
\begin{proof}
Choose two distinct occurrences $v,w$ supported in $Y$, allowing
their words to coincide. Let $z$ list every variable of $Y$ once.
The second law in \eqref{f57:eq:absorb} inserts $vwz$. Its initial
variable is repeated because $z$ contains every variable of $Y$.
Its standard word is consequently $b_Y$. Repeat the insertion.
Every new word is active only on supersets of $Y$, all of which
were already saturated. Thus the profile is unchanged, and the
procedure can be performed for every saturated support.
\end{proof}

\subsubsection{A basis for the saturated algebra}
\begin{theorem}\label{f57:thm:saturation}
A complete basis of thirteen identities for $A$ consists of $\SR$,
\eqref{f57:eq:wordlaws}, \eqref{f57:eq:absorb}, and
\begin{equation}\label{f57:eq:capbasis}
3x\eq2x,\qquad 2x\eq2x^2,\qquad 2xy\eq2yx.
\end{equation}
For every polynomial $p$, a canonical form is
\begin{equation}\label{f57:eq:capnormal}
N_A(p)=\sum_{w\in\Res_p}w+\sum_{Y\in\Sat_p}2b_Y.
\end{equation}
Only present summands are written; an empty component is omitted.
\end{theorem}
\begin{proof}
All displayed laws hold in $A$. For \eqref{f57:eq:capbasis}, a doubled
word has value $\bot$ exactly when some used variable is $\bot$,
and has value $o$ otherwise. For \eqref{f57:eq:absorb}, the added
product vanishes if either of its first two factors is $\bot$;
if neither is $\bot$, their sum is already $o$.

The doubled identities imply
\begin{equation}\label{f57:eq:doublecontent}
2v\eq2w\quad\text{whenever }c(v)=c(w).
\end{equation}
To see this equationally, multiply $2xy\eq2yx$ by any needed left
and right word contexts and distribute. This permits adjacent
transpositions inside a doubled word. Likewise $2x^2\eq2x$ deletes
adjacent repetitions in such contexts. Sort and shorten both words
to the same fixed squarefree word. No empty context is required as
a term; an absent context is simply omitted.

Apply Lemma~\ref{f57:lem:insert}. Every original word $w$ with
saturated support $D$ can now be removed while keeping its marker:
\[
2b_D+w\eq2w+w\eq2w\eq2b_D.
\]
Every remaining original occurrence has residual support and can
be replaced by its standard word using Lemma~\ref{f57:lem:words}.
This derives \eqref{f57:eq:capnormal} solely from the proposed basis.
The expression is nonempty because $p$ has at least one occurrence.

For uniqueness, assign $e$ on $Y$ and $\bot$ outside $Y$.
The values for profile entries zero, one and two are respectively
$\bot,e,o$. Hence the term function determines $b_p$, and therefore
$\Sat_p$ and the residual supports. At a residual support $D$,
put $\bot$ outside $D$ and allow arbitrary values inside $D$.
Exactly one occurrence is supported there, so these assignments
recover its word function. Lemma~\ref{f57:lem:words} recovers its
standard word uniquely. Completeness follows from
Lemma~\ref{lem:normal-form-completeness}.
\end{proof}
Although the nonzero multiplication in \eqref{f57:eq:a} resembles a
zero extension, the subset $\{e,a\}$ is not a subsemigroup:
$a^2=o$. The proof above does not invoke a transfer theorem that
requires this subset to be closed.

\subsubsection{Boundary initials and a second basis}
For the following definitions, replace each word occurrence of $p$
by its standard word; this preserves its support and multiplicity.
For each $x\in X$, consider the supports $D$ of marked occurrences
$u_{x,D}$ in $p$ such that
\begin{equation}\label{f57:eq:boundary}
b_p(D\setminus\{x\})=0.
\end{equation}
Let $\Bd_{p,x}$ be the inclusion-minimal supports in this collection.
This is an antichain; repeated occurrences do not create repeated
members. A residual marked word necessarily belongs to this
antichain. Write $\Bd_p$ for all the marked words $u_{x,D}$ with
$D\in\Bd_{p,x}$. This set includes residual marked words.

\begin{theorem}\label{f57:thm:boundary}
A complete basis of fifteen identities for $B$ consists of $\SR$,
\eqref{f57:eq:wordlaws}, \eqref{f57:eq:absorb}, and
\begin{equation}\label{f57:eq:boundarybasis}
\begin{aligned}
3x&\eq2x,&2x&\eq x+2x^2,\\
xy+y&\eq x^2y+y,&xyz+y&\eq x^2yz+y,\\
x+xy&\eq x+2x^2y^2.
\end{aligned}
\end{equation}
A canonical form is
\begin{equation}\label{f57:eq:boundarynormal}
N_B(p)=\sum_{w\in\Res_p}w+
       \sum_{Y\in\Sat_p}2b_Y+
       \sum_{\substack{u_{x,D}\in\Bd_p\\D\in\Sat_p}}u_{x,D}.
\end{equation}
In particular, a residual marked word occurs only once in this form.
\end{theorem}
\begin{proof}
First we check soundness. The table gives $3x=2x=x+2x^2$.
In a projection law, squaring $x$ changes its value only when
$x=a$. Then $x^2=0$ and the other product is either $0$ or $a$.
For it to be $a$, the value of $y$ must be $t$ or $e$, each of
which absorbs $a$ under addition. This proves both projection laws.
For $x+xy=x+2x^2y^2$, the cases $x=0,t,a$ follow directly from
\eqref{f57:eq:badd}--\eqref{f57:eq:bmul}. If $x=e$, a value $y=0$ or $a$
makes the second summand on both sides zero; a value $y=t$ or $e$
makes the whole sum $t$. For \eqref{f57:eq:absorb}, if one of $x,y$
is zero the added product is zero. If $x+y=t$, it absorbs the
product. The remaining nonzero cases have $x+y=e$ or $a$;
then a product $xy$ or $xyz$ is respectively $0$ or $a$, or is
zero when $x=y=a$. It is again absorbed. The word laws were
checked in Lemma~\ref{f57:lem:words}.

For derivability of the normal form, standardize the original words
and insert all markers using Lemma~\ref{f57:lem:insert}. Keep two copies
of each $b_Y$ permanently available. The law $3v\eq2v$ removes
any excess square copies on saturated supports. We describe the
remaining reductions in detail.

\emph{Repeated marked words.} If a marked word $w=u_{x,D}$ occurs
at least twice, then $D\in\Sat_p$. The law $2w\eq w+2w^2$
and \eqref{f57:eq:squares} reduce the number of copies by one and add
two copies of $b_D$. Repetition of this step leaves a single marked
copy. The added square copies are absorbed into the retained marker.

\emph{Nonboundary marked words.} Suppose $u_{x,D}$ remains and
$b_p(D\setminus\{x\})>0$. Some original occurrence has support
$E\subseteq D\setminus\{x\}$. A word with that support is still
available: an unsaturated original occurrence remains, and at a
saturated original support there is the permanently retained square
marker. Denote any such current word by $v$. It need not have the
same initial variable as the original donor. Let $z$ list every
variable of $D\setminus\{x\}$ once. This is a nonempty word
because $E$ is nonempty. The word laws give
\[
xvz\eq u_{x,D},\qquad x^2vz\eq b_D.
\]
The instance $xvz+v\eq x^2vz+v$ of the contextual projection law
replaces the marked word by $b_D$, retaining $v$. The support $D$
is saturated, because it contains the marked occurrence and the
distinct donor occurrence. The new square copy can thus be absorbed.
This step never removes the last available word at an original
support, so the donor argument remains valid throughout the process.

\emph{Nonminimal boundary supports.} All marked words now satisfy
\eqref{f57:eq:boundary}. If $E\subsetneq D$ are two such supports
with the same initial $x$, choose an inclusion-minimal one $E$
below $D$, and let $v=u_{x,E}$. Let $z$ list $D\setminus E$.
Then $vz\eq u_{x,D}$ and $v^2z^2\eq b_D$. The instance
\[
v+vz\eq v+2v^2z^2
\]
removes the larger marked word and retains the smaller one. Again
$D$ is saturated, since two original marked supports are contained
in it. Its extra square copies are absorbed. Applying this to every
nonminimal boundary support leaves exactly $\Bd_p$. A residual
marked word cannot be projected, duplicated, or have a smaller
marked support, because any such situation would make its support
saturated. We have derived exactly \eqref{f57:eq:boundarynormal}.

It remains to show uniqueness. Assign $e$ on $Y$ and $0$ outside
$Y$. Profile entries zero, one and two yield $0,e,t$, respectively.
Hence the term function recovers $b_p$. At a residual support,
assignments vanishing outside it isolate the unique supported word,
whose label is recovered by Lemma~\ref{f57:lem:words}.

Fix $x\in X$ and $Y\subseteq X\setminus\{x\}$ with $b_p(Y)=0$.
Assign $a$ to $x$, $e$ on $Y$, and $0$ to every other variable.
Words not containing $x$ vanish, because none is supported inside
$Y$. Square words, and marked words with $x$ in a suffix, also
vanish by \eqref{f57:eq:bmul}. A marked word with initial $x$ has
value $a$ exactly when its support is contained in $Y\cup\{x\}$.
Any positive number of these values has sum $a$, since $a+a=a$.
Thus the result is $a$ precisely when some member of $\Bd_{p,x}$
has its suffix support contained in $Y$; otherwise it is $0$.
The suffixes $D\setminus\{x\}$ for $D\in\Bd_{p,x}$ form an
antichain, and each is an admissible test by \eqref{f57:eq:boundary}.
Lemma~\ref{lem:finite-set-recovery}(ii) therefore recovers them as
the minimal accepted sets. All canonical data are determined, so
Lemma~\ref{lem:normal-form-completeness} proves completeness.
\end{proof}

\begin{corollary}\label{f57:cor:opposites}
The generated varieties satisfy
\[
\Var(A)\subsetneq\Var(B),\qquad
\Var(A^{\rm op})\subsetneq\Var(B^{\rm op}).
\]
\end{corollary}
\begin{proof}
The table of $A$
satisfies every identity of the basis in Theorem~\ref{f57:thm:boundary}.
By its completeness, $A\in\Var(B)$. The identity $2x\eq2x^2$
holds in $A$ but fails in $B$ at $x=a$, where its two values
are $a$ and $0$. This makes the inclusion strict. Corollary~\ref{cor:catalogue-transfers} proves the second inclusion.
\end{proof}

\endgroup
\section{Further support, prefix and endpoint forms}\label{sec:further-forms}
The remaining positive constructions combine bounded witnesses
with support and endpoint information. The bounds concern words and
the variables needed in a derivation. The completeness arguments still
apply to identities with arbitrarily many variables.

\subsection{Parity coefficients and initial variables}
The sources $A_{675}$ and $A_{692}$ have the same multiplication.
A nonlinear word is 1 if its head is 1, is 2 if all its variables
are 2, and is 0 otherwise. They contain the field $\FF$ on
$\{0,2\}$, with field zero 0 and field one 2; their additive
identity is the element labelled 1. The element 3 distinguishes
a variable from its square. Their tables are in
Appendix~\ref{app:tables}.

Take $\SR$ and the six common identities
\begin{align}
x^2y&\eq xy,&xyz&\eq xzy,&xyx&\eq xy,\label{eq:parity-words}\\
2xy&\eq2x^2,&x+xy&\eq x^2+xy,\label{eq:parity-heads}\\
2x^2+2y^2+xy&\eq2x^2+2y^2+yx.\label{eq:parity-transfer}
\end{align}
For $A_{675}$ add $3x\eq x^2$ and $2x\eq2x^2$; denote the
resulting thirteen-identity set by $\Sigma_T$.
For $A_{692}$ instead add $3x\eq x$; denote this
twelve-identity set by $\Sigma_P$.

For a nonempty support $D$ and $h\in D$, let
$W(D,h)=h\prod_{x\in D}x$, with the product in a fixed order.
The word laws reduce every nonlinear word with support $D$ and
head $h$ to $W(D,h)$. To check tail repetition deletion, derive
\[
xy^2\eq(xyx)y\eq x(yxy)\eq x(yx)\eq xyx\eq xy.
\]
Tail permutation then deletes all repetitions; the head can be
inserted into the tail using $xyx\eq xy$ backwards. Squares
cover singleton supports. Every nonlinear word is multiplicatively
idempotent modulo these laws.

For a polynomial, let $H$ contain its nonlinear heads. In the
threshold case, include every variable having at least two linear
occurrences. Let $L$ be the linearly occurring variables outside
$H$. In the threshold case their coefficients are one. Let
$\mathcal R$ record supports with odd coefficient, omitting those
single linear occurrences in the threshold case and counting all
occurrences in the periodic case. In the threshold case every
$D\in\mathcal R$ meets $H$. In the periodic case this holds
except for allowed singletons $\{x\}$ with $x\in L$.
Choose $h(D)$ as the least member of $D\cap H$ whenever this
intersection is nonempty. Define
\begin{align}
N_T(L,H,\mathcal R)
&=\sum_{x\in L}x+\sum_{h\in H}2h^2
  +\sum_{D\in\mathcal R}W(D,h(D)),\label{eq:parity-normal-T}\\
N_P(L,H,\mathcal R)
&=\sum_{x\in L}2x+\sum_{h\in H}2h^2
  +\sum_{\substack{D\in\mathcal R\\D\cap H\ne\varnothing}}W(D,h(D))
  +\sum_{\substack{x\in L\\\{x\}\in\mathcal R}}x.
\label{eq:parity-normal-P}
\end{align}
Here $L\cap H=\varnothing$ and $L\cup H\ne\varnothing$.
Empty components are omitted, so both expressions are terms.

\begin{theorem}\label{thm:parity}
$\Sigma_T$ is a basis for $A_{675}$, and $\Sigma_P$ is a basis
for $A_{692}$. Their polynomial functions on a finite variable set
are uniquely represented by the allowed forms
\eqref{eq:parity-normal-T} and \eqref{eq:parity-normal-P}, respectively.
The threshold basis also defines $\Var(A_{1641})$.
\end{theorem}
\begin{proof}
The tables verify the displayed laws. First normalize all words.
In the threshold case, pairs of linear $x$ become $2x^2$.
If a further linear copy remains, $2x^2+x\eq3x\eq x^2$;
iterating leaves only nonlinear copies at that head. In both
cases, a linear $x$ coexisting with a nonlinear word starting in
$x$ is promoted using $x+xy\eq x^2+xy$. These operations preserve
the parity of every support and retain the relevant head.

For a nonlinear word $w$, multiplicative idempotence and the scalar
law give $3w\eq w$. Thus two extra copies may be inserted.
If its head is $h$, $2xy\eq2x^2$ replaces these copies by
$2h^2$, yielding the full background $\sum_{h\in H}2h^2$.
Each background term is additively idempotent and absorbs a pair
of nonlinear words at that head. In the periodic case the same
argument with $3x\eq x$ inserts the background $2x$ for every
remaining linear variable.

The common background permits transfer between heads in the same
support. For $x,y\in D\cap H$, represent the standard word as
$xyv$, with $v$ a nonempty word containing the whole support;
the word laws allow the redundant copies. Multiply
\eqref{eq:parity-transfer} by $v$ and distribute. The doubled
terms $2x^2v,2y^2v$ reduce to $2x^2,2y^2$, and the other
two terms are the standard words with heads $x$ and $y$.
Move all remaining odd copies to the chosen head and absorb
the even copies into the background. This gives exactly the
two displayed normal forms.

For uniqueness, assign a selected variable $x$ the value 3 and
all others the additive identity 1. The result is 3 for $x\in L$,
0 for $x\in H$, and 1 otherwise. This recovers $L$ and $H$.
On the field subsemiring $\{0,2\}$ the background vanishes and
each word is the squarefree monomial on its support. After
subtracting the recovered linear part, Lemma~\ref{lem:finite-set-recovery}(i)
recovers $\mathcal R$ from these evaluations.
Consequently equal functions have the same normal form, proving
completeness in every rank.

For $A_{1641}$ a nonlinear word is 2 if its head is 2,
otherwise 0 when all variables lie in $\{0,3\}$ and 1 in
the other cases. Its additive identity is 2, and $\{1,0\}$
is a field copy with zero 1. The same equations normalize to
\eqref{eq:parity-normal-T}; assigning $x=3$ and the others 2
distinguishes $L$, $H$, and absence, and the field evaluations
recover parity. Hence it has the same complete basis.
\end{proof}

Reversal gives the bases for $A_{1153},A_{1170},A_{1742}$,
respectively. If $n$ is the
number of variables, the threshold free algebra has size
\[
2^n-1+\sum_{h=1}^n\binom nh 2^{n-h}
                  2^{\,2^n-2^{n-h}},
\]
and the periodic free algebra has size
\[
\sum_{\substack{h,l\geq0\\1\leq h+l\leq n}}
\binom nh\binom{n-h}{l}
2^{\,2^n-2^{n-h}+l}.
\]
The first three values are $3,27,711$ and $4,40,1162$.
These formulas follow by choosing $H$, then $L$, and then the allowed
parity supports in the unique normal forms.
\subsection{Support closures with prescribed initial variables}
We use the nonlinear standard words $W(D,h)$ introduced in the
preceding section. The common word identities are
\begin{equation}\label{eq:closure-word}
 x^2y\eq xy,\qquad xyz\eq xzy,\qquad xyx\eq xy.
\end{equation}
The derivation of $xy^2\eq xy$ given there shows that these laws
reduce every nonlinear word to $W(D,h)$, with nonempty support $D$
and initial variable $h\in D$. A linear word $x$ remains distinct
from $W(\{x\},x)=x^2$. In particular, $w^2\eq w$ for every
nonlinear word $w$.

\subsubsection{Separate and paired closure rules}
\begin{proposition}\label{thm:head-closure}
For $A_{1274}$ and $A_{1278}$, a basis consists of $\SR$,
\eqref{eq:closure-word}, and
\begin{equation}\label{eq:paired-closure}
 2x\eq x^2,\quad 3x\eq2x,\quad
 xy+xyz\eq xy,\quad x+y+xy\eq x+y.
\end{equation}
Deleting the last identity gives a basis for $A_{1282}$.
Replacing the first two identities in that shorter list by
$2x\eq x$ and $x+x^2\eq x^2$ gives a basis for $A_{1295}$.
\end{proposition}
\begin{proof}
All four generators have multiplication
$0100/1111/2122/0100$. A nonlinear word is 1 if a used variable
is assigned 1; otherwise it is 2 when its initial is 2, and is 0
when its initial is 0 or 3. The element 1 is an additive identity.
The nonlinear values form the chain $1<2<0$, except in $A_{1278}$,
where the chain is $1<0<2$. In the two unpaired cases, $3+2=3$.
These rules verify the identities. They also show why the extra
paired law is available in the first two sources. Algebraically,
every nonlinear word is additively idempotent, and $x+x^2\eq x^2$
follows in the nonidempotent cases from $3x\eq2x$ and $2x\eq x^2$.

Fix a finite alphabet $X$ and expand a term $p$, retaining
multiplicities. Form a set $C(p)$ of standard words by the following
rules. Include each original word; include $x^2$ at a repeated
linear $x$ in the nonidempotent cases. From $W(E,h)$ include
$W(D,h)$ whenever $E\subseteq D\subseteq X$. From $x^2$ include
$x$. In the paired cases only, include the standard product of
any two distinct retained words. A nonlinear word may be paired
with itself because it can be duplicated. Take the least set closed
under these rules. There are only
$|X|+|X|2^{|X|-1}$ standard words, so this is a finite set.

Each closure step is derivable from the proposed basis. The extension
of a support is an instance of $xy+xyz\eq xy$, after normalizing
the right multiple. The square rule follows from square domination.
The paired rule is the last identity of \eqref{eq:paired-closure}.
Repeated linear summands first become squares, while nonlinear
duplicates can be deleted. It follows that
\begin{equation}\label{eq:closure-normal}
 p\eq N(p):=\sum_{w\in C(p)}w
\end{equation}
is derivable. Temporary nonlinear summands may be retained and
subsequently removed using Lemma~\ref{lem:temporary}; its proof
also works without global additive idempotence whenever the inserted
summands themselves are additively idempotent.

We prove uniqueness of \eqref{eq:closure-normal}. Let $Q$ be the
variables whose singleton squares are included, and let $M$ be the
included linear variables. A singleton square comes only from an
original singleton nonlinear word or a repeated linear variable
in a nonidempotent source. Thus $M$ consists of the original linear
variables together with $Q$. In an unpaired case,
$W(D,h)\in C(p)$ precisely when $h\in Q$, or an original nonlinear
word with initial $h$ has support contained in $D$. In a paired
case there is one further possibility: an original linear $h$
occurs once, and some other original occurrence is supported in $D$.
To justify this description, a product beginning in a nonlinear
word introduces only a support extension. A product beginning in
a linear word uses exactly the additional possibility just stated.
Further products cannot create a smaller support or a new initial.

Assign 3 to $x$ and 1 to every other variable. The value of $p$
is 0 for $x\in Q$, 3 for $x\in M\setminus Q$, and 1 otherwise.
Hence its function determines $M$ and $Q$. Suppose now that
$W(D,h)\notin C(p)$, and assign 1 outside $D$. Write $a$ for the
larger and $b$ for the smaller of the two nonlinear values 0 and 2.
If $h\notin M$, assign $a$ to $h$ and $b$ to the other variables
of $D$. No original active word has initial $h$, so $p$ is 1 or
$b$, whereas $W(D,h)$ is $a$. It is not absorbed.

If $h\in M$, then $h\notin Q$ and an original single linear $h$
is present. Assign 3 to $h$. In a paired case assign $b$ to the
other variables of $D$. No other original occurrence can be supported
there, by the membership criterion, so $p=3$. In an unpaired case
assign 2 to those variables. Every other active word then gives 1
or 2, so again $p=3$, since $3+2=3$. The missing nonlinear word
has value 0, and $3+0=0$. Thus it is not absorbed in either case.

The nonlinear members of $C(p)$ are therefore exactly the nonlinear
standard words absorbed by $p$. Equal term functions have equal
$M,Q$ and equal nonlinear absorption sets; hence they have equal
$C$-sets and equal normal forms. Applying \eqref{eq:closure-normal}
to both sides derives every identity and proves all four assertions.
\end{proof}

\subsubsection{A shared antichain of supports}
\begin{proposition}\label{thm:antichain}
A basis for $A_{1636}$ consists of $\SR$, \eqref{eq:closure-word},
and
\begin{gather}
 xy+xyz\eq xy,\qquad xu+yz\eq xu+yz+yxu,\label{eq:antichain-share}\\
 2x\eq x^2,\qquad3x\eq2x,\qquad x+y+xy\eq x+y.
\label{eq:antichain-nonai}
\end{gather}
A basis for $A_{1648}$ is obtained by replacing
\eqref{eq:antichain-nonai} with
\begin{equation}\label{eq:antichain-ai}
 2x\eq x,\qquad x+x^2\eq x^2,\qquad x+yz\eq x+yz+xyz.
\end{equation}
\end{proposition}
\begin{proof}
Both multiplications are $0110/1111/2222/0110$. A nonlinear word
is 2 when its initial is 2; otherwise it is 0 if all its variables
are in $\{0,3\}$, and is 1 in every other case. Addition on the
nonlinear values is maximum for $2<1<0$, with
$3+1=3$, $3+0=0$, and $3+3=0$ or 3, respectively.
For the sharing identity, a newly added value 0 of $yxu$ is
already supplied by $xu$; otherwise the added value is 1 or 2,
and a retained nonlinear word absorbs it. For the linear sharing
law in \eqref{eq:antichain-ai}, a value 0 of $xyz$ is already
supplied by $yz$. The paired law in the other source also uses
$3+3=0$. The scalar and word laws follow from the same value rules.

Let $H$ be all original initials, including linear ones. Start a
family $\mathcal A$ with the supports of the nonlinear occurrences.
For $A_{1636}$, also include $\{x\}$ at a repeated linear $x$, and
$\{x,y\}$ for distinct linearly occurring variables. Only its upward
closure matters, so $\mathcal A$ may be replaced by its antichain
of minimal members. Put $Q=\{x:\{x\}\in\mathcal A\}$, and let
$M$ be the original linear variables together with $Q$. Define
\[
 C(p)=\{x:x\in M\}\ \cup\
 \{W(D,h):h\in H\cap D,\ E\subseteq D
                    \text{ for some }E\in\mathcal A\}.
\]

Every support in $\mathcal A$ is realized by an inserted nonlinear
word: use the paired law or $2x\eq x^2$ for the added supports.
If $v$ realizes $E$ and $h\in H$ is linear, linear sharing inserts
$hv$. If $h$ is the initial of a nonlinear word, the second law of
\eqref{eq:antichain-share}, using $v$ and that word, again inserts
$hv$. Support extension then produces every prescribed $W(D,h)$.
Square domination supplies the linears in $Q$. No operation changes
$H$ or introduces a support containing no member of $\mathcal A$.
Deleting excess occurrences proves $p\eq\sum C(p)$.

To recover the parameters, assign 3 to $x$ and 2 elsewhere. The
four possible outputs are 0, 3, 1, 2 according as
$x\in Q$, $x\in M\setminus Q$, $x\in H\setminus M$, or
$x\notin H$. Suppose $W(D,h)$ is missing. If $h\notin H$, assign
$h=1$ and all others 2; the polynomial is 2 and the target is 1.
Otherwise no member of $\mathcal A$ is contained in $D$. Assign 3
on $D$ and 2 elsewhere. No original nonlinear word gives 0.
In $A_{1636}$ at most one linear occurrence is active, since two
would supply a support in $\mathcal A$ contained in $D$.
In $A_{1648}$ several active linears still sum to 3. Thus $p$ is
never 0, while $W(D,h)=0$. This again prevents absorption.

It follows that the term function recovers the linear part and the
entire nonlinear part of $C(p)$. Equal functions have the same
derived normal form, proving completeness in both cases.
\end{proof}
\subsection{Capped support coefficients}
If $p$ is a polynomial with occurrences counted, define
\begin{equation}\label{eq:capped-profile}
 b_p(Y)=\min\bigl(2,\#\{\text{occurrences of }p
                                  \text{ supported in }Y\}\bigr).
\end{equation}
A set $Y$ is \emph{saturated} if $b_p(Y)=2$.
The minimal positive supports $D$ with $b_p(D)=1$ are the
\emph{residual supports}. Each contains exactly one original
occurrence, whose support is $D$. Conversely, the support of any
original occurrence not lying at a saturated support is residual.
These assertions follow directly by counting occurrences contained
in the given support.

\subsubsection{Initials supported on the same set}
\begin{proposition}\label{thm:capped-local}
The varieties generated by $A_{2143}$ and $A_{2146}$ have the
common basis $\SR$ together with
\begin{equation}\label{eq:capped-local-laws}
 x^2\eq x,\qquad xyz\eq xzy,\qquad
 x+y+xy\eq x+y,\qquad x+y+xyz\eq x+y.
\end{equation}
\end{proposition}
\begin{proof}
The word laws reduce a word to $V(D,h)$: its initial $h$ followed
by the other support variables once in a fixed order. In these two
tables 1 is both an additive identity and a multiplicative zero.
A word involving 1 is 1; otherwise it is 2 when its initial is 2,
is 3 when every variable is 3, and is 0 otherwise. On the
idempotent elements the additive chains are $1<2<0$ and $1<0<2$,
respectively, and $3+3=3+0=0$. The pair and triple laws hold
because a factor 1 makes the added product neutral; when neither
of the first two variables is 1, their sum is 0 or 2 and already
absorbs that product. Taking $y=x$ in the pair law gives $3x\eq2x$.

Put $H_p(Y)=\{\head(v):v\text{ is an occurrence supported in }Y\}$.
At every saturated $Y$ and each $h\in H_p(Y)$, choose an original
occurrence $u$ with initial $h$, and a distinct occurrence $v$, both
supported in $Y$. Substituting $u,v$ and a word listing $Y$ into
the triple law inserts $V(Y,h)$. Apply it twice to insert two
copies. The occurrences may be equal words, but their two copies
are retained. These insertions change neither $b_p$ nor $H_p$.

Let $E_p$ be the sum of $2V(Y,h)$ over all these pairs $(Y,h)$.
An original occurrence at a saturated support is removed against
its marker by $3v\eq2v$. The remaining occurrences are exactly
the unique words on the residual supports. Thus we derive the
normal form consisting of $E_p$ and those residual words.

Assign 3 on $Y$ and 1 outside it. The outputs 1, 3, 0 recover
$b_p(Y)=0,1,2$. To recover $H_p(Y)$ in $A_{2143}$, assign 0 to
a selected $h\in Y$, 2 to the other variables in $Y$, and 1
outside. The result is 0 precisely when $h\in H_p(Y)$.
In $A_{2146}$ interchange 0 and 2; the detecting value is now 2.
The profile therefore recovers the initial of each residual word
as well as the entire doubled part. Equal term functions have
the same normal form, proving the common basis.
\end{proof}

\subsubsection{Initials shared across supports}
\begin{proposition}\label{thm:capped-global}
A basis for $A_{2166}$ is $\SR$, \eqref{eq:capped-local-laws}, and
\begin{equation}\label{eq:capped-global-laws}
 2x+yz\eq2x+yz+yx,\qquad
 xu+yv+xyz\eq xu+yv+yxz.
\end{equation}
\end{proposition}
\begin{proof}
The additive identity is 2; the other rules include
$0+3=3$, $3+3=1$, and 1 absorbs all values. A word is 2 if its
initial is 2, is 3 if all variables are 3, is 1 if all variables
are in $\{1,3\}$ with some 1, and is 0 otherwise. Consequently
its polynomial function is determined by the set $H$ of all
initials and the capped profile $b_p$. Indeed, take the
$\{1,3\}$-fibre $G$ of an assignment. If $b_p(G)=2$, the result
is 1. If $b_p(G)=1$, its unique active support determines whether
the active word is 3 or 1. If $b_p(G)=0$, the result is 0 when
an initial is outside the 2-fibre, and is 2 otherwise.
Each displayed identity preserves these data, proving its validity.
In particular the first new law adds a support only where the
two retained $x$ occurrences already saturate the profile.

The preceding proof first inserts doubled words on each saturated
$Y$ with initials originally active there. Let $w$ be such a word.
For any $h\in H\cap Y$, retain an original word $hv$, with $v$
nonempty; a linear word can be written $h^2$. The first law of
\eqref{eq:capped-global-laws} inserts $hw$ beside $2w+hv$.
Repeating supplies $2V(Y,h)$. Hence we obtain the background
\[
 E_p=\sum_{\substack{Y:\ b_p(Y)=2\\h\in H\cap Y}}2V(Y,h).
\]
Remove original words at saturated supports using $3x\eq2x$.
For each residual $D$ let $h(D)$ be the least member of $H\cap D$.
We now change its residual word to $V(D,h(D))$.

If $p$ had at least two occurrences, its full support $X$ is
saturated. The background contains $2V(X,h)$ for every $h\in H$.
To change a residual initial $x$ to $y\in H\cap D$, use one
copy of the full-support words with initials $x$ and $y$ as $xu$
and $yv$ in the second new law. Represent the residual as $xyz$,
where $z$ is a nonempty word listing $D$; the word laws permit
these extra repetitions. This replaces it by $V(D,y)$, retaining
the background. If the original polynomial had one occurrence,
$H$ is its single initial and no change is needed. Thus
\[
 p\eq E_p+\sum_{D\text{ residual}}V(D,h(D)).
\]
Empty parts of this expression are omitted.

Assign 0 to $h$ and 2 elsewhere to detect $h\in H$ by value 0.
Assign 3 on $Y$ and 2 elsewhere to detect $b_p(Y)=1,2$ by values
3, 1; the other values indicate zero. Hence the function recovers
all parameters of the derived normal form. Equal functions reduce
to the same expression, establishing completeness.
\end{proof}
\subsection{First-occurrence sequences and parity}
The use of ordered bands in semiring identity theory goes back to
\cite{OrderedI,OrderedII}. We give a direct proof for the particular
three-element algebra needed here, so that the subsequent coefficient
arguments have an explicit equational foundation.

\subsubsection{The three-element prefix algebra}
Let $B=\{l,e,t\}$ have addition given by maximum in $l<e<t$.
Multiplication has identity $e$ and left zeros $l,t$. Thus a word
takes its first value other than $e$, or $e$ when all its variables
are assigned $e$. Let $i(w)$ be the sequence of first occurrences
of the variables in $w$. The identities
\begin{equation}\label{eq:lrb}
 x^2\eq x,\qquad xyx\eq xy
\end{equation}
reduce $w$ to $i(w)$: a later occurrence of $x$ is deleted by the
second identity if an intervening word is present, and by the first
identity for adjacent occurrences.

For a polynomial $p$ on a finite alphabet $X$, put
\begin{align*}
 A_p(Y)&\Longleftrightarrow
       \text{some summand of }p\text{ is supported in }Y,\\
 H_p(Y)&=\{\text{first letters outside }Y\text{ of summands of }p
                                      \text{ not supported in }Y\}.
\end{align*}
If $Y$ is the $e$-fibre of an assignment and $T$ its $t$-fibre,
the value is $t$ when $H_p(Y)\cap T\ne\varnothing$, and otherwise
is $e$ or $l$ according as $A_p(Y)$ holds or fails. Conversely,
assigning $e$ on $Y$ and $l$ elsewhere recovers $A_p(Y)$.
Assigning $t$ to $h\notin Y$, $e$ on $Y$, and $l$ elsewhere
recovers $h\in H_p(Y)$ by a value $t$.

Let $C(p)$ be the repetition-free nonempty words absorbed by $p$
in $B$. For $w=x_1\cdots x_k$ and
$P_j=\{x_1,\ldots,x_j\}$, the preceding description gives
\begin{equation}\label{eq:prefix-membership}
 w\in C(p)\quad\Longleftrightarrow\quad
 A_p(P_k)\text{ and }x_j\in H_p(P_{j-1})\ (1\leq j\leq k).
\end{equation}
Here $P_0=\varnothing$. To see sufficiency for an arbitrary $Y$,
let $x_j$ be the first letter of $w$ outside $Y$. Its preceding
prefix lies in $Y$. A summand whose first letter outside
$P_{j-1}$ is $x_j$ also has first letter outside $Y$ equal to
$x_j$. If the whole target is supported in $Y$, the supported
summand supplied by $A_p(P_k)$ gives the necessary value.
The recovery assignments prove necessity.

\begin{lemma}\label{lem:prefix-B}
A basis for $B$ consists of $\SR$, additive idempotence,
\eqref{eq:lrb}, and
\begin{equation}\label{eq:B-pair-trim}
 x+y\eq x+y+xy,\qquad xy+z\eq xy+z+xz.
\end{equation}
Moreover, for each fixed support $D$, the words of $C(p)$ on $D$
are connected by adjacent interchanges staying in $C(p)$.
\end{lemma}
\begin{proof}
The two absorption laws are valid because $xy$ is either $x$ or
$y$, and $xz$ is either $xy$ or $z$. To derive an arbitrary valid
absorption, write $w=x_1\cdots x_k\in C(p)$ and choose a summand
$u$ supported in $P_k$. Retain $r_0=u$. Suppose that
$r_{j-1}=x_1\cdots x_{j-1}u$ has been inserted. A witness summand
has the form $a x_j b$ with $\supp(a)\subseteq P_{j-1}$. The
prefix $a$ may be absent. If $b$ is absent, write the witness as
its own square to obtain a nonempty suffix. The trimming law
inserts $a x_j r_{j-1}$. For $j=1$ this is $r_1$; for $j>1$
apply the same law using the split $r_{j-1}=Pu$, where
$P=x_1\cdots x_{j-1}$, to insert
\[
 P a x_j P u\eq P x_j u=r_j.
\]
The equality uses \eqref{eq:lrb}. At the end $r_k\eq w$ because
$u$ is supported in $P_k$. The absorption lemmas remove temporary
words and derive every identity of $B$.

For connectivity, transform a word on $D$ into a target word from
left to right. If their common prefix has support $P$ and the next
target letter is $h$, move $h$ leftwards by adjacent interchanges.
The allowedness condition is monotone:
\[
 h\in H_p(P),\quad P\subseteq Q,\quad h\notin Q
       \quad\Longrightarrow\quad h\in H_p(Q).
\]
Indeed, the witness for $h$ outside $P$ remains a witness outside
$Q$. Thus $h$ is allowed at each intermediate position. A crossed
letter has acquired $h$ in its prefix and remains allowed by the
same implication. All intermediate words satisfy
\eqref{eq:prefix-membership}, proving the assertion.
\end{proof}

\subsubsection{A periodic coefficient extension}
\begin{proposition}\label{thm:periodic-prefix}
The algebras $A_{2027}$ and $A_{2288}$ have the common basis $\SR$,
\eqref{eq:lrb}, and
\begin{gather}
 3x\eq x,\label{eq:periodic-scalar}\\
 2x+2y\eq2x+2y+2xy,\qquad
 2xy+2z\eq2xy+2z+2xz,\label{eq:periodic-lift}\\
 xy+2yx\eq2xy+yx.\label{eq:periodic-swap}
\end{gather}
\end{proposition}
\begin{proof}
Both contain $B$ and $\FF$. In $A_{2027}$ their element labels are
$(l,e,t)=(1,2,0)$ and $(0_F,1_F)=(0,3)$; in $A_{2288}$ they are
$(0,2,1)$ and $(2,3)$, respectively. Doubling maps each algebra
onto its displayed $B$. A word has value 3 exactly when all its
variables are assigned 3. In $A_{2027}$ an odd number of such
values makes the polynomial 3; with even parity its value is that
of the doubled polynomial in $B$. In $A_{2288}$ the same rule
holds except that the $B$-value $t=1$ absorbs the odd value 3.
This follows by addition of the four possible word values.
Thus the $B$-function together with the parity at every support
determines the full function. The proposed laws agree on both
these data: \eqref{eq:periodic-lift} gives the laws of
\eqref{eq:B-pair-trim} after doubling, and
\eqref{eq:periodic-swap} changes no support parity. They are valid.

Equationally, $d(x)=2x$ is an idempotent endomorphism by
\eqref{eq:periodic-scalar}, since $4x\eq2x$ and
$(2x)(2y)\eq4xy\eq2xy$. Its image satisfies the displayed
basis of $B$. The derivation in Lemma~\ref{lem:prefix-B}, applied
to doubled terms, gives
\[
 2p\eq E_p:=\sum_{w\in C(p)}2w,
 \qquad p\eq p+E_p.
\]
For each support $D$ occurring in $C(p)$ choose its least word
$r_D$ in a fixed lexicographic order. An adjacent interchange
$w=a x y b$, $w'=a y x b$ in $C(p)$ gives, by
\eqref{eq:periodic-swap} in its contexts,
$w+2w'\eq2w+w'$. Empty contexts mean no multiplication on that
side. Since $E_p$ contains both doubled words and $3w\eq w$,
\[
 E_p+w\eq E_p+w'.
\]
Connectivity therefore moves every original occurrence to $r_D$
at its support. Pairs of residual copies are absorbed by their
background, giving the normal form
\begin{equation}\label{eq:periodic-normal}
 E_p+\sum_{\varepsilon_p(D)=1}r_D,
 \qquad
 \varepsilon_p(D)=\#\{\text{occurrences on }D\}\pmod2.
\end{equation}

An identity in either source restricts to an identity in $B$, so
its $C$-sets agree. Restriction to $\FF$ gives, on an indicator
assignment for $Y$,
$f(Y)=\sum_{D\subseteq Y}\varepsilon_p(D)$ in $\FF$.
Lemma~\ref{lem:finite-set-recovery}(i) recovers the coefficients,
hence the normal form \eqref{eq:periodic-normal}.
The derivations prove completeness and equality of the two varieties.
\end{proof}
\subsection{Capped coefficients for first-occurrence words}
\begin{proposition}\label{thm:capped-prefix}
A basis for $A_{2287}$ is $\SR$, \eqref{eq:lrb}, and
\begin{gather}
 3x\eq2x,\qquad x+y\eq x+y+xy,\label{eq:cappair}\\
 xy+2z\eq xy+2z+2xz,\label{eq:captrim}\\
 xu+yv+xyz\eq xu+yv+yxz.\label{eq:capswap}
\end{gather}
\end{proposition}
\begin{proof}
Multiplication has identity 3. On $\{0,2,1\}$ it restricts to
$B$ with $(l,e,t)=(0,2,1)$. The additional sums are
$0+3=3$, $2+3=2$, $1+3=1$, and $3+3=2$.
Normalize words by \eqref{eq:lrb}. Define $H_p,A_p,C(p)$ as in
the preceding section and the capped profile $b_p$ by
\eqref{eq:capped-profile}. In particular $A_p(Y)$ means
$b_p(Y)>0$.

These parameters determine evaluations. If $Y$ is the
$\{2,3\}$-fibre, a word not supported in $Y$ is 0 or 1 according
to its first letter outside $Y$. The presence of value 1 forces
the sum to 1. Otherwise, no supported word gives 0; at least two
supported occurrences give 2; exactly one, with support $D$, gives
3 if all of $D$ is assigned 3, and 2 otherwise. The unique active
support is the minimal positive subset of $Y$ in the capped
profile. Conversely the embedded $B$ recovers $H_p$, and assigning
3 on $Y$ and 0 elsewhere recovers $b_p(Y)=0,1,2$ by values 0,3,2.

The identities preserve $H_p$ and $b_p$. The product added by
the pair law is active only when both retained occurrences are
active. The two added occurrences in \eqref{eq:captrim} are
active only when the two retained copies of $z$ are active.
The swap law preserves the changed support, and the two retained
words preserve both possible first letters outside a set.
This proves validity, including all coincidences of variables.

Let $C_2(p)$ consist of the words $w\in C(p)$ such that
$b_p(\supp(w))=2$, and put $E_p=\sum_{w\in C_2(p)}2w$.
We first derive $p\eq p+E_p$. Fix $w=x_1\cdots x_k\in C_2(p)$
with support $D$. Two original occurrences $u,v$ are supported
in $D$. Applying the pair law twice while retaining them inserts
$2uv$. Put $r_0=i(uv)$. The prefix construction in
Lemma~\ref{lem:prefix-B} can now be repeated with doubled words:
\eqref{eq:captrim} first inserts $2i(a x_j r_{j-1})$ and then
$2r_j$, where $r_j=i(x_1\cdots x_jr_0)$. A terminal witness is
written as its own square to keep the substituted suffix nonempty.
At the final step $r_k=w$. Every doubled marker is additively
idempotent. If $p\eq p+K$ and $p+K\eq p+K+m$, then
$p+m\eq p+K+m\eq p+K\eq p$; hence temporary markers can be
eliminated without additive cancellation.

Each original occurrence at a saturated support is removed against
its doubled copy by $3w\eq2w$. There remains one occurrence on
each residual support. Choose the least $r_D\in C(p)$ at each
such support. We show how to replace the residual by $r_D$.
If the original polynomial has one occurrence, its $C$-set on
that support contains only that word, by
\eqref{eq:prefix-membership}. Otherwise its whole support $X$
is saturated. Concatenating all original words and deleting
repetitions gives a full-support member $v$ of $C(p)$.

Consider a permitted adjacent interchange $PxyQ\mapsto PyxQ$
along a path from the residual to $r_D$. There are full-support
words in $C(p)$ beginning with $Px$ and $Py$: append $v$ to
each prefix, and use monotonicity of allowedness as in
Lemma~\ref{lem:prefix-B}. Their doubled markers lie in $E_p$.
Write them as $Pxu$ and $Pyv'$ with nonempty suffixes, using
idempotence if necessary. Multiply \eqref{eq:capswap} on the
left by $P$, omitting the context if absent, and substitute
$z=Q$ when $Q$ is nonempty. If $Q$ is absent, use $z=xy$;
the added repetitions disappear by \eqref{eq:lrb}. One copy of
each marker supplies the retained summands. The identity changes
the residual and leaves the doubled background in place. If the
markers coincide, its two copies still supply the two summands.
Connectivity completes the replacements.

We have derived
\[
 p\eq E_p+\sum_{D\text{ residual}}r_D,
\]
with absent parts omitted. The function recovers $H_p,b_p$, which
determine $C_2(p)$ and the chosen residuals. Equal functions have
the same derived normal form, proving completeness.
\end{proof}
\subsection{A neutral exceptional element}
We now keep linear words distinct from nonlinear words. The common
word identities are
\begin{equation}\label{eq:stable-word}
 x^2y\eq xy,\qquad xy^2\eq xy,\qquad xyx\eq xy.
\end{equation}
They reduce a nonlinear word with at least two distinct variables
to its first-occurrence sequence. A nonlinear word of singleton
support reduces to $x^2$. Denote this nonlinear representative of
a sequence $s$ by $s^{\mathrm n}$. In particular $w^2\eq w$
for nonlinear $w$.

\begin{proposition}\label{thm:stable-prefix}
A basis for $A_{1819}$ consists of $\SR$, \eqref{eq:stable-word},
and
\begin{gather}
 2x\eq x^2,\quad3x\eq2x,\quad x+y\eq x+y+xy,
 \label{eq:stable-scalar}\\
 x+yz\eq x+yz+xyz,\qquad
 xy+zu\eq xy+zu+xzu.\label{eq:stable-trims}
\end{gather}
A basis for $A_{1822}$ is obtained by replacing
\eqref{eq:stable-scalar} with $2x\eq x$ and $x+x^2\eq x^2$.
\end{proposition}
\begin{proof}
In both sources $\{1,0,2\}$ is $B$ with $(l,e,t)=(1,0,2)$.
The exceptional element 3 acts as 0 in products of length at
least two, with $3^2=0$. Thus a nonlinear word takes its first
value outside $\{0,3\}$, or 0 if there is none. Linear words
retain their assigned values. The sources differ in $3+3$: it is
0 in $A_{1819}$ and 3 in $A_{1822}$.
The word laws and the displayed scalar laws are valid. For linear
sharing and nonlinear trimming, the added word has a retained
nonlinear seed with support contained in its support. Its first
exceptional value is already represented by a retained summand.
The pair law additionally uses $3+3=0$. These observations verify
all proposed identities.

Every nonlinear word is additively idempotent, and
$x+x^2\eq x^2$ holds in both bases. For a polynomial $p$ define
$H_p(Y)$ using all its first letters outside $Y$, including those
of linear words. Let $S_p(D)$ hold when a nonlinear summand is
supported in $D$; for $A_{1819}$ also allow two linear occurrences
supported in $D$, with multiplicity counted. Put
$Q=\{x:S_p(\{x\})\}$ and $M=L\cup Q$, where $L$ is the set
of linear variables. Define $C_p$ to consist of the nonlinear
standard words $w$ for which $S_p(\supp(w))$ holds and every next
distinct letter $h$ is in $H_p(P)$ for its preceding prefix
support $P$. The intended normal form is
\begin{equation}\label{eq:stable-normal}
 N(p)=\sum_{w\in C_p}w+\sum_{x\in M\setminus Q}x.
\end{equation}

Fix $w\in C_p$. Choose a nonlinear seed $r_0$ supported in
$\supp(w)$. It is an original nonlinear word, or is inserted from
two linear occurrences using the pair law or $2x\eq x^2$.
Write the distinct-letter sequence of $w$ as $x_1\cdots x_k$.
We construct $r_j=(x_1\cdots x_jr_0)^{\mathrm n}$ while retaining
the original polynomial. At step $j$, a linear witness must be
$x_j$ itself; linear sharing in \eqref{eq:stable-trims} inserts
$x_jr_{j-1}$, since $r_{j-1}$ splits into two nonempty factors.
A nonlinear witness has form $a x_j b$, with
$\supp(a)\subseteq\{x_1,\ldots,x_{j-1}\}$.
If $b$ is absent, replace the witness by its own square. Nonlinear
trimming inserts $a x_jr_{j-1}$. For $j>1$, use that law again
with the split $r_{j-1}=Pr_0$ to insert
\[
 P a x_j P r_0\eq P x_jr_0=r_j,
 \qquad P=x_1\cdots x_{j-1}.
\]
The reductions use \eqref{eq:stable-word}; the result is nonlinear
even when its support is a singleton. At the last step we obtain
$w$. All inserted terms are additively idempotent, so the temporary
summands are removed by the absorption calculation already used
above. Repeating inserts all of $C_p$. Original nonlinear words
are absorbed against their copies. Repeated linears become squares
in the nonidempotent source, and square domination removes exactly
the linears in $Q$. This proves $p\eq N(p)$.

To recover $H_p(Y)$, assign 2 to a selected variable outside $Y$,
0 on $Y$, and 1 elsewhere, using the embedded $B$.
Assigning 3 on $D$ and 1 elsewhere gives 0 precisely when
$S_p(D)$ holds. Otherwise the value is 3 or 1. For $D=\{x\}$,
these outputs distinguish $x\in Q$, $x\in M\setminus Q$, and
absence. Thus the term function determines $H_p,S_p,M$, hence
the normal form \eqref{eq:stable-normal}. Equal functions derive
to the same polynomial, proving both bases.
\end{proof}
\subsection{An exceptional element below the additive top}
The next multiplication restricts to $B$ on $\{1,2,0\}$ with
$(l,e,t)=(1,2,0)$. The fourth element 3 acts like 0 in every
nonlinear product. Its additive behaviour must still be retained
for linear terms.

\begin{proposition}\label{thm:exceptional-prefix}
A basis for $A_{793}$ consists of $\SR$, \eqref{eq:stable-word},
and
\begin{gather}
 2x\eq x^2,\qquad3x\eq2x,\label{eq:exception-scalar}\\
 x+xy\eq x^2+xy,\label{eq:exception-promote}\\
 x+yxz\eq x+yxz+yx^2,\label{eq:exception-seed}\\
 xy+zu\eq xy+zu+xzu.\label{eq:exception-trim}
\end{gather}
For $A_{799}$ replace \eqref{eq:exception-scalar} by
$2x\eq x$ and $x+x^2\eq x^2$.
\end{proposition}
\begin{proof}
In $A_{793}$, $3+3=0$; in $A_{799}$, $3+3=3$.
Both satisfy $1+3=2+3=3$ and $0+3=0$.
The word laws are valid by the first-exception rule. The promotion
law replaces a linear value by its square in the presence of a
nonlinear word beginning there. For the seed law, if the prefix
$y$ already exits the neutral fibre its value is retained in
$yxz$; otherwise $x$ and $yx^2$ have the needed common
exceptional value after squaring. The trimming law acts on
nonlinear values in $B$ and has the same validity argument as
\eqref{eq:B-pair-trim}. Together with the scalar rules this proves
soundness of the list.

Preprocess a polynomial so that linears occur once: in $A_{793}$
replace repeated linears by their squares, and in $A_{799}$ use
additive idempotence. All nonlinear words are additively
idempotent and $x+x^2\eq x^2$ is available. Let $L$ be the
remaining linear variables, $N$ the nonlinear standard words,
$A(Y)$ assert the presence of any summand supported in $Y$, and
$H_N(Y)$ record only the first letters outside $Y$ of words in
$N$. Define $C$ to be the nonlinear standard words $w$ with
$A(\supp(w))$ and with each next distinct letter in $H_N(P)$
at its preceding prefix $P$. Put
\[
 M=\{x:A(\{x\})\},\quad Q=M\cap H_N(\varnothing),\quad
 N(p)=\sum_{w\in C}w+\sum_{x\in M\setminus Q}x.
\]

For $w\in C$ with support $D$, a nonlinear supported summand is
already a suitable seed. If none exists, $A(D)$ supplies a linear
$x\in D$. The transition at $x$ in $w$ supplies a nonlinear
witness $a x b$ with $\supp(a)\subseteq D$. When $a,b$ are both
nonempty, \eqref{eq:exception-seed} inserts $a x^2$. When only
$a$ is nonempty, the witness itself is supported in $D$. When
$a$ is absent, write the witness as $xy$ with nonempty $y$.
The promotion law gives
$x+xy\eq x^2+xy\eq x+x^2+xy$, using square domination;
thus it inserts the seed $x^2$ while retaining the linear word.

All transition witnesses are nonlinear. The two-step prefix
construction of the preceding proposition, now using
\eqref{eq:exception-trim}, successively inserts
$(x_1\cdots x_jr_0)^{\mathrm n}$ and finally $w$.
Nonlinear witnesses with empty tails are again squared before
substitution. Retaining and then removing temporary nonlinear
words proves absorption of the whole set $C$. Each original
nonlinear word is in $C$; remove its duplicate. A singleton
square is in $C$ exactly for variables in $Q$, so square
domination removes exactly those linears. This derives $N(p)$.

Assign 2 on $Y$ and 1 outside. The result is 2 precisely when
$A(Y)$ holds. For $h\notin Y$, assign $h=3$, 2 on $Y$, and
1 elsewhere. A nonlinear word gives 0 precisely when its first
letter outside $Y$ is $h$. Other words give 1 or 2, and a
remaining linear $h$ gives 3. There is at most one such linear
after preprocessing; its sum with 1 or 2 stays 3. Thus the
polynomial is 0 exactly when $h\in H_N(Y)$.
These assignments recover $A,H_N$, and hence $M,Q,C$.
The derived normal form is unique as a term function, establishing
completeness for both generators.
\end{proof}
\subsection{A prefix profile that includes linear summands}
\begin{proposition}\label{thm:completed-prefix}
A basis for $A_{790}$ is the basis displayed for $A_{793}$ in
Proposition~\ref{thm:exceptional-prefix}, together with
\begin{equation}\label{eq:append-linear}
 x+yz\eq x+yz+yzx.
\end{equation}
\end{proposition}
\begin{proof}
Multiplication is the same as for $A_{793}$. The exceptional sums
are now $1+3=3$, $2+3=0$, $3+3=0$, and $0+3=0$.
The earlier laws remain valid. For \eqref{eq:append-linear}, if
$yz$ has already taken a nonneutral value in $B$, then
$yzx=yz$. Otherwise $yz=2$, and $2+x$ absorbs its nonlinear
product with $x$, including $x=3$ because $2+3=0$.

Use the same preprocessing and notation $L,N,A,H_N$ as in the
preceding proof. The observable profile is
\begin{equation}\label{eq:completed-profile}
 H^*(Y)=H_N(Y)\cup
 \begin{cases}L\setminus Y,&A(Y)\text{ holds},\\
               \varnothing,&A(Y)\text{ fails}.
 \end{cases}
\end{equation}
In particular a supported linear word can supply the second
clause. Let $C$ be defined by the same support and successive
letter conditions, using $H^*$ instead of $H_N$. Put
$M=\{x:A(\{x\})\}$ and $Q=M\cap H^*(\varnothing)$.
We derive the normal form $\sum_{w\in C}w+\sum_{x\in M\setminus Q}x$.

Let $C_0$ use $H_N$ alone. The preceding seed and saturation
argument uses only laws in the present basis, so every member
of $C_0$ can already be inserted. For a word in $C\setminus C_0$,
take its first transition $h$ supplied only by the second clause
of \eqref{eq:completed-profile}. Its preceding distinct-letter
prefix $P$ is nonempty, because $A(\varnothing)$ is false.
Every transition of $P$ uses $H_N$, and $A(\supp(P))$ holds.
Thus $P^{\mathrm n}\in C_0$ is available. The retained linear
$h$ and \eqref{eq:append-linear} insert $P^{\mathrm n}h$.

Continue along the target. A transition with an original nonlinear
witness is supplied by the two-step trimming construction, taking
the current nonlinear prefix as seed; it adds no unwanted suffix.
A transition using only the second clause is appended from its
retained linear word by \eqref{eq:append-linear}. At the end
the target is present. All intermediate words are nonlinear and
additively idempotent, so they can be removed by the absorption
calculation. Inserting the finite set $C$, deleting nonlinear
duplicates and applying square domination gives the claimed
normal form. Since $H^*(\varnothing)=H_N(\varnothing)$, its
remaining linears are exactly $M\setminus Q$.

The assignment 2 on $Y$, 1 elsewhere, recovers $A(Y)$.
For $h\notin Y$, assign $h=3$, 2 on $Y$, and 1 elsewhere.
The output is 0 in exactly two situations: a nonlinear word has
first letter outside $Y$ equal to $h$, or an arbitrary summand
is supported in $Y$ and the remaining linear $h$ is present.
The latter gives $2+3=0$ and includes a supported linear
summand. In all other cases the values 1,2 and at most one
linear 3 cannot yield 0. This recovers precisely $H^*(Y)$.
Hence equal term functions have the same canonical form, and
the proposed identities form a basis.
\end{proof}
\subsection{Factoring terminal variables}
\begin{proposition}\label{thm:factorization}
A basis for $A_{464}$ consists of $\SR$ and the following identities:
\begin{subequations}\label{eq:factor-laws}
\begin{align}
 x+x&\eq x,\label{eq:f1}\\
 xyz&\eq xy^2z,&xyzu&\eq xzyu,\label{eq:finterior}\\
 xyx&\eq x^2y^2,&x^3&\eq x^2,\label{eq:fsquare}\\
 x^2+y^2&\eq x^2y^2,\label{eq:f6}\\
 xyz&\eq xyz+y^2z,\label{eq:f7}\\
 x^2+y&\eq x^2+yx^2,\label{eq:f8}\\
 xy+zu&\eq xy+zu+xu,\label{eq:f9}\\
 xy+yz&\eq xy+yz+y^2,\label{eq:f10}\\
 xy+zyu&\eq xy+zyu+y^2,\label{eq:f11}\\
 x+yxz&\eq x+yxz+x^2,\label{eq:f12}\\
 x+xy&\eq x+xy+y^2,\label{eq:f13}\\
 x^2+yz&\eq x^2+yz+z^2.\label{eq:f14}
\end{align}
\end{subequations}
\end{proposition}
\begin{proof}
Write $(0,e,l,r)=(0,1,2,3)$. The nonzero products are
$e^2=e$, $le=l$ and $er=r$. A word is $e$ when all letters
are $e$, is $l$ when its unique initial occurrence is assigned
$l$ and all other occurrences are $e$, and is $r$ under the
corresponding unique-terminal condition. It is zero otherwise.
Addition has top 0, with $e+l=l$ and $e+r=l+r=0$.
These four word-value alternatives verify \eqref{eq:factor-laws}.
For example, a nonzero square has value $e$, so adjoining a square
to $r$ already gives 0; the extra squares in
\eqref{eq:f10}--\eqref{eq:f14} cannot introduce a new nonzero
exception without a retained summand already forcing 0.

For nonempty $D$ put $b_D=\prod_{x\in D}x^2$ in a fixed order.
Square factors commute by \eqref{eq:f6}, and
\[
 (uv)^2=(uvu)v\eq u^2v^3\eq u^2v^2.
\]
Thus the square of a word with support $D$ is $b_D$, and
$b_D+b_E\eq b_{D\cup E}$. Also $x+x^2\eq x^2$ follows
from \eqref{eq:f8} with $y=x$. Interior letters may be squared
and permuted using \eqref{eq:finterior}. Repeated endpoints are
collected using these laws and \eqref{eq:fsquare}. Consequently,
if all noninitial occurrences of a word $v$ belong to $D$, then
\begin{equation}\label{eq:factor-reduction}
 vb_D\eq
 \begin{cases}\head(v)b_D,&\head(v)\notin D,\\
 b_D,&\head(v)\in D.
 \end{cases}
\end{equation}
Indeed, square the interior, move equal letters together, and
commute square factors. In the second case an appended occurrence
of the initial makes it square as well. The cubic law deletes
all remaining repeated square factors.

For a polynomial $p$ let $X$ be its support, and let $H$ be the
variables occurring only as a unique initial letter in every
summand containing them. If every summand has a unique terminal
variable and no such terminal variable occurs nonfinally in any
summand, let $T$ be the set of these variables; otherwise set
$T=\bot$. The word-value rules give the following complete
description. The value of $p$ is $e$ exactly when all of $X$ is
assigned $e$. It is $l$ exactly when variables outside $H$ are
$e$, variables in $H$ are in $\{e,l\}$, and at least one is $l$.
It is $r$ exactly when $T\ne\bot$, all of $T$ is assigned $r$,
and all of $X\setminus T$ is assigned $e$. Other assignments
give 0. In particular two variables of $H$ cannot occur together
in a word. A value $r$ requires every summand to be $r$, since
$e+r=l+r=0$.

\emph{Case 1: $T=\bot$.} Put $D=X\setminus H$. We first
insert a nonempty square seed supported in $D$. If a summand has
a repeated terminal variable, its terminal segment $xux$ becomes
a square by \eqref{eq:fsquare}; an adjacent $xx$ already is one.
If a prefix remains, use \eqref{eq:f7}, writing the square suffix
$s$ as $s^2$, to insert $s^3\eq s$. All variables of the suffix
occur noninitially, so its support lies in $D$.
Otherwise some unique terminal $x$ occurs nonfinally in another
summand. Write the two summands as $Ax$ and $BxC$, with $C$
nonempty. If $A,B$ are both nonempty, \eqref{eq:f11} inserts
$x^2$. If only $A$ is nonempty use \eqref{eq:f10}; if only $B$
is nonempty use \eqref{eq:f12}. If both are absent,
\eqref{eq:f13} inserts $C^2$. This supplies the seed in every case.

Retain it and apply \eqref{eq:f14} to each original nonlinear
summand split after its first letter. The inserted suffix squares
have union of supports exactly $D$. Combining them by
\eqref{eq:f6} gives $b_D$, so $D$ is nonempty and $p\eq p+b_D$.
Now \eqref{eq:f8} and distributivity give
$p\eq b_D+pb_D$. Applying \eqref{eq:factor-reduction} to each
summand yields
\begin{equation}\label{eq:factor-absent}
 p\eq b_D+\sum_{h\in H}hb_D.
\end{equation}

\emph{Case 2: $T\ne\bot$.} Put $I=H\cap T$,
$H'=H\setminus T$, $T'=T\setminus H$, and
$D=X\setminus(H\cup T)$. Each variable of $I$ occurs only
as an isolated linear summand. Remove those summands temporarily.
If nothing remains, $\sum_{i\in I}i$ is the normal form.
Otherwise every remaining summand is $u_it_i$, with $u_i$
nonempty and $t_i\in T'$, and no variable of $T$ occurs in a
prefix $u_i$. The cross-product law \eqref{eq:f9} inserts all
$u_it_j$. Distributivity gives the product of the sum of the
prefixes and $\sum_{t\in T'}t$.

If $D$ is empty, each prefix is a single variable of $H'$, so
the normal form is
$\sum_{i\in I}i+(\sum_{h\in H'}h)(\sum_{t\in T'}t)$.
If $D$ is nonempty, every member of $D$ occurs noninitially in
some prefix. Split that prefix after its first letter and use
\eqref{eq:f7} to insert its square suffix times $t_i$.
Use \eqref{eq:f9} to supply every terminal factor. Combining
the squares with a common terminal factor by \eqref{eq:f6}
produces $b_Dt$ for each $t\in T'$. Then \eqref{eq:f8},
multiplied on the right by $t$, replaces
$b_Dt+u_it$ by $b_Dt+u_ib_Dt$. Applying
\eqref{eq:factor-reduction} gives the normal form
\begin{equation}\label{eq:factor-present}
 \sum_{i\in I}i+
 \left(b_D+\sum_{h\in H'}hb_D\right)\left(\sum_{t\in T'}t\right).
\end{equation}
The empty-$D$ and purely linear cases were treated separately,
so all products here are nonempty terms.

Finally, assign 0 to a selected variable and $e$ elsewhere to
recover $X$. Assign $l$ to it and $e$ elsewhere to recover $H$.
Among assignments giving $r$ on a subset of $X$ and $e$ on its
complement, a value $r$ occurs for exactly the subset $T$ when
$T\ne\bot$, and never otherwise. These tests recover all
parameters in the derived normal forms. Equal term functions
therefore give the same form, proving the basis.
\end{proof}
\subsection{Normal forms consisting of one- and two-letter words}
The three source multiplications in this section coincide. Elements
0,1,2 are left zeros; $3\cdot1=2$ and $3\cdot a=0$ for
$a\ne1$. Every binary product is in $\{0,1,2\}$, so
$xyz\eq xy$. The additive orders are $0<3<2<1$ for $A_{1371}$,
$2<3<0$ and $1<0$, with 1 incomparable to 2 and 3, for
$A_{1372}$, and $1<2<3<0$ for $A_{1373}$.

\begin{proposition}\label{thm:two-letter}
All three bases contain $\SR$, $x+x\eq x$, and $xyz\eq xy$.
For $A_{1371}$ adjoin
\begin{equation}\label{eq:graph-low}
 x+x^2\eq x,\quad xy\eq xy+x^2,\quad
 xy+yz\eq x^2+yz,\quad xy+y\eq x^2+y.
\end{equation}
For $A_{1372}$ adjoin
\begin{equation}\label{eq:graph-separated}
 x+x^2\eq x^2,\quad x^2+xy\eq x^2,\quad
 xy+yz\eq x^2+yz,\quad xy+y\eq x^2+y.
\end{equation}
For $A_{1373}$ only the first two identities of
\eqref{eq:graph-separated} are required.
\end{proposition}
\begin{proof}
The displayed multiplication and additive orders verify these
identities. Expand and shorten a term to a nonempty set of
linears $h$ and pairs $hk$. Let $L$ be its linear variables,
$E$ its directed pairs, including loops, and $H$ its set of
heads. Every head occurs in $L$ or as a first coordinate of $E$.

For $A_{1371}$ put $F=\{(h,k)\in E:k\notin H\}$. The first
two laws of \eqref{eq:graph-low} insert a square at each head.
Whenever an edge $hk$ points into $H$, use
$hk+k^2\eq h^2+k^2$ to remove it while retaining the head
markers. Remove square markers at linear heads using
$h+h^2\eq h$. This derives
\begin{equation}\label{eq:lowgraph-normal}
 N(H,L,F)=\sum_{h\in L}h+\sum_{h\in H\setminus L}h^2
                      +\sum_{(h,k)\in F}hk.
\end{equation}
Assign $h=1$, others 0, to detect $h\in H$ by value 1.
For $h\in H$, assign $h=3$, others 0, to detect $h\in L$
by value 3. For $k\notin H$, assign $h=3,k=1$, others 0;
the result is 2 precisely when $(h,k)\in F$. This recovers
every parameter of \eqref{eq:lowgraph-normal}.

For $A_{1373}$ let $D=\{h:(h,h)\in E\}$,
$L'=L\setminus D$, and $F=\{(h,k)\in E:h\notin D\}$.
At each original loop, square domination removes the linear and
all outgoing pairs. The derived normal form is
\begin{equation}\label{eq:othergraph-normal}
 N(H,L',D,F)=\sum_{h\in L'}h+\sum_{h\in D}h^2
                           +\sum_{(h,k)\in F}hk.
\end{equation}
Assign $h=0$, others 1, to detect $H$ by 0. For $h\in H$,
assign $h=3$, others 1; the value is 0 for $h\in D$, 3 for
$h\in L'$, and 2 otherwise. Finally, for $h\in H\setminus D$
and $k\ne h$, assign $h=3,k=2$, others 1. The value is 0
precisely for $(h,k)\in F$. A head assigned 2 contributes only
2, so it does not create a spurious zero.

For $A_{1372}$ instead put
$D=\{h:(h,k)\in E\text{ for some }k\in H\}$ and define
$L',F$ as above. If an edge $hk$ points into $H$, a summand
at head $k$ is either $k$ or $ku$. The last two identities of
\eqref{eq:graph-separated} replace $hk$ by $h^2$, retaining
that summand. Square domination then removes the other outgoing
terms at $h$. Heads never disappear: a processed head remains
as its square and can still serve as a witness. Thus each member
of the original $D$ is processed, while no head outside $D$
changes. This derives \eqref{eq:othergraph-normal}, now with
all targets of $F$ outside $H$.

Recover $H$ by $h=0$, others 1. For a selected $h\in H$,
assign $h=3$, other heads 2, and variables outside $H$ the
value 1. Values 0,3,2 distinguish $D,L'$, and the remaining
heads. For $h\in H\setminus D$ and $k\notin H$, change
$k$ from 1 to 2; the result becomes 0 precisely when
$(h,k)\in F$. All other heads contribute 2 independently
of their tails. These tests recover the entire normal form.
In all three cases equal term functions therefore reduce to the
same polynomial, proving completeness.
\end{proof}
\subsection{Global terminal variables}
\begin{proposition}\label{thm:global-tails}
The algebras $A_{967},A_{969},A_{1018}$ have the common basis
$\SR$ together with
\begin{subequations}\label{eq:global-tail-laws}
\begin{align}
 x+x&\eq x,\label{eq:gtadd}\\
 x^2y&\eq xy,&xyz&\eq xy^2z,\label{eq:gtword1}\\
 xyzu&\eq xzyu,&xyx&\eq xy^2,\label{eq:gtword2}\\
 xy&\eq xy+x^2,\label{eq:gtseed}\\
 xy+zu&\eq xy+zu+xzu,\label{eq:gttrim}\\
 x^2+yz&\eq x^2+yx^2z,\label{eq:gtmerge}\\
 x+xy&\eq x^2+xy,&x+yxz&\eq x^2+yxz,\label{eq:gtpromote}\\
 x^2+y&\eq x^2+y+x^2y.\label{eq:gtappend}
\end{align}
\end{subequations}
\end{proposition}
\begin{proof}
For $A_{967},A_{969}$ a nonlinear word is 2 when its head is 2;
otherwise it is 1 when every letter is 1, is 3 for a sole
exceptional occurrence in the final position assigned 3, and is
0 otherwise. Their orders are respectively $1<3<0$, $2<0$
with 2 incomparable to 1,3, and $1<3<0<2$.
For $A_{1018}$ a word on $\{1,2\}$ takes its head value; a
sole final 3 with earlier letters in $\{1,2\}$ gives 3;
all other nonlinear cases give 0. Its order is $2<1<3<0$.
Linears retain their values in all cases. These rules verify
the displayed laws: the word laws preserve head, support and
the possibility of a unique exceptional terminal; the additive
laws add or replace a word only when its exceptional values are
already forced by a retained summand.

Write $b_D=\prod_{x\in D}x^2$ for nonempty $D$, with fixed
ordering. The word laws square and permute the interior and
cap repetitions, since $x^3\eq x^2$. A nonlinear word with
head $h$ and support $C$ consequently reduces to
$hb_{C\setminus\{t\}}t$ if its terminal variable $t$ is unique,
and to $hb_C$ otherwise. The first expression has a nonempty
middle support because $h\ne t$. To obtain the second, collect
an internal occurrence of the repeated terminal next to its
last occurrence, or use $huh\eq hu^2$ if it repeats only
the head. These reductions introduce no empty substitutions.
For $h\in D$ put $R(h,D)=hb_D$. It is multiplicatively
idempotent, and products of such repeated-terminal words
normalize by uniting supports and retaining the first head.

For a polynomial $p$ let $X$ be its support, $H$ its heads,
and $T$ the variables that occur only as a unique terminal
in every summand containing them. Set
$I=H\cap T$, $K=H\setminus T$, and $D=X\setminus T$.
If $K$ is empty, every summand is linear. Otherwise $D$ is
nonempty, contains $K$, and is the union of the supports of
the nonfinal parts of the original words. We derive
\begin{equation}\label{eq:global-tail-normal}
 N(X,H,T)=\sum_{i\in I}i+\sum_{h\in K}R(h,D)
                  +\sum_{\substack{h\in K\\t\in T}}R(h,D)t.
\end{equation}
If $K=\varnothing$, only the first sum is used.

First replace every linear $h\notin T$ by $h^2$. A nonfinal
occurrence of $h$ in another summand supplies one of the two
contexts in \eqref{eq:gtpromote}. The remaining linears are
precisely those in $I$. By \eqref{eq:gtseed} each nonlinear
head gives an absorbed square marker $h^2$.
Fix $h\in K$, and split an original nonlinear summand as
$uv$ with $v$ its last letter. Two uses of \eqref{eq:gttrim}
insert $uh^2$ and then $huh^2$. The second normalizes to
$R(h,\{h\}\cup\supp(u))$. The supports of these prefixes
$u$ cover $D$. Two markers $r,s$ with the same head are
multiplicatively idempotent; \eqref{eq:gtmerge} gives
$r+s\eq r+srs$. Adding $s$ shows that the union marker
$srs$ may be inserted while both old markers are retained.
Successive unions therefore insert $R(h,D)$ for every $h\in K$.

Keep the full markers. For each nonlinear $w=yz$ with head
$h$, use \eqref{eq:gtmerge} with $x=R(h,D)$ to replace $w$
by $yR(h,D)z$. The word normal form is $R(h,D)$ if the
original terminal is in $D$, and $R(h,D)t$ otherwise.
For a missing pair with $t$ linear, \eqref{eq:gtappend}
inserts $R(h,D)t$. Otherwise some retained word is $R(g,D)t$;
\eqref{eq:gttrim}, using $R(h,D)=hb_D$, inserts
$hR(g,D)t\eq R(h,D)t$. Thus every summand of
\eqref{eq:global-tail-normal} has been supplied. Removing
duplicates and temporary words derives that normal form.

Assign 0 to a selected variable and 1 elsewhere to recover $X$
by value 0. Assign 3 to it and 1 elsewhere; the value is 3
for membership in $T$, 1 for absence, and 0 otherwise.
For $A_{969}$, $h=2$ and others 1 detects heads by value 2.
For $A_{967}$, $h=1$ and others 2 gives a value different
from 2 precisely when $h$ is a head. The same assignment in
$A_{1018}$ detects heads by value 1. These tests recover
$X,H,T$ in each source. Consequently equal functions give
identical derived normal forms, proving the common basis.
\end{proof}
\subsection{Bounded witnesses for polynomial absorption}\label{sec:witnesses}
\subsubsection{A general finite basis principle}
For a fixed alphabet $Y$ of size $k$, let
\[
W(k,\ell)=\bigcup_{j=1}^{\ell}Y^j.
\]
For an ai-semiring $A$, define the finite set
\begin{equation}\label{eq:bounded-set}
\begin{split}
\Bnd_A(k,\ell,d)=\bigl\{
&\textstyle\sum_{u\in U}u+w\eq\sum_{u\in U}u:\quad
U\subseteq W(k,\ell),\quad 1\leq |U|\leq d,\\
&w\in W(k,\ell),\quad
A\models\textstyle\sum_{u\in U}u+w\eq\sum_{u\in U}u
\bigr\}.
\end{split}
\end{equation}
Choose a fixed ordering of words when writing a sum. If
$M=\sum_{j=1}^{\ell}k^j$, then
\[
|\Bnd_A(k,\ell,d)|\leq M\sum_{j=1}^{d}\binom{M}{j}.
\]
For finite $A$ this is an effective finite specification: membership
is decided by evaluating the displayed identity on all assignments.
It is not necessary to assume that $A$ is finitely based in order to
define or decide membership in this set.

\begin{theorem}[Short words and bounded witness words]\label{thm:bounded}
Let $A$ be a finite ai-semiring and let $\Gamma\subseteq\Id(A)$ be
finite and contain $\SR$ and additive idempotence. Suppose that:
\begin{enumerate}[label=(\roman*)]
\item $\Gamma$ reduces every term to a sum of words of length at most
$\ell$;
\item whenever such a polynomial $p$ absorbs a word $w$ of length at
most $\ell$, a nonempty subpolynomial of $p$ with at most $d$ summands
already absorbs $w$;
\item these witness words together with $w$ involve at most $k$ variables.
\end{enumerate}
Then $\Gamma\cup\Bnd_A(k,\ell,d)$ is a finite basis for $A$.
\end{theorem}
\begin{proof}
Let $p\eq q$ hold in $A$, and shorten both sides using $\Gamma$.
For each summand $w$ of $q$, choose the witness words provided by (ii).
Injectively rename the at most $k$ variables of their absorption into
$Y$. The result belongs to \eqref{eq:bounded-set}. Substituting back
and adding the remaining summands derives $p+w\eq p$. Apply
Lemma~\ref{lem:mutual}.
\end{proof}

In our applications, shortening follows from the single identity
\begin{equation}\label{eq:shortening}
xyzu\eq xyz+xyu+xzu.
\end{equation}
For a word of length $m\geq4$, substitute its first three letters
for $x,y,z$ and its remaining nonempty suffix for $u$.
The resulting words have lengths $3,m-1,m-1$. Induction on length
reduces every polynomial to words of length at most three.

\subsubsection{Undirected occurrence graphs}
For a polynomial $p$, let $X_p$ be its support, $H_p$ its set of heads,
and
\[
\begin{aligned}
Q_p&=\bigcup_{v\in p}\rep(v),\\
E_p&=\bigl\{\{x,y\}:x\ne y,\ x,y\notin Q_p,\\
 &\qquad\text{some summand contains both }x,y\bigr\}.
\end{aligned}
\]
Thus $E_p$ ignores edges incident with a variable repeated in a
summand. Put $J_p=H_p\setminus Q_p$.

\begin{proposition}\label{prop:undirected}
For each of $A_{1504},A_{1506},A_{1960}$ the polynomial function of
$p$ is determined, and uniquely distinguished, by
$(X_p,H_p,Q_p,E_p)$. For $A_{1400}$ it is determined and uniquely
distinguished by $(X_p,J_p,Q_p,E_p)$.
\end{proposition}
\begin{proof}
The tables are in Appendix~\ref{app:tables}. In $A_{1504}$ and
$A_{1506}$ a word is 2 when its head is assigned 2; otherwise it is
1 when all letters are 1, is 3 when exactly one occurrence is 3
and the others are 1, and is 0 in all other cases. Their additive
orders are, respectively,
\[
1<3<0,\quad 2<0\quad(2\text{ incomparable with }1,3),
\qquad 1<3<0<2.
\]
In $A_{1960}$ an all-$\{1,2\}$ word takes its head value. A sole
occurrence of 3 with the other occurrences in $\{1,2\}$ gives 3;
a zero or two occurrences of 3 give 0. Its additive order is
$2<1<3<0$. In $A_{1400}$ an all-1 word gives 1, a single exceptional
occurrence from $\{2,3\}$ gives 3 if it is an initial 3 and gives
2 otherwise, and all other cases give 0. Its additive order is
$1<2<3<0$. These descriptions include linear words and follow by
successive multiplication.

These rules prove that the stated parameters determine evaluations.
In particular, repeated exceptional variables are detected by $Q_p$,
and pairs of exceptional variables occurring together are detected
by $E_p$. The distinction involving heads is retained by $H_p$;
for $A_{1400}$ it is immaterial on $Q_p$, and $J_p$ suffices.

They can also be recovered. Assigning $h=0$ and all other variables
1 detects $h\in X_p$ by a zero result. Assigning $h=3$ and the
others 1 detects $h\in Q_p$ by a zero result. For distinct
$h,k\notin Q_p$, assign $h=k=3$ and the others 1; the result is
zero exactly when $\{h,k\}\in E_p$. In $A_{1506}$, assigning
$h=2$ and the others 1 gives 2 exactly when $h\in H_p$.
In $A_{1504}$, assigning $h=1$ and all others 2 gives a value
different from 2 exactly when $h\in H_p$; in $A_{1960}$ the same
assignment gives 1 exactly in this case. Finally, for
$h\in X_p\setminus Q_p$ in $A_{1400}$, assigning $h=3$ and the
others 1 gives 3 for $h\in J_p$ and 2 otherwise.
\end{proof}

It follows that $p$ absorbs $w$ in a full-head source precisely when
\begin{equation}\label{eq:graph-absorption}
\begin{gathered}
\supp(w)\subseteq X_p,\qquad \head(w)\in H_p,\qquad
\rep(w)\subseteq Q_p,\\
\{x,y\}\in E_p\quad
\text{for all distinct }x,y\in\supp(w)\setminus Q_p.
\end{gathered}
\end{equation}
For $A_{1400}$ replace the head condition by
$\head(w)\in J_p\cup Q_p$. These are exactly the conditions that
adjoining $w$ leaves the parameters unchanged.

\begin{lemma}[Ten witness words]\label{lem:ten}
If $p$ is a sum of words of length at most three and absorbs a word
$w$ of length at most three in one of the four sources, at most ten
summands of $p$ absorb $w$, and they contain all variables of $w$.
\end{lemma}
\begin{proof}
Choose one summand containing each variable of $w$, at most three
choices. Choose a head witness word when the head condition requires one,
at most one further choice. For each member of
$\supp(w)\cap Q_p$, choose a repetition witness, at most three.
Finally choose a summand witnessing each required pair in
\eqref{eq:graph-absorption}, at most three. Let $q$ be their distinct
sum. Its repeated set agrees with $Q_p$ on $\supp(w)$: the selected
witnesses give one inclusion, and $q$ being a subpolynomial gives
the other. Every required edge and head is witnessed. Hence the
same criterion shows that $q$ absorbs $w$. The variables of $w$
already occur in $q$, so at most $10\cdot3=30$ variables are used.
\end{proof}

\begin{corollary}\label{cor:graph-basis}
For $A=A_{1504},A_{1506},A_{1960},A_{1400}$, a finite basis is
\[
\SR\cup\{x+x\eq x,\ \eqref{eq:shortening}\}
\cup\Bnd_A(30,3,10).
\]
The reversed bases apply to $A_{1857},A_{1859},A_{1955},A_{1370}$,
respectively.
\end{corollary}
\begin{proof}
The word rules in Proposition~\ref{prop:undirected} verify
\eqref{eq:shortening}; each triple on the right is a subword of the word on the left;
the repeated-variable and co-occurrence conditions coincide. Every
head is retained. These are precisely the parameters of
Proposition~\ref{prop:undirected}. Apply Theorem~\ref{thm:bounded} and
Lemma~\ref{lem:ten}, and then reverse products.
\end{proof}

\subsubsection{Directed occurrence graphs}
The source $A_{1513}$ has additive order $1<2<3<0$.
Its multiplication has zero 0, identity 1, $2^2=2$,
$2\cdot3=3\cdot1=1\cdot3=3$, and $3\cdot2=3^2=0$.
A word has value 3 precisely when there is exactly one occurrence
assigned 3, every preceding letter lies in $\{1,2\}$, and every
subsequent letter is 1. Its other values are
1 when all letters are 1, 2 when all are in $\{1,2\}$ with some 2,
and 0 otherwise.

Define $X_p,Q_p$ as before and put
\[
\begin{aligned}
R_p&=\bigl\{(x,y):x\ne y,\ x\notin Q_p,\\
 &\qquad\text{an occurrence of }x\text{ precedes one of }y\\
 &\qquad\text{in a summand}\bigr\}.
\end{aligned}
\]
Only the source of a pair is required to be outside $Q_p$.
An evaluation of $p$ is zero exactly when a variable of $X_p$ is
zero, a variable of $Q_p$ is 3, or a pair $(x,y)\in R_p$ has
$x=3$ and $y\in\{2,3\}$. If none occurs, the result is 3 when
some variable is 3, is 2 when some is 2, and is 1 otherwise.
Thus $(X_p,Q_p,R_p)$ determines the function. Assignments
$h=0$, $h=3$, and $(h,k)=(3,2)$, with all other variables 1,
recover the three parameters in turn, using $h\notin Q_p$ for
the last test.

Consequently, $p$ absorbs $w$ exactly when its support is present,
its repeated variables belong to $Q_p$, and all its directed pairs
with source outside $Q_p$ belong to $R_p$. After shortening by
\eqref{eq:shortening}, choose at most three support witnesses,
three repetition witnesses, and six directed-pair witnesses.
Their sum still satisfies this criterion and uses at most 36
variables. We obtain the following result by
Theorem~\ref{thm:bounded}.

\begin{corollary}\label{cor:directed}
The set
$\SR\cup\{x+x\eq x,\ \eqref{eq:shortening}\}
\cup\Bnd_{A_{1513}}(36,3,12)$ is a finite basis for $A_{1513}$.
Its reverse is a basis for $A_{1932}$.
\end{corollary}

These finite sets can be large: $|W(30,3)|=27930$ and
$|W(36,3)|=47988$. Their sizes do not affect the completeness
argument. The bounded-set definitions are part of the bases;
the seven seed identities alone are not asserted to be complete.
\subsection{Compression of long words}\label{sec:compression}
Shortening to a fixed word length is not always available. In this
section, the number of witness words is controlled first; variable
compression then produces a finite set of generating absorptions.

\subsubsection{Multiplicity types with a fixed head}
Consider the word identities
\begin{equation}\label{eq:tail-laws}
xyz\eq xzy,\qquad xyx\eq x^2y,\qquad x^3\eq x^2.
\end{equation}
They put a word in a form specified by its support $D$, head $h$,
and singleton set $F$. The head is placed first, the remaining
letters are ordered, and the multiplicity of each variable is one
on $F$ and two on $D\setminus F$. Indeed, the first law permutes
the tail, the second collects occurrences of the head, and the
third caps multiplicities at two. Write $r_G$ for the resulting
word with fixed $D,h$ and singleton set $G\subseteq F$.

\begin{lemma}[Multiplicity compression]\label{lem:multiplicity}
Suppose a finite ai-semiring $A$ satisfies \eqref{eq:tail-laws} and
has a right identity element. Every valid absorption involving at
most two witness words follows from
\[
\SR\cup\{x+x\eq x,\ \eqref{eq:tail-laws}\}
\cup\Bnd_A(30,60,2).
\]
\end{lemma}
\begin{proof}
Normalize the at most three words involved. Keep every variable
which is the head of one of them as a separate symbol; there are
at most three. Each other variable has a multiplicity vector in
$\{0,1,2\}^3\setminus\{(0,0,0)\}$, with unused coordinates
omitted if fewer words occur. Group variables with the same vector.
There are at most 26 groups, and therefore at most 29 symbols in
all. Replacing each group by one symbol, with its common
multiplicity, gives words of length at most 58. The stated bounds
30 and 60 are thus valid.

To see that the compressed absorption holds in $A$, assign one
representative of each group its compressed value and all other
members the right identity. None of these variables is a head.
Erasing the inserted right identities and sorting tails gives the
compressed word values. Validity of the original absorption therefore
implies validity of the compressed one.

Conversely, substitute the nonempty product of all members of a
group for its compressed symbol. Tail permutation restores each
member once or twice as prescribed, and \eqref{eq:tail-laws}
normalize back to the original words. The compressed identity is
an instance of the bounded finite set. The right identity was used
only in assignments and is not a constant in these substitutions.
\end{proof}

We now apply the lemma to five sources. Their tables are in
Appendix~\ref{app:tables}. All have a right identity labelled 1
and satisfy \eqref{eq:tail-laws}. For a supported witness word $v\in p|_D$
and $x\in F$, abbreviate the conditions
\[
\begin{aligned}
O_x(v)&: x\notin\supp(v)\text{ or }x\in\sing(v),\\
I_x(v)&: x\notin\supp(v)\text{ or }(x\in\sing(v),\ \head(v)=x),\\
N_x(v)&: x\notin\supp(v)\text{ or }(x\in\sing(v),\ \head(v)\ne x).
\end{aligned}
\]
Each existence assertion in the next proposition is independent;
different conditions may use different witness words.

\begin{proposition}\label{prop:tail-criteria}
For a standard word $w$ with parameters $(D,h,F)$, the exact
conditions for $p$ to absorb $w$ are as follows.
\begin{enumerate}[label=(\alph*)]
\item In $A_{1389}$ and $A_{1399}$, $p|_D$ is nonempty;
every $x\in F$ has a supported witness word satisfying $O_x$;
if $h\in F$, a supported witness word satisfies $I_h$.
\item In $A_{1395}$, $p|_D$ is nonempty; a singleton $h$
has a supported witness word satisfying $O_h$, and each
$x\in F\setminus\{h\}$ has one satisfying $N_x$.
\item In $A_{1507}$, $h$ is the head of some summand of $p$;
$p|_D$ is nonempty; and every $x\in F$ has a supported
witness word satisfying $O_x$.
\item In $A_{1959}$, a summand in $p|_D$ has head $h$;
every $x\in F$ has a supported witness word satisfying $O_x$.
\end{enumerate}
\end{proposition}
\begin{proof}
In $A_{1389}$, a word gives 1 when all letters are 1, gives 2
for a single occurrence assigned 2 and all others 1, and gives 3
for a single initial occurrence assigned 3 and all others 1;
otherwise it gives 0. Its additive order is $0<3<2<1$.
In $A_{1399}$ a sole noninitial 3 instead gives 2, and the additive
order is $0<2<3<1$. In $A_{1395}$ this same multiplication is
combined with order $0<3<2<1$.

For these three sources, assign 0 outside $D$. Assigning 1
throughout $D$ forces a supported witness word. Changing just $x$ to 2
detects the single-occurrence condition. Changing it to 3
distinguishes its initial and noninitial positions. In (a), an
initial target 3 can be absorbed only by omission of that variable
or by an initial single occurrence. In (b), a noninitial target 3
has value 2 and requires a noninitial single occurrence or omission,
whereas an initial one has value 3 and permits either position.
These tests prove necessity. The same cases exhaust the nonzero
word values, and the chosen witness word absorbs the target in each,
which proves sufficiency.

In $A_{1507}$, the additive order is $2<0<3<1$. A head assigned
2 makes the word 2. Otherwise an all-1 word is 1, exactly one
occurrence of 3 with the rest 1 gives 3, and every other word is
0. Assigning $h=1$ and all other variables 2 forces a global head
witness word. Assigning 1 on $D$ and 2 outside forces a supported witness word;
changing $x\in F$ to 3 forces $O_x$. Conversely the global head
witness word covers a target value 0, a supported witness word covers 1, and an
$O_x$ witness word covers 3. The value 2 is least.

Finally $A_{1959}$ has additive order $0<3<2<1$. A word on
$\{1,2\}$ takes its head value; a sole occurrence of 3 with
the others in $\{1,2\}$ gives 3, and all other cases give 0.
Assigning $h=1$, the other variables in $D$ the value 2 and
those outside $D$ the value 0 forces a supported head witness word.
Assigning $x=3$ and the other variables in $D$ the value 2
forces $O_x$. Those witness words cover, respectively, the values 1 or 2
and the value 3, so the conditions are sufficient.
\end{proof}

\begin{theorem}\label{thm:tail-bases}
For each of $A_{1389},A_{1399},A_{1395},A_{1507},A_{1959}$, a
finite basis is the set in Lemma~\ref{lem:multiplicity}.
The reversed bases apply to $A_{473},A_{1369},A_{1365},A_{1860},
A_{1954}$, respectively.
\end{theorem}
\begin{proof}
Consider a valid absorption $p+w\eq p$, where the standard word $w$
has support $D$, head $h$ and singleton set $F$. For $G\subseteq F$,
let $r_G$ be the standard word with support $D$, head $h$ and
singleton set $G$.

In cases (a) and (b) of Proposition~\ref{prop:tail-criteria}, choose
any supported witness $v_0\in p|_D$. It absorbs $r_{\varnothing}$.
In case (a), for $x\in F$ choose $v_x\in p|_D$ satisfying $I_h$
when $x=h$, and satisfying $O_x$ when $x\ne h$. The first choice
is available by the additional head condition in the proposition,
and $I_h$ implies $O_h$. In case (b), choose $v_h$ satisfying $O_h$
when $h\in F$, and choose $v_x$ satisfying $N_x$ for each
$x\in F\setminus\{h\}$. The corresponding absorption criterion,
applied to the polynomial consisting of the single word $v_x$,
now gives $v_x+r_{\{x\}}\eq v_x$ in every case.

In case (c), fix a summand with head $h$. Pair it with any
supported witness to absorb $r_{\varnothing}$, and with a supported
$O_x$ witness to absorb $r_{\{x\}}$. In case (d), fix instead a
supported summand with head $h$. This word absorbs
$r_{\varnothing}$; pairing it with a supported $O_x$ witness
absorbs $r_{\{x\}}$. Thus each of these initial insertions uses
at most two witness words. Coincident witnesses require only one
word because addition is idempotent.

For $G\subseteq F$ and $x\in F\setminus G$, the polynomial
$r_G+r_{\{x\}}$ absorbs $r_{G\cup\{x\}}$. Both available words
have support $D$ and head $h$. A singleton in $G$ is witnessed by
$r_G$, and the singleton $x$ is witnessed by $r_{\{x\}}$.
In case (a), a singleton head is initial in the relevant witness,
so its $I_h$ condition holds. In case (b), every singleton other
than $h$ is noninitial in that witness, so its $N_x$ condition
holds. The supported-head and global-head conditions in cases
(c) and (d) hold as well. Proposition~\ref{prop:tail-criteria}
therefore justifies every merging step.

Lemma~\ref{lem:multiplicity} derives all these unary and binary
absorptions from the stated finite set. Add the remaining summands
of $p$ and retain each inserted word. Successively adjoining the
members of $F$ yields $r_F=w$; when $F$ is empty, the initial
insertion already supplies $w$. Lemma~\ref{lem:temporary} removes
the intermediate words. Every valid absorption is consequently
derivable. Lemma~\ref{lem:mutual} proves completeness, and
Lemma~\ref{lem:opposite} gives the five reversed bases.
\end{proof}

\subsubsection{Retaining both endpoints}
Consider instead
\begin{equation}\label{eq:interior-laws}
xyz\eq xy^2z,\qquad xyzu\eq xzyu.
\end{equation}
These identities put each nonlinear word in the form consisting
of its head, an ordered set of interior variables, and its tail.
The endpoints are retained even if they agree or also occur in
the interior. The second law permutes adjacent interior factors,
and the first inserts or deletes repeated interior letters.

\begin{lemma}[Endpoint compression]\label{lem:endpoint-compression}
Suppose $A$ is a finite ai-semiring satisfying
\eqref{eq:interior-laws}, and has an element $e$ with $aeb=ab$
for all $a,b\in A$. Every valid absorption involving at most two
witness words follows from
\[
\SR\cup\{x+x\eq x,\ \eqref{eq:interior-laws}\}
\cup\Bnd_A(16,20,2).
\]
\end{lemma}
\begin{proof}
Normalize the at most three words, keeping the at most six endpoint
variables separate. Every other variable has an interior presence
vector in $\{0,1\}^3\setminus\{(0,0,0)\}$. Group equal vectors,
giving at most seven more symbols. A compressed nonlinear word
has two endpoints and at most thirteen interior symbols, hence
length at most fifteen. Linear words remain linear.

To prove validity of the compressed identity, give one representative
of each group its compressed value and give the others $e$.
Each inserted occurrence is strictly interior and is erased by
$aeb=ab$. Conversely substitute the nonempty product of each
group's members for its symbol. Interior permutation and deletion
recover the original words. This proves the assertion with the
stated, slightly relaxed, bounds.
\end{proof}

For a word with support $D$, head $h$ and tail $t$, call a witness word
\emph{terminally suitable} if it omits $t$ or has a unique occurrence
of $t$ at its tail. Suitability is required below only if $t$ occurs
once in the target.

\begin{theorem}\label{thm:endpoints}
The finite set in Lemma~\ref{lem:endpoint-compression} is a basis
for each of $A_{970},A_{1013},A_{1017}$. Their opposites are
$A_{1851},A_{1927},A_{1919}$, respectively.
\end{theorem}
\begin{proof}
The sources satisfy \eqref{eq:interior-laws} and $a1b=ab$.
The exact absorption criteria are:
\begin{enumerate}[label=(\roman*)]
\item for $A_{970}$, a global witness word has head $h$, $p|_D$ is
nonempty, and a terminally suitable supported witness word exists when
required;
\item for $A_{1013}$, a witness word in $p|_D$ has tail $t$, and a
terminally suitable supported witness word exists when required;
\item for $A_{1017}$, a witness word in $p|_D$ has head $h$, and a
terminally suitable supported witness word exists when required.
\end{enumerate}
In $A_{970}$ the additive order is
$2<0<3<1$. A head assigned 2 gives 2; otherwise a word is 1 if
all letters are 1, is 3 for a sole exceptional occurrence which
is a final 3, and is 0 otherwise. The tests $h=1$, others 2;
1 on $D$, 2 outside; and $t=3$, others in $D$ equal to 1,
force the three conditions in (i). They also suffice: the global
head witness word covers 0, the supported witness word covers 1, and the
suitable witness word covers 3.

The other two sources have additive order $0<3<2<1$. On
$\{1,2\}$, a word takes its tail value in $A_{1013}$ and its
head value in $A_{1017}$. A sole final 3 with all previous
letters in $\{1,2\}$ gives 3; zero or a nonfinal 3 gives 0.
Assigning 0 outside $D$, the distinguished endpoint 1 and the
other variables 2 forces the endpoint witness word. Assigning $t=3$
and the other variables in $D$ the value 2 forces suitability.
The same two witness words cover all nonzero target values.

At most two witness words suffice for each criterion: choose the global head
and the suitable supported witness word in (i), or the supported endpoint
and the suitable witness word in (ii)--(iii). When suitability is not
required, an arbitrary supported witness word suffices in (i).
Lemma~\ref{lem:endpoint-compression} derives their absorption from
the finite set. Add the other summands and use
Lemma~\ref{lem:mutual}.
\end{proof}
\subsection{Prefix saturation with a terminal exception}\label{sec:prefixes}
The two sources in this section have the same multiplication:
1 is a left identity, 0 and 2 are left zeros, and every product
beginning with 3 and having length at least two is 0. A word
therefore takes the value of its first occurrence outside $\{1\}$,
except that a first exceptional 3 gives 3 only in the final
position and gives 0 otherwise. All-1 words give 1. The additive
orders are $0<3<1<2$ in $A_{971}$ and $2<1<3<0$ in $A_{973}$.

\subsubsection{Word reduction and prefix data}
Both sources satisfy
\begin{equation}\label{eq:prefix-word}
x^2y\eq xy,\qquad xyxz\eq xyz,\qquad xyx\eq xy^2.
\end{equation}
Let $\sigma(v)$ be the sequence of first occurrences in a word
$v$. For a nonempty repetition-free sequence $s=x_1\cdots x_k$,
write $R(s)=s x_k$.

\begin{lemma}\label{lem:prefix-word}
The identities \eqref{eq:prefix-word} reduce $v$ to $\sigma(v)$
if its final variable occurs once, and to $R(\sigma(v))$ otherwise.
In particular, $v^2$ reduces to $R(\sigma(v))$.
\end{lemma}
\begin{proof}
The first two laws delete repeated letters at nonfinal positions,
leaving at most the final occurrence as a repetition. If that
occurrence repeats an earlier $x$ with nonempty intervening word
$u$, the last law replaces $xux$ by $xu^2$. Deleting nonfinal
repetitions now leaves the last variable of $\sigma(v)$ doubled.
An adjacent final $xx$ is already in this form. This includes
the one-variable words $x$ and $x^2$.
\end{proof}

For a polynomial $p$ and a set $Y$ of variables, define:
\begin{align*}
A_p(Y)&:\quad p|_Y\ne\varnothing,\\
H_p(Y)&=\bigl\{\text{first variables outside }Y\\
 &\qquad\text{in summands not supported in }Y\bigr\},\\
G_p(Y)&=\bigl\{\text{such first variables occurring}\\
 &\qquad\text{at nonfinal positions}\bigr\},\\
T_p(Y)&=\bigl\{t\notin Y:\text{ a summand ends uniquely in }t\\
 &\qquad\text{and all preceding variables lie in }Y\bigr\}.
\end{align*}
A linear summand $t$ contributes to $T_p(Y)$ when $t\notin Y$.
This is a separate allowed case, not an empty prefix term.

For a standard word with first-occurrence sequence
$s=x_1\cdots x_k$, put $D=\supp(s)$,
$P_j=\{x_1,\ldots,x_j\}$ and $P_0=\varnothing$.

\begin{proposition}\label{prop:prefix-criteria}
In $A_{971}$, $p$ absorbs a standard word $w$ exactly when
\begin{equation}\label{eq:prefix-low}
\begin{gathered}
A_p(D),\qquad x_j\in H_p(P_{j-1})\quad(1\leq j\leq k),\\
\text{if }w=s\text{ and }t=x_k,\quad
A_p(D\setminus\{t\})\text{ or }t\in T_p(D\setminus\{t\}).
\end{gathered}
\end{equation}
In $A_{973}$, the exact criterion is
\begin{equation}\label{eq:prefix-upper}
A_p(D),\qquad x_j\in G_p(P_{j-1})\quad(1\leq j\leq k),
\end{equation}
except that $H_p(P_{k-1})$ suffices at the last position when
$w=s$ has a unique terminal variable.
\end{proposition}
\begin{proof}
For $A_{971}$, assigning 1 on $D$ and 0 outside forces $A_p(D)$.
Assigning 1 on $P_{j-1}$, 2 to $x_j$, and 0 elsewhere forces
the $H_p$ condition, since the target is 2. When the terminal
variable is unique, assign it 3, assign 1 on $D\setminus\{t\}$,
and 0 elsewhere. A witness word can absorb the target 3 only by
avoiding $t$ and being supported there, or by ending uniquely
in $t$ with the rest in that support. Conversely, the supported
witness word covers target 1, the relevant first-exit witness word covers 2,
and the terminal condition covers 3. Zero is least.

For $A_{973}$, assign 1 on $D$ and 2 outside to force
$A_p(D)$. To detect the next first-exit variable, assign 1 on
the preceding prefix, 0 to that variable and 2 elsewhere.
The target is 0 and requires $H_p$. If the target occurrence
is nonfinal, assign 3 instead of 0. The target is again 0;
only a witness word with its first exceptional 3 nonfinal gives 0.
This forces $G_p$. For sufficiency, target 1 is covered by
$A_p(D)$, target 3 by the terminal $H_p$ condition, and target
0 by $H_p$ or $G_p$ according to whether its first exceptional
value is 0 or a nonfinal 3. The value 2 is least.
\end{proof}

\subsubsection{A common prefix construction}
The identity
\begin{equation}\label{eq:left-trim}
xy+z\eq xy+z+xz
\end{equation}
allows a nonempty terminal part of a witness word to be replaced by
another available word. Suppose a seed $r_0$ is available, and
for each desired next prefix variable $x_j$ a witness word has the form
$a x_j b$ with $\supp(a)\subseteq P_{j-1}$ and $b$ nonempty;
the prefix $a$ may be absent. We can successively insert
\begin{equation}\label{eq:prefix-construction}
r_j=x_1\cdots x_j r_0.
\end{equation}
Indeed, apply \eqref{eq:left-trim} to $(a x_j)b$ and
$r_{j-1}$ to insert $a x_j r_{j-1}$. For $j=1$, $a$ is absent.
For $j>1$, write $P=x_1\cdots x_{j-1}$ and apply the same law
to $P r_0$ and the newly inserted word. This inserts
$P a x_j P r_0$. The repeated letters of $a$ and the second
$P$ are nonfinal, so \eqref{eq:prefix-word} reduce it to
$P x_j r_0$. Every substitution is nonempty.

\begin{theorem}\label{thm:prefix-low}
A complete basis for $A_{971}$ consists of $\SR$ and the seven
identities
\begin{gather}
x+x\eq x,\qquad \eqref{eq:prefix-word},\qquad x+x^2\eq x,
\label{eq:prefix-low-laws}\\
xy+z\eq xy+z+xz,\qquad xy+z\eq xy+z+zx.
\label{eq:two-trims}
\end{gather}
This gives twelve identities in total. Reversing them gives a
basis for $A_{1861}$.
\end{theorem}
\begin{proof}
The displayed word-value rule and additive order verify all
the identities. We derive any valid absorption using
\eqref{eq:prefix-low}.
First let $w=R(s)$. The low-square law permits the square of
every original witness word to be retained. Choose a supported witness word
$v\in p|_D$ and use $r_0=v^2$ as seed. A summand witnessing
$x_j\in H_p(P_{j-1})$, after squaring, has a nonempty suffix
after the chosen first-exit occurrence. Therefore the construction
\eqref{eq:prefix-construction} applies. The final word $s r_0$
has all its first occurrences in $s$, and its last variable
repeats, so it reduces to $R(s)$.

Now suppose $w=s$ ends uniquely in $t$. If $t\in
T_p(D\setminus\{t\})$, use that witness word as a seed ending uniquely
in $t$. Otherwise choose a witness word $v$ supported in
$D\setminus\{t\}$ and a witness word $a t b$ witnessing
$t\in H_p(D\setminus\{t\})$. If $b$ is absent, this witness word
is a suitable seed. If $b$ is nonempty, the second trim law
in \eqref{eq:two-trims} inserts $v a t$, also a suitable seed.
Here an absent $a$ is omitted. Construct the prefix
$x_1,\ldots,x_{k-1}$ using squared witnesses as before.
No inserted prefix contains $t$. Deleting only nonfinal
repetitions therefore gives $s$, retaining its unique terminal
variable. For $k=1$, the terminal condition already requires
the summand $t$. Lemmas~\ref{lem:temporary} and
\ref{lem:mutual} prove completeness.
\end{proof}

\begin{theorem}\label{thm:prefix-upper}
A complete basis for $A_{973}$ consists of $\SR$,
$x+x\eq x$, the three identities \eqref{eq:prefix-word},
$x+x^2\eq x^2$, and \eqref{eq:left-trim}.
It has eleven members. Its reverse is a basis for $A_{1863}$.
\end{theorem}
\begin{proof}
All displayed laws hold by the value rule. If every first-exit
condition in \eqref{eq:prefix-upper} uses $G_p$, take any original
supported witness word as seed. The nonempty suffixes needed for
\eqref{eq:prefix-construction} are supplied directly by $G_p$.
The resulting word reduces to $R(s)$.

For a target $s$ ending uniquely in $t$, if
$t\in G_p(D\setminus\{t\})$, this construction gives $R(s)$.
It is equal to $s^2$, and the upper-square law inserts $s$.
Otherwise $t\in H_p(D\setminus\{t\})\setminus
G_p(D\setminus\{t\})$. Its witness is then already a word
ending uniquely in $t$, all earlier variables belonging to
$D\setminus\{t\}$. Use this as seed and construct the preceding
prefix using the available $G_p$ witnesses. Nonfinal deletion
gives $s$ without repeating $t$. This also covers a linear seed;
when $k=1$ no prefix is required. Apply the two absorption lemmas.
\end{proof}
\Needspace{4\baselineskip}
\subsection{Reconstructing single occurrences}\label{sec:singletons}
The final construction uses the multiplication of $A_{1513}$ but
the reverse ordering of its exceptional additive levels.
The source $A=A_{1512}$ has addition given by maximum on
$0<3<2<1$ and multiplication
\[
\begin{array}{c|rrrr}
\cdot&0&1&2&3\\\hline
0&0&0&0&0\\
1&0&1&2&3\\
2&0&2&2&3\\
3&0&3&0&0
\end{array}.
\]
Thus 1 is a multiplicative identity and 0 a multiplicative zero.
The word values are those described for $A_{1513}$ in
Section~\ref{sec:witnesses}: the value is 3 exactly when a unique
occurrence is 3, all preceding letters lie in $\{1,2\}$, and all
subsequent letters are 1. In particular, any 0 or a 2 following 3
makes the product zero.

Let $\Gamma_A$ consist of $\SR$ and the seven identities
\begin{align}
x+x&\eq x,&x+x^2&\eq x,\label{eq:singleton-add}\\
xyx&\eq yx^2,&x^3&\eq x^2,&x^2y^2&\eq y^2x^2,\label{eq:singleton-word}\\
x^2y+xy^2&\eq x^2y+xy^2+xy,\label{eq:merge-two}\\
xy^2z+x^2yz^2&\eq xy^2z+x^2yz^2+xyz.\label{eq:merge-three}
\end{align}

\begin{theorem}\label{thm:singletons}
An exact finite basis for $A_{1512}$ is
\[
\Gamma_A\cup\Bnd_A(12,24,1).
\]
Its reverse is a basis for $A_{1931}$. In particular both algebras
are finitely based.
\end{theorem}

We prove the theorem through three lemmas. The bound specifies
a finite collection of valid unary word absorptions. It does not
claim completeness of the twelve seed identities by themselves.

\subsubsection{Word normal forms and absorption}
The first word law in \eqref{eq:singleton-word} moves an earlier
occurrence of a repeated variable next to its last occurrence.
Iterating, and capping powers at two, gives the sequence of last
occurrences, each variable written once or twice according as it
originally occurred once or at least twice. Distinct squared blocks
commute. Consequently the square of a word of support $D$ reduces
to
\[
b_D=\prod_{x\in D}x^2,
\]
in any chosen order. All displayed laws hold by the word-value rule: the word laws
preserve support and the suffix of each single occurrence, while
the two additive laws provide, for each surviving single occurrence,
a summand with the required suffix.

For a standard word $w$ with support $D$ and singleton set $F$,
let $A_w(x)$ be the set of variables strictly after the unique
occurrence of $x$, for $x\in F$.

\begin{lemma}\label{lem:singleton-criterion}
A polynomial $p$ absorbs $w$ in $A$ if and only if
\begin{enumerate}[label=(\roman*)]
\item $p|_D$ is nonempty;
\item for every $x\in F$, some $v\in p|_D$ either omits $x$,
or has $x\in\sing(v)$ and $A_v(x)\subseteq A_w(x)$.
\end{enumerate}
\end{lemma}
\begin{proof}
Assign 1 on $D$ and 0 outside. The target is 1, which forces
(i). For (ii), assign $x=3$, assign 1 on $A_w(x)$, assign 2
to the remaining variables of $D$, and 0 outside. The target is
3. Unsupported witness words vanish. A supported witness word omitting $x$
gives 1 or 2; one containing $x$ is nonzero only if $x$ occurs
once and all variables after it belong to $A_w(x)$. These are
exactly the ways to absorb the target 3.

Conversely, target 0 is automatically absorbed. A supported witness word
covers target 1 or 2. If the target is 3, let $x$ be its sole
variable assigned 3. A witness word from (ii) omitting $x$ gives 1 or
2; otherwise its suffix variables all receive 1 and it gives 3.
This proves sufficiency.
\end{proof}

\subsubsection{Isolating one singleton and compressing}
Fix an absorption $p+w\eq p$ and retain the order of the
last-occurrence blocks in $w$. For $G\subseteq F$, let $r_G$
have this block order, multiplicity one on $G$, and multiplicity
two on $D\setminus G$. Then $r_F=w$ and
$r_{\varnothing}\eq b_D$ by the word laws.

Choose a supported witness word $v$. Contextual use of $x+x^2\eq x$
and distributivity allows any selected occurrence in $v$ to be
squared while retaining $v$ in the polynomial. Squaring all of
them and normalizing inserts $b_{\supp(v)}$. By
Lemma~\ref{lem:singleton-criterion}, the valid unary absorption
\begin{equation}\label{eq:singleton-required}
b_{\supp(v)}+r_{\varnothing}\eq b_{\supp(v)}
\end{equation}
would insert the all-repeated target.

For each $x\in F$, choose a witness word $v_x$ from part (ii) of the
criterion. Square all its variables except $x$ and normalize,
giving $u_x$. It either omits $x$, in which case it is a square
word, or has $x$ as its only singleton. The set of variables
following $x$ is unchanged. Thus
\begin{equation}\label{eq:singleton-unary}
u_x+r_{\{x\}}\eq u_x
\end{equation}
is valid. We next derive these identities from the fixed finite
set, rather than assume an unbounded scheme.

\begin{lemma}\label{lem:unary-compression}
Every valid unary absorption between words having at most one
singleton, the same distinguished variable when present, follows
from \eqref{eq:singleton-word} and $\Bnd_A(12,24,1)$.
\end{lemma}
\begin{proof}
Normalize each word and order its squared blocks separately before
and after the distinguished variable $x$. Each other variable
has one of three positions in each word: absent, before $x$, or
after $x$. If the word omits $x$, place its single square block
in the before position. There are at most $3^2-1=8$ nonempty
position classes across the two words. Replace each class by one
symbol and keep $x$ separate. The compressed words use at most nine variables. Each word
contains at most six class symbols and at most two occurrences of
$x$, so its length is at most fourteen. If neither word
has a singleton, three support classes suffice.

For any assignment to the compressed symbols, assign one member
of each class its given value and every other member 1. Since 1
is a multiplicative identity and the squared blocks commute,
the original pair then has the compressed values. The compressed
absorption is therefore valid. Conversely, replace a class symbol
by the nonempty product of its members. The square of this product
reduces to their square word by \eqref{eq:singleton-word}.
Ordering the blocks on either side of $x$ recovers the original
pair. The compressed absorption belongs to $\Bnd_A(12,24,1)$,
which proves the lemma.
\end{proof}

The lemma applies to \eqref{eq:singleton-required} and
\eqref{eq:singleton-unary}. We can therefore insert
$r_{\varnothing}$ and every $r_{\{x\}}$ using the proposed basis.

\subsubsection{Merging single occurrences}
\begin{lemma}\label{lem:merge-singletons}
If $r_G$ and $r_{\{x\}}$ have been inserted, where
$x\in F\setminus G$, then the seed identities insert
$r_{G\cup\{x\}}$.
\end{lemma}
\begin{proof}
Write $r_{G\cup\{x\}}=P x Q$, with either $P$ or $Q$ possibly
absent. The supports of $P,Q$ and $\{x\}$ are disjoint because
each variable appears in a single last-occurrence block. Modulo
the word laws the two available words are
\[
r_G=P x^2Q,\qquad r_{\{x\}}=P^2xQ^2.
\]
If both $P$ and $Q$ are nonempty, substitute them for $x,z$
and the distinguished variable for $y$ in
\eqref{eq:merge-three}. It inserts $P x Q$. If just one of
$P,Q$ is absent, \eqref{eq:merge-two} gives the insertion,
interchanging its summands when needed. If both are absent,
the target is the already available $x$. No empty substitution
is used.
\end{proof}

\begin{proof}[Proof of Theorem~\ref{thm:singletons}]
The preceding steps first insert the all-repeated word and all
one-singleton words. Apply Lemma~\ref{lem:merge-singletons}
successively along $F$ to obtain $r_F=w$. Every step retains
its old summands, and Lemma~\ref{lem:temporary} eliminates the
intermediates. Hence every valid word absorption follows from
the fixed finite basis. Lemma~\ref{lem:mutual} proves completeness
for an arbitrary finite variable set. Reversal proves the second
case.
\end{proof}

\section{Nonfinite basis obstructions}\label{sec:obstructions}
\subsection{The three-element obstruction and term images}
Let $S_7=\{z,e,a\}$ have flat addition, meaning $u+u=u$ and
$u+v=z$ for distinct elements, and multiplication with zero $z$,
identity $e$, and $a^2=z$. This algebra is nonfinitely based
\cite{JRZ}. Its name refers to the customary three-element
classification and does not indicate its cardinality.

\begin{proposition}\label{prop:nonai}
For every semiring $A$ of order at most four, the term $p(x)=12x$
induces an idempotent endomorphism onto the additive idempotents
$E^+(A)$. If $E^+(A)\cong S_7$, then $A$ is nonfinitely based.
\end{proposition}
\begin{proof}
The additive multiples of any element have index at most four and
period dividing 12. Hence $24x=12x$ and $144x=12x$.
Commutativity gives $12(x+y)=12x+12y$, and distributivity gives
$(12x)(12y)=144(xy)=12(xy)$. These equalities also give the
idempotence of $p$ and show that its values are additively idempotent.
Every additive idempotent is fixed. The last assertion follows from
Theorem~\ref{thm:retract}.
\end{proof}

Exactly twelve nonidempotent four-element representatives satisfy the
hypothesis in Proposition~\ref{prop:nonai}. Their local indices are
\[
\begin{gathered}
436,\ 438,\ 442,\ 446,\ 1346,\ 1351,\ 1435,\ 1439,\ 1759,\\
1760,\ 1935,\ 1936
\end{gathered}
\]
The image is obtained directly by taking twelve additive copies of each element. For instance, on $A_{436}$ the term
$12x$ has values $(0,1,2,0)$, and the image on $\{0,1,2\}$ is $S_7$.

Another imported obstruction is the following criterion
\cite[Proposition 2.2]{PartI}, based on \cite{JRZ}: a finite
ai-semiring $A$ is nonfinitely based if its noncyclic elements form an
order ideal and $S_7\in\Var(A)$. Here an element is cyclic if
$a^m=a$ for some $m>1$. We use the
stronger condition $S_7\in\mathsf{HS}(A)$ by a subalgebra and a
surjective homomorphism. The required order ideals and surjections are given below.

\subsection{A symmetrized hypergraph criterion}
Flat-algebra constructions provide a general connection between
universal Horn conditions and identities \cite{JacksonFlat}. Here the
specific ingredient is the following semiring locality construction.
We use one established locality construction in its following form
\cite[proof of Theorem 4.9]{JRZ}, as explicitly recalled in
\cite[proof of Theorem 2.2]{GaoS7}.
Let $\Sc$ be the flat semiring whose seven nonzero elements are the
nonempty squarefree divisors of $abc$, with disjoint union as
multiplication and all repeated-letter products zero.
For every $n\geq2$ there is a non-2-colourable 3-uniform hypergraph
$H_n$, and a commutative flat semiring $S_{H_n}$, such that every
subalgebra generated by at most $n$ elements belongs to $\Var(\Sc)$.
Under the vertex assignment, all edge products have the same nonzero
value $c$ with $c^2=z$. We use this locality statement from the cited constructions for
arbitrary $n$.

For completeness, $\Sc\in\HSP(S_7)$ has a small explicit witness.
In $S_7^3$, remove $(e,e,e)$ and collapse all tuples having a $z$
coordinate to zero. The remaining seven tuples, whose coordinates
are $e$ or $a$, correspond to their nonempty sets of $a$ coordinates.
The 26-element domain is a subalgebra and the collapsed set is an
ideal for both operations. Its quotient is exactly $\Sc$.

Define the symmetric term
\[
\tau(x,y,z)=xyz+xzy+yxz+yzx+zxy+zyx.
\]

\begin{theorem}\label{thm:hypergraph}
Suppose an ai-semiring $A$ contains the displayed $S_7$ and satisfies,
for all $r,s,t\in A$,
\begin{equation}\label{eq:tau}
\begin{gathered}
\tau(r,s,t)\in\{z,e,a\},\\
\tau(r,s,t)=a\ \Longrightarrow\ (r,s,t)\text{ is a permutation of }(e,e,a).
\end{gathered}
\end{equation}
Then $A$ is nonfinitely based.
\end{theorem}
\begin{proof}
For a non-2-colourable 3-uniform hypergraph $H=(V,E)$, put
\[
T_H=\sum_{\{u,v,w\}\in E}\tau(x_u,x_v,x_w).
\]
If an evaluation in $A$ gave $T_H=a$, every summand in the flat
copy of $S_7$ would be $a$. By \eqref{eq:tau}, colouring the
vertices assigned $a$ red and the others blue would give exactly
one red vertex in each edge, a proper 2-colouring. Therefore
$T_H\in\{z,e\}$ under every assignment, and
\begin{equation}\label{eq:hyper-identity}
A\models T_H\eq T_H+T_H^2.
\end{equation}
Suppose a finite basis for $A$ uses at most $n\geq2$ variables.
Every assignment of a basis identity into $S_{H_n}$ lies in an
$n$-generated subalgebra belonging to
$\Var(\Sc)\subseteq\Var(A)$. Thus $S_{H_n}$ satisfies the basis.
On the vertex generators, however, commutativity gives $T_{H_n}=c$,
whereas $T_{H_n}+T_{H_n}^2=c+z=z\ne c$. This contradicts
\eqref{eq:hyper-identity}.
\end{proof}

The theorem applies to $A_{452},A_{471},A_{1387}$. In their tables
in Appendix~\ref{app:tables}, $(z,e,a)=(0,1,2)$ and the fourth element
is $b=3$. A triple involving $b$ has a zero product in at least one
of its six orderings: all products involving $b$ vanish in $A_{452}$,
one can put $b$ first in $A_{471}$, and last in $A_{1387}$.
Since $z$ absorbs addition, $\tau=z$. Triples on $\{e,a\}$ have
$\tau=e$ only for $(e,e,e)$ and $\tau=a$ exactly for the three
permutations of $(e,e,a)$. This proves \eqref{eq:tau} directly.

\subsection{The remaining negative variety equalities}
For brevity write $V(X)=\Var(X)$ and use
\[
\begin{aligned}
 V(B_{36})&=V(B_{75})=V(E_4)\vee V(E_7),\\
 V(B_{96})&=V(B_{119})=V(E_7)\vee V(E_8).
\end{aligned}
 \tag{*}\label{negative:eq:small}
\]
Here the tables of the named small algebras are
\begin{center}\small
\begin{tabular}{lll}
\toprule Algebra&Addition&Multiplication\\\midrule
$E_4$&\texttt{0111}&\texttt{0000}\\
$E_7$&\texttt{0001}&\texttt{0001}\\
$E_8$&\texttt{0111}&\texttt{0001}\\
$B_{36}$&\texttt{002012222}&\texttt{000010000}\\
$B_{75}$&\texttt{012111212}&\texttt{010111010}\\
$B_{96}$&\texttt{010111012}&\texttt{000010002}\\
$B_{119}$&\texttt{000012022}&\texttt{000011012}\\
\bottomrule
\end{tabular}
\end{center}
For $B_{36}$ the quotient maps $(0,0,1)$ and $(0,1,0)$ have targets
isomorphic to $E_4,E_7$; for $B_{96}$ the same two maps have targets
isomorphic to $E_7,E_8$. They jointly separate points.
The algebra $B_{75}$ contains copies of $E_4,E_7$ on $\{0,2\}$ and
$\{0,1\}$, and $B_{119}$ contains copies of $E_7,E_8$ on $\{0,1\}$
and $\{1,2\}$, respectively.

To prove the remaining inclusions in (*), the following auxiliary algebras
embed in squares of $B_{36}$ and $B_{96}$ and have the indicated quotients:
\begin{center}\small
\begin{tabular}{llll}
\toprule Algebra&Addition&Multiplication&Embedding coordinates\\\midrule
$A_{592}$&\texttt{0033012332233333}&\texttt{0000011001100000}
 &$(0,0,2,2),(0,1,1,0)$\\
$A_{2226}$&\texttt{0220212122220123}&\texttt{0000012302200303}
 &$(0,1,1,0),(0,2,0,2)$\\
\bottomrule
\end{tabular}
\end{center}
The coordinate maps in the first row are homomorphisms to $B_{36}$;
those in the second are homomorphisms to $B_{96}$. They jointly separate
the four elements in each row. Quotient maps $(0,1,2,0)$ from $A_{592}$
and $(0,1,0,2)$ from $A_{2226}$ have targets isomorphic to $B_{75}$ and
$B_{119}$, respectively. Thus $B_{75}\in HSP(B_{36})$ and
$B_{119}\in HSP(B_{96})$, completing the proof of (*).

\begin{proposition}\label{negative:prop:negative}
The semirings $A_{1442},A_{1770},A_{1886}$ are nonfinitely based.
\end{proposition}
\begin{proof}
The chain case $A_{1442}$ is also covered directly by
Theorem~\ref{thm:chain-classification}. The following variety equalities
transfer its negative status to the non-chain algebra
$A_{1886}$. The following three quotient maps jointly separate points
in each row.
Targets are specified up to simultaneous isomorphism of both operations.
\begin{center}\small
\begin{tabular}{lll}
\toprule Source&Quotient maps&Targets\\\midrule
$A_{1442}$&$(0,0,1,0),(0,1,0,0),(0,1,0,2)$&$E_8,E_7,B_{89}$\\
$A_{1886}$&$(0,0,1,0),(0,1,0,0),(0,1,0,2)$&$E_8,E_7,B_{89}$\\
$A_{1950}$&$(0,1,1,0),(0,1,1,2),(0,1,2,0)$&$E_7,B_{89},B_{119}$\\
$A_{1770}$&$(0,1,1,0),(0,1,1,2),(0,1,2,0)$&$E_7,B_{75},B_{87}$\\
$A_{460}$&$(0,0,0,1),(0,1,0,0),(0,1,2,0)$&$E_4,E_7,B_{87}$\\
\bottomrule
\end{tabular}
\end{center}
Consequently (*) and the separating-quotient argument give
\begin{align*}
 V(A_{1442})=V(A_{1886})=V(A_{1950})
   &=V(E_7)\vee V(E_8)\vee V(B_{89}),\\
 V(A_{1770})=V(A_{460})
   &=V(E_4)\vee V(E_7)\vee V(B_{87}).
\end{align*}
The operation tables are
\begin{center}\small
\begin{tabular}{lll}
\toprule Algebra&Addition&Multiplication\\\midrule
$B_{87}$&\texttt{000010002}&\texttt{000012020}\\
$B_{89}$&\texttt{000012022}&\texttt{000012020}\\
$A_{1442}$&\texttt{0020012322220323}&\texttt{0000010300200300}\\
$A_{1770}$&\texttt{0113111111213113}&\texttt{0120111121120120}\\
$A_{1886}$&\texttt{0000010300200303}&\texttt{0020012322220320}\\
$A_{1950}$&\texttt{0000012302230333}&\texttt{0000011301230330}\\
$A_{460}$&\texttt{0003010300233333}&\texttt{0000012002000000}\\
\bottomrule
\end{tabular}
\end{center}
Table~\ref{tab:literature} identifies $A_{1950}$ with $S_{(4,545)}$
of \cite[Theorem 3.10]{YueTwo}, and $A_{460}$ with $S_{(4,282)}$ of
\cite[Theorem 1.1]{PartIII}. Both are nonfinitely based by those theorems,
using the restricted cancellation justification in
Section~\ref{sec:restricted-cancellation} for the former source.
Equality of the generated varieties proves all three assertions.
\end{proof}

$B_{87}$ is the
nonfinitely based semiring $S_7$, whereas $E_7,E_8,B_{89}$ are finitely
based; the latter three nevertheless generate the nonfinitely based join
displayed above.

\begin{proposition}\label{prop:negative1407}
The semiring $A_{1407}$ generates the same variety as $A_{1411}$,
and is nonfinitely based.
\end{proposition}
\begin{proof}
For $A_{1407}$, the maps $0100$, $0102$, $0120$ are onto
$E_7$, $B_{87}$, $B_{89}$, respectively. For $A_{1411}$ the maps
$0100$, $0102$, $0122$ have these same targets. The strings list the
values on $0,1,2,3$ in order. Direct use of the operation tables
shows that all six maps preserve both operations. The three maps
in either row separate every pair of elements. Lemma~\ref{lem:separating-maps}
therefore identifies both varieties with
$\Var(E_7)\vee\Var(B_{87})\vee\Var(B_{89})$.
The negative status of $A_{1411}$ follows from its literature
correspondence in Table~\ref{tab:literature}.
\end{proof}
\subsection{All negative representatives}
The order-ideal applications outside the cited table correspondences
are listed below. The letters in the image column use the order
$(e,a,z)$ of $S_7$. Each map is defined on the indicated subalgebra
in its increasing local order. The set $N$ is the set of noncyclic
elements and is downward closed in the additive order.
\begin{center}\small
\begin{tabular}{rlll}\toprule
$i$ & $N$ & Subalgebra & Images in $(e,a,z)$\\\midrule
456 & $\{2,3\}$ & $\{0,1,2\}$ & $(z,e,a)$\\
474 & $\{2,3\}$ & $\{0,1,2\}$ & $(z,e,a)$\\
1354 & $\{2,3\}$ & $\{0,1,2\}$ & $(z,e,a)$\\
1359 & $\{2,3\}$ & $\{0,1,2\}$ & $(z,e,a)$\\
1364 & $\{2,3\}$ & $\{0,1,2\}$ & $(z,e,a)$\\
1368 & $\{2,3\}$ & $\{0,1,2\}$ & $(z,e,a)$\\
1390 & $\{2,3\}$ & $\{0,1,2\}$ & $(z,e,a)$\\
1394 & $\{2,3\}$ & $\{0,1,2\}$ & $(z,e,a)$\\
1398 & $\{2,3\}$ & $\{0,1,2\}$ & $(z,e,a)$\\
1409 & $\{2,3\}$ & $\{0,1,2\}$ & $(z,e,a)$\\
1466 & $\{3\}$ & $\{0,1,3\}$ & $(z,e,a)$\\
1469 & $\{3\}$ & $\{0,1,2,3\}$ & $(z,e,z,a)$\\
1511 & $\{3\}$ & $\{0,1,3\}$ & $(z,e,a)$\\
1769 & $\{2,3\}$ & $\{0,1,2\}$ & $(e,z,a)$\\
1930 & $\{3\}$ & $\{0,1,3\}$ & $(z,e,a)$\\
1937 & $\{3\}$ & $\{0,1,3\}$ & $(z,e,a)$\\
1938 & $\{3\}$ & $\{0,1,3\}$ & $(z,e,a)$\\
1953 & $\{3\}$ & $\{0,1,3\}$ & $(z,e,a)$\\
1958 & $\{3\}$ & $\{0,1,3\}$ & $(z,e,a)$\\
\bottomrule\end{tabular}\end{center}
The homomorphism property follows by substituting every pair of
elements in each displayed subalgebra; its three image values are
distinct. The powers in the multiplication table identify $N$, and
the addition table verifies downward closure. Thus every premise of
the cited order-ideal criterion holds for these nineteen algebras.

\begin{theorem}\label{thm:negative-list}
The following 57 representatives are nonfinitely based: the
twelve nonidempotent algebras listed after Proposition~\ref{prop:nonai}
and the following 45 ai-semirings:
\[
\begin{gathered}
450,\ 452,\ 456,\ 460,\ 469,\ 471,\ 474,\ 1354,\ 1359,\\
1364,\ 1368,\ 1385,\ 1387,\ 1390,\ 1394,\ 1398,\ 1406,\ 1407,\\
1409,\ 1411,\ 1436,\ 1440,\ 1442,\ 1466,\ 1469,\ 1475,\ 1476,\\
1503,\ 1505,\ 1511,\ 1769,\ 1770,\ 1856,\ 1858,\ 1885,\ 1886,\\
1887,\ 1890,\ 1892,\ 1930,\ 1937,\ 1938,\ 1950,\ 1953,\ 1958
\end{gathered}
\]
\end{theorem}
\begin{proof}
The nineteen negative entries of Table~\ref{tab:literature}, the
nineteen order-ideal applications above, the four additional variety
equalities, and the three symmetrized-hypergraph examples form disjoint sets and
give $19+19+4+3=45$ ai-semirings. The additive-idempotent-image
argument gives twelve further algebras, all with nonidempotent
addition. Their union is the stated set of 57 representatives.
\end{proof}
\section{Completion of the classification}\label{sec:assembly}
We combine the structural and normal-form results with the finite
enumeration and coverage specified below.

\subsection{Compatible multiplications and isomorphism types}
Let $X=\{0,1,2,3\}$ and let $(X,p)$ be a commutative semigroup.
Write $a+b=p(a,b)$ and give $\operatorname{End}(X,p)$ pointwise
addition. Let $\mathcal L(p)$ consist of the additive homomorphisms
\[
 \lambda:(X,p)\longrightarrow
       \bigl(\operatorname{End}(X,p),+\bigr)
\]
satisfying
\begin{equation}\label{eq:lambda-associativity}
 \lambda\bigl(\lambda(a)(b)\bigr)=\lambda(a)\circ\lambda(b)
 \qquad(a,b\in X).
\end{equation}

\begin{lemma}[Compatible multiplications]\label{lem:compatible-multiplications}
The rule $a\cdot b=\lambda(a)(b)$ gives a bijection between
$\mathcal L(p)$ and the associative multiplications distributive
over $p$. If $G_p=\operatorname{Aut}(X,p)$, its action is
\[
 (g\lambda)(a)=g\circ\lambda(g^{-1}(a))\circ g^{-1}.
\]
The isomorphism types with additive reduct isomorphic to $(X,p)$
are exactly the $G_p$-orbits in $\mathcal L(p)$.
\end{lemma}
\begin{proof}
For endomorphisms $f,g$ of the commutative semigroup,
\[
 (f+g)(a+b)=f(a)+f(b)+g(a)+g(b)
          =(f+g)(a)+(f+g)(b),
\]
so pointwise addition is defined on the endomorphism set. The fact
that each $\lambda(a)$ is an endomorphism is left distributivity.
The homomorphism equation for $\lambda$ is right distributivity.
Evaluating \eqref{eq:lambda-associativity} at $c$ gives
$(a\cdot b)\cdot c=a\cdot(b\cdot c)$, and proves the bijection.
A change of labels preserving $p$ transforms the multiplication by
the displayed conjugation action. Conversely, an isomorphism
between two such semirings preserves $p$ and hence belongs to $G_p$.
\end{proof}

\begin{corollary}[Orbit formula]\label{cor:orbit-formula}
Choose one representative $p$ of each four-element commutative
semigroup. The number of four-element semiring isomorphism types is
\[
 \sum_{[p]}|\mathcal L(p)/G_p|
 =\sum_{[p]}\frac1{|G_p|}
       \sum_{g\in G_p}|\operatorname{Fix}_{\mathcal L(p)}(g)|.
\]
\end{corollary}
\begin{proof}
Lemma~\ref{lem:compatible-multiplications} gives the first sum.
For fixed $p$, count pairs $(g,\lambda)$ with $g\lambda=\lambda$.
An orbit with stabilizer of order $s$ contributes
$(|G_p|/s)s=|G_p|$ such pairs. Division by $|G_p|$ gives its
contribution one and proves the formula.
\end{proof}

\begin{proposition}[Enumeration]\label{prop:enumeration}
There are 2341 four-element semiring isomorphism types: 866 with
idempotent addition and 1475 with nonidempotent addition.
\end{proposition}
\begin{proof}
The ten entries $p(i,j)$, $0\leq i\leq j\leq3$, determine a
commutative operation. Imposing
$p(p(i,j),k)=p(i,p(j,k))$ and identifying simultaneous relabellings
gives the 58 additive representatives in Table~\ref{tab:additive}.
For each such $p$, the finite endomorphism conditions above determine
\[
 m_+=|\mathcal L(p)|,\qquad
 f_+=\sum_{\substack{g\in G_p\\g\ne1}}|\operatorname{Fix}(g)|,
\]
where $1$ denotes the identity permutation. The table records these
numbers and the orbit count $(m_++f_+)/|G_p|$. For example, the
constant addition has $|G_p|=6$, $m_+=442$ and $f_+=98$, so it
admits $90$ semiring isomorphism types. The operation tables are specified in Appendix~\ref{app:catalogue};
the finite conditions determine the rows of Table~\ref{tab:additive}.

Summing the orbit counts gives $2341$. The corresponding labelled
counts are
\[
 \sum_{[p]}\frac{24}{|G_p|}=1140,\qquad
 \sum_{[p]}\frac{24m_+}{|G_p|}=53592.
\]
Equivalently, the full simultaneous relabelling action has nine
orbits of size four, 201 of size twelve and 2131 of size twenty-four:
\[
 9+201+2131=2341,\qquad
 4\cdot9+12\cdot201+24\cdot2131=53592.
\]
The condition $p(x,x)=x$ selects 866 of the representatives; the
remaining 1475 have nonidempotent addition. This action preserves
the order of multiplication, so it identifies opposite semirings
only when an actual simultaneous isomorphism exists.
\end{proof}
\subsection{Assignment of the mathematical cases}
The classification table assigns every algebra a basis construction,
a cited additive-type theorem, or a transfer to another specified
algebra. The structural rows use term definitions, retractions, or
separating homomorphisms. The explicit retraction and linear maps
are recorded in the appendices. The finite-variety rows have the
meaning of Theorem~\ref{thm:finite-membership}, with the full algebra
of term functions as their domain.

\begin{proposition}\label{prop:coverage}
Every row of Table~\ref{tab:all} satisfies the hypotheses of its
indicated mathematical case. In particular, all 386 chain rows are
covered by Theorem~\ref{thm:chain-classification}. The transfer dependencies terminate in
direct basis constructions or cited results. There are 2284 positive
rows and 57 negative rows.
\end{proposition}
\begin{proof}
First consider the additively idempotent representatives with chain
additive reducts. Table~\ref{tab:chain-correspondence} identifies them
with all the tables of \cite{PartIV}, without repetition. By
Theorem~\ref{thm:chain-classification}, these contribute 383 positive
and three negative rows. This part of the classification uses the
cited theorem directly.

For the remaining representatives, the term-definition rows satisfy
the displayed definitions of the second operation. The retraction rows satisfy idempotence, both
endomorphism equations, and the ideal or depth hypotheses of the
indicated theorem. The maps in Table~\ref{tab:linear} are
endomorphisms of the retained commutative reduct and give the stated
interpretation of the other operation.

The hereditary rows satisfy their cited defining identities.
The permutations in Table~\ref{tab:literature} preserve both
operations. These previously classified additive types are disjoint
from the chain class, except for the explicitly repeated negative
source attributions. The separating quotient maps have trivial joint kernel,
so Lemma~\ref{lem:separating-maps} applies. For the forty-two
designated variety equalities, Appendix~\ref{app:finite-witnesses}
gives finite closure constructions whose evaluation maps are
surjective in both directions, by Proposition~\ref{prop:closure-certificate}.

The other positive rows are exactly the tables in the normal-form
sections. Their basis identities are valid under every assignment;
the corresponding normal-form, absorption and recovery proofs apply
without a bound on the alphabet. The additional source/opposite
table accounts for its forty-three pairs. Reversal of each identity
gives the indicated opposite basis.

The negative rows are accounted for by the literature correspondence,
the order-ideal criterion, the twelve additive-idempotent images, and
the symmetrized hypergraph construction, with equality of varieties
where indicated. The finite hypotheses are those stated in
Section~\ref{sec:obstructions}. Counting the rows in these cases
gives 57; none is assigned a positive construction.

Following a transfer in the table reaches a direct case after at most
two further four-element transfers. No circular transfer is used.
Outside the chain class, the rows comprise 1901 positive and
54 negative cases. Together with its 383 positive and three negative
members, the disjoint totals are 2284 and 57, respectively. This proves
the assertions.
\end{proof}

\begin{proof}[Proof of Theorem~\ref{thm:classification}]
Proposition~\ref{prop:enumeration} gives the exhaustive list.
Proposition~\ref{prop:coverage} assigns a valid basis or obstruction
to each row. The constructive proofs and cited theorems therefore
give the asserted binary classification. Restricting the same rows
to idempotent addition yields $821+45=866$; the other rows yield
$1463+12=1475$. These are subdivisions of the same exhaustive list.
\end{proof}

\paragraph{Further questions.} The constructions do not determine minimum
basis size or hereditary finite basability in all cases. An intrinsic
structural characterization of the exhaustive finite coverage would
also complement the classification.

\section*{AI Use Disclosure}
The source manuscript reports the use of generative AI in mathematical
exploration, including the development and refinement of constructions
and proofs. For this condensed version, OpenAI Codex additionally
assisted with exposition, proof reorganization and condensation,
LaTeX preparation, and the development of finite verification scripts.
The finite checks supplement the mathematical arguments and do not
replace proofs for arbitrary terms and variable sets.

\appendix
\section{The additional support and endpoint generators}\label{app:tables}
The rows below specify the source generators of the support, prefix and
endpoint constructions. All operation strings are read row by row on
$\{0,1,2,3\}$. The displayed permutation $\pi$ defines the opposite
partner by
\[
\pi(a)+\pi(b)=\pi(a+b),\qquad
\pi(a)\mathbin{\cdot_{\mathrm{op}}}\pi(b)=\pi(ba).
\]
\begingroup\footnotesize\setlength{\tabcolsep}{2pt}
\begin{longtable}{@{}rrllll@{}}
\caption{The generators and their multiplicative opposites.}\label{tab:generators}\\
\toprule $i$ & $j$ & Addition of $A_i$ & Multiplication of $A_i$ & $\pi$ & Result\\\midrule\endfirsthead
\toprule $i$ & $j$ & Addition of $A_i$ & Multiplication of $A_i$ & $\pi$ & Result\\\midrule\endhead
\midrule\multicolumn{6}{r}{Continued on the next page}\\\endfoot
\bottomrule\endlastfoot
675 & 1153 & \texttt{0020012322020320} & \texttt{0000111100200000} & \texttt{0123} & \ref{thm:parity} \\
692 & 1170 & \texttt{0020012322020323} & \texttt{0000111100200000} & \texttt{0123} & \ref{thm:parity} \\
1641 & 1742 & \texttt{1001011001231031} & \texttt{0110111122220110} & \texttt{0123} & \ref{thm:parity} \\
1274 & 1705 & \texttt{0000012302200300} & \texttt{0100111121220100} & \texttt{0123} & \ref{thm:head-closure} \\
1278 & 1709 & \texttt{0020012322220320} & \texttt{0100111121220100} & \texttt{0123} & \ref{thm:head-closure} \\
1282 & 1713 & \texttt{0000012302230330} & \texttt{0100111121220100} & \texttt{0123} & \ref{thm:head-closure} \\
1295 & 1726 & \texttt{0000012302230333} & \texttt{0100111121220100} & \texttt{0123} & \ref{thm:head-closure} \\
1636 & 1737 & \texttt{0000011301230330} & \texttt{0110111122220110} & \texttt{0123} & \ref{thm:antichain} \\
1648 & 1749 & \texttt{0000011301230333} & \texttt{0110111122220110} & \texttt{0123} & \ref{thm:antichain} \\
2143 & 2187 & \texttt{0000012302200300} & \texttt{0100111121220103} & \texttt{0132} & \ref{thm:capped-local} \\
2146 & 2197 & \texttt{0020012322220320} & \texttt{0100111121220103} & \texttt{0132} & \ref{thm:capped-local} \\
2166 & 2081 & \texttt{0103111101233131} & \texttt{0000010122220103} & \texttt{0132} & \ref{thm:capped-global} \\
2027 & 2102 & \texttt{0003012302233330} & \texttt{0000111101200003} & \texttt{0123} & \ref{thm:periodic-prefix} \\
2288 & 2282 & \texttt{0123111121233132} & \texttt{0000111101220123} & \texttt{0123} & \ref{thm:periodic-prefix} \\
2287 & 2281 & \texttt{0123111121223122} & \texttt{0000111101220123} & \texttt{0123} & \ref{thm:capped-prefix} \\
1819 & 1813 & \texttt{0020012322220320} & \texttt{0120111122220120} & \texttt{0123} & \ref{thm:stable-prefix} \\
1822 & 1816 & \texttt{0020012322220323} & \texttt{0120111122220120} & \texttt{0123} & \ref{thm:stable-prefix} \\
793 & 1183 & \texttt{0000012302230330} & \texttt{0000111101200000} & \texttt{0123} & \ref{thm:exceptional-prefix} \\
799 & 1189 & \texttt{0000012302230333} & \texttt{0000111101200000} & \texttt{0123} & \ref{thm:exceptional-prefix} \\
790 & 1180 & \texttt{0000012302200300} & \texttt{0000111101200000} & \texttt{0123} & \ref{thm:completed-prefix} \\
464 & 465 & \texttt{0000012002200003} & \texttt{0000010302000000} & \texttt{0132} & \ref{thm:factorization} \\
1371 & 1728 & \texttt{0123111121223123} & \texttt{0000111122220200} & \texttt{0123} & \ref{thm:two-letter} \\
1372 & 1729 & \texttt{0000010000230033} & \texttt{0000111122220200} & \texttt{0123} & \ref{thm:two-letter} \\
1373 & 1730 & \texttt{0000012302230333} & \texttt{0000111122220200} & \texttt{0123} & \ref{thm:two-letter} \\
967 & 1848 & \texttt{0000010300200303} & \texttt{0000010322220000} & \texttt{0123} & \ref{thm:global-tails} \\
969 & 1850 & \texttt{0020012322220323} & \texttt{0000010322220000} & \texttt{0123} & \ref{thm:global-tails} \\
1018 & 1920 & \texttt{0000011301230333} & \texttt{0000011302230000} & \texttt{0123} & \ref{thm:global-tails} \\
1504 & 1857 & \texttt{0000010300200303} & \texttt{0000010322220300} & \texttt{0123} & \ref{cor:graph-basis} \\
1506 & 1859 & \texttt{0020012322220323} & \texttt{0000010322220300} & \texttt{0123} & \ref{cor:graph-basis} \\
1960 & 1955 & \texttt{0000011301230333} & \texttt{0000011302230330} & \texttt{0123} & \ref{cor:graph-basis} \\
1400 & 1370 & \texttt{0000012302230333} & \texttt{0000012202000300} & \texttt{0123} & \ref{cor:graph-basis} \\
1513 & 1932 & \texttt{0000012302230333} & \texttt{0000012302230300} & \texttt{0123} & \ref{cor:directed} \\
1389 & 473 & \texttt{0123111121223123} & \texttt{0000012002000300} & \texttt{0123} & \ref{thm:tail-bases} \\
1399 & 1369 & \texttt{0123111121233133} & \texttt{0000012202000300} & \texttt{0123} & \ref{thm:tail-bases} \\
1395 & 1365 & \texttt{0123111121223123} & \texttt{0000012202000300} & \texttt{0123} & \ref{thm:tail-bases} \\
1507 & 1860 & \texttt{0103111101233133} & \texttt{0000010322220300} & \texttt{0123} & \ref{thm:tail-bases} \\
1959 & 1954 & \texttt{0123111121223123} & \texttt{0000011302230330} & \texttt{0123} & \ref{thm:tail-bases} \\
970 & 1851 & \texttt{0103111101233133} & \texttt{0000010322220000} & \texttt{0123} & \ref{thm:endpoints} \\
1013 & 1927 & \texttt{0123111121223123} & \texttt{0000012301230000} & \texttt{0123} & \ref{thm:endpoints} \\
1017 & 1919 & \texttt{0123111121223123} & \texttt{0000011302230000} & \texttt{0123} & \ref{thm:endpoints} \\
971 & 1861 & \texttt{0123112122223123} & \texttt{0000012322220000} & \texttt{0123} & \ref{thm:prefix-low} \\
973 & 1863 & \texttt{0000011301230333} & \texttt{0000012322220000} & \texttt{0123} & \ref{thm:prefix-upper} \\
1512 & 1931 & \texttt{0123111121223123} & \texttt{0000012302230300} & \texttt{0123} & \ref{thm:singletons} \\
\end{longtable}
\endgroup
\section{Small source algebras}\label{app:small}
The following table fixes all three-element labels used in reductions.
The unique negative case is $B_{87}\cong S_7$. The sources and their
basis constructions are described in Sections~\ref{sec:small-results}
and~\ref{sec:family4}.
\begin{center}\small
\begin{tabular}{rll}\toprule
$j$ in $E_j$ & Addition & Multiplication\\\midrule
1 & \texttt{0000} & \texttt{0000}\\
2 & \texttt{0110} & \texttt{0000}\\
3 & \texttt{0001} & \texttt{0000}\\
4 & \texttt{0111} & \texttt{0000}\\
5 & \texttt{0000} & \texttt{0001}\\
6 & \texttt{0110} & \texttt{0001}\\
7 & \texttt{0001} & \texttt{0001}\\
8 & \texttt{0111} & \texttt{0001}\\
9 & \texttt{0001} & \texttt{0101}\\
10 & \texttt{0001} & \texttt{0011}\\
\bottomrule\end{tabular}\end{center}
\begingroup\footnotesize
\setlength{\tabcolsep}{2pt}
\begin{longtable}{@{}rlll@{\hspace{1.4em}}rlll@{\hspace{1.4em}}rlll@{}}
\caption{Three-element source algebras.}\label{tab:small}\\
\toprule $j$ & $+$ & $\cdot$ & Status & $j$ & $+$ & $\cdot$ & Status & $j$ & $+$ & $\cdot$ & Status\\
\midrule\endfirsthead
\toprule $j$ & $+$ & $\cdot$ & Status & $j$ & $+$ & $\cdot$ & Status & $j$ & $+$ & $\cdot$ & Status\\
\midrule\endhead
\bottomrule\endlastfoot
1 & \texttt{000000000} & \texttt{000000000} & FB & 2 & \texttt{000010000} & \texttt{000000000} & FB & 3 & \texttt{000020000} & \texttt{000000000} & FB\\
4 & \texttt{011100100} & \texttt{000000000} & FB & 5 & \texttt{010101010} & \texttt{000000000} & FB & 6 & \texttt{010111010} & \texttt{000000000} & FB\\
7 & \texttt{012111210} & \texttt{000000000} & FB & 8 & \texttt{000012020} & \texttt{000000000} & FB & 9 & \texttt{002012220} & \texttt{000000000} & FB\\
10 & \texttt{012120201} & \texttt{000000000} & FB & 11 & \texttt{000011011} & \texttt{000000000} & FB & 12 & \texttt{011111111} & \texttt{000000000} & FB\\
13 & \texttt{012111211} & \texttt{000000000} & FB & 14 & \texttt{000012021} & \texttt{000000000} & FB & 15 & \texttt{012112221} & \texttt{000000000} & FB\\
16 & \texttt{000010002} & \texttt{000000000} & FB & 17 & \texttt{000011012} & \texttt{000000000} & FB & 18 & \texttt{010111012} & \texttt{000000000} & FB\\
19 & \texttt{011111112} & \texttt{000000000} & FB & 20 & \texttt{012111212} & \texttt{000000000} & FB & 21 & \texttt{000000000} & \texttt{000010000} & FB\\
22 & \texttt{000010000} & \texttt{000010000} & FB & 23 & \texttt{010101010} & \texttt{000010000} & FB & 24 & \texttt{010111010} & \texttt{000010000} & FB\\
25 & \texttt{012111210} & \texttt{000010000} & FB & 26 & \texttt{000012020} & \texttt{000010000} & FB & 27 & \texttt{002002220} & \texttt{000010000} & FB\\
28 & \texttt{002012220} & \texttt{000010000} & FB & 29 & \texttt{000000002} & \texttt{000010000} & FB & 30 & \texttt{000010002} & \texttt{000010000} & FB\\
31 & \texttt{010101012} & \texttt{000010000} & FB & 32 & \texttt{010111012} & \texttt{000010000} & FB & 33 & \texttt{012111212} & \texttt{000010000} & FB\\
34 & \texttt{000012022} & \texttt{000010000} & FB & 35 & \texttt{002002222} & \texttt{000010000} & FB & 36 & \texttt{002012222} & \texttt{000010000} & FB\\
37 & \texttt{000000000} & \texttt{000020000} & FB & 38 & \texttt{000020000} & \texttt{000020000} & FB & 39 & \texttt{000010002} & \texttt{000020000} & FB\\
40 & \texttt{012111212} & \texttt{000020000} & FB & 41 & \texttt{000012022} & \texttt{000020000} & FB & 42 & \texttt{000010000} & \texttt{000111000} & FB\\
43 & \texttt{010111010} & \texttt{000111000} & FB & 44 & \texttt{012111210} & \texttt{000111000} & FB & 45 & \texttt{000012020} & \texttt{000111000} & FB\\
46 & \texttt{002012220} & \texttt{000111000} & FB & 47 & \texttt{000010002} & \texttt{000111000} & FB & 48 & \texttt{010111012} & \texttt{000111000} & FB\\
49 & \texttt{012111212} & \texttt{000111000} & FB & 50 & \texttt{000012022} & \texttt{000111000} & FB & 51 & \texttt{002012222} & \texttt{000111000} & FB\\
52 & \texttt{000000000} & \texttt{000012000} & FB & 53 & \texttt{000010002} & \texttt{000012000} & FB & 54 & \texttt{012111212} & \texttt{000012000} & FB\\
55 & \texttt{000012022} & \texttt{000012000} & FB & 56 & \texttt{000010000} & \texttt{010010010} & FB & 57 & \texttt{010111010} & \texttt{010010010} & FB\\
58 & \texttt{012111210} & \texttt{010010010} & FB & 59 & \texttt{000012020} & \texttt{010010010} & FB & 60 & \texttt{002012220} & \texttt{010010010} & FB\\
61 & \texttt{000010002} & \texttt{010010010} & FB & 62 & \texttt{010111012} & \texttt{010010010} & FB & 63 & \texttt{012111212} & \texttt{010010010} & FB\\
64 & \texttt{000012022} & \texttt{010010010} & FB & 65 & \texttt{002012222} & \texttt{010010010} & FB & 66 & \texttt{000010000} & \texttt{010111010} & FB\\
67 & \texttt{010111010} & \texttt{010111010} & FB & 68 & \texttt{012111210} & \texttt{010111010} & FB & 69 & \texttt{000012020} & \texttt{010111010} & FB\\
70 & \texttt{002012220} & \texttt{010111010} & FB & 71 & \texttt{101010101} & \texttt{010111010} & FB & 72 & \texttt{111111111} & \texttt{010111010} & FB\\
73 & \texttt{000010002} & \texttt{010111010} & FB & 74 & \texttt{010111012} & \texttt{010111010} & FB & 75 & \texttt{012111212} & \texttt{010111010} & FB\\
76 & \texttt{000012022} & \texttt{010111010} & FB & 77 & \texttt{002012222} & \texttt{010111010} & FB & 78 & \texttt{201012120} & \texttt{012111210} & FB\\
79 & \texttt{111111111} & \texttt{012111210} & FB & 80 & \texttt{011111112} & \texttt{012111210} & FB & 81 & \texttt{000000000} & \texttt{000010020} & FB\\
82 & \texttt{000010002} & \texttt{000010020} & FB & 83 & \texttt{012111212} & \texttt{000010020} & FB & 84 & \texttt{000012022} & \texttt{000010020} & FB\\
85 & \texttt{000000000} & \texttt{000012020} & FB & 86 & \texttt{000020000} & \texttt{000012020} & FB & 87 & \texttt{000010002} & \texttt{000012020} & NFB\\
88 & \texttt{012111212} & \texttt{000012020} & FB & 89 & \texttt{000012022} & \texttt{000012020} & FB & 90 & \texttt{000000000} & \texttt{000010002} & FB\\
91 & \texttt{000010000} & \texttt{000010002} & FB & 92 & \texttt{010101010} & \texttt{000010002} & FB & 93 & \texttt{010111010} & \texttt{000010002} & FB\\
94 & \texttt{002012220} & \texttt{000010002} & FB & 95 & \texttt{000010002} & \texttt{000010002} & FB & 96 & \texttt{010111012} & \texttt{000010002} & FB\\
97 & \texttt{000010000} & \texttt{000111002} & FB & 98 & \texttt{010111010} & \texttt{000111002} & FB & 99 & \texttt{002012220} & \texttt{000111002} & FB\\
100 & \texttt{000010002} & \texttt{000111002} & FB & 101 & \texttt{010111012} & \texttt{000111002} & FB & 102 & \texttt{002012222} & \texttt{000111002} & FB\\
103 & \texttt{000000002} & \texttt{002012002} & FB & 104 & \texttt{000010002} & \texttt{002012002} & FB & 105 & \texttt{010101012} & \texttt{002012002} & FB\\
106 & \texttt{010111012} & \texttt{002012002} & FB & 107 & \texttt{002002222} & \texttt{002012002} & FB & 108 & \texttt{002012222} & \texttt{002012002} & FB\\
109 & \texttt{000000000} & \texttt{000011012} & FB & 110 & \texttt{011100100} & \texttt{000011012} & FB & 111 & \texttt{000011011} & \texttt{000011012} & FB\\
112 & \texttt{011111111} & \texttt{000011012} & FB & 113 & \texttt{012111211} & \texttt{000011012} & FB & 114 & \texttt{000012021} & \texttt{000011012} & FB\\
115 & \texttt{012112221} & \texttt{000011012} & FB & 116 & \texttt{000011012} & \texttt{000011012} & FB & 117 & \texttt{011111112} & \texttt{000011012} & FB\\
118 & \texttt{012111212} & \texttt{000011012} & FB & 119 & \texttt{000012022} & \texttt{000011012} & FB & 120 & \texttt{012112222} & \texttt{000011012} & FB\\
121 & \texttt{012111212} & \texttt{010011012} & FB & 122 & \texttt{012111212} & \texttt{000111012} & FB & 123 & \texttt{111111111} & \texttt{012111012} & FB\\
124 & \texttt{010111012} & \texttt{012111012} & FB & 125 & \texttt{000012022} & \texttt{012111012} & FB & 126 & \texttt{000010002} & \texttt{012012012} & FB\\
127 & \texttt{000011012} & \texttt{012012012} & FB & 128 & \texttt{111111111} & \texttt{010111212} & FB & 129 & \texttt{010111012} & \texttt{010111212} & FB\\
130 & \texttt{000012022} & \texttt{010111212} & FB & 131 & \texttt{000010002} & \texttt{000111222} & FB & 132 & \texttt{000011012} & \texttt{000111222} & FB\\
\end{longtable}
\endgroup
\section{The variety generated by all two-element semirings}\label{app:twoelement}
This appendix establishes the structural results used in
Theorem~\ref{thm:twohfb} and describes the full subvariety lattice.
The specialized external inputs are the identity basis and subvariety
classification of Shao and Ren \cite[Corollary 2.3 and Theorem 3.8]{ShaoRen}
for the six two-element ai-semirings. The full variety is also the subject of \cite{WWLY}; the
decomposition, completeness and separation statements required here
are proved below, so its subvariety enumeration is not used as an input.
\begingroup
\providecommand{\Var}{}\renewcommand{\Var}{\operatorname{Var}}
\providecommand{\HSP}{}\renewcommand{\HSP}{\operatorname{HSP}}
\providecommand{\HH}{}\renewcommand{\HH}{\operatorname{H}}
\providecommand{\Sub}{}\renewcommand{\Sub}{\operatorname{S}}
\providecommand{\cl}{}\renewcommand{\cl}{\operatorname{cl}}
\providecommand{\Lat}{}\renewcommand{\Lat}{\mathcal L}
\providecommand{\K}{}\renewcommand{\K}{\mathcal K}
\providecommand{\E}{}\renewcommand{\E}{\mathcal E}
\providecommand{\V}{}\renewcommand{\V}{\mathcal V}
\providecommand{\W}{}\renewcommand{\W}{\mathcal W}
\providecommand{\A}{}\renewcommand{\A}{\mathcal A}
\providecommand{\Sr}{}\renewcommand{\Sr}{\mathcal V_2}
\providecommand{\C}{}\renewcommand{\C}{\mathcal C}
\providecommand{\T}{}\renewcommand{\T}{\mathcal T}
\providecommand{\powerset}{}\renewcommand{\powerset}{\mathcal P}
\providecommand{\eq}{}\renewcommand{\eq}{\approx}
\providecommand{\mat}{}\renewcommand{\mat}[4]{\begin{psmallmatrix}#1&#2\\#3&#4\end{psmallmatrix}}
\[
\begin{gathered}
\K=\{L_2,R_2,M_2,D_2,N_2,T_2,Z_2,W_2,Z_7,Z_8\},\\
\Sr=\Var(\K),\qquad\E=\K\cup\{Q,L_0,R_0\}.
\end{gathered}
\]
The three additional generators will be defined below as members of
$\Sr$. They do not change the ambient variety.
\subsection{Small semirings and identity tests}\label{tw:sec:tests}

All two-element tables used in this section are reproduced in Table~\ref{tw:tab:two}. Rows and columns are ordered $0,1$. These fix the notation used in this section. Throughout, $2u=u+u$ and $3u=u+u+u$.

\begin{table}[htbp]
\centering
\caption{The ten two-element semirings.}\label{tw:tab:two}
\small
\renewcommand{\arraystretch}{1.65}
\begin{tabular}{c cc@{\qquad}c cc}
\toprule
Semiring&Addition&Multiplication&Semiring&Addition&Multiplication\\
\midrule
$L_2$&$\mat0111$&$\mat0011$&$T_2$&$\mat0111$&$\mat1111$\\
$R_2$&$\mat0111$&$\mat0101$&$Z_2$&$\mat0000$&$\mat0000$\\
$M_2$&$\mat0111$&$\mat0111$&$W_2$&$\mat0000$&$\mat0001$\\
$D_2$&$\mat0111$&$\mat0001$&$Z_7$&$\mat0110$&$\mat0000$\\
$N_2$&$\mat0111$&$\mat0000$&$Z_8$&$\mat0110$&$\mat0001$\\
\bottomrule
\end{tabular}
\end{table}

By distribution, every semiring term has a polynomial expression
\[
 u=w_1+\cdots+w_m\qquad(m\ge1),
\]
where each $w_i$ is a nonempty word in variables. Summands are counted with multiplicity. For a word $w$, let $h(w),t(w),c(w)$ and $|w|$ be its first letter, last letter, content set and length. For a polynomial $u$ define
\begin{align*}
 H(u)&=\{h(w_i):1\le i\le m\},\\
 T(u)&=\{t(w_i):1\le i\le m\},\\
 C(u)&=\bigcup_i c(w_i),\\
 A(u)&=\{x:w_i=x\text{ for some }i\},\\
 B(u)&=\begin{cases}1,&\text{some }|w_i|\ge2,\\0,&\text{otherwise},\end{cases}\\
 P(u)&=\{x:\#\{i:w_i=x\}\text{ is odd}\},\\
 M(u)&=\min_{\subseteq}\{c(w_i):1\le i\le m\},\\
 O(u)&=\{U\ne\varnothing:\#\{i:c(w_i)=U\}\text{ is odd}\}.
\end{align*}
Here $M(u)$ is the antichain of inclusion-minimal content sets. Introduce three tagged values
\begin{align*}
 \tau(u)&=\begin{cases}\star,&B(u)=1,\\A(u),&B(u)=0,\end{cases}\\
 \zeta(u)&=\begin{cases}x,&m=1\text{ and }w_1=x,\\\star,&\text{otherwise},\end{cases}\\
 \omega(u)&=\begin{cases}c(w_1),&m=1,\\\star,&m\ge2.\end{cases}
\end{align*}
In each definition $\star$ is a formal marker, distinct from all other possible values.

For $L_2,R_2,M_2,D_2,N_2$, the following criteria restate
\cite[Lemma 1.1(i)--(v)]{ShaoRen} in our notation. The $T_2$
criterion also covers linear identities and follows from the absorption
formulation in \cite[Lemma 1.2(6)]{PartI}. Direct evaluations are
included for all ten semirings.

\begin{proposition}[Complete identity test]\label{tw:prop:test}
For arbitrary polynomial terms $u,v$, Table~\ref{tw:tab:tests} gives a necessary and sufficient condition for $S\models u\eq v$. Consequently, $\Sr\models u\eq v$ if and only if all ten conditions hold.
\end{proposition}

\begin{table}[htbp]
\centering
\caption{Identity tests for the ten two-element semirings.}\label{tw:tab:tests}
\begin{tabular}{c l@{\qquad}c l}
\toprule
$S$&Condition&$S$&Condition\\\midrule
$L_2$&$H(u)=H(v)$&$T_2$&$\tau(u)=\tau(v)$\\
$R_2$&$T(u)=T(v)$&$Z_2$&$\zeta(u)=\zeta(v)$\\
$M_2$&$C(u)=C(v)$&$W_2$&$\omega(u)=\omega(v)$\\
$D_2$&$M(u)=M(v)$&$Z_7$&$P(u)=P(v)$\\
$N_2$&$A(u)=A(v)$&$Z_8$&$O(u)=O(v)$\\
\bottomrule
\end{tabular}
\end{table}

\begin{proof}
In $L_2,R_2,M_2$, evaluation of $u$ is the Boolean join of the variables in $H(u),T(u),C(u)$ respectively. In $N_2$, all words of length at least two evaluate to zero, so the result is the join of the variables in $A(u)$; an empty join here denotes the constant zero function, not a nullary symbol in the language. In $T_2$, the occurrence of any such word makes the value constantly one. Otherwise the value is the join over $A(u)$. These observations prove the criteria for $L_2,R_2,M_2,N_2,T_2$.

For $D_2$, a polynomial gives the monotone Boolean function whose minimal true supports are exactly $M(u)$. These minimal supports uniquely determine the function. In $Z_2$, a lone variable evaluates as itself and every composite term evaluates constantly to zero. In $W_2$, a single word is the conjunction of its content, and any sum of two or more words evaluates constantly to zero. This proves the corresponding tagged criteria.

In $Z_7$, the value is the sum modulo two of the variables in $P(u)$. Finally, in $Z_8$ the word $w_i$ gives the square-free monomial with support $c(w_i)$, and repeated equal supports cancel in pairs. Evaluation at the characteristic vector of $V$ gives $\sum_{U\subseteq V}a_U$, so Lemma~\ref{lem:finite-set-recovery}(i) recovers the square-free coefficients. Thus $O(u)$ uniquely determines the term function. The last assertion follows from the definition of $\Sr$.
\end{proof}

\begin{corollary}\label{tw:cor:free}
The identity problem of $\Sr$ is decidable, and $\Sr$ is locally finite. Its free algebra on $n$ generators can be constructed by closing the $n$ projections under pointwise operations in
\[
 \prod_{S\in\K}S^{\{0,1\}^n}.
\]
\end{corollary}
\begin{proof}
The proposition is a finite test. Alternatively, equality of the displayed term-function tuples is precisely equality in the free algebra. The ambient product is finite, with at most $2^{10\cdot2^n}$ elements.
\end{proof}

\begin{theorem}[Shao--Ren \cite{ShaoRen}]\label{tw:lem:ai-known}
The variety
\[
\A_0=\Var(L_2,R_2,M_2,D_2,N_2,T_2)
\]
is defined, among semirings, by
\begin{equation}\label{tw:eq:ai-basis}
 2x\eq x,\quad (xy)^2\eq xy,\quad xyzt\eq xzyt,\quad
 x+yz\eq x+yz+xz+yx.
\end{equation}
Every subvariety of $\A_0$ is generated by a unique subset of the six displayed semirings. Thus $\Lat(\A_0)$ is the 64-element Boolean lattice.
\end{theorem}

The reduced identity list in \eqref{tw:eq:ai-basis} is the
presentation of \cite[Corollary 2.3]{ShaoRen}, together with its
preceding results. The assertion about subvarieties is
\cite[Theorem 3.8]{ShaoRen}. The hereditary finite basability of
$\A_0$ was also established in \cite[Theorem 3.7]{ShaoRen}.
These results concern the ai-semiring subvariety of $\Sr$.

The two assertions of this theorem have distinct uses below. The finite basis ensures that the doubling image of an arbitrary model of our proposed identities belongs to $\A_0$. The subvariety classification then identifies the small generators of that image.

\subsection{Three homomorphic images in the original variety}\label{tw:sec:three}

Define $Q$ on $\{e,a,\infty\}$ by making $\{e,a\}$ the group of order two under addition, adjoining an additive absorbing element $\infty$, and making every product equal to $\infty$. Define $L_0$ on $\{0,a,b\}$ by constant-zero addition and by adjoining a multiplicative zero to the two-element left-zero semigroup. Define $R_0$ by reversing the multiplication of $L_0$. Their matrices are
\begin{align}
 Q_+&=\begin{pmatrix}e&a&\infty\\a&e&\infty\\\infty&\infty&\infty\end{pmatrix},\qquad
 Q_\cdot=\begin{pmatrix}\infty&\infty&\infty\\\infty&\infty&\infty\\\infty&\infty&\infty\end{pmatrix},\label{tw:eq:Qtable}\\
 (L_0)_+=(R_0)_+&=\begin{pmatrix}0&0&0\\0&0&0\\0&0&0\end{pmatrix},\notag\\
 (L_0)_\cdot&=\begin{pmatrix}0&0&0\\0&a&a\\0&b&b\end{pmatrix},\qquad
 (R_0)_\cdot=\begin{pmatrix}0&0&0\\0&a&b\\0&a&b\end{pmatrix}.\label{tw:eq:LRtable}
\end{align}
The orders of rows and columns are $e,a,\infty$ and $0,a,b$, respectively. Their semiring axioms follow immediately from the definitions: all sums in $L_0,R_0$ are zero and zero is multiplicatively absorbing; in $Q$, both sides of either distributive law are $\infty$.

\begin{lemma}\label{tw:lem:quotients}
The following inclusions hold:
\begin{align*}
 T_2&\in\Sub(Q),&W_2&\in\Sub(L_0)\cap\Sub(R_0),\\
 Q&\in\HH(T_2\times Z_7)\cap\HH(T_2\times Z_8),\\
 L_0&\in\HH(L_2\times W_2),&R_0&\in\HH(R_2\times W_2).
\end{align*}
Furthermore, $Z_2$ is a homomorphic image of each of $N_2\times W_2$, $T_2\times W_2$, and $Z_7\times W_2$. In particular, $Q,L_0,R_0\in\Sr$.
\end{lemma}
\begin{proof}
The subsets $\{e,\infty\}\subseteq Q$ and $\{0,a\}\subseteq L_0,R_0$ give the first line. For either $G=Z_7$ or $G=Z_8$, define
\[
 q:T_2\times G\longrightarrow Q,\qquad
 q(0,0)=e,\quad q(0,1)=a,\quad q(1,i)=\infty.
\]
On the layer with first coordinate zero, addition is the group operation of order two; any first coordinate one forces the absorbing layer. Every product has first coordinate one. Thus $q$ is a surjective homomorphism.

For $L_2\times W_2$, send $(i,0)$ to $0$ and send $(0,1),(1,1)$ to $a,b$. Sums have second coordinate zero. Products with both second coordinates one retain the first coordinate of the left factor; the other products map to zero. This verifies the map to $L_0$. Reversing multiplication gives the map to $R_0$.

For the three remaining products, map one point to $1\in Z_2$ and all other points to $0$. Choose the distinguished point $(1,1)$ for $N_2\times W_2$, $(0,1)$ for $T_2\times W_2$, and $(1,1)$ for $Z_7\times W_2$. No sum or product in the relevant direct product equals the distinguished point: sums have second coordinate zero, and products have first coordinate $0,1,0$, respectively. All three maps are therefore surjective homomorphisms.
\end{proof}

\subsection{Consequences of the defining identities}\label{tw:sec:system}

Use the thirty-one identities $\Sigma$ and the six term abbreviations
from Section~\ref{sec:twoelement}.

Set $\Upsilon=\operatorname{Mod}(\mathrm{SR}\cup\Sigma)$. Direct substitution in Table~\ref{tw:tab:two} verifies every entry of $\Sigma$, so
\begin{equation}\label{tw:eq:easy-basis}
 \Sr\subseteq\Upsilon.
\end{equation}
The next three lemmas are valid in $\Upsilon$, without assuming membership in $\Sr$.

\begin{lemma}\label{tw:lem:double-general}
For every $A\in\Upsilon$, $d(A)\in\A_0$.
\end{lemma}
\begin{proof}
Equation~\eqref{tw:eq:endos} makes $d$ an idempotent endomorphism, and its image is additively idempotent. Identities $(xy)^2\eq xy$ and $xyzt\eq xzyt$ hold on this image by \eqref{tw:eq:squares} and \eqref{tw:eq:normal}. Substitute preimages into \eqref{tw:eq:ai-lift} and use $d(yz)=d(y)d(z)$ to obtain the final identity of \eqref{tw:eq:ai-basis} on $d(A)$. Theorem~\ref{tw:lem:ai-known} now gives the assertion.
\end{proof}

\begin{lemma}[Normal-band reduction]\label{tw:lem:normal-reduction}
In a band satisfying $xyzt\eq xzyt$, two nonempty words with the same head, content and tail have equal values. Consequently, in every $A\in\Upsilon$, two words of length at least two with these three data equal also have equal values.
\end{lemma}
\begin{proof}
In a band, duplicate the first and last letters to retain them as endpoints. The identity $xyzt\eq xzyt$ interchanges adjacent interior letters, with the nonempty prefix and suffix substituted for $x$ and $t$. Permute the interior letters and delete adjacent repetitions using idempotence. A word is thereby reduced to its endpoints and one interior occurrence of each distinct content letter in a fixed order. This proves the first assertion, also for a one-letter word by first duplicating that letter.

For the second assertion, $e(A)$ is a band because $e(e(x))=e(x)$, and it satisfies the same normal identity. The identities $x^2y\eq xy$ and $xy^2\eq xy$ permit replacement of every letter of a word of length at least two by its square. Apply the first assertion to these squared-letter words in $e(A)$.
\end{proof}

\begin{lemma}[Balanced-content transfer]\label{tw:lem:content-transfer}
If $w,v$ have length at least two and $c(w)=c(v)$, then every $A\in\Upsilon$ satisfies
\begin{equation}\label{tw:eq:content-transfer}
 2w+2v+w\eq2w+2v+v.
\end{equation}
\end{lemma}
\begin{proof}
Lemma~\ref{tw:lem:normal-reduction} gives $wvw\eq w$ and $vwv\eq v$. Substitute $x=w,y=v$ into the single identity \eqref{tw:eq:bridge} and make these two reductions. Thus the entire family \eqref{tw:eq:content-transfer} follows from the finite system; no infinite identity scheme is being added to $\Sigma$.
\end{proof}

\subsection{Structural decomposition and small generators}\label{tw:sec:exhaustive}

We now work in $\Upsilon$. The argument proves both the converse of \eqref{tw:eq:easy-basis} and hereditary generation by the thirteen semirings. In particular, no step below presupposes that an arbitrary model of $\Sigma$ already lies in $\Sr$.

\paragraph{Retractions onto ideals}

An ideal of a semiring here means a nonempty subset $I$ such that $a+b,ab,ba\in I$ whenever $a\in I$ and $b\in S$. The Rees congruence $\rho_I$ identifies all elements of $I$ and leaves every other element in a singleton class.

\begin{lemma}[Retraction decomposition]\label{tw:lem:retract}
Suppose $h:S\to S$ is an idempotent endomorphism and $I=h(S)$ is an ideal. Then
\[
 S\hookrightarrow I\times S/\rho_I,\qquad
 a\longmapsto(h(a),[a]_{\rho_I})
\]
is a subdirect embedding. Consequently,
\[
 \Var(S)=\Var\bigl(I,S/\rho_I\bigr).
\]
Both factors belong to $\operatorname{HS}(S)$.
\end{lemma}
\begin{proof}
Lemma~\ref{lem:ideal-decomposition} gives the subdirect embedding.
Its factors are a subalgebra and a quotient of $S$, respectively,
so they generate $\Var(S)$.
\end{proof}

By \eqref{tw:eq:r-sum}, the image $r(A)$ equals $A+A$: it contains all sums, and $3a=a+2a$ is itself a sum. Together with \eqref{tw:eq:r-ideal}, this proves that $r(A)$ is an ideal. Lemma~\ref{tw:lem:retract} splits $A$ into an additively regular factor $r(A)$, satisfying $3x\eq x$, and a factor with constant addition. We describe the two factors separately.

\paragraph{The factor with constant addition}

\begin{lemma}\label{tw:lem:zero-add}
If $A\in\Upsilon$ has constant addition, then $\Var(A)$ is generated by a subset of $\{Z_2,W_2,L_0,R_0\}$, each member of which belongs to $\Var(A)$.
\end{lemma}
\begin{proof}
Write $0$ for the constant sum. Distributivity gives $a0=0a=0$. By \eqref{tw:eq:squares}--\eqref{tw:eq:delete}, $e(A)=A^2$ is a multiplicative band, and it is an ideal of $A$: every product lies in $A^2$ and every sum is $0\in A^2$. Lemma~\ref{tw:lem:retract} applies to $e$. In $A/\rho_{A^2}$ both operations are constantly zero. This factor is trivial or generates $\Var(Z_2)$; in the latter case a two-element subalgebra is $Z_2$.

It remains to classify $B=A^2$. Its multiplication is a band with zero and satisfies \eqref{tw:eq:normal}; its addition is zero. Lemma~\ref{tw:lem:normal-reduction} shows that the value of a nonempty word is determined by its head, content and tail.

If $B$ is nontrivial, $\{0,b\}\cong W_2$ for every $b\ne0$. If $aba\ne ba$, put $u=aba$ and $v=ba$. The normal-band identities give
\[
 u^2=u,\quad v^2=v,\quad uv=u,\quad vu=v.
\]
Neither $u$ nor $v$ is zero: either possibility would contradict the displayed products and $u\ne v$. Thus $\{0,u,v\}\cong L_0$. Dually, failure of $aba\eq ab$ produces $R_0$.

Take $W_2$ when $B$ is nontrivial, take $L_0$ exactly when $xyx\eq yx$ fails in $B$, and take $R_0$ exactly when $xyx\eq xy$ fails. These chosen algebras are subalgebras of $B$. Conversely they generate a variety containing $B$, as follows. Every polynomial with at least two summands is zero. For two single words, the identity tests of the chosen algebras require equal contents and, respectively, equal heads when $L_0$ is chosen and equal tails when $R_0$ is chosen. If $L_0$ is absent, $xyx\eq yx$ makes the head irrelevant. Indeed, if $a$ occurs in a word $w=uav$, then $aw=auav=uav=w$, using $aua=ua$ when $u$ is nonempty and idempotence when it is empty. Thus two words with different heads can be given a common head without changing their values. Dually, absence of $R_0$ permits a common tail. Lemma~\ref{tw:lem:normal-reduction} then compares the words. Empty prefixes or suffixes here describe word positions, not constants in the language. Thus every identity of the chosen generators holds in $B$. A single word cannot equal a polynomial with two or more summands in $W_2$, so no additional mixed case occurs. This proves the assertion for $B$ and then for $A$.
\end{proof}

\paragraph{The additively regular factor}

Call an algebra additively regular here if it satisfies $3x\eq x$. This terminology refers only to the identity just stated. In particular, each element $a$ lies in the additive group $\{a,2a\}$ when $a\ne2a$.

\begin{lemma}\label{tw:lem:top-null}
Suppose $J\in\Upsilon$ satisfies $3x\eq x$, has constant multiplication $xy=c$, and $x+c=c$ for every $x$. Then its generated variety is one of
\[
 \T,\qquad\Var(T_2),\qquad\Var(Q).
\]
In each case it is generated by an algebra belonging to $\operatorname{HS}(J)$, except that the trivial case needs no generator.
\end{lemma}
\begin{proof}
Distributivity gives $c+c=c$. Any polynomial containing a word of length at least two evaluates to $c$. In a purely linear polynomial, the identity $3x\eq x$ reduces each positive coefficient to one or two according to its parity.

An identity between two purely linear polynomials holds in $Q$ precisely when both their supports and their coefficient parities agree. Support is detected by assigning $\infty$ to a variable present on only one side; parity is detected in the two-element additive group. Identities between two polynomials containing a long word hold automatically, and no identity between these two types of polynomial holds in $Q$. It follows that every identity of $Q$ holds in $J$, so $J\in\Var(Q)$.

If $J$ is additively idempotent, the same argument using only supports gives $J\in\Var(T_2)$; every nontrivial such $J$ contains $\{a,c\}\cong T_2$ for some $a\ne c$. Otherwise choose $a\ne2a$. The three elements $a,2a,c$ are distinct, since $c$ is an additive absorbing idempotent. They form a subalgebra isomorphic to $Q$. The three alternatives follow.
\end{proof}

\begin{lemma}\label{tw:lem:p-factor}
Let $P\in\Upsilon$ satisfy $3x\eq x$ and $x\eq x+x^2$. If $P$ is additively idempotent, then $P=d(P)$. Otherwise
\[
 \Var(P)=\Var\bigl(d(P),Z_7\bigr),\qquad Z_7\in\operatorname{HS}(P).
\]
\end{lemma}
\begin{proof}
Every value of a word of length at least two is multiplicatively idempotent by \eqref{tw:eq:squares}--\eqref{tw:eq:delete}; the identity $x\eq x+x^2$ makes such a value additively idempotent as well. For a polynomial $u$, let $P(u)$ be the parity set of its linear summands from Section~\ref{tw:sec:tests}. Coefficient reduction gives, on $P$,
\begin{equation}\label{tw:eq:linear-normal}
 u=2u+\sum_{x\in P(u)}x.
\end{equation}
When $P(u)$ is empty the right side means $2u$, with no empty sum in the signature. If $u\eq v$ holds in $d(P)$ and $Z_7$, then $2u=2v$ on $P$ and $P(u)=P(v)$. Equation~\eqref{tw:eq:linear-normal} proves $u=v$ on $P$. Hence $P\in\Var(d(P),Z_7)$.

If $a\ne2a$, put $b=a^2$ and $e=2a$. The element $b$ is additively idempotent, $a+b=a$, and $e+b=e$. All products in $\{a,e,b\}$ are $b$, by \eqref{tw:eq:squares}--\eqref{tw:eq:delete}. The map taking $a$ to $1$ and $e,b$ to $0$ is a surjective homomorphism onto $Z_7$; it remains valid when $e=b$. Since $d(P)$ is also a subalgebra, the reverse variety inclusion follows.
\end{proof}

\begin{lemma}\label{tw:lem:q-factor}
Let $B\in\Upsilon$ satisfy $3x\eq x$ and $x\eq2x+x^2$. If $B$ is additively idempotent, then $B=d(B)$. Otherwise
\[
 \Var(B)=\Var\bigl(d(B),Z_8\bigr),\qquad Z_8\in\Sub(B).
\]
\end{lemma}
\begin{proof}
Write $u=\sum_i w_i$ and $v=\sum_j v_j$, where the summands are nonempty words. Suppose $u\eq v$ holds in $d(B)$ and $Z_8$, and fix any assignment in $B$. Write $b=2u=2v$. Since squaring is an endomorphism,
\[
 u=b+u^2=b+\sum_i w_i^2,
 \qquad
 v=b+v^2=b+\sum_j v_j^2.
\]
The identity $x\eq2x+x^2$ implies $2x=2x+2x^2$. Thus the additively idempotent element $b$ absorbs every $2w_i^2$ and every $2v_j^2$: absorption of their sum implies absorption of each summand in the semilattice of additive idempotents.

Group all these squared words according to their content. Within a fixed content class, Lemma~\ref{tw:lem:content-transfer} permits replacement of any representative $w$ by any other representative $v$, in the presence of the background $b$, since $b+2w+2v=b$. After these replacements, pairs of equal representatives are absorbed by $b$. The remaining representative occurs exactly when the original number of words in that content class was odd. The identity in $Z_8$ says exactly that these parities agree for $u$ and $v$. The two expressions are therefore equal. This proves $B\in\Var(d(B),Z_8)$.

If $a\ne2a$, then $a^2$ cannot be additively idempotent, since $a=2a+a^2$ would then be additively idempotent. Put $b=a^2$. We have $b^2=b$ and $b\ne2b$. The two-element subalgebra $\{2b,b\}$ is $Z_8$: addition is the group of order two, $b^2=b$, $b(2b)=(2b)b=2b$, and $(2b)^2=2b$. This gives the opposite variety inclusion.
\end{proof}

\begin{proposition}\label{tw:prop:regular-generation}
If $R\in\Upsilon$ satisfies $3x\eq x$, then $\Var(R)$ is generated by a subset of
\[
 \{L_2,R_2,M_2,D_2,N_2,T_2,Z_7,Z_8,Q\}
\]
whose members all belong to $\Var(R)$.
\end{proposition}
\begin{proof}
By \eqref{tw:eq:f-ideal}, $I=f(R)$ is an ideal: it is closed under addition with an arbitrary element, and every product in $R$ belongs to $I$, since $f(xy)=3xy=xy$. Lemma~\ref{tw:lem:retract} gives
\[
 \Var(R)=\Var(I,J),\qquad J=R/\rho_I.
\]
The quotient $J$ has constant multiplication, and its collapsed element is additively absorbing. Lemma~\ref{tw:lem:top-null} classifies $\Var(J)$.

On $I$, both $3x=x$ and $f(x)=x$ hold. Equation~\eqref{tw:eq:split-pq} therefore gives $x=p(x)+q(x)$. It follows that
\[
 I\longrightarrow p(I)\times q(I),\qquad x\longmapsto(p(x),q(x))
\]
is an injective homomorphism, surjective onto each coordinate. Hence $\Var(I)=\Var(p(I),q(I))$. The two images satisfy $x=x+x^2$ and $x=2x+x^2$, respectively, because $p$ and $q$ are idempotent. Apply Lemmas~\ref{tw:lem:p-factor} and~\ref{tw:lem:q-factor}. Finally, their doubling images lie in $\A_0$ by Lemma~\ref{tw:lem:double-general}, so Theorem~\ref{tw:lem:ai-known} generates each of their varieties by its two-element members. All generators obtained lie in $\Var(R)$, as required.
\end{proof}

\begin{theorem}[Completeness and hereditary small generation]\label{tw:thm:exhaustive}\label{tw:thm:basis}
We have $\Upsilon=\Sr$. For every $A\in\Sr$,
\begin{equation}\label{tw:eq:hereditary-generation}
 \Var(A)=\Var\bigl(\E\cap\Var(A)\bigr).
\end{equation}
Consequently, every $\W\le\Sr$ satisfies $\W=\Var(\E\cap\W)$.
\end{theorem}
\begin{proof}
Let $A\in\Upsilon$. Apply Lemma~\ref{tw:lem:retract} to $r:A\to A$. Its image is additively regular and its Rees quotient has constant addition. Proposition~\ref{tw:prop:regular-generation} and Lemma~\ref{tw:lem:zero-add} express the variety of each factor using algebras in $\E\cap\Var(A)$. Since $\E\subseteq\Sr$ by Lemma~\ref{tw:lem:quotients}, both factors lie in $\Sr$, and their subdirect product $A$ does too. This proves $\Upsilon\subseteq\Sr$; \eqref{tw:eq:easy-basis} gives equality. The same argument proves \eqref{tw:eq:hereditary-generation}, since each small generator obtained belongs to $\Var(A)$.

For $\W\le\Sr$, every $A\in\W$ lies in $\Var(\E\cap\Var(A))\subseteq\Var(\E\cap\W)$. This proves $\W\subseteq\Var(\E\cap\W)$, and the reverse inclusion is immediate. No bound on the cardinality or the number of generators of $A$ has been imposed.
\end{proof}

\begin{corollary}\label{tw:cor:absolute-basis}
The list $\mathrm{SR}\cup\Sigma$ is an identity basis for $\Sr$. It has at most 36 entries, each using at most four variables.
\end{corollary}

\subsection{The complete subvariety lattice}\label{tw:sec:closure}

Let $\E=\K\cup\{Q,L_0,R_0\}$ and $\mathfrak F=\{\Var(B):B\subseteq\E\}$. By Theorem~\ref{tw:thm:exhaustive}, $\mathfrak F=\Lat(\Sr)$. On subsets of $\E$, define $\cl(B)$ to be the least superset of $B$ satisfying the following ten implications:
\begin{equation}\label{tw:eq:rules}
\begin{array}{rcl@{\qquad}rcl}
 Q&\Longrightarrow&T_2,&L_0&\Longrightarrow&W_2,\\
 R_0&\Longrightarrow&W_2,&T_2,Z_7&\Longrightarrow&Q,\\
 T_2,Z_8&\Longrightarrow&Q,&L_2,W_2&\Longrightarrow&L_0,\\
 R_2,W_2&\Longrightarrow&R_0,&N_2,W_2&\Longrightarrow&Z_2,\\
 T_2,W_2&\Longrightarrow&Z_2,&Z_7,W_2&\Longrightarrow&Z_2.
\end{array}
\end{equation}
Let $\C=\{C\subseteq\E:\cl(C)=C\}$. Lemma~\ref{tw:lem:quotients} proves that every implication is valid for membership in a variety. Thus
\begin{equation}\label{tw:eq:oneinc}
 \cl(B)\subseteq\E\cap\Var(B),\qquad
 \Var(\cl(B))=\Var(B).
\end{equation}
The reverse inclusion follows from the explicit separating identities in Table~\ref{tw:tab:separators}.

\begin{table}[htbp]
\centering
\caption{Separating identities. The last column lists \emph{all} members of $\E$ in which the identity fails; it holds in every unlisted member.}\label{tw:tab:separators}
\renewcommand{\arraystretch}{1.45}
\begin{tabular}{c l l}
\toprule
Label&Identity&Failure set\\\midrule
$s_1$&$x+yxy\eq x+xy$&$\{L_2\}$\\
$s_2$&$x+yxy\eq x+yx$&$\{R_2\}$\\
$s_3$&$x+2xyx\eq x+2x^2$&$\{M_2\}$\\
$s_4$&$xy+y+yx\eq2x^2+y$&$\{D_2\}$\\
$s_5$&$2x+x^2+y\eq x^2+y$&$\{N_2\}$\\
$s_6$&$x+2x^2\eq3x$&$\{T_2,Q\}$\\
$s_7$&$3x\eq x$&$\{Z_2,W_2,L_0,R_0\}$\\
$s_8$&$x^2\eq x$&$\{N_2,T_2,Z_2,Z_7,Q\}$\\
$s_9$&$3x^2\eq x^2$&$\{W_2,L_0,R_0\}$\\
$s_{10}$&$x+2x^2\eq2x+x^2$&$\{Z_7\}$\\
$s_{11}$&$x+2x^2\eq x+x^2$&$\{Z_8\}$\\
$s_{12}$&$3x\eq2x$&$\{Z_7,Z_8,Q\}$\\
$s_{13}$&$xyx\eq yx$&$\{L_2,L_0\}$\\
$s_{14}$&$xyx\eq xy$&$\{R_2,R_0\}$\\
\bottomrule
\end{tabular}
\end{table}

\begin{lemma}\label{tw:lem:separators}
Table~\ref{tw:tab:separators} is correct.
\end{lemma}
\begin{proof}
For the two-element members, apply the tests in Proposition~\ref{tw:prop:test}, or substitute the entries of Table~\ref{tw:tab:two}. For $L_0,R_0$, every term containing addition is constantly zero, so all identities with sums on both sides hold. The remaining identities are checked in the left- and right-zero bands with an absorbing zero. For $Q$, any polynomial with a word of length at least two is constantly $\infty$. A sum of copies of $x$ has its group value on $\{e,a\}$ and value $\infty$ at $x=\infty$. These observations give precisely the failure sets listed. Each failure set can thus be determined directly from these evaluations.
\end{proof}

\begin{lemma}\label{tw:lem:separateclosed}
If $C\in\C$ and $S\in\E\setminus C$, some identity in Table~\ref{tw:tab:separators} holds in every member of $C$ and fails in $S$.
\end{lemma}
\begin{proof}
For $S=L_2,R_2,M_2,D_2,N_2,Z_7,Z_8$, use the corresponding singleton failure set. The remaining cases are as follows.
\begin{enumerate}[label=\textup{(\arabic*)},leftmargin=*]
\item If $T_2\notin C$, then $Q\notin C$, and $s_6$ applies.
\item If $W_2\notin C$, then $L_0,R_0\notin C$, and $s_9$ applies.
\item Suppose $Z_2\notin C$. If $W_2\notin C$, also $L_0,R_0\notin C$, and $s_7$ applies. If $W_2\in C$, the last three implications force $N_2,T_2,Z_7\notin C$, and $Q\notin C$ follows from $Q\Rightarrow T_2$. Now $s_8$ applies.
\item Suppose $Q\notin C$. If $T_2\notin C$, use $s_6$. Otherwise both $Z_7,Z_8$ are absent by \eqref{tw:eq:rules}, and $s_{12}$ applies.
\item Suppose $L_0\notin C$. If $W_2\notin C$, use $s_9$. Otherwise $L_2\notin C$, and $s_{13}$ applies. The dual argument uses $s_{14}$ when $R_0\notin C$.
\end{enumerate}
Each chosen failure set is disjoint from $C$ and contains $S$.
\end{proof}

\begin{theorem}[Complete classification of all subvarieties]\label{tw:thm:closure}
For every $B\subseteq\E$,
\begin{equation}\label{tw:eq:closureexact}
 \E\cap\Var(B)=\cl(B).
\end{equation}
The map
\[
 \Psi:\C\longrightarrow\Lat(\Sr),\qquad C\longmapsto\Var(C)
\]
is an order isomorphism. For $C,D\in\C$,
\begin{align}
 \Psi(C)\vee_{\Lat(\Sr)}\Psi(D)&=\Psi\bigl(\cl(C\cup D)\bigr),\label{tw:eq:join}\\
 \Psi(C)\wedge_{\Lat(\Sr)}\Psi(D)&=\Psi(C\cap D).\label{tw:eq:internalmeet}
\end{align}
These are the meet and join operations of the whole subvariety lattice.
\end{theorem}
\begin{proof}
Put $C=\cl(B)$. By \eqref{tw:eq:oneinc}, $\Var(C)=\Var(B)$ and $C\subseteq\E\cap\Var(B)$. If $S\in\E\setminus C$, Lemma~\ref{tw:lem:separateclosed} gives an identity satisfied by $C$ and hence by $\Var(C)$, but not by $S$. This proves \eqref{tw:eq:closureexact}.

By Theorem~\ref{tw:thm:exhaustive}, every subvariety of $\Sr$ equals $\Var(C)$ for a closed $C$, so $\Psi$ is surjective. Equation~\eqref{tw:eq:closureexact} gives injectivity and shows that $\Var(C)\subseteq\Var(D)$ implies $C\subseteq D$. The converse is immediate. Closed subsets form a lattice whose meet is intersection and whose join is the closure of union. Finally, $\Var(C)\vee\Var(D)=\Var(C\cup D)=\Var(\cl(C\cup D))$, which proves agreement of joins.
\end{proof}

\begin{corollary}[Relative identity bases]\label{tw:cor:relative-bases}
Let $F_i$ denote the failure set of $s_i$ in Table~\ref{tw:tab:separators}. For each closed $C\subseteq\E$,
\begin{equation}\label{tw:eq:relative-basis}
 \Var(C)=\{S\in\Sr:S\models s_i\text{ whenever }F_i\cap C=\varnothing\}.
\end{equation}
Thus every subvariety has a relative basis consisting of at most fourteen displayed identities.
\end{corollary}
\begin{proof}
Let $\W$ be the variety on the right. It contains $\Var(C)$. Lemma~\ref{tw:lem:separateclosed} shows that every member of $\E\setminus C$ fails an identity imposed on $\W$. Hence $\E\cap\W=C$. Theorem~\ref{tw:thm:exhaustive} gives $\W=\Var(C)$.
\end{proof}

\endgroup
\section{Complete representatives and their basis cases}\label{app:catalogue}
Each operation is written in row order on $\{0,1,2,3\}$. To save space,
two successive entries $a,b$ are encoded by the hexadecimal digit $4a+b$.
Thus each eight-digit string specifies all sixteen entries without loss.
A starred index denotes a nonfinitely based algebra. The last column
gives a theorem, section or cited classification, followed by a source algebra
when a transfer is used. Every isomorphism type occurs once.
\begingroup\fontsize{7.5}{8.5}\selectfont\setlength{\tabcolsep}{1.5pt}
\begin{longtable}{@{}rlll@{\hspace{1em}}rlll@{\hspace{1em}}rlll@{}}
\caption{All four-element semirings and their basis cases.}\label{tab:all}\\
\toprule $i$&$+$&$\cdot$&Case&$i$&$+$&$\cdot$&Case&$i$&$+$&$\cdot$&Case\\\midrule\endfirsthead
\toprule $i$&$+$&$\cdot$&Case&$i$&$+$&$\cdot$&Case&$i$&$+$&$\cdot$&Case\\\midrule\endhead
\bottomrule\endlastfoot
$1$ & \texttt{00000000} & \texttt{00000000} & \ref{prop:constant-operation} & $2$ & \texttt{00100000} & \texttt{00000000} & \ref{prop:constant-operation} & $3$ & \texttt{00200000} & \texttt{00000000} & \ref{prop:constant-operation}\\
$4$ & \texttt{00182000} & \texttt{00000000} & \ref{prop:constant-operation} & $5$ & \texttt{000C3000} & \texttt{00000000} & \ref{prop:constant-operation} & $6$ & \texttt{002C3000} & \texttt{00000000} & \ref{prop:constant-operation}\\
$7$ & \texttt{003C3000} & \texttt{00000000} & \ref{prop:constant-operation} & $8$ & \texttt{00141400} & \texttt{00000000} & \ref{prop:constant-operation} & $9$ & \texttt{00182400} & \texttt{00000000} & \ref{prop:constant-operation}\\
$10$ & \texttt{00100800} & \texttt{00000000} & \ref{prop:constant-operation} & $11$ & \texttt{00300800} & \texttt{00000000} & \ref{prop:constant-operation} & $12$ & \texttt{00141800} & \texttt{00000000} & \ref{prop:constant-operation}\\
$13$ & \texttt{00300C00} & \texttt{00000000} & \ref{prop:constant-operation} & $14$ & \texttt{003C3C00} & \texttt{00000000} & \ref{prop:constant-operation} & $15$ & \texttt{15404040} & \texttt{00000000} & \ref{prop:term-expansion}\\
$16$ & \texttt{10451010} & \texttt{00000000} & \ref{prop:term-expansion} & $17$ & \texttt{10551010} & \texttt{00000000} & \ref{thm:depth}; $E_{4}$ & $18$ & \texttt{10451810} & \texttt{00000000} & \ref{f3:cor:linear}\\
$19$ & \texttt{10551810} & \texttt{00000000} & \ref{thm:depth}; $B_{18}$ & $20$ & \texttt{10451C10} & \texttt{00000000} & \ref{prop:term-expansion} & $21$ & \texttt{10551C10} & \texttt{00000000} & \ref{thm:depth}; $E_{4}$\\
$22$ & \texttt{1A559090} & \texttt{00000000} & \ref{f3:cor:linear} & $23$ & \texttt{001A2020} & \texttt{00000000} & \ref{prop:constant-operation} & $24$ & \texttt{0A0AA0A0} & \texttt{00000000} & \ref{prop:term-expansion}\\
\end{longtable}
\endgroup
\section{Additive orbit counts}\label{app:orbits}
\begingroup\footnotesize
\setlength{\tabcolsep}{3pt}
\begin{longtable}{rlrrrrrr}
\caption{The 58 additive types and their classification. Here $m_+=|\mathcal L(p)|$, $f_+$ counts fixed multiplications over nonidentity additive automorphisms, and Total is $(m_++f_+)/|\mathrm{Aut}(+)|$.}\label{tab:additive}\\
\toprule Type & Addition & $|\mathrm{Aut}(+)|$ & $m_+$ & $f_+$ & FB & NFB & Total\\\midrule\endfirsthead
\toprule Type & Addition & $|\mathrm{Aut}(+)|$ & $m_+$ & $f_+$ & FB & NFB & Total\\\midrule\endhead
\midrule\multicolumn{8}{r}{Continued}\\\endfoot
\bottomrule\endlastfoot
1 & \texttt{0000000000000000} & 6 & 442 & 98 & 90 & 0 & 90\\
2 & \texttt{0000000000000001} & 1 & 31 & 0 & 31 & 0 & 31\\
3 & \texttt{0000000000000003} & 2 & 101 & 21 & 61 & 0 & 61\\
4 & \texttt{0000000000010010} & 2 & 41 & 7 & 24 & 0 & 24\\
5 & \texttt{0000000000010011} & 1 & 24 & 0 & 24 & 0 & 24\\
6 & \texttt{0000000000010012} & 1 & 4 & 0 & 4 & 0 & 4\\
7 & \texttt{0000000000020023} & 1 & 17 & 0 & 17 & 0 & 17\\
8 & \texttt{0000000000100001} & 2 & 18 & 4 & 11 & 0 & 11\\
\end{longtable}
\endgroup

\section{Correspondence and structural maps}\label{app:maps}
\begingroup\footnotesize
\begin{longtable}{@{}rllll@{\hspace{1.4em}}rllll@{}}
\caption{Correspondence with cited four-element tables. The permutation maps source labels $0,1,2,3$ to local labels.}\label{tab:literature}\\
\toprule $i$ & $k$ in $S_{(4,k)}$ & Permutation & Status & Source & $i$ & $k$ in $S_{(4,k)}$ & Permutation & Status & Source\\
\midrule\endfirsthead
\toprule $i$ & $k$ in $S_{(4,k)}$ & Permutation & Status & Source & $i$ & $k$ in $S_{(4,k)}$ & Permutation & Status & Source\\
\midrule\endhead
\bottomrule\endlastfoot
103 & 1 & \texttt{0123} & FB & \cite{PartI} & 105 & 276 & \texttt{0123} & FB & \cite{PartIII}\\
106 & 388 & \texttt{0123} & FB & \cite{PartII} & 107 & 371 & \texttt{1023} & FB & \cite{PartIII}\\
109 & 58 & \texttt{1023} & FB & \cite{PartI} & 111 & 470 & \texttt{1023} & FB & \cite{PartII}\\
114 & 480 & \texttt{1302} & FB & \cite{PartII} & 117 & 387 & \texttt{1203} & FB & \cite{PartIII}\\
190 & 3 & \texttt{0231} & FB & \cite{PartI} & 193 & 287 & \texttt{0123} & FB & \cite{PartIII}\\
194 & 389 & \texttt{0123} & FB & \cite{PartII} & 198 & 278 & \texttt{0231} & FB & \cite{PartIII}\\
202 & 386 & \texttt{2103} & FB & \cite{PartIII} & 204 & 373 & \texttt{2031} & FB & \cite{PartIII}\\
209 & 479 & \texttt{2103} & FB & \cite{PartII} & 249 & 2 & \texttt{0231} & FB & \cite{PartI}\\
\end{longtable}
\endgroup
\subsection*{The chain additive reducts}
Table~\ref{tab:chain-correspondence} gives the complete correspondence
with Table 28 of \cite{PartIV}. Each permutation lists the local
images of the source elements $1,2,3,4$, in that order. Unlike the
preceding table, the source labels here start at 1. The source addition
is $r+s=\min(r,s)$ on these labels. Thus both operations, including
the direction of multiplication, are specified by this correspondence.
\begingroup\footnotesize
\begin{longtable}{@{}rlllllll@{\hspace{1.4em}}rlllllll@{}}
\caption{Chain representatives: source index, local index and isomorphism.}
\label{tab:chain-correspondence}\\
\toprule $k$ & $i$ & $\pi$ & Status &
$k$ & $i$ & $\pi$ & Status & $k$ & $i$ & $\pi$ & Status &
$k$ & $i$ & $\pi$ & Status\\
\midrule\endfirsthead
\toprule $k$ & $i$ & $\pi$ & Status &
$k$ & $i$ & $\pi$ & Status & $k$ & $i$ & $\pi$ & Status &
$k$ & $i$ & $\pi$ & Status\\
\midrule\endhead
\bottomrule\endlastfoot
481 & 112 & \texttt{0123} & FB & 482 & 262 & \texttt{0231} & FB & 483 & 266 & \texttt{0321} & FB & 484 & 207 & \texttt{0231} & FB\\
485 & 487 & \texttt{0321} & FB & 486 & 430 & \texttt{0321} & FB & 487 & 1340 & \texttt{0231} & FB & 488 & 1382 & \texttt{0231} & FB\\
489 & 1343 & \texttt{0321} & FB & 490 & 1077 & \texttt{0231} & FB & 491 & 321 & \texttt{0321} & FB & 492 & 516 & \texttt{0321} & FB\\
493 & 504 & \texttt{0321} & FB & 494 & 294 & \texttt{0321} & FB & 495 & 457 & \texttt{0321} & FB & 496 & 1391 & \texttt{0321} & FB\\
497 & 875 & \texttt{0312} & FB & 498 & 585 & \texttt{0312} & FB & 499 & 766 & \texttt{0312} & FB & 500 & 786 & \texttt{0312} & FB\\
501 & 1455 & \texttt{0321} & FB & 502 & 814 & \texttt{0312} & FB & 503 & 767 & \texttt{0321} & FB & 504 & 1189 & \texttt{0321} & FB\\
505 & 586 & \texttt{0321} & FB & 506 & 1626 & \texttt{1203} & FB & 507 & 1915 & \texttt{0312} & FB & 508 & 1920 & \texttt{0312} & FB\\
509 & 1928 & \texttt{0312} & FB & 510 & 1916 & \texttt{0321} & FB & 511 & 1863 & \texttt{0312} & FB & 512 & 1627 & \texttt{1230} & FB\\
\end{longtable}
\endgroup

\subsection{Structural maps}
\begingroup\footnotesize
\setlength{\tabcolsep}{2pt}
\begin{longtable}{@{}rllll@{\hspace{1.4em}}rllll@{}}
\caption{Unary term retractions used in the primary reductions. A dash means that no depth bound is used.}\label{tab:retractions}\\
\toprule $i$ & $h(x)$ & $h(0)h(1)h(2)h(3)$ & Image & $d$ & $i$ & $h(x)$ & $h(0)h(1)h(2)h(3)$ & Image & $d$\\
\midrule\endfirsthead
\toprule $i$ & $h(x)$ & $h(0)h(1)h(2)h(3)$ & Image & $d$ & $i$ & $h(x)$ & $h(0)h(1)h(2)h(3)$ & Image & $d$\\
\midrule\endhead
\bottomrule\endlastfoot
17 & $(x+x)$ & \texttt{0100} & $E_{4}$ & 1 & 19 & $(x+x)$ & \texttt{0120} & $B_{18}$ & 1\\
21 & $(x+(xx))$ & \texttt{0100} & $E_{4}$ & 2 & 33 & $(x+x)$ & \texttt{0110} & $E_{4}$ & 1\\
34 & $(x+x)$ & \texttt{0120} & $B_{19}$ & 1 & 44 & $(x+x)$ & \texttt{0120} & $B_{20}$ & 1\\
64 & $(x+x)$ & \texttt{0111} & $E_{4}$ & 1 & 66 & $(x+x)$ & \texttt{0121} & $B_{18}$ & 1\\
67 & $(x+x)$ & \texttt{0121} & $B_{19}$ & 1 & 68 & $(x+x)$ & \texttt{0121} & $B_{20}$ & 1\\
69 & $(x+(xx))$ & \texttt{0111} & $E_{4}$ & 2 & 81 & $(x+x)$ & \texttt{0121} & $B_{18}$ & 1\\
83 & $(x+x)$ & \texttt{0121} & $B_{20}$ & 1 & 84 & $(x+x)$ & \texttt{0121} & $B_{19}$ & 1\\
124 & $(x+x)$ & \texttt{0120} & $B_{30}$ & 1 & 125 & $(x+x)$ & \texttt{0120} & $B_{34}$ & 1\\
\end{longtable}
\endgroup
\begingroup\footnotesize
\setlength{\tabcolsep}{2pt}
\begin{longtable}{@{}p{.58\textwidth}llll@{}}
\caption{Linear interpretations of the other operation. A dash means an absent contribution.}\label{tab:linear}\\
\toprule Indices $i$ & Retained & $f(x)$ & $g(x)$ & Constant term\\
\midrule\endfirsthead
\toprule Indices $i$ & Retained & $f(x)$ & $g(x)$ & Constant term\\
\midrule\endhead
\bottomrule\endlastfoot
18, 22, 25, 30, 31, 35, 36, 38, 40, 42, 43, 45, 46, 47, 48, 49, 50, 54, 55, 56, 57, 65, 70, 71, 72, 73, 74, 76, 77, 78, 82, 85, 86, 88, 89, 90, 91, 93, 94, 95, 96, 98, 100, 101, 102 & + & $-$ & $-$ & $(xx)$\\[2pt]
140, 141, 152, 154, 161, 163, 164, 177, 178, 180, 186, 188, 223, 224, 227, 230, 233, 234, 243, 244, 246, 272, 525, 547, 548, 552, 553, 556, 584, 586, 589, 658, 863, 874, 1101, 1107, 1118, 1126, 1270 & + & $(xx)$ & $(xx)$ & $-$\\[2pt]
156, 182, 544, 1106, 1125 & + & $x$ & $(xx)$ & $-$\\[2pt]
278, 282, 285, 389, 392, 988 & + & $(xx)$ & $(x+x)$ & $-$\\[2pt]
334, 335, 342, 344, 345, 346, 347, 348, 349, 350, 351, 357, 360, 361, 362, 363, 364, 365, 366, 367, 880, 883, 890, 891, 892, 893, 894, 895, 896, 903 & + & $(xx)$ & $-$ & $-$\\[2pt]
414, 418, 421, 1332, 1335, 1618 & + & $(x+x)$ & $(xx)$ & $-$\\[2pt]
609, 716, 1195, 1429, 1460, 1578, 1594, 1662, 1795, 1839, 1879, 1962, 1963, 1982, 1984, 1998, 1999, 2110, 2111, 2218, 2245, 2252 & $\cdot$ & $x$ & $(x+x)$ & $-$\\[2pt]
1031, 1032, 1039, 1041, 1042, 1043, 1044, 1045, 1046, 1047, 1048, 1054, 1057, 1058, 1059, 1060, 1061, 1062, 1063, 1064, 1520, 1523, 1530, 1531, 1532, 1533, 1534, 1535, 1536, 1543 & + & $-$ & $(xx)$ & $-$\\[2pt]
\end{longtable}
\endgroup

\subsection{Opposite correspondences}\label{app:opposite-maps}
The source and target lists in each row are paired in order. The permutation
gives the images of $(0,1,2,3)$ in an isomorphism $A_i^{\mathrm{op}}\to A_j$.
The cited sections provide the source bases; Corollary~\ref{cor:catalogue-transfers}
provides all reversed bases and their variety relations.
\begingroup\small
\begin{longtable}{@{}c p{.31\textwidth} p{.31\textwidth} c@{}}
\caption{Opposite correspondences for the normal-form families.}\label{tab:opposite-transfers}\\
\toprule Section&Source indices $i$&Target indices $j$&Permutation\\\midrule\endfirsthead
\toprule Section&Source indices $i$&Target indices $j$&Permutation\\\midrule\endhead
\bottomrule\endlastfoot
\ref{sec:family31} & 277, 281, 289, 293, 297, 400 & 413, 417, 425, 429, 433, 1343 & $(0,1,2,3)$\\
\ref{sec:family31} & 990 & 1627 & $(1,0,3,2)$\\
\ref{sec:family32} & 387, 390, 395, 398 & 1330, 1333, 1338, 1341 & $(0,1,2,3)$\\
\ref{sec:family40} & 851 & 1253 & $(0,1,2,3)$\\
\ref{sec:family42} & 1303 & 1487 & $(0,2,1,3)$\\
\ref{sec:family42} & 1167, 1573, 1724 & 689, 933, 1293 & $(0,1,2,3)$\\
\ref{sec:family43} & 1150, 1558, 1711 & 672, 918, 1280 & $(0,1,2,3)$\\
\ref{sec:family44} & 1149, 1557 & 671, 917 & $(0,1,2,3)$\\
\ref{sec:family45} & 1299, 1300 & 1481, 1482 & $(0,2,1,3)$\\
\ref{sec:family46} & 1156, 1157, 775, 778, 1556 & 678, 679, 803, 806, 916 & $(0,1,2,3)$\\
\ref{sec:family52} & 819, 822, 827 & 1447, 1450, 1455 & $(0,1,2,3)$\\
\ref{sec:family53} & 1491, 1493, 1498 & 1652, 1653, 1658 & $(0,2,1,3)$\\
\ref{sec:family54} & 1421, 1904, 1423, 1914 & 709, 994, 711, 1004 & $(0,1,2,3)$\\
\ref{sec:family55} & 1424, 1916 & 712, 1006 & $(0,1,2,3)$\\
\ref{sec:family55} & 1870 & 1318 & $(0,2,1,3)$\\
\ref{sec:family57} & 1866 & 1310 & $(0,2,1,3)$\\
\ref{sec:family57} & 1909 & 999 & $(0,1,2,3)$\\
\end{longtable}
\endgroup

\section{Finite term-function certificates}\label{app:finite-witnesses}

Use $R_{D,C}$ and $T_D$ from Proposition~\ref{prop:closure-certificate}.
Their finite closures are specified recursively.
Set $R_0=\{(\pi_j,c_j):1\leq j\leq k\}$ and
\[
 R_{r+1}=R_r\cup\{u+v,uv:u,v\in R_r\}.
\]
The sequence stabilizes in the finite algebra $D^{D^k}\times C$;
its stable value is $R_{D,C}$. Projection gives $T_D$. Using the
operation tables in Appendices~\ref{app:catalogue} and \ref{app:small},
the resulting cardinalities in both directions are given in
Table~\ref{tab:finite-witnesses}. Each displayed tuple generates its
target, and each pair of cardinalities is equal. Thus
Proposition~\ref{prop:closure-certificate} proves both inclusions
for every row. The complete bases of the direct endpoints are proved
in the body of the paper.

\begingroup\small\setlength{\tabcolsep}{3pt}
\begin{longtable}{cllcrcr}
\caption{Finite closure certificates for the designated variety equalities.}
\label{tab:finite-witnesses}\\
\toprule
Witness & $C$ & $D$ & $\mathbf c$ & $|T_D|=|R_{D,C}|$ & $\mathbf d$ & $|T_C|=|R_{C,D}|$\\
\midrule
\endfirsthead
\multicolumn{7}{c}{\tablename~\thetable\ (continued)}\\
\toprule
Witness & $C$ & $D$ & $\mathbf c$ & $|T_D|=|R_{D,C}|$ & $\mathbf d$ & $|T_C|=|R_{C,D}|$\\
\midrule
\endhead
\midrule\multicolumn{7}{r}{Continued on the next page}\\
\endfoot
\bottomrule
\endlastfoot
W01 & $A_{292}$ & $B_{54}$ & $(1,2,3)$ & 272 & $(1,2)$ & 14\\
W02 & $A_{294}$ & $A_{397}$ & $(1,2,3)$ & 44 & $(1,3)$ & 10\\
W03 & $A_{428}$ & $B_{83}$ & $(1,2,3)$ & 272 & $(1,2)$ & 14\\
W04 & $A_{430}$ & $A_{1340}$ & $(1,2,3)$ & 44 & $(1,3)$ & 10\\
W05 & $A_{455}$ & $B_{88}$ & $(1,2,3)$ & 272 & $(1,2)$ & 14\\
W06 & $A_{457}$ & $A_{1358}$ & $(1,2,3)$ & 57 & $(1,3)$ & 11\\
W07 & $A_{475}$ & $A_{1367}$ & $(1,2,3)$ & 90 & $(1,3)$ & 13\\
W08 & $A_{616}$ & $A_{735}$ & $(1,2,3)$ & 33 & $(2,3)$ & 11\\
W09 & $A_{619}$ & $A_{736}$ & $(1,2,3)$ & 207 & $(2,3)$ & 19\\
W10 & $A_{620}$ & $A_{720}$ & $(1,2,3)$ & 61 & $(1,2,3)$ & 61\\
W11 & $A_{635}$ & $A_{765}$ & $(1,2,3)$ & 280 & $(2,3)$ & 28\\
W12 & $A_{730}$ & $A_{1215}$ & $(2,3)$ & 30 & $(1,2,3)$ & 299\\
W13 & $A_{759}$ & $A_{1244}$ & $(2,3)$ & 30 & $(1,2,3)$ & 299\\
W14 & $A_{798}$ & $A_{800}$ & $(1,2,3)$ & 307 & $(1,2,3)$ & 307\\
W15 & $A_{825}$ & $B_{54}$ & $(1,2,3)$ & 272 & $(1,2)$ & 14\\
W16 & $A_{990}$ & $A_{297}$ & $(2,3)$ & 19 & $(1,2,3)$ & 105\\
W17 & $A_{1188}$ & $A_{1190}$ & $(1,2,3)$ & 307 & $(1,2,3)$ & 307\\
W18 & $A_{1209}$ & $A_{720}$ & $(1,2,3)$ & 61 & $(1,2,3)$ & 61\\
W19 & $A_{1292}$ & $A_{813}$ & $(1,2,3)$ & 342 & $(1,2,3)$ & 342\\
W20 & $A_{1391}$ & $A_{1397}$ & $(1,2,3)$ & 90 & $(1,3)$ & 13\\
W21 & $A_{1453}$ & $B_{83}$ & $(1,2,3)$ & 272 & $(1,2)$ & 14\\
W22 & $A_{1468}$ & $B_{88}$ & $(1,2,3)$ & 272 & $(1,2)$ & 14\\
W23 & $A_{1586}$ & $A_{720}$ & $(2,3)$ & 11 & $(1,2,3)$ & 61\\
W24 & $A_{1587}$ & $A_{1197}$ & $(2,3)$ & 14 & $(1,2,3)$ & 83\\
W25 & $A_{1588}$ & $A_{1199}$ & $(2,3)$ & 13 & $(1,2,3)$ & 106\\
W26 & $A_{1597}$ & $A_{736}$ & $(2,3)$ & 19 & $(2,3)$ & 19\\
W27 & $A_{1610}$ & $A_{1228}$ & $(2,3)$ & 14 & $(1,2,3)$ & 132\\
W28 & $A_{1723}$ & $A_{785}$ & $(1,2,3)$ & 342 & $(1,2,3)$ & 342\\
W29 & $A_{1756}$ & $A_{607}$ & $(2,3)$ & 33 & $(2,3)$ & 33\\
W30 & $A_{1783}$ & $A_{1200}$ & $(1,2,3)$ & 145 & $(1,2,3)$ & 145\\
W31 & $A_{1785}$ & $A_{1199}$ & $(1,2,3)$ & 106 & $(1,2,3)$ & 106\\
W32 & $A_{1793}$ & $A_{1206}$ & $(1,3)$ & 34 & $(1,2,3)$ & 1383\\
W33 & $A_{1797}$ & $A_{1198}$ & $(1,3)$ & 21 & $(1,2,3)$ & 220\\
W34 & $A_{1807}$ & $A_{761}$ & $(1,2,3)$ & 145 & $(1,2,3)$ & 145\\
W35 & $A_{1808}$ & $A_{1228}$ & $(1,2,3)$ & 132 & $(1,2,3)$ & 132\\
W36 & $A_{1817}$ & $A_{1190}$ & $(1,2,3)$ & 307 & $(1,2,3)$ & 307\\
W37 & $A_{1823}$ & $A_{800}$ & $(1,2,3)$ & 307 & $(1,2,3)$ & 307\\
W38 & $A_{2006}$ & $A_{2113}$ & $(2,3)$ & 14 & $(1,2,3)$ & 144\\
W39 & $A_{2251}$ & $A_{2119}$ & $(0,3)$ & 35 & $(1,2,3)$ & 1563\\
W40 & $A_{2255}$ & $A_{2114}$ & $(0,1,3)$ & 192 & $(1,2,3)$ & 192\\
W41 & $A_{2259}$ & $A_{2113}$ & $(0,1,3)$ & 144 & $(1,2,3)$ & 144\\
W42 & $A_{2263}$ & $A_{2123}$ & $(0,1,3)$ & 998 & $(1,2,3)$ & 998\\
\end{longtable}
\endgroup

\par\bigskip
\Needspace{6\baselineskip}
\begingroup
\raggedright
{\small\noindent\paperaffiliation\par}
\smallskip
{\footnotesize\noindent\paperemails\par}
\endgroup
\end{document}